\documentclass[12pt] {article}
\usepackage{amsmath,amssymb}

\title{The Beltrami equation with parameters and uniformization of holomorphic foliations 
with hyperbolic leaves\thanks{The work was supported by grant RFBR 16-01-00748.}}

\author{A.A. Shcherbakov\thanks{ E-mail: arsshcher@mail.ru}}

\begin{document}

\maketitle

\newtheorem{theorem}{Theorem}
\newtheorem{lemma}{Lemma}
\newtheorem{corollary}{Corollary}
\newtheorem{proposition}{Proposition}
\newtheorem{definition}{Definition}

\newcommand{\bbP}{\mathbb{P}}
\newcommand{\bbQ}{\mathbb{Q}}
\newcommand{\bbZ}{\mathbb{Z}}
\newcommand{\bbC}{\mathbb{C}}
\newcommand{\mF}{\mathcal{F}}
\newcommand{\mA}{\mathcal{A}}
\newcommand{\mB}{\mathcal{B}}
\newcommand{\mC}{\mathcal{C}}
\newcommand{\mS}{\mathcal{S}}
\newcommand{\pd}{\partial}
\newcommand{\vpf}{\varphi}



\bigskip                                                                                                                                                                                                                                                                               

\section*{Introduction}

Let $X$ be a complex manifold. We say that a foliation with singularities is defined on $X$ if
there exists an analytic subset $\Sigma$ of codimension at least two and a foliation
of its complement by analytic curves that cannot be extended to a neighborhood of any point of $\Sigma$. A foliation  
can be locally defined by polynomial vector fields.  In generic case the singular set $\Sigma$ consists
of isolated points.

A covering manifold of leaves of a foliation was defined in [Il1], [Il2]. Let $\mF$ be a foliation with singularities on a
complex manifold $X$ and let $B$ be
transversal cross-section. Let $\vpf_p$ be a leaf passing through a point $p\in B$ and let $\hat\vpf_p$ be the universal covering over this
leaf with the marked point $p$. Define $M=\bigcup\limits_{p\in B} \hat\vpf_p$. It is shown in [Il1], [Il2] that at least in affine case or, in more general
Stein case, a topology and a complex structure on this union can be defined so that it is a complex manifold with locally
biholomorphic projection $\tilde\pi :M\to X$ and a holomorphic section $B\to M$ right inverse to the 
holomorphic retraction $\pi :M\to B$. For any leaf $\vpf_p$ the restriction
of $\tilde\pi$ to $\hat\vpf_p$ is the universal covering map over $\vpf_p$. For a foliation of a compact manifold the manifold of
universal coverings can be non-Hausdorff but in a generic case it is Hausdorff (see [Br1], [Br2]). It is possible to define
a Hausdorff universal covering for general foliations of compact Kahler manifolds if we include the singular points in the leaves
in some not generic cases ([Br1], [Br3]) but here we don't consider such situations.

Let $T_{\mF}$ be the linear bundle tangent to the leaves. If this bundle is negative, then there exists 
an hermitean metric on $X\setminus\Sigma$ and restriction of this metric on each leaf has
a negative curvature. For generic such foliation each leaf is hyperbolic ([Gl1] or [LN]).
In particular, it's true for a generic foliation of $\bbC\bbP^n$. The uniformizing map of every leaf is
unique modulo authomorphisms of the disk, and after some normalization (to get uniqueness) we may ask: how the uniformizing map of
$\vpf_p$ depends on the point $p$? Equivalently, we may put on every leaf its Poincare metric, i. e., the unique complete hermitian
metric of curvature -1 and ask about dependence of this metric
on the point $p$. It is known that the Poincare metric is continuous [V] and even Holder-continuous [DNS].
The simultaneous uniformization conjecture states that there exists an analytic in $p$ biholomorphizm 
of $\vpf_p$ onto an appropriate $p$-depending domain on the Riemann sphere. It is known 
that this conjecture is wrong for general foliation in dimension of more than two or even for foliations of general
two-dimensional manifolds [Gl2]. It is not known is this conjecture true or not for generic foliations of $\bbC^2$ or $\bbC\bbP^2$. 

One of the main problems of the theory of holomorphic foliations is the problem of analytic continuation of the Poincare map
defined on a transversal to the leaves and the related problem of the persistence of cycles. It was shown in [Il3] that these
problems have the positive solution if there exists an analytic simultaneous uniformization and the image domains 
continuously depend on the initial conditions. In the absence of an analytic simultaneous uniformization we can consider the
results below as its more feeble version. We hope that these results can be useful in following attempts to clear the situation
with the persistence problems. Though the situation can not be simple. There exists examples of the non-extendability of the
Poincare map though for rather special cases [CDFG]. 

 There were shown in [Sh1], [Sh2] that for generic foliation with negative $T_{\mF}$ we can define the complex structure
on $M$ as an almost complex structure on the product $B\times D$ ($D$ is the unit disk) and
this almost complex structure can be defined by forms of type (1,0)
$$
  dz_i ,\,\,i=1,...,n ,               \eqno (0.1)
$$  
$$
  dw+\mu d\bar w +<c,d\bar z> ,         \eqno (0.2)
$$
where $n+1$ is the dimension of $X$, $z,w$ are charts on $B$ and $D$ correspondingly, $c$ is a smooth
vector-function, $<c,d\bar z>=c_1 d\bar z_1 +...+c_n d\bar z_n$, $\mu$ is a smooth function
satisfying the estimate $|\mu |\le d<1$ for some non-negative $d$. I.e., we can say that
this almost complex structure is quasiconformal on each fiber. 

The conclusion made in [Sh1] that the Poincare metric smoothly depends on a base point isn't
correct because there isn't satisfied the sufficient condition: uniform boundedness along the fibers of
derivatives with respect to the parameters, i.e., to the coordinates on the base (see [AhB]).
It was shown by B. Deroin that the Poincare metric isn't smooth for the foliation of a neighborhood of 
a hyperbolic singular point.

However in a generic case there exists a finitely smooth map holomorphic on the fibers and
mapping each fiber on a bounded domain in $\bbC$ continuously depending on a base point.
Moreover, there exist estimates for derivatives of this map similar to the estimates of
$\mu$ and $c$ obtained in [Sh2]. Now we formulate our main result. We denote by $\pd^{(k)}f$
any derivative with respect to the variable $w,\bar w ,z_i, \bar z_i$ of the total order $|(k)|$.
 
\begin{theorem}  Suppose $X$ is a compact complex manifold of dimension $n+1$ and $\mF$ is a holomorphic foliation of 
$X$ with negative with $T_{\mF}$.  Suppose that the singular set $\Sigma$ is finite and in some neighborhood
of each singular point the vector field locally defining the foliation is analytically
linearizable and the linear part is diagonalizable. Let $M$ be a manifold of universal
covering with simply connected base $B$ and let the complex structure on $M$ be defined
by forms (0.1), (0.2). Then for every integer $p\ge 0$ there
exists a fiberwise map $f:M\to B\times\bbC$ differentiable up to the order $p$, holomorphic on the fibers,
continuously depending on a base point in $C^0 (\bbC )$ and satisfying the estimates 
$$
  |f_w (w,z)|\le C(1-|w|)^{-\alpha},\,|f_{\bar w}(w,z)|\le C(1-|w|)^{-\alpha},\,0\le\alpha <1 ,  (0.3)
$$
$$
  |f_{w^m \bar w^l}(w,z)|/|f_w (w,z)|\le C(1-|w|)^{-(m+l-1)},\,m+l\le p ,          (0.4)
$$
$$
   |\pd^{(k)}f (w.z)|\le C(1-|w|)^{-N(p)},\, |(k)|\le p .   (0.5)
$$
Here $C$ is some uniform constant and $N(p)$ is a constant depending on $p$.
\end{theorem}

Since our complex structure is defined by forms (0.1), (0.2) the assertion that $f$ is holomorphic on the fibers
means that $f$ satisfies the Beltrami equation $f_{\bar z}=\mu f_z$. The proof of the theorem
is based on the estimates obtained in [Sh2] and the present article can be considered as a continuation
of that work. In fact, all that follows is a study of the Beltrami equation with a coefficient 
depending on parameters and satisfying the estimates of [Sh2].    

\section{Preliminary notes and the sketch of the proof}

In what follows $D$ is the unit disk, $D_r$ is the disk or radii $r$ centered at zero,
$D_{a,r}$ is the disk or radii $r$ centered at $a$. Suppose $f$ is a function of a vector variable
$z$ of dimension $n$ and of a scalar variable $w$, and $(k)$ is a multi-index 
$(k)=\{k_0 ,k_{\bar 0}, k_1 ,...,k_n ,k_{\bar 1},...,k_{\bar n}\}$. 
We denote by $\pd^{(k)}f$ or $f_{(k)}$ the derivative 
$f_{w^{k_0}{\bar w}^{k_{\bar 0}}z_1^{k_1}...z_n^{k_n}{\bar z}_1^{k_{\bar 1}}...{\bar z}_n^{k_{\bar n}}}$
and define $|(k)|=k_0 +k_{\bar 0}+k_1 +...+k_n +k_{\bar 1}...+k_{\bar n}$.
Also, sometimes we shall use double multi-indexes $(k,l)=\{\{k_1 ,k_{\bar 1}\},\{l_1 ,...,l_n ,l_{\bar 1},...,l_{\bar n}\}\}$
and denote by $\pd^{(k,l)}f$ or $f_{(k,l)}$ the derivatives 
$f_{w^{k_1}{\bar w}^{k_{\bar 1}}z_1^{l_1}...z_n^{l_n}{\bar z}_1^{l_{\bar 1}}...{\bar z}_n^{l_{\bar n}}}$.
In this case we define $|(k)|=k_1 +k_{\bar 1}, |(l)|=l_1 +...+l_n +l_{\bar 1}...+l_{\bar n}$.
The main result of [Sh2] is the theorem about almost complex structures on manifolds of universal
coverings:

{\bf Theorem ACS}.
{\it Suppose $X$ is a compact complex manifold of dimension $n+1$ and $\mF$ is a holomorphic foliation of 
$X$ with negative $T_{\mF}$. Suppose that the singular set $\Sigma$ is finite and in some neighborhood
of each singular point the vector field locally defining the foliation is analytically
linearizable and the linear part is diagonalizable. Let $M$ be a manifold of universal
covering with a simply connected base $B$. Then the complex structure
on $M$ can be defined as an almost complex structure on the product $B\times D$ ($D$ is the unit disk) and
this almost complex structure can be defined by forms (0.1), (0,2) of type (1,0).
There $c$ is a smooth vector-function, $\mu$ is a smooth function and we
have the estimates 
$$
   |\mu |\le d<1 ,              
$$
$$
   |\mu_{w^k ,\bar w^l}(w,z)| \le C(1-|w|)^{-(k+l)} ,   \eqno (1.1)
$$   
$$
  |\pd^{(k)}\mu (w,z)|\le C(1-|w|)^{-A|(k)|^4},\,|\pd^{(k)}c_i (w,z)|\le C(1-|w|)^{-A(|(k)|^4 +1)} \eqno (1.2) 
$$
for any pair $k,l$ and multi-index $(k)$. The constant $C$ in these estimates depends on $k+l$ or on $|(k)|$,
the constant $A$ doesn't depend.}

In fact, we don't need in the exact exponent $-|(k)|^4$ in estimate (1.2). It is enough only
to know that this exponent is negative and depends only on $|(k)|$. From the other hand, exact
estimate (1.1) is essential.

Applying Theorem ACS, we can reformulate Theorem 1.

{\bf Theorem 1} (second formulation).  
{\it  Suppose $X$ is a compact complex manifold of dimension $n+1$ and $\mF$ is a holomorphic foliation of 
$X$ with negative  $T_{\mF}$. Suppose that the singular set $\Sigma$ is finite and in some neighborhood
of each singular point the vector field locally defining the foliation is analytically
linearizable and the linear part is diagonalizable. Let $M$ be a manifold of universal
coverings with a simply connected base $B$. Then for any $p\ge 1$ the manifold $M$ is 
diffeomorphic by a $p$-smooth fiberwise diffeomorphism to a domain $\tilde M\subset B\times\bbC$
having continuous boundary and fibered by topological disks $K_z$. The domain $\tilde M$ is
an image of $B\times D$ under the diffeomorphism $f$ satisfying estimates (0.3) - (0.5). As a complex manifold $M$ is
biholomorphic to the manifold $\tilde M$ supplied with an almost complex structure defined by the forms
$$
  dz_i ,\,\,i=1,...,n ,                 
$$  
$$
  dw+<c,d\bar z> ,        
$$ 
where $z_i ,w,c$ have the same sense as in (0.1), (0.2) and the vector-function $c=\{c_1 ,...,c_n$
satisfies the estimates
$$
  |\pd^{(k)}c_i (w,z)|\le C[1-{\rm dist}(w,\pd K_z)]^{-N}
$$
for $|(k)|\le p$ with the constants $C$ and $N\ge 0$ depending only on $p$.} 
 
Theorem 1 in either formulation reduces to the next theorem about the Beltrami equation with
parameters:

\begin{theorem} Suppose $\mu$ is a $p$-smooth function of a variable $z\in D$ and a vector
variable $t=\{t_1 ,...,t_n\}$ belonging to some domain $B\subset\bbC^n$. Let $\mu$ satisfies
the estimates  
$$
  |\mu |\le d<1 ,              \eqno (1.3)
$$
$$
   |\mu_{z^k ,\bar z^l}(z,t)| \le C(1-|z|)^{-(k+l)} ,  \eqno (1.4)
$$   
$$
  |\pd^{(k)}\mu (z,t)|\le C(1-|z|)^{-N}     ,    \eqno (1.5)
$$
for $k+l \le p,|(k)|\le p$ with constants $C$ and $N\ge 0$ depending only on $p$. Then there
exists a solution $f$ to the Beltrami equation
$$
  f_{\bar z}=\mu f_z  
$$  
that is continuous in $C_0 (\bbC )$ as a function of $t$, is $p$-smooth with respect to 
all variables, at every $t$ maps $D$ homeomorphically onto some bounded subdomain of $\bbC$,
and satisfies the estimates
$$
  |f_z (z,t)|\le C(1-|z|)^{-\alpha},\,|f_{\bar z}(z,t)|\le C(1-|z|)^{-\alpha},\,0\le\alpha <1 , \eqno (1.6)
$$
$$
  |f_{z^k \bar z^l}()|/|f_z (z,t)|\le C(1-|z|)^{-(k+l-1)},\,k+l\le p ,  \eqno (1.7)
$$
$$
   |\pd^{(k)}f (z,t)|\le C(1-|z|)^{-N},\, k\le p .   \eqno (1.8)
$$
The constants $C$ and $N\ge 0$ depend only on $p$. 
\end{theorem}

The proof of this theorem starts in Section 2. Now we present some motivations for the
below considerations and outline main steps of the proof.

Remind at first the classical construction of homeomorphic solutions to the Beltramy
equation $f_{\bar z}=\mu f_z$ for a compactly supported $\mu$ (see, for example [Al] or [As]). 
Recall the definition of the classical integral operators acting on functions 
$f\in C_0^{\infty}(\bbC )$: the Cauchy  transform
$$
  \mC f(z)=\frac{1}{\pi}\int\frac{f(\tau )}{z-\tau}dS_{\tau} 
$$
and the Beorling transform
$$
  \mS f(z)=-\frac{1}{\pi}\int\frac{f(\tau )}{(z-\tau )^2}dS_{\tau}  .
$$
Here $dS_{\tau}$ is the usual measure on the $\tau$-plane and the second integral we
understand in terms of its principal value.  
The Cauchy transform is right inverse to the Cauchy-Riemann operator 
$$
  \frac{\pd}{\pd\bar z}\mC f=f 
$$
and
$$
  \mS f= \frac{\pd}{\pd z}\mC f   .
$$  
If $1<q<2<p<\infty$ is a Holder conjugate pair, then the Cauchy transform extends to a bounded
linear mapping from $L^p (\bbC )\cap L^q (\bbC $ into $C_0 (\bbC )$. The Beorling transform
extends to a continuous operator from $L^p (\bbC )$ to $L^p (\bbC )$ for all $1<p<\infty$.
The norm of this operator tends to 1 as $p\to 2$ (the Kalderon-Zygmund inequality).

Suppose $\mu$ has a compact support and $|\mu (z)|\le k<1$. For every $\vpf\in L^p (bbC )$
with a compact support there exists a unique solution $\sigma$ to the inhomogeneous Beltrami equation
$$
  \sigma_{\bar z}=\mu \sigma_z +\vpf    \eqno (1.9)
$$
with derivatives in $L^p (\bbC )$ and decay $f(z)=O(1/z)$ at infinity. We obtain this solution
in the following way. The operator $({\rm Id}-\mu \mS )^{-1}$ defined by the Neumann series
$$ 
  ({\rm Id}-\mu \mS )^{-1}={\rm Id}+\mu \mS +\mu \mS\mu \mS+...   \eqno (1.10)
$$  
is bounded in $L^p (\bbC )$ for $p$ close enough to 2. It is easy to see that
$$
  \sigma=\mC ({\rm Id}-\mu \mS )^{-1}\vpf   \eqno (1.11)
$$  
is a solution to (1.9) and this solution has the required properties. We obtain a solution to the Beltrami 
equation if we put $\vpf =\mu$ in (1.11) and set $f(z)=z+\sigma (z)$. It is an unique
solution to the Beltramy equation with $f_z$ belonging to $L^p_{loc}(\bbC )$ for $2\ge p$ close enough to 2
and normalized by the condition $f(z)=z+O(1/z)$ as $z\to\infty$. Such solution is called
{\it principal solution}. We have
$$
  f(z)=z+\mC (\mu +\mu \mS\mu +\mu \mS\mu \mS\mu +...)(z) . \eqno (1.12)
$$

In fact, the principal solution is a homeomorphism of the complex plane. At first suppose
that $\mu\in C_0^{\infty}(\bbC )$. There is the solution $\sigma\in L^p (\bbC)$ to equation
(1.9) with $\vpf=\mu_z$. We put
$$
  F(z)=z+\mC (\mu e^{\sigma})(z)  .     \eqno (1.13)
$$  
Since $\sigma (z)=O(/z)$ near $\infty$, it follows that $e^{\sigma}-1 \in L^p(\bbC )$ and $\mu e^{\sigma}\in L^p(\bbC )$.
Hence, $F_z =e^{\sigma}$ belongs to $L^p_{loc}(\bbC )$ and $F(z) -z=O(/z)$ near $\infty$. 
$F$ satisfies the Beltrami equation $F_{\bar z}=\mu F_z$ and, by uniqueness of the principal solution, we have 
$F=f$. Further, $F_z =e^{\sigma}$ and $F$ is a local homeomorphism. Since we can extend $f$ to
$\hat\bbC$ by setting $f(\infty )=\infty$, we find that $f$ is a local homeomorphism
$\hat\bbC \to\hat\bbC$ and, hence, $f$ is a global homeomorphism by the monodromy theorem.
For compactly supported measurable $\mu$ we approximate by convolutions $\mu_{\varepsilon}\to\mu$
in $L^p (\bbC )$ for proper $p$ and obtain the principal solution as a limit of smooth
conformal mappings. 

We can remove the restriction that $\mu$ is compactly supported and find a homeomorphism
satisfying the Beltramy equation as a composition of solutions to the equations with the coefficients
having supports in $D$ and in the closer of $\hat\bbC\setminus D$. 

Now suppose ${\rm supp}\mu \subset D$. Applying the extension of $\mu$ by symmetry
we obtain a unique $\mu$-quasiconformal homeomorphism $f:D\to D$ normalized by the 
conditions
$$
   f(0)=0 ,\,\,\,f(1)=1 .
$$
We call this map a {\it normal solution} or a {\it normal mapping}.   

If $\mu$ smoothly depends on a parameter $t$, then the principal solution and the normal
solution aren't necessarily $t$-differentiable. The normal solution has $t$-derivative only when
$\mu$ has uniformly bounded $t$-derivative ([AlB] or [Al]). Indeed, in general case we
can't differentiate, for example, series (1.10). The integrals $\mC\mu_t$ and $\mS\mu_t$ aren't
defined if $\mu_t$ grows sufficiently rapidly near the boundary of $D$.

However, when derivatives of $\mu$ satisfy estimates (1.4), (1.5), we can attempt to find 
$t$-differentiable solutions to the Beltramy equation if we replace the transforms $C$ and $S$
by integral operators with counter-items. Suppose ${\rm supp}f\subset D$. We define
$$
  {\mC}_m f(z)=\frac{1}{\pi}\int_D f(\zeta )\left[
  \frac{1}{z-\zeta}-\frac{1}{z-\bar\zeta^{-1}}-...-\frac{(\zeta-\bar\zeta^{-1})^{m-1}}{(z-\bar\zeta^{-1})^m}\right]dS_{\zeta}=
$$  
$$
  =\frac{1}{\pi}\int_D \frac{f(\zeta)}{z-\zeta}\left(\frac{\zeta-\bar\zeta^{-1}}{z-\bar\zeta^{-1}}\right)^m dS_{\zeta}=
  \frac{1}{\pi}\int_D \frac{f(\zeta)}{z-\zeta}\left(\frac{1-|\zeta |^2}{1-\bar\zeta z}\right)^m dS_{\zeta}.
$$
Here we used the identity
$$
  \left(\frac{\zeta-\bar\zeta^{-1}}{z-\bar\zeta^{-1}}\right)^{k-1}\left(\frac{1}{z-\zeta}-\frac{1}{z-\bar\zeta^{-1}}\right)=
  \frac{1}{z-\zeta}\left(\frac{\zeta-\bar\zeta^{-1}}{z-\bar\zeta^{-1}}\right)^k  .  
$$
Define also
$$
  {\mS}_m f(z)=-\frac{1}{\pi}\int_D f(\zeta )\left[\frac{1}{(z-\zeta )^2}-\frac{1}{(z-\bar\zeta^{-1})^2}-...-\frac{(\zeta-\bar\zeta^{-1})^{m-1}}{(z-\bar\zeta^{-1})^{m+1}}\right]dS_{\zeta}=
$$
$$
  =-\frac{1}{\pi}\int_D\frac{f(\zeta)}{z-\zeta}\left(\frac{1-|\zeta |^2}{1-\bar\zeta z}\right)^m 
 \left[\frac{1}{z-\zeta}-\frac{m\bar\zeta}{1-\bar\zeta z}\right]dS_{\zeta}.  
$$
Again ${\mC}_m$ is right-inverse to the Cauchy-Riemann operator on $D$ and
${\mS}_m f (z)=({\mC}_m f)_z (z)$.
  
\begin{definition} We say that a function $f$ on $D$ belongs to $L^p_s (D)$, $0\ge s<\infty$ if the function $f(z)(1-|z|)^s$
belongs to $L^p (D)$. We denote by $\|f\|_{p,s}$ the $L^p$-norm of the function $f(z)(1-|z|)^s$. A function $f$ belongs to
$C^0_s$ if $f(z)(1-|z|)^s$ is uniformly bounded. We denote by $\|f\|_{0,s}$ the $C^0$-norm of the function
$f(z)(1-|z|)^s$.   
\end{definition}

If $2<p<\infty$, $m\ge s$, then ${\mC}_m$ is a bounded mapping from $L^p_s (D)$ into $C_s^0 (D)$. The transform
${\mS}_m$ acts as continuous operator from $L^p_s (D)$ to $L^p_s (D)$ for all $1<p<\infty$, $m\ge s$. In what
follows we shall prove these estimates in more general setting.

If $\mu$ satisfies estimates (1.4), (1.5) we can seek solutions to the Beltrami equation
analogous to (1.12) or (1.13) replacing the operators $\mC$ and $\mS$ by ${\mC}_m$ and
${\mS}_m$ correspondingly. That is, we can write
$$
  f(z)=z+{\mC}_m (\mu +\mu {\mS}_m \mu +\mu {\mS}_m \mu {\mS}_m \mu +...)(z) .  \eqno (1.14)
$$
For any $l>0$ $t$-derivatives or mixed derivatives up to the order $l$ of the items of
series (1.14) will be defined if $m$ is large enough. But there appear two difficulties.  

First, though the operator ${\mS}_m$ is bounded in $L^p_s (D)$, its norm isn't close to 1
if $m>0$. It implies that estimate (1.3) isn't enough for convergence of series (1.14). 
The constant $d$ in (1.3) must be small enough.   

Second, even if we shall find a locally homeomorphic solution analogously to (1.13), we can't
extend this solution to $\hat\bbC$ and apply the topological argument to prove its global
univalence.

We apply some results of theory of univalent functions to overcome this obstacle. A function
$h$ holomorphic on $D$ and mapping 0 to 0 is univalent if
$$
  |h''(z)/h'(z)|\le (1-|z|)^{-1} .    \eqno (1.15)
$$  
(See [Pom]). If $f$ is a solution to the equation $f_{\bar z}=\mu f_z$, then $g=f_{zz}/f_z$
satisfies the equation
$$
  g_{\bar z}=\mu g_z +\mu_z g_z +\mu_{zz}  .   \eqno (1.16)
$$  
Conversely, if we find a solution to this equation, then we can find a solution to the Beltrami
equation by integration. It appears, we can find a solution to equation (1.16) with the estimate
$|g(z)|\le b$ for any $b>0$ if the constant $C$ in the right side of (1.4) is small enough
for several first $k+l$. Further, for any $\mu$-quasiholomorhic function $f$ on $D$ we have the
decomposition into the product $f=h\circ f_{\mu}$, where $f_{\mu}$ is the normal mapping and $h$ is holomorphic.
It implies that if we have sufficiently good estimates for the derivatives of $f_{\mu}$ and 
for $|g|=|f_{zz}/f_z |$, then $h$ satisfies estimate (1.15). Thus $h$ will be univalent
and $f$ homeomorphic.

Now we outline the main steps of the proof of Theorem 2. First, in Section 2 we obtain estimates for derivatives
of the normal solutions when $\mu$ satisfies estimates (1.3), (1.4). In Section 3 we obtain estimates
for the differences of these derivatives when we have the normal solutions with the complex dilatations $\mu_1$ and $\mu_2$. 
For family $\mu (z,t)$ satisfying
estimates (1.5) we obtain Heolder estimates for differences of $z$-derivatives. In application
to foliations we can consider this result as some generalization of the result of [DNS]
on existence of Heolder estimates for the Poincare metrics on the leaves.

In Section 4, applying the obtained estimates, we approximate the family of normal mappings $f_{\mu}$
with $\mu (z,t)$ satisfying estimates (1.4), (1.5) by a 
finitely smooth family $f_t$ of $\mu_t$-quasiconformal homeomorphisms with $\mu_t$ approximating
$\mu (.,t)$ in terms of $\|.\|_{0,k+l}$-norms. The family $f_t$ maps $B\times D$ onto some
domain $\Omega\subset B\times\bbC$ fibered by topological disks $\Omega_t$. For functions
on $\Omega_t$ we define the spaces $L^p_s$ and $C^0_s$ as in Definition 1 replacing the
difference $1-|z|$ by ${\rm dist}(z,\pd\Omega_t )$.
The mappings $f_{\mu}(.,t)$ decompose into the products $g_t \circ f_t$,
where $g_t$ are $\tilde\mu_t$-quasiconformal mappings with $\tilde\mu_t$ small and having
derivatives up to some finite order small in terms of $\|.\|_{0,k+l}$-norms. As a result,
we reduce Theorem 2 to the analogous theorem with $\mu$ defined
on $\Omega_t$ and small with derivatives in terms of $\|.\|_{0,k+l}$-norms. 

In Section 5 we define integral operators analogous to ${\mC}_m$ and ${\mS}_m$ on the domains
$\Omega_t$. In Sections 6 and 7 we obtain estimates for these operators in appropriate norms. 
 
To obtain univalence we must find a solution $f$ to equation (1.16) on $\Omega_t$ satisfying the uniform
estimate $|f(z,t)|\le c{\rm dist}(z,\pd\Omega_t )^{-1}$ with sufficiently small $c$. 
It appears possible if $\mu_t$ is small in terms of $\|.\|_{0,k+l}$-norms. It is essential
that appearing singular integral operators are bounded in appropriate Holder norms.

\section{Estimates for derivatives of normal mappings}

In the following estimates we shall often use the expression "uniform constant"
in a sense that we shall specify in each case. In what follows $c$ or $C$ often means an indeterminate 
uniform constant. For example, in
inequalities of the type $|(.)|\le C(..)\le C(...)$ in the right side $C$ in two cases 
isn't necessary the same. 

In this section $(k)$ is a multi-index of type $(k)=(k_0 ,k_{\bar 0})$ and $f_{(k)}$ is
the derivative $f_{z^{k_0}\bar z^{k_{\bar 0}}}$.

\begin{lemma} Suppose $f$ is a $\mu$-quasiconformal normal mapping and for $w\in D$ 
we have the estimates
$$
      |\mu (w)|\le d<1 ,
$$	  
$$
    |\mu_{(k)}(w)|\le \frac{b_{|(k)|}}{(1-|w|)^{|(k)|}} .  
$$
Then we have the estimates
$$    
    a\le |f_z (w)|\le A    \eqno (2.1)
$$
with uniform $a,A$, and we can put these constants tending to 1 as $d$ and $b_1$ tend to 0.
$$ 
   |f_{(k)}(w)|\le \frac{B}{(1-|w|)^{|(k)|-1}}     \eqno (2.2)
$$
with $B$ depending only on $d$ and $b_{|(l)|}$, $|(l)|\le |(k)|$ .
If $d, b_1 ,...,b_k$ are small enough, then 
$$
  |f_{(k)}(w)|\le C\frac{d+b_1 +...+b_k}{(1-|w|)^{k-1}}    \eqno (2.3)
$$
for $|(k)|\ge 2$ with some uniform $C$ independent of $d,b_1 ,...,b_k$.
\end{lemma} 

In this section we say that an estimate or a constant is uniform if it depends only on 
the constants $d, b_1 ,...,b_k$ of Lemma 1. 

{\bf Proof of Lemma 1}. 1){\it Reduction to the case $w=0$, $\mu (0)=0$.} 

The map
$$
  \vpf_w (z)=\frac{w-z}{1-z\bar w}
$$
maps $D\to D$ and the point $w$ to zero. Note that
$$
  \vpf_w^{-1}(z)=\frac{w-z}{1-z\bar w}
$$
also.

Define $w'=f(w)$, $f_w =\vpf_{w'}\circ f\circ \vpf_w^{-1}$.
Then $f_w$ maps zero to zero. There is the useful inequality
$$
  c_1 (1-|w|)\le 1-|w'|\le c_2 (1-|w|)     \eqno (2.4)
$$  
for some uniform $c_1$, $c_2$ depending only on $d$.  Indeed, if 
$(1+|\mu |)/(1-|\mu |)\le K$, then, by distortion theorems for quasiconformal mappings, 
(see, for example, [L])
$$
   |w|^K \le |w'|\le |w|^{1/K} 
$$
and we can put in (2.4) $c_1 =1/2K$, $c_2 =2K$. 
Notice that we can put the bounds $c_1 ,c_2$ independent of $d$ for $d\le B$ for any $B<1$.

We need in estimates for derivatives of $\mu_{f_w}$. 

\begin{proposition} We have
$$
  |(\mu_{f_w})_{(k)}(z)|\le\frac{C(d+...+b_{|(k)|}}{(1-|z|)^{|(k)|}(1-|\bar w z|)^{|(k)|}}    \eqno (2.5)
$$  
with $C$ independent of $d, b_1 ,...,b_{|(k)|}$. In particular,
$$
  |(\mu_{f_w})_{|(k)|}|\le C(b+...+b_{|(k)|})   \eqno (2.6)
$$
on every disk $|z|\le a<1$, where $C$ depends on $a$ but not depends on $w$.
\end{proposition}

{\bf Proof}. We have
$$
  \mu_{f_w}(z)=\mu\circ \vpf_w^{-1}(z)\frac{(\vpf_w )_z}{\overline{(\vpf_w )_z}}=\mu\circ \vpf_w^{-1}(z)\left(\frac{1-\bar z w}{1-z\bar w}\right)^2 . \eqno (2.7)
$$   
Denote $\mu_0 = \mu_{f_w}(0)=\mu (w)$. We have 
$$
  (\mu_{f_w})_z =\mu_z \circ\vpf_w^{-1}(z)\left(\frac{1-w\bar z}{1-\bar w z}\right)^2 \frac{1-|w|^2}{(1-\bar w z)^2} +
   2\mu\circ\vpf_w^{-1}(z)\frac{\bar w(1-w\bar z)^2}{(1-\bar w z)^3} .   
$$
Since 
$$
  |\mu_z \circ\vpf_w^{-1}(z)|\le b_1 (1-|\vpf_w^{-1}(z)|)^{-1}=b_1 \frac{|1-\bar w z|}{||1-\bar w z|-|w-z||},  \eqno (2.8)
$$
we obtain
$$
  |(\mu_{f_w})_z |\le (d+b_1 )\left(\frac{1-|w|^2}{||1-\bar w z|-|w-z|||1-\bar w z|}+\frac{|w|}{|1-\bar w z|}\right)  
$$ 
  
We need in an estimate of the difference  $||1-\bar w z|-|w-z||$ from above. Note that
$$
  ||1-\bar w z|-|w-z||\ge ||1-\bar w z|^2 -|w-z|^2|/2 .
$$  
Let $\theta$ be the angle between $z$ and $w$, $\rho =|w|$, $r=|z|$. We have
$$
  |1-\bar w z|^2 =(1-\rho r\cos\theta )^2 +\rho^2 r^2 \sin^2 \theta =1-2\rho r\cos\theta +\rho^2 r^2 ,
$$
$$
  |w-z|^2 =(\rho -r\cos\theta )^2 +r^2 \sin^2 \theta =\rho^2 -2\rho r\cos\theta +r^2 .
$$
Hence,
$$
  |1-\bar w z|^2 -|w-z|^2 =(1-\rho^2 )(1-r^2 )      
$$ 
and
$$
  ||1-\bar w z|-|w-z||\ge (1-\rho^2 )(1-r^2 )/2     \eqno (2.9)
$$
From this estimate and (2.8) we obtain
$$
  |(\mu_{f_w})_z | \le (d+b_1 )\left(\frac{1-|w|^2}{(1-|w|)(1-|z|)|1-\bar wz|}+\frac{|w|}{|1-\bar w|z|}\right)\le
  \frac{2(d+b_1 )}{(1-|z|)1-\bar w z)} . 
$$  
We get an analogous estimate for $(\mu_{f_w})_{\bar z}$. We obtained estimate (2.5) for the
derivatives of first order.

Now we shall get estimates for derivatives of higher orders. 
We have
$$
  (\mu_{f_w})_{(k)}(z)=\mu_{(k)}\circ\vpf_w^{-1}(z))\left(\frac{1-w\bar z}{1-\bar w z}\right)^2 \left(\frac{1-|w|^2}{(1-\bar w z )^2}\right)^{|(k)|}+
$$
$$
  +\sum\mu_{(l)}\circ\vpf_w^{-1}(z))(\vpf_w^{-1})_{z^{l_1}}(z)...(\vpf_w^{-1})_{z^{l_p}}(z)
  \overline{(\vpf_w^{-1})_{z^{\bar l_1}}(z)}...\overline{(\vpf_w^{-1})_{z^{\bar l_q}}(z)}\left[\left(\frac{1-w\bar z}{1-\bar w z}\right)^2 \right]_{(m)} . 
$$
Here in the second line we have the sum of items with $(|(l)|<|(k)|$, $p=l_0$, $q=l_{\bar 0}$,  
$l_1 +...+l_p +\bar l_1 +...\bar l_q +|(m)|=|(k)|$. We estimate these items by induction.
Denote by $O(1)$ terms having uniform estimates. 

The multiple $\mu_{(l)}\circ\vpf_w^{-1}(z))$ has the estimate (we apply (2.8) and (2.9))
$$
  b_{|(l)|}(1-|\vpf_w^{-1}(z)|)^{-|(l)|}\le cb_{|(l)|}\frac{|1-\bar w z|^{|(l)|}}{(1- |z|^2 )^{|(l)|}(1-|w|^2 )^{|(l)|}} .
$$  
with $c$ depending only on $|(l)|$. After differentiation we obtain the multiple with the estimate
$$
  cb_{|(l)|+1}\frac{|1-\bar w z|^{|(l)|+1}}{(1- |z|^2 )^{|(l)|+1}(1-|w|^2 )^{|(l)|+1}}\frac{1-|w|^2}{|1-\bar w z|^2}=
$$
$$
  =cb_{|(l)|+1}\frac{|1-\bar w z|^{|(l)|}}{(1- |z|^2 )^{|(l)|}(1-|w|^2 )^{|(l)|}}\frac{1}{(1- |z|^2 )|1-\bar w z|} 
$$
Any multiple of type $(\vpf_w^{-1})_{z^l}(z)$ has the estimate
$$
  \frac{1-|w|^2}{|1-\bar w z|^{l+1}}O(1)
$$
and the multiple $[(1-w\bar z )^2 /(1-\bar w z )^2 ]_{(m)}$ is of the type
$$
  |1-\bar w z |^{-|(m)|}O(1) .
$$  
Thus at every differentiation there appears either the multiple $|1-\bar w z |^{-1}$ either
the multiple $|1-\bar w z |^{-1}(1-|z|^2 )^{-1}$, and in the estimates we must replace $b_{|(l)|}\leadsto b_{|(l)|+1}$.  
It finishes the proof of estimate (2.5).  
$\Box$
 
Let $g_w$ be the map $z\mapsto z+\mu_0 \bar z$. Then 
$$
  g_w^{-1}(z)=\frac{z-\mu_0 \bar z}{1-|\mu_0 |^2}
$$
The map $g_w$ maps $D$ onto some ellipsis. Let $H_w$ be the conformal map mapping this ellipsis onto $D$ and having
real derivative at zero. Let $Z_w$ be the composition $Z_w =H_w \circ g_w$. 
We have
$$
  \mu_{Z_w}(z)=\mu_{g_w}(z)=\mu_0 .
$$    
Define the map
$$
  h_w =f_w \circ Z_w^{-1}  .
$$
The complex dilatation $\mu_{h_w}$ is
$$
  \mu_{h_w} =\frac{\mu_{f_w}-\mu_0}{1-\bar\mu_0 \mu_{f_w}}\frac{Z_{w,z}}{\overline{Z_{w,z}}}\circ Z_w^{-1} . \eqno (2.10)
$$
In particular, $\mu_{h_w}(0)=0$. The map $Z_w$ and its inverse $Z_w^{-1}$ have all derivatives bounded uniformly 
with respect to $w$. It follows that for derivatives of $\mu_{h_w}$ we have 
estimates analogous to the estimates of Proposition 1.

In what follows we adopt the notation $Z=Z_w^{-1}(z)$.
\begin{proposition}.  We have the estimates  
$$
  |(\mu_{h_w})_{(k)}(z)|\le \frac{C(d+...+b_{|(k)|})}{(1-|Z|)^{|(k)|}(1-|\bar w Z|)^{|(k)|}}  .  \eqno (2.11) 
$$
Here $C$ doesn't depend on $d,...b_{|(k)|}$ if $d\le B$ for any $B<1$, for example, if $d \le 1/2$. 
In particular,
$$
  |(\mu_{h_w})_{(k)}|\le C(d+...+b_{|(k)|}) .  \eqno (2.12) 
$$
on any disk $|z|\le a<1$. 
\end{proposition}

{\bf Proof} Almost all assertions were already proved. The assertion about constant $C$ in (2.11)
holds because $Z_w$ is quasiconformal with the complex dilatation $\mu_0$, $|(Z_w )_z|$ 
is bounded from below and from above by constants depending only
on $d$, and $(Z_w )_z$ tends to 1 uniformly as $d$ tends to zero, all other derivatives of
$Z_w$ are also bounded by constants depending only on $d$, and we can put these constants 
arbitrary small as $d$ tends to zero.
$\Box$
  
Now we return to the original map $f$. 
\begin{proposition} a) We have the estimates 
$$
  a|(h_w )_z (0)|\le |f_z (w)|\le A|(h_w )_z (0)|   \eqno (2.13)
$$
with some uniform $a, A$. These constants depend on $d$ but we can put them independent of $d$ for $d\le B$ for any $B<1$.
$$
   |f_{(k)}(w)| \le C\frac{1-|w'|}{(1-|w|)^{|(k)}}\sum_{s_1 |(l_1)|+...s_j |(l_j)|\le |(k)|}|((h_w )_{(l_1)}(0))^{s_1}...((h_w )_{(l_j)}(0))^{s_j}|    \eqno (2.14)
$$
with some uniform $C$.
\end{proposition}

{\bf Proof}. We have
$$
  f=\vpf_{w'}^{-1}\circ h_w  \circ  Z_w \circ\vpf_w .    \eqno (2.15)
$$
We adopt the notations: $\pd$ is the derivative of a function with respect to its analytic argument and $\bar\pd$
is the derivative with respect to the conjugate variable. 
We have
$$
   f_z =\pd (\vpf_{w'}^{-1})\circ h_w \circ  Z_w \circ \vpf_w [\pd (h_w )\circ  Z_w \circ\vpf_w \cdot (\pd Z_w )\circ\vpf_w \cdot\pd\vpf_w +  
$$
$$    
   +\bar\pd (h_w )\circ  Z_w \circ\vpf_w \cdot\overline{\pd Z_w}\circ\vpf_w \cdot\pd\vpf_w ] . \eqno (2.16)
$$
In particular, recalling that $\bar\pd h_w (0)=0$, we obtain
$$
  f_z (w)=\frac{1-|w'|^2}{1-|w|^2}\pd h_w (0)\pd Z_w (0).  
$$ 
The fraction $\frac{1-|w'|^2}{1-|w|^2}$ is uniformly bounded from below and from above by (2.4)
and we can set the bounds independent of $d$ for $d\le B$ for any $B<1$. 
The same holds for the value $\pd Z_w (0)$. We obtain (2.13).  

Now differentiating (2.16) we see that
the derivative $\frac{\pd^2 f}{\pd z^2}(w)$ is the sum of terms of the types
$$
  \pd\pd\vpf_{w'}^{-1}(0)[\pd h_w (0)\pd Z_w (0)\pd\vpf_w (w)]^2=[\pd h_w (0)]^2 O((1-|w|^2 )^{-1}) , 
$$
$$
  \pd\vpf_{w'}^{-1}(0) \pd\pd h_w (0)[\pd Z_w (0)\pd\vpf_w (w)]^2 =\pd\pd h_w (0)O((1-|w|^2 )^{-1}) , 
$$
$$
  \pd\vpf_{w'}^{-1}(0)\pd h_w (0)\pd\pd Z_w (0)[\pd\vpf_w (w)]^2 =\pd h_w (0)O((1-|w|^2 )^{-1})  
$$ 
$$
  \pd\vpf_{w'}^{-1}(0)\pd h_w (0)\pd Z_w (0)\pd\pd\vpf_w (w) =\pd h_w (0)O((1-|w|^2 )^{-1})   
$$  
and of analogous terms containing $\bar\pd$ derivatives of $h_w$ and $Z_w^{-1}$. All these terms have
estimates $O((1-|w|^2 )^{-1})$.

At further differentiations we obtain each time either the multiple $\pd\vpf_w (w)=1/(1-|w|^2 )$ either the term, where the
multiple  $\pd^k \vpf_w (w)$ is replaced by the multiple $\pd^{k+1}\vpf_w (w)$, which also results in multiplication by
$(1-|w|^2 )^{-1}$. Also, if we had some item containing a product $((h_w )_{(l_1)}(0))^{s_1}...((h_w )_{(l_j)}(0))^{s_j}$, 
then, after differentiation, we obtain the term of the same type with the sum $s_1 |(l_1)|+...s_j |(l_j)|$ 
higher no more, than by 1. By induction we see that the derivative $f_{(k)}(w)$ at $|(k)|\ge 1$
can be represented as a sum of items of the types 
$$
  \frac{1-|w'|^2}{(1-|w|^2 )^m}((h_w )_{(k_1 )}(0))^{s_1}...((h_w )_{(k_j )}(0))^{s_j}
 ((Z_w )_{(q_1)}(0))^{r_1}...((Z_w )_{(q_i)}(0))^{r_i}O(1)  ,
$$
where  $m\le |(k)|$,$s_1 |(k_1)|+...s_j |(k_j)| +r_1 |(q_1)|+...r_i |(q_i)|\le|(k)|$
and $O(1)$ is a multiple uniformly bounded and independent of $t$ in conditions of point b).
Estimate (2.14) follows immediately.$\Box$

We see that to prove Lemma 1 it is enough to show that $(h_w )_z (0)$ is uniformly bounded 
from below and from above and that derivatives $(h_w )_{(k)}(0)$, $|(k)|\ge 2$ are uniformly 
bounded. 

2){\it Transition to the logarithmic chart}

Define the logarithmic coordinates
$$
  \zeta=\log z=\xi +i\vpf ,\,\omega=\log h_w (z)=\eta +i\theta .
$$
We shall use the notation $F_w$ for $h_w$ represented in the logarithmic coordinates
$\omega =F_w (\zeta )$. The complex dilatation $\mu_{F_w}$ is
$$
  \mu_{F_w}(\zeta )=\frac{e^{\bar\zeta}}{e^{\zeta}}\mu_{h_w}(e^{\zeta})  .  \eqno (2.17)
$$

Consider the first terms of the power series decomposition for $h_w$
$$
  h_w (z)=a_0 z + b_{20} z^2 +b_{1\bar 1}z\bar z+b_{0\bar 2}\bar z^2 +...  \eqno (2.18)
$$
In the logarithmic coordinates we have the decomposition at $\xi =-\infty$ (we assume $z=re^{i\vpf}$, $\xi=\log r$)    
$$
  \omega =F_w (\zeta )=\eta +i\theta =\log [a_0 re^{i\vpf}+r^2 (b_{20}e^{2i\vpf}+b_{1\bar 1}+b_{0\bar 2}e^{-2i\vpf})+...]=
$$
$$
  =\xi +i\vpf +\log a_0 +\frac{e^{\xi}}{a_0} (b_{20} e^{i\vpf}+b_{1\bar 1} e^{-i\vpf}+b_{0\bar 2}e^{-3i\vpf})+O(e^{2\xi}) .
 \eqno (2.19)
$$
Consider the higher terms. If the term of order $k$ in the decomposition of $h_w$ is
$$
  b_{k0} z^k +b_{(k-1)\bar 1} z^{k-1}\bar z+...+b_{0\bar k}\bar z^k  ,
$$
then the term of order $e^{(k-1)\xi}$ in the decomposition of $F_w$ is
$$  
  a_0^{-1}e^{(k-1)\xi}(b_{k0}e^{i(k-1)\vpf}+b_{(k-1)\bar 1} e^{i(k-3)\vpf}+b_{0\bar k}e^{-(k+1)i\vpf})+S) , \eqno (2.20)
$$
where $S$ is a sum of products of multiples containing only the coefficients of the decomposition $h_w$ of degrees
less than $k$. Thus uniform estimates for coefficients of decomposition (2.18) follow from uniform estimates
for coefficients of decomposition (2.20). Note also that we can write expression (2.20) as
$$
  a_0^{-1}e^{(k-1)\xi}(c_{k0}e^{i(k-1)\vpf}+c_{(k-1)\bar 1} e^{i(k-3)\vpf}+c_{0\bar k}e^{-(k+1)i\vpf}) . \eqno (2.21)
$$
Indeed, each term $z^{k-1-p}\bar z^p$ in the decomposition of $\log (h_w (z)/z)$ yields the term $e^{(k-1)\xi}e^{i(k-1-2p)\vpf}$.
We see also that if we obtain uniform estimates for the coefficients of decomposition (2.21), 
then we also obtain uniform estimates for the coefficients $b_{k0},b_{(k-1)\bar 1},...$.
 
3) {\it Integral operators in the logarithmic chart}.

Let $H$ be the stripe $-\pi \le \vpf \le \pi$ and $L_H^p$ be a space of $\vpf$-periodic functions on $H$ with usual
$L^p$-norm. We adopt the notations $L^p_{H_-}$ for the space of functions belonging to $L_H^p$ with
a support in the left half-stripe $H\cap\{\xi\le 0\}$ and $L^p_{H_a}$ for the space of functions
supported in the set $H\cap\{\xi\le a\}$.

We define an integral transform $P_H$ acting on $\vpf$-periodic functions on $H$
$$
  P_H f(\zeta)=\frac{1}{2\pi i}\int_H f(\tau )\left(\frac{e^{\tau}}{e^{\tau}-e^{\zeta}}-1\right)d\tau d\bar\tau=
   \frac{e^{\zeta}}{2\pi i}\int_H \frac{f(\tau )}{e^{\tau}-e^{\zeta}}d\tau d\bar\tau .
$$ 

\begin{proposition} The transform $P_H$ is right inverse to the operator $\pd/\pd\bar\zeta$ 
for smooth functions belonging to $L^p_{H_a}$ at $a<\infty$. Suppose $e^{-\xi}f\in L^p_{H_a}$, $2< p\le\infty$.
Then
$$
  \|P_H f\|_C \le C_{ap}\|e^{-\xi}f\|_p    \eqno (2.22)
$$
with $C_{ap}$ depending on $a$ and $p$.
\end{proposition} 

{\bf Proof.} At first we prove that the operator $P_H$ is right inverse to the Cauchy-Riemann operator. We have 
$$
  P_H f(\zeta)=\frac{1}{2\pi i}\int_H f(\tau )\left(\frac{e^{\tau}}{e^{\tau}-e^{\zeta}}-1\right)d\tau  d\bar\tau=
$$
$$
=\frac{1}{2\pi i}\int_H \frac{f(\tau )e^{-\bar\tau}}{e^{\tau}-e^{\zeta}}d(e^{\tau}) d(e^{\bar\tau})+h(\zeta )=
  \frac{1}{2\pi i}\int_D \frac{f(\log t)}{\bar t(t-e^{\zeta})}dtd\bar t +h(\zeta ) .
$$
where $h$ is a holomorphic function. Hence, we have
$$
  \frac{\pd P_H f}{\pd (e^{\bar\zeta})}(\zeta)=\frac{f(\zeta )}{e^{\bar\zeta}} ,
$$
and
$$
  \frac{\pd P_H f}{\pd \bar\zeta}(\zeta)=f(\zeta ).
$$ 

To proof estimate (2.22) it is enough to show that 
$e^{\zeta}\frac{\theta_a (\tau )e^{\tau}}{e^{\tau}-e^{\zeta}}$ considered as a function of $\tau$
has a uniform estimate with respect to $\zeta$ in $L_H^q$ for $1\le q<2$. Here $\theta_a (\tau)$ is the "step": $\theta (\tau) =1$
if ${\rm Re}\tau \le a$ and $\theta (\tau) =0$ if ${\rm Re}\tau >a$. 

Suppose at first $\xi ={\rm Re}\zeta \le 0$. The integral
$$
  \int_{|\tau -\zeta |\le 1}\frac{|e^{q\tau}|}{|e^{\tau}-e^{\zeta}|^q}d\tau d\bar\tau =\int_{|\tau -\zeta |\le 1}\frac{d\tau d\bar\tau }{|1-e^{\zeta -\tau}|^q}
$$
is uniformly bounded  and the integral
$$
  \int_{\{|\tau -\zeta |>1,{\rm Re}\tau \le a\}}\frac{|e^{q\tau}|}{|e^{\tau}-e^{\zeta}|^q}d\tau d\bar\tau 
$$
has the estimate $c_a (e^{-q\xi})$ with some $c_a$ depending only on $a$. Hence we obtain an uniform estimate for the kernel.

If $\xi >0$, then the integral over the domain $|\tau -\zeta |\le 1$ is obviously uniformly bounded (we should consider only the integral
over the intersection of this disk with the semi-plane $\{\xi\le a\}$). 
From the other hand, the integral
$$
  \int_{H\cap \{|\tau -\zeta |>1,{\rm Re}\tau \le a\}}\frac{|e^{q(\zeta +\tau )}|}{|e^{\tau}-e^{\zeta}|^q}d\tau d\bar\tau 
$$
at $\xi \ge 1$ also has an uniform estimate because the function $\frac{|e^{\zeta}|}{|e^{\tau}-e^{\zeta}|}$ is
uniformly bounded.
$\Box$

Now we define
$$
  T_H f=\frac{\pd P_H f}{\pd f}=\frac{1}{2\pi i}e^{\zeta}\int_H \frac{f(\tau )e^{\tau}}{(e^{\tau}-e^{\zeta})^2}d\tau d\bar\tau  ,
$$
where the integral is defined in terms of its principal value. 

\begin{proposition} The operator $T_H$ is bounded in $L_H^p$ at $1<p<\infty$ and its $L_p$-norm tends to 1 when $p$ tends to 2.
\end{proposition}

{\bf Proof.} If $f$ is a $\vpf$-periodic function with support in some vertical stripe, then $T_H f$ at large positive or negative $\xi$
decreases as $e^{-|\xi |}$. In the following standard calculations all possible
integrals over the boundary equal zero because $P_H$ transforms periodic functions
into periodic ones.
$$
  \int_H |T_H f|^2 dS_{zeta}=-\frac{1}{2i}\int_H (P_H f)_{\zeta}\overline{(P_H f)}_{\bar\zeta} dS_{\zeta}=
 \frac{1}{2i}\int_H P_H f \overline{(P_H f)}_{\bar\zeta \zeta} dS_{\zeta}=
$$
$$
  =\frac{1}{2i}\int_H (P_H f)\bar f_{\bar\zeta} dS_{\zeta}=-\frac{1}{2i}\int_H \bar f(P_H f)_{\bar\zeta} dS_{\zeta}=
\int_H |f|^2 dS_{\zeta}
$$
Now we prove the $L_H^p$- boundedness. We have
$$
  T_H f (\zeta )=\frac{1}{2\pi i}\int_H \frac{f(\tau )e^{\tau-\zeta}}{(e^{\tau-\zeta}-1)^2}d\tau d\bar\tau=
  \frac{1}{2\pi i}\int_H \frac{f(\tau )d\tau d\bar\tau}{(\tau-\zeta )^2}+
$$
$$
  +\frac{1}{2\pi i}\int_{|\tau -\zeta |\le 1}f(\tau )\left[\frac{e^{\tau-\zeta}}{(e^{\tau-\zeta}-1)^2}-\frac{1}{(\tau -\zeta )^2}\right]d\tau d\bar\tau +
$$
$$
  +\frac{1}{2\pi i}\int_{H\setminus\{|\tau -\zeta |\le 1\}}f(\tau )\left[\frac{e^{\tau-\zeta}}{(e^{\tau-\zeta}-1)^2}-\frac{1}{(\tau -\zeta )^2}\right]d\tau d\bar\tau 
\eqno (2.23)
$$ 
The first integral is the usual Beurling transform of a function $f$ with support in $H$ and has an estimate 
in $L^p (\bbC)$ by the Calderon-Zigmund theorem. The second integral is a convolution of the function $f$ with the function
$$
  \left(\frac{e^t}{(e^t -1)^2}-\frac{1}{t^2}\right)\chi_D (t) ,
$$
where $\chi_D$ is the characteristic function of the unit disk. This function belongs to $L^1$ and the integral has
an $L^p (\bbC)$-estimate by the Yung inequality. The $L_H^p$-norms of these integrals are no greater then its 
$L^p (\bbC)$-norms. The last integral can be considered as a bounded operator on $L_H^1$ and on $L_H^{\infty}$.
Indeed, if we denote by $K(\tau ,\zeta )$ the kernel (i.e., the function in the square brackets), than we can easy see that
$$
  \int_{H\setminus\{|\tau -\zeta \le 1\}}K(\tau ,\zeta )dS_{\tau}
$$
has an uniform estimate as a function of $\zeta$ and, by a symmetry, the integral
$$
  \int_{H\setminus\{|\tau -\zeta \le 1\}}K(\tau ,\zeta )dS_{\zeta}
$$
also has an uniform estimate as a function of $\tau$. We obtain an $L^1$-estimate of the last integral in (2.23)
$$
  \int_H \left| \int_{H\setminus\{|\tau -\zeta \le 1\}}f(\tau )K(\tau ,\zeta )dS_{\tau}\right|dS_{\zeta} \le
$$
$$
  \le\int_H |f(\tau )| \sup_{\tau} \int_{H\setminus\{|\tau -\zeta \le 1\}}|K(\tau ,\zeta )|dS_{\zeta}dS_{\tau}\le
C\|f\|_{L_H^1} .
$$
for some $C$. Analogously,
$$
  \left| \int_{H\setminus\{|\tau -\zeta |\le 1\}}f(\tau )K(\tau ,\zeta )dS_{\tau}\right|\le 
 \sup_{\tau}|f(\tau )| \sup_{\zeta}\int_{H\setminus\{|\tau -\zeta \le 1\}}|K(\tau ,\zeta ))dS_{\tau}\le C|f\|_{L_H^{\infty}}.
$$  
We obtain $L_H^p$-boundedness of the last integral in (2.23) for $1\le p< \infty$ by the Riesz-Thorin interpolation theorem.

The $L_H^p$-norm of the operator $T_H$ tends to 1 as $p$ tends to 2 also by the Riesz-Thorin theorem.
$\Box$

4) {\it Solution to the Beltramy equation in the logarithmic chart}.
 
Suppose $\nu (\zeta )$ is some $\vpf$-periodic function with the properties:

a) $|\nu (\zeta)|\le d<1$ for some $d$, 

b) $\nu$ has the support in some domain $\xi \le a$, 

c) $|e^{-\xi}\nu (\xi )|\le c$ for some $c<\infty$.
 
Since the function $e^{\xi}$ is bounded in $L^p_{H_a}$, the series 
$$
  h=(Id - T_H \nu)^{-1}(1) =1 +T_H \nu +T_H \nu T_H \nu +...  ,
$$
converge in $L_H^p$ for $p>2$ sufficiently close to 2, as in the classical case, and we have for the $L^{(p)}_H$-norm
of the function $h-1$ the estimate $C_a c/(1-d)$ with some $C_a$ depending only on $a$. 
By Proposition 4, $|P_H (\nu h)|\le C_a c^2 /(1-d)$ with, possibly, some new $C_a$. 
We see that
$$
  f_{\nu}(\zeta )=\zeta +P_H [\nu (Id -T_H \nu)^{-1}(1)](\zeta )
$$
is a solution to the Beltrami equation with the Beltrami coefficient $\nu$.
We call this solution {\it the principal logarithmic solution}.
 
\begin{proposition} The principal logarithmic solution is a homeomorphism of the plane 
satisfying the estimate
$$
  |f_{\nu}(\zeta )-\zeta |\le \frac{C_a c^2}{1-d} .   \eqno (2.24)
$$  
\end{proposition}

{\bf Proof} Estimate (2.24) follows from Propositions 4 and 5. Prove that $f_{\nu}$
is a homeomorphism.

The map $\tilde f_{\nu}=\exp\circ f_{\nu}\circ\log$ is a quasiholomorphic function on 
the punctured plane with the compactly supported complex dilatation. From (2.24) follows
that the estimate 
$$
  c|z|\le |\tilde f_{\nu}(z)|\le C|z|,\, 0<c,C<\infty     \eqno (2.25)
$$ 
holds at zero and at infinity. From the other hand, $\tilde f_{\nu}$ extends as a quasiholomorphic function to zero.
Indeed, after the change of the chart by a quasiconformal homeomorphism we obtain a holomorphic
function with a removable singularity. Further, since $\tilde f_{\nu}$ is holomorphic
outside some disk, we obtain from (2.25) that at infinity this map has the asymptotics 
$cz+b+O(1/z)$, $c\ne 0$. Subtracting $b$ and dividing by $c$ we obtain a quasiholomorphic map
with a compactly supported complex dilatation and with the asymptotics $z+O(1/z)$ at infinity.
It implies that it is the principal solution to the corresponding Beltrami equation, which is
unique and homeomorphic. It follows that $f_{\nu}$ also is a homeomorphism. $\Box$

\begin{proposition}
The coefficient $\mu_{F_w}$ satisfy conditions a) - c) for $a=0$ with estimates uniform with respect to $w$. 
The operator $Id - T_H \mu_{F_w}$ is invertible in $L^p_H$ if $2\le p\le P_{d+b_1}$ for some $P_{d+b_1}$,
which can be made arbitrary large if $d,b_1$ are small enough. Here $d$ and $b_1$ are the constants from
the formulation of Lemma 1.
Moreover, we have the estimate
$$
  |f_{\mu_{F_w}}(\zeta )-\zeta |\le C\frac{(d+b_1 )^2}{1-d} , \eqno (2.26) 
$$
where $C$ doesn't depend $d,b_1$.
\end{proposition}

{\bf Proof}. If $\xi \le -\log 2$ we have, applying (2.17) and estimate (2.12)
$$
  |\mu_{F_w}(\zeta )|=|\mu_{h_w}(e^{\zeta})\frac{e^{\bar\zeta}}{e^{\zeta}}|\le
$$
$$
  \le e^{\xi}\sup_{D_{1/2}}(|(\mu_{h_w})_z |+|(\mu_{h_w})_{\bar z}|) \le C(d+b_1 )e^{\xi}  .  
$$  
Also, 
$$
  |\mu_{F_w}(\zeta )|\le Cd
$$
with some uniform $C$ if $\xi > -\log 2$. We see that $\|\mu_{F_w}\|_p \le C(d+b_1 )$ and 
$$
\|(Id - T_H \mu_{F_w})^{-1}(1)-1\|_p \le C\frac{d+b_1}{1-d} .
$$
Analogously, we can see that the function $e^{\zeta}\frac{\theta (\tau )|\mu_{F_w}(\tau )|}{e^{\tau}-e^{\zeta}}$
has $L_H^q$-norm no greater, than $C(d+b_1 )$. By definition of $f_{\mu_{F_w}}$ and applying
estimate (2.24), we obtain (2.26).  
$\Box$

5) {\it The proof of estimate (2.1)}
We shall show that the map $F_w$ has a representation in terms of principal logarithmic
solutions. 

Define the map  
$$
  \tilde f_{\mu_{F_w}}(\zeta )=-\overline{f_{\mu_{F_w}}(-\bar\zeta )}  .
$$
This map is "symmetrical" to  $f_{\mu_{F_w}}$ with respect to the imaginary axis. Its complex dilatation is
$\tilde\mu_{F_w}(\zeta )=\overline{\mu_{F_w}(-\bar\zeta )}$. Now we define the function
$$
  \lambda_{F_w}=\mu_{F_w}\frac{\tilde f_{\mu_{F_w},\zeta}}{\overline{\tilde f_{\mu_{F_w},\zeta}}}\circ\tilde f^{-1}_{\mu_{F_w}}= 
 \mu_{F_w}\circ \tilde f^{-1}_{\mu_{F_w}}\frac{\overline{f_{\mu_{F_w},\zeta}}}{f_{\mu_{F_w},\zeta}}\circ (-\overline{\tilde f^{-1}_{\mu_{F_w}}}) .  
 \eqno (2.27))
$$    

\begin{proposition} 
The coefficient $\lambda_{F_w}$ satisfies conditions a) - c) with uniform constants. 
More, we have the estimate
$$
  |f_{\lambda_{F_w}}(\zeta )-\zeta |\le C\frac{(d+b_1 )^2}{1-d} ,
$$
where $d,b_1$ are from the formulation of Lemma 1 and $C$ doesn't depend on $d,b_1$ if $d\le 1/2$.
\end{proposition}

{\bf Proof}. Obviously $|\lambda_{F_w}|=|\mu_{F_w}\circ \tilde f^{-1}_{\mu_{F_w}}|$. We see that 
$\lambda_{F_w}$ has a support in some half-plane $\xi\le a$ with an uniform $a$
and, since for 
${\rm id}-\tilde f^{-1}_{\mu_{F_w}}$ we have the estimate of Proposition 6, we have
the estimate $|\lambda_{F_w}(\zeta )|\le C(b+b_1 )e^{\xi}$ at $\xi \le -\log 2$ with some uniform $C$. 
We finish the proof analogously to the proof of Proposition 7.
$\Box$

Consider the map
$$
  \tilde F_w =f_{\lambda_{F_w}}\circ \tilde f_{\mu_{F_w}}
$$
In the chart $z$ it is the map
$$
  \tilde h_w =\exp \circ\tilde F_w \circ\log  .
$$

\begin{proposition} 
$$
  h_w (z)=\frac{\vpf_{w'}(1)}{\tilde h_w \circ Z_w^{-1}\circ\vpf_w (1)}\tilde h_w (z) ,   \eqno (2.28)
$$
for $z\in D$.
\end{proposition}
   
{\bf Proof}. 
$\tilde F_w$ is a quasiconformal map with the complex dilatation $\mu_{F_w}+\tilde\mu_{F_w}$. 
The map $\tilde h_w$, obtained after transition to the chart $z=e^{\zeta}$, 
is a solution to the Beltrami equation on the punctured plane with the Beltrami
coefficient symmetrical with respect to the unit circle and equal to $\mu_{h_w}$ if $|z|\le 1$. 
Also, this solution satisfies the uniform estimate $b|z|\le |\tilde h_w (z) |\le B|z|$, 
where for $|\log |b||$ and $|\log |B||$ we have the estimates
$$
  C\frac{(d+b_1 )^2}{1-d}   \eqno (2.29)
$$
according to Propositions 6,7 and 8.
As in the proof of Proposition 6 we can extend $\tilde h_w$ to zero and obtain the quasiconformal
map fixing zero and infinity. Suppose $\tilde h_w (1)=a_w$. Then $|a_w |$ has uniform estimates
from below and from above. Dividing by $a_w$ we obtain a map fixing 0, $\infty$ and 1,
i.e., the normal map. This map is unique and symmetrical with respect to the unit circle.
It means that on $D$ it coincides with $h_w$ up to, possibly, some rotation ($h_w$ doesn't map necessarily
1 to 1).

We see that the map $h_w$ differs from the map $\tilde h_w$ at $|z|\le 1$ only by a constant multiple.
We can find this multiple from the condition 
$$
  1=f(1)=\vpf_{w'}^{-1}\circ h_w \circ Z_w^{-1}\circ\vpf_w (1)  ,
$$
i.e., $h_w \circ Z_w^{-1}\circ\vpf_w (1)=\vpf_{w'}(1)$. We obtain
$$
  h_w (z)=\frac{h_w \circ Z_w^{-1}\circ\vpf_w (1)}{\tilde h_w \circ Z_w^{-1}\circ\vpf_w (1)}\tilde h_w (z)=\frac{\vpf_{w'}(1)}{\tilde h_w \circ Z_w^{-1}\circ\vpf_w (1)}\tilde h_w (z) .
$$  
$\Box$

{\bf Corollary}. {\it If $d$ and $b_1$ are small enough, then 
$$
  |h_w (z) -(h_w )_z (0)z|\le C(d+b_1 )^2
$$
with some uniform $C$.}

{\bf Proof}. It follows from Propositions 7 and 8 and (2.28). $\Box$
    
{\bf Proof of estimate (2.1)}.  Since $\vpf_{w'}(1)$ and $\tilde h_w \circ Z_w^{-1}\circ\vpf_w (1)$
equal to 1 by modulus, and for $\tilde h_w$ we have estimate $b|z|\le |\tilde h_w (z) |\le B|z|$, 
where $|\log |b||$ and $|\log |B||$ satisfy estimate (2.29), we obtain the estimate
$$
  \exp\frac{c(d+b_1 )^2}{1-d}\le |(h_w )_z (0)|\le \exp\frac{C(d+b_1 )^2}{1-d}   \eqno (2.30)
$$
with uniform $c,C$ don't depending on $d,b_1$. Now estimate (2.1) follows from estimate (2.15)
of Proposition 3 because $(1-|w'|)/(1-|w|)$ tends to 1 as $d\to 0$. $\Box$

6) {\it Estimates for derivatives of the mappings $h_w$, $f_{\mu_{F_w}}$ and $f_{\lambda_{F_w}}$.}
For convenience we place here these estimates, which we shall use in the next section.

\begin{proposition} a)
$$
  c\le |(h_w )_z (z)| \le \frac {C}{|1-\bar w Z_w^{-1}(z)|^2}  , \eqno (2.31)
$$
with some uniform $c,C$, and we have for $|(h_w )_{\bar z}|$ an estimate analogous to the 
right inequality.

b)
$$
  c|\xi | \le |f_{\mu_{F_w},\zeta}(\zeta )|\le C|\xi |^{-3}.  \eqno (2.32)
$$
with some uniform $c,C$, and $f_{\mu_{F_w},\zeta}$ is uniformly bounded at $\xi < -1$.
We have the same estimates for $(\tilde f_{\mu_{F_w}})_{\zeta}$.

c)
$$
  |(f_{\lambda_{F_w}})_{\zeta}(\zeta )|\le C|{\rm Re}\tilde f^{-1}_{\mu_{F_w}}(\zeta )|^{-3}    \eqno (2.33)
$$
with some uniform $c,C$ at $|{\rm Re}\tilde f^{-1}_{\mu_{F_w}}(\zeta )|\le 1$ and $(f_{\lambda_{F_w}})_{\zeta}(\zeta )$
is uniformly bounded at ${\rm Re}\tilde f^{-1}_{\mu_{F_w}}(\zeta )<-1$. 
\end{proposition}

{\bf Proof}. a) Remind that $h_w =f_w \circ Z_w^{-1}$, where $f_w =\vpf_{w'}\circ f\circ \vpf_w^{-1}$.
The map $Z_w$ is quasiconformal with the complex dilatation $\mu_0$, $|(Z_w )_z|$ 
is bounded from below and from above by constants depending only on $d$. Thus it is enough to
prove estimate analogous to (2.31) for $f_w$ if we replace $Z_w^{-1}(z)$ by $z$. We have
$$
  (f_w )_z (z) =(\vpf_{w'})_z \circ f\circ\vpf_w^{-1}(z) f_z \circ\vpf_w^{-1}(z)(\vpf_w^{-1})_z (z)  =
$$
$$
  =\frac{1-|w'|^2}{(1-\bar w' f\circ\vpf_w^{-1}(z))^2 }f_z \circ\vpf_w^{-1}(z)\frac{1-|w|^2}{(1-\bar w z)^2} .
$$
Applying (2.4) we see that 
$$
  |1-\bar w' f\circ\vpf_w^{-1}(z)|\ge c (1-|w|) .
$$
Thus,
$$
  c|f_z \circ\vpf_w^{-1}(z)|\le |(f_w )_z (z)| \le \frac {C}{|1-\bar w z|^2}|f_z \circ\vpf_w^{-1}(z)|  .  
$$
Since $|f_z |$ is uniformly bounded from below and from above, we obtain (2.31). 
The estimate for $\bar z$-derivative is analogous.

b) The map $f_{\mu_{F_w}}$ written in the chart $z$ is the map $\hat g_w=\exp\circ f_{\mu_{F_w}}\circ\log$.
The function $f_{\mu_{F_w}}(\zeta )-\zeta$ is uniformly bounded at $\xi\to -\infty$. It means that
$\hat g_w$ is a quasiconformal homeomorphism of $\bbC$ with the complex dilatation $\mu_{h_w}$ and  
with derivatives at zero uniformly bounded from below and from above.
We can represent it on $D$ as $\hat h_w \circ h_w$, where $\hat h_w$ is a holomorphic univalent 
uniformly bounded function with the derivative at zero uniformly bounded from below and from above. 
For such functions we have estimates ([Pom])
$$
  c(1-|z|)\le |\hat h'_w (z)|\le \frac{C}{1-|z|} .   \eqno (2.34)
$$ 
Also, for a variable $w$, $h_w$ is a family of quasiconformal mappings of $D$ onto itself with
uniformly bounded dilatations and, hence, there are the estimates
$$
  c(1-|z|)\le 1-|h_w (z)| \le C(1-|z|)   \eqno (2.35)
$$  .
with uniform $c,C$.

Further, at $\xi \le 0$
$$
  (f_{\mu_{F_w}})_{\zeta}(\zeta) = (\hat h'_w \circ h_w \circ\exp\zeta )^{-1}\cdot (h_w )_z \circ\exp\zeta 
$$
By left inequality (2.31), estimate (2.35), and right inequality (2.34), we obtain
$$
    |f_{\mu_{F_w},\zeta}(\zeta )|\ge c|\xi |  
$$
at $|\xi |\le 1$. Analogously, right inequality (2.31) 
and left inequality (2.34) together with (2.35) yield right estimate (2.32). 

On $\bbC\setminus D$ the functions $\hat g_w$ are univalent holomorphic, uniformly bounded away from zero and having uniformly bounded derivatives at infinity.
Thus there is also the uniform estimates $c(|z|-|1)|\hat g_w (z)_z )\le C (|z|-1)^{-1}$. 
We obtain estimates (2.32) also at $0<|\xi |\le 1$.

The case of $(\tilde f_{\mu_{F_w}})_{\zeta}$ is analogous.

c) From (2.28) follows that we have for the derivatives of $\tilde h_w$ the estimates
analogous to (2.31). Since $\tilde h_w =\exp \circ f_{\lambda_{F_w}}\circ \tilde f_{\mu_{F_w}}\circ\log$ 
we have
$$
  |(f_{\lambda_{F_w}})_{\zeta}\circ \tilde f_{\mu_{F_w}}(\zeta )||(\tilde f_{\mu_{F_w}})_{\zeta}(\zeta )|\le C|\xi |^{-2} ,
$$
and by left estimate (2.32) applied to $\tilde f_{\mu_{F_w}}$, we obtain 
$$
  |(f_{\lambda_{F_w}})_{\zeta}(\zeta )|\le C|{\rm Re}\tilde f^{-1}_{\mu_{F_w}}(\zeta )|^{-3}    
$$
at $|{\rm Re}\tilde f^{-1}_{\mu_{F_w}}(\zeta )|\le 1$. The rest of the proposition is obvious.
$\Box$

7) {\it Inductive change of variables}.
Now we pass to estimates of derivatives of higher orders. 
Consider the decomposition of $F_w$ at $\xi\to -\infty$
$$
  F_w (\zeta )=\zeta +\log a_0 +a_0^{-1}e^{\xi}R_1 (\vpf )+...+a_0^{-1}e^{k\xi}R_k (\vpf )+...    (2.36)
$$
Let $\chi (\xi )$ be a function equal to 1 at $\xi <0$ and to 0 at $\xi >2$ with the derivative less
then 1 by modulus. We define some successive change of variables.
Put $\chi_0 (\xi )=\chi ((4|\log a_0 |)^{-1}(\xi -c_0))$ if $|\log a_0 |\ge 1/4$ and $\chi_0 (\xi )=\chi (\xi -c_0)$
if $|\log a_0 |< 1/4$ with some $c\le -2$ large enough by modulus in both cases and define the new variable $\zeta_0$
$$
  \zeta =\sigma_0 (\zeta_0 )=s_0 (\xi_0 +i\vpf_0 )=\zeta_0 -\log a_0 \chi_0 (\xi_0 ) .  \eqno (2.37)
$$
The derivative of the function $|\log a_0 |\chi_0$ is less then 1/4, $\log a_0$ has uniform bounds independent of $w$, 
and we can set $c_0$ also independently of $w$. We obtain the estimates $|(\sigma_0 )_{\zeta_0}-1|<1/4$, 
$|(\sigma_0 )_{\bar\zeta_0}|<1/4$. Hence, $\sigma_0$ is a homeomorphism of the left half-plain onto itself.  
We get the new map $F_{0w}$ with the asymptotics
$$
  F_{0w}(\zeta_0 )=F_w (\sigma_0 (\zeta_0 ))=\zeta_0 +a_0^{-1}e^{\xi_0}R_1 (\vpf_0 )+...+a_0^{-k}e^{k\xi_0}R_k (\vpf_0 )+...=
$$
$$
  =\zeta_0 +e^{\xi_0}R_{01}(\vpf_0 )+...+e^{k\xi_0}R_{0k}(\vpf_0 )+...  .
$$

Suppose we have uniform estimates for the coefficients of the expansion of $h_w$ at zero
up to order $k\ge 1$.
Define the variable $\zeta_1$
$$
  \zeta_0 =\sigma_1 (\zeta_1 )=\sigma_1 (\xi_1 +i\vpf_1 )=\zeta_1 -e^{\xi_1}R_{01}(\vpf_1 )\chi (\xi_1 -c_1 ) .  \eqno (2.38)
$$
We get the map
$$
  F_{1w}(\zeta_1 )=F_{0w}(\sigma_1 (\zeta_1 ))=\zeta_1 +e^{2\xi_1}R_{12}(\vpf_1 )+...  .
$$ 
Analogously, we define the variable $\zeta_l$, $l\ge 2$
$$
  \zeta_{l-1} =\sigma_l (\zeta_l )=\sigma_l (\xi_l +i\vpf_l )=\zeta_l -e^{2\xi_l}R_{(l-1)l}(\vpf_l )\chi (\xi_l -c_l ) .  
$$
After $k$ successive changes of variables we obtain the map $F_{(k-1)w}=F_w \circ\sigma_0 \circ\sigma_1\circ...\circ\sigma_{k-1}$
$$
  F_{(k-1)w}(\zeta )=\zeta +e^{k\xi}R_{(k-1)k}(\vpf )+...  .
$$
Here we returned to the notation $\zeta$ for the variable. 

In the chart $z$ we get the transformations
$$
  z=s_0 (z_0 )=\exp (\sigma_0 (\log z_0 ))=z_0 \exp [-\log a_0 \chi_0 (\log |z_0 |)] ,  \eqno (2.39)
$$  
...
$$
  z_{l-1} =\exp (\sigma_l (\log z_l )) ,    \eqno (2.40)
$$
...

and the resulting mapping $h_w^{(k)}=h_w \circ s_0...\circ s_{k-1}$, 
$$
  h_w^{(k)}(z)=z+c_{(k+1)\bar 0}z^{k+1}+...+c_{0\overline{(k+1)}}\bar z^{k+1}+...=z+c_{(k+1)\bar 0}z^{k+1}+P_k (z)+... .
$$

By the inductive assumption, we have uniform estimates for the coefficients of the expansion of $h_w$ at zero
up to order $k$.
It is easy to see that to obtain estimates for the derivatives at zero of order $k+1$ 
it is enough to estimate the coefficients $c_{(k+1)\bar 0},...,c_{0\overline{(k+1)}}$. 
Indeed, if we have decomposition (2.18) for $h_w$,
then $c_{(k+1)\bar 0}$ (and analogously other coefficients at the terms of order $k+1$)
has the representation 
$$
 c_{(k+1)\bar 0}= b_{(k+1)\bar 0}/a_0^{k+1}+p_{(k+1)\bar 0} ,
$$
where $p_{(k+1)\bar 0}$ is a polynomial from $b_{i\bar j}/a_0^m$, $i+j=m$ with $m,k+1$. 

Let $s$ be the composition $s=s_0 \circ s_1...\circ s_{k-1}$. The complex dilatation of the map $h_w^{(k)}$ is
$$
  \mu_{h_w^{(k)}}(z)=\frac{\mu_{h_w}-\mu_s}{1-\bar\mu_s \mu_{h_w}}\cdot\frac{s_z}{\overline{s_z}}\circ s^{-1}(z )=
$$
$$
  =H_{k+1}(z)+|z|^{k+2}h(z)  ,
$$
where
$$
  H_{k+1}(z)=c_{k\bar 1}z^k +...+(k+1)c_{0\overline{(k+1)}}\bar z^k  .
$$ 
The coefficients $c_{k\bar 1}/a_0$,...,$c_{0\overline{(k+1)}}/a_0$ can be represented
as polynomials in the variables $a_0^{-1}$, the coefficients of of the expansion of $h_w$ at zero
up to order $k$, and the coefficients of order $k+1$ of $\mu_{h_w}$. All these coefficients
are uniformly bounded either by the inductive assumption either by Proposition 2. 
It follows that $c_{k\bar 1}/a_0$,...,$c_{0\overline{(k+1)}}/a_0$ are uniformly bounded. Thus it is enough
only to estimate the coefficient $c_{(k+1)\bar 0}$. Define the new transformation
$$
  z=\tilde s_k (z_1 )=z_1 [1-z_1^{-1}P_k (z_1 )\chi (\log |z_1 | -c_k )]  .  \eqno (2.41)
$$
Again adopting the notation $z$ for the chart we obtain the map
$$
  g_w^{(k)}(z)=h_w^{(k)}\circ \tilde s_k (z)=z+c_{(k+1)\bar 0}z^{k+1}+O(|z|^{k+2})
$$
with the complex dilatation $\mu_{g_w^{(k)}}(z)=O(|z|^{k+1})$.    

Note that we obtained the relations
$$
  (g_w^{(k)})_{z^{k+1}}(0)=(h_w )_{z^{k+1}}(0)+p_k (\{(h_w )_{(l)}(0)\}) ,  \eqno (2.42)
$$
where $p_k$ is a polynomial in $(h_w )_{z^l}(0)$ with $l\le k$ and in 
$(\mu_{h_w})_{(r)}(0)$, $|r|\le k$. From (2.42) follow the inverse relations
$$
  (h_w )_{(m)}(0)=q_{(m)}(\{(g_w^{(l)})_{z^{l+1}}(0),(\mu_{h_w})_{(r)}(0)\}) .  \eqno (2.43)
$$
Here $|m|=k+1$ and $q_{(m)}$ is a polynomial in $(g_w^{(l)})_{z^{l+1}}(0),(\mu_{h_w})_{(r)}(0)$ with $l\le k$,
$|r|\le k$.

\begin{proposition} Suppose we have uniform estimates for the coefficients of the expansion of $h_w$ at zero
up to order $k\ge 1$. Then we can put $c_o ,c_1 ,...,c_k$ such that $s_0 ,s_1 ,...\tilde s_k$
will be homeomorphisms with uniformly bounded derivatives of order up to $k+1$. The first derivatives $(g_w^{(k)})_z$, $(g_w^{(k)})_{\bar z}$ 
and all derivatives $(\mu_{g_w^{(k)}})_{(l)}$,
 $2\le |(l)|\le k+1$ will be uniformly bounded on the disk $D_{1/2}$ . 

If the coefficients of the expansion of $h_w$ at zero of order $l$, $l\le k$ have the estimates
$C(d+b_1 +...+b_l )$ with some uniform $C$ and $d, b_1 ,...,b_k$ are small enough with some uniform estimates, we can 
put $c_o ,c_1 ,...,c_k$ such that
$$
  |(g_w^{(k)})_z -1|\le C(d+b_1 +...b_k ),\,|(g_w^{(k)})_{\bar z}|\le C(d+b_1 +...b_k )  \eqno (2.44)` 
$$
$$
  |(\mu_{g_w^{(k)}})_{(l)}|\le C(d+b_1 +...b_k ),2\le |(l)|\le k+1   \eqno (2.45)
$$
on $D_{1/2}$ with uniform $C$ independent of $d,b_1 ,...,b_k$.
\end{proposition}

{\bf Proof}. We already proved that $s_0$ is homeomorphic and here we give only a more
explicit estimate. By (2.39), (2.38),
$$
  (s_0 )_z (z_0 )= \exp \{-\log a_0 \chi ((4|\log a_0 |)^{-1}(\log |z_0 | -c_0 ))\}\times 
$$
$$  
  \{1-\frac{z_0}{4} \chi' ((4|\log a_0 |)^{-1}(\log |z_0 | -c_0 ))\bar z_0^{1/2}/(2z_0^{3/2})\}  \eqno (2.46) 
$$
if $|\log a_0 |\ge 1/4$, and
$$
  (s_0 )_z (z_0 )= \exp \{-\log a_0 \chi (\log |z_0 | -c_0 )\} 
  \{1-z_0 \log a_0 \chi' (\log |z_0 | -c_0 )\bar z_0^{1/2}/(2z_0^{3/2})\}  \eqno (2.47) 
$$
if $|\log a_0 |\le 1/4$. Since $|\chi' |<1$, we can see that
$$
  \exp (-|\log a_0|)/2  \le |(s_0 )_z (z_0 )|\le 3/2\max\{1,\exp (|\log a_0|)\} .
$$
Also, applying estimate (2.30), we see that
$$
   |(s_0 )_z (z_0 )-1|\le C(d+b_1 )^2 ]   \eqno (2.48)
$$  
with some uniform $C$ for $d+b_1$ small enough. Analogously,
$$
  (s_0 )_{\bar z}(z_0 )=\exp\{-\log a_0 \chi ((4|\log a_0 |)^{-1}(\log |z_0 | -c_0 ))\}\times
$$
$$  
  \frac{z_0}{4}\chi' ((4|\log a_0 |)^{-1}(\log |z_0 | -c_0 ))z_0^{1/2}/(2\bar z_0^{3/2})   \eqno (2.46')
$$
if $|\log a_0 |\ge 1/4$, and
$$
  (s_0 )_{\bar z}(z_0 )= \exp\{-\log a_0 \chi (\log |z_0 | -c_0 )\}z_0 \log a_0 \chi' (\log |z_0 | -c_0 )z_0^{1/2}/(2\bar z_0^{3/2})   \eqno (2.47')
$$   
if $|\log a_0 |\le 1/4$. We obtain the estimate
$$
  |(s_0 )_{\bar z}(z_0 )|\le C[(d+b_1 )^2 ]\exp[C(d+b_1 )^2 ]\le 2C(d+b_1 )^2   \eqno (2.48')
$$
with some uniform $C$ at $d+b_1$ small enough.

Note that we can chose $c_0$ independent of $d$ and $b_1$, for example, we can
put $c_0 =-5$.

Now, by (2.39) for $l\le k-1$,
$$
  (s_l )_z (z_l )=1-(P_l )_{z_1}(z_1 )\chi (\log |z_l | -c_l )- P_l (z_l )\chi' (\log |z_l | -c_l )\bar z_l^{1/2}/(2z_l^{3/2})  .  \eqno (2.49)
$$
Here $P_l$ is a homogeneous polynomial in the variables $z_1,\bar z_1$ of order $l+1$. Its coefficients
have uniform estimates by the inductive assumption. Moreover, these coefficients are polynomials
in $(h_w )_{(j)}(0)$ with $2\le |(j)|\le l$ without a term of zero order. By the second inductive assumption,
we can estimate them as $C(d+b_1 +...+b_l )$ if $d,b_1 ...,b_l$ are small enough. Remind that
$\chi (\log |z_l | -c_l )\ne 0$ only if $\log |z_l | \le c_l$. Suppose the coefficients of $P_l$
are bounded by the constant $M$. Then $(|P_l )_{z_1}(z_l )|\le M(l+2)^2 |z_l|^l$, 
$|P_l (z_l )\bar z_l^{1/2}/(2z_l^{3/2})|\le M(l+2)|z_l|^{l}$. If we put $c_l <\log [-8M(l+2)^2]$,
then we obtain $|(s_l )_z (z_l )-1|\le 1/4$. If $d,b_1 ...,b_l$ are small enough, then we can put $c_l =-5$
and obtain the estimate
$$
  |(s_l )_z (z_l )-1|\le C(d+b_1 +...+b_l )  .  \eqno(2.50)
$$
Also,
$$
  (s_l )_{\bar z}(z_l )= (P_l )_{\bar z}(z_l )\chi (\log |z_l | -c_l )+ |P_l (z_l )\chi' (\log |z_l | -c_l )z_l^{1/2}/(2\bar z_l^{3/2})  \eqno (2.49')
$$    
and we obtain the estimate $|(s_l )_{\bar z}(z_l )|\le 1/4$ in the general case and the estimate 
$$
  |(s_l )_{\bar z}(z_l )|\le C(d+b_1 +...+b_l )   \eqno (2.50')
$$  
if $d,b_1 ...,b_l$ are small enough.

Consider at last the derivative 
$$
  (\tilde s_k )_z (z_k )=1-(\tilde P_k )_{z_k}(z_k )\chi (\log |z_k | -c_k )- \tilde P_k (z_k )\chi' (\log |z_k | -c_k )\bar z_k^{1/2}/(2z_k^{3/2})  .  \eqno (2.51)
$$
The polynomial $\tilde P_k$ is the sum $\tilde c_{z^k \bar 1}(0)z^k \bar z +...+\tilde c_{0\overline{k+1}}\bar z^{k+1}$
and the coefficients have the representation 
$c_{l\overline{k+1-l}}=\alpha_{l\overline{k+1-l}}(\mu_{h_w})_{z^l \bar z^{k-l}}(0)+p_{l\overline{k+1-l}}$,
where $p_{l\overline{k+1-l}}$ is a polynomial in $(h_w )_{(j)}(0)$ with $2\le |(j)|\le k$ without a term of zero order.
Applying Proposition 2 and the inductive assumptions we obtain for these coefficients an uniform estimate
in the general case and the estimate $C(d+b_1 +...+b_k )$ if $d,b_1 ...,b_k$ are small enough.
As above, we obtain the estimate $|(\tilde s_k )_z (z_k )-1|\le 1/4$ at an appropriate $c_k$.
If $d,b_1 ...,b_k$ are small enough, we can put $c_k =-5$ and obtain
$$
  |(\tilde s_k )_z (z_k )-1|\le C(d+b_1 +...+b_l )  .  \eqno(2.52)
$$
Analogously,
$$
  (\tilde s_k )_{\bar z}(z_k )= (\tilde P_k )_{\bar z}(z_k )\chi (\log |z_k | -c_k )+ |\tilde P_k (z_k )\chi' (\log |z_k | -c_k )z_k^{1/2}/(2\bar z_k^{3/2})  \eqno (2.51')
$$    
and we obtain an uniform estimate in the general case and the estimate 
$$
  |(\tilde s_k )_{\bar z}(z_k )|\le C(d+b_1 +...+b_k )    \eqno(2.52')
$$ 
if $d,b_1 ...,b_k$ are small enough.

We proved that $s_0 ,...,\tilde s_k$ are homeomorphic if we set $c_0,...,c_k$ as above.
Put $S^{(k)}=s_0 \circ s_1 \circ...\circ\tilde s_k$. By definition, $g_w^{(k)}=h_w \circ S^{(k)}$, and 
we obtain estimates (2.44) from (2.48), (2.48'), (2.50), (2.50'), (2.52), (2.52') and corollary
of Proposition 9. 

Now consider derivatives of the right parts of (2.46), (2.47), (2.46') and (2.47'). 
In the first derivatives there appear the terms of the types
$O(1)\log a_0 \chi' (\log |z_0 | -c_0 )|z_0 |^{-1}$ and $O(1)\log a_0 \chi''(\log |z_0 | -c_0 )|z_0 |^{-1}$.
Since derivatives of the function $\chi$ don't equal zero only at $0< \log (|z_0 |) -c_0 <2$, 
we obtain the estimates  
$$
  |(s_0 )_{(l)}(z_0 )|\le C(d+b_1 )^2 e^{|c_0 |},\, |(l)|=2,\,|z_0 |\le 1/2 .
$$  
Derivatives of order $l$ are sums of items containing multiples of the types
$O(1)\log a_0$, $\chi^{(j)}(\log [(|z_0 |) -c_0 ),1\le |(j)|\le l$ and $|z_0 |^{-s}$ with some integer $s$.
At each differentiation there appears no more than one multiple $|z_0 |^{-1}$,
and we see that we have the estimate 
$$
  |(s_0 )_{(j)}(z_0 )|\le C(d+b_1 )^2 e^{(|(j)|-1)|c_0 |}  ,\,|z_0 |\le 1/2
$$  
with some uniform $C$. Remind that we can put $c_0 =-5$. At small $d$ and $b_1$ we 
obtain the estimate $C(d+b_1 )^2$ with $C$ independent of $d$ and $b_1$.

We estimate derivatives of $s_l$ and $(\tilde s_k )$ analogously. 
Differentiating the right parts of (2.49), (2.49'), (2.51) and (2.51') we obtain terms of order 
$O(1)(d+b_1 +b_l)e^{(|(j)|-l)|c_l |}$ for derivatives of order $j$. At our choice of
$c_l$ we have uniform estimates in the general case and the estimate $C(d+b_1 +..+b_l )$
if $d$, $b_1$, ...,$b_l$ are small.

Now we have
$$
  \mu_{g_w^{(1)}}=\frac{\mu_{h_w}-\mu_{S^{(k)}}}{1-\bar\mu_{S^{(k)}} \mu_{h_w}}\cdot\frac{S^{(k)}_z}{\overline{S^{(k)}_z}}\circ (S^{(k)})^{-1} .
$$
Now the proposition follows from the estimates for derivatives of $S^{(k)}$ and Proposition 2. 
$\Box$

8) {\it Estimates for higher derivatives}. From the next proposition we obtain by induction
estimate (2.2) and estimate (2.3). 

\begin{proposition} a) Suppose that we have uniform bounds for the coefficients of the expansion of $h_w$ at zero
up to order $k\ge 1$. Then the derivative $(g_w^{(k)})_{z^{k+1}}$ is uniformly bounded.

b) If the coefficients of the expansion of $h_w$ at zero of order $l$, $l\le k$ have the estimates
$C(d+b_1 +...+b_l )$ with some uniform $C$ and $d, b_1 ,...,b_k$ are small enough with some uniform estimate,
then  $|(g_w^{(k)})_{z^{k+1}}|\le C(d+b_1 +...+b_{k+1})$.
\end{proposition}
  
{\bf Proof}. 
a) Denote by $\mu_w^{(k)}$ the Beltrami coefficient of the map $g_w^{(k)}$. We have
$$
  g_w^{(k)}(z)=\frac{1}{2\pi i}\int_{\pd D_{1/2}}\frac{g_w^{(k)}(t)}{t-z}dt+\frac{1}{\pi i}\int_{D_{1/2}}\frac{\mu_w^{(k)}(t)g_{w,t}^{(k)}(t)}{t-z}dS_t  . 
\eqno (2.53)
$$
at $|z|\le 1/2$. All derivatives in zero of the first integral in (2.53) are uniformly bounded. 
Now we have
$$
  \frac{1}{t-z}=\frac{1}{t}+\frac{z}{t^2}+...+\frac{z^{k+1}}{t^{k+2}}+\frac{z^{k+2}}{t^{k+2}(t-z)} .
$$
The first derivatives of the function $g_w^{(k)}$
and derivatives of $\mu_w^{(k)}$ of order $k+1$ are uniformly bounded on $D_{1/2}$ by Proposition 11.

Now $\mu_w^{(k)}(z)=|z|^{k+1}\gamma (z)$, where 
$$
  |\gamma (z)|\le \frac{1}{(k+1)!}\max_{|z|\le 1/2}\|\pd^{k+1}\mu_w^{(k)}\| ,
$$
where $\|\pd^{k+1}\mu_w^{(k)}\|$ is the sum of modulus of derivatives of order $k+1$. Hence, $|\gamma |$ is uniformly
bounded on $D_{1/2}$. 

The integrals
$$
  \int_{D_{1/2}}\frac{|t|^{k+1}\gamma (t)g_{w,t}^{(k)}(t)}{t^l}dS_t ,\,l\le k+2   \eqno (2.54)
$$ 
are all uniformly bounded. We get
$$
  \int_{D_{1/2}}\frac{\mu_w^{(k)}(t)g_{w,t}^{(k)}(t)}{t-z}dS_t =c_0 +zc_1 +...+z^{k+1}c_{k+1}+ 
z^{k+2}\int_{D_{1/2}}\frac{|t|^{k+1}\gamma (t)g_{w,t}^{(k)}(t)}{t^{k+2}(t-z)}dS_t .  \eqno (2.55)
$$
From the other hand, we have the estimate
$$
  \int_{D_{1/2}}\frac{dS_t}{|t||t-z|}\le C\|\log |z||
$$
with some uniform $C$. We see that integral (2.55) has $k+1$ derivatives at zero and its $z$-derivative 
of order $k+1$ is equal to $(k+1)!c_{k+1}$. 

b) Remind that $g_w^{(k)}(z)=h_w (z)$ at $|z|=1/2$. Applying corollary of Proposition 9, we see
that all derivatives at zero of order higher than 1 of the first integral in (2.53)
have the estimate $C(d+b_1 )^2$. To estimate the derivatives of the second integral in (2.52)
we need to estimate integrals (2.54). We obtain the estimates by Proposition 11.
$\Box$

Now we finish the proof of Lemma 1.

{\bf Proof of estimate (2.3)}. Applying inductive relations (2.43) we obtain the estimates for the derivatives 
at zero of $h_w$ from the estimates for $g_w^{(k)}$. The corollary follows now from 
Proposition 3. $\Box$

We shall use in the next section the following estimate:

\begin{proposition}
$$
  |(f_{\mu_{F_w}})_{(k)}(\zeta )|\le Ce^{\xi}(1+|\xi |^{-6}),\,|(k)|=2 
$$
with some uniform $C$ at $\xi \le 1$.
\end{proposition}

{\bf Proof}. Estimate at first the second derivatives of the map $h_w$. Analogously to the
proof of Proposition 10 a) it is enough to obtain the estimates for $f_w$. We have
$$
  (f_w )_{zz}(z)=\frac{\pd}{\pd z}\left[\frac{1-|w'|^2}{(1-\bar w' f\circ\vpf_w^{-1}(z))^2} f_z \circ\vpf_w^{-1}(z)\frac{1-|w|^2}{(1-\bar w z)^2}\right]=
$$
$$
  =\frac{2w'(1-|w'|^2 )}{(1-\bar w' f\circ\vpf_w^{-1}(z))^3}(f_z \circ\vpf_w^{-1}(z))^2\left[\frac{1-|w|^2}{(1-\bar w z)^2}\right]^2 +
$$
$$  
  +\frac{1-|w'|^2}{(1-\bar w' f\circ\vpf_w^{-1}(z))^2}\left[f_{zz}\circ\vpf_w^{-1}(z)\left(\frac{1-|w|^2}{(1-\bar w z)^2}\right)^2 +
  f_z \circ\vpf_w^{-1}(z)\frac{2w(1-|w|^2 )}{(1-\bar w z)^3}\right] .
$$  
Applying (2.4) we see that
$$
  |(f_w )_{zz}(z)|\le C\left[\frac {1-|w|}{(1-\bar w z)^4}|f_{zz}\circ\vpf_w^{-1}(z)|+\frac {1}{(1-\bar w z)^4}|f_z \circ\vpf_w^{-1}(z)|^2 \right] 
$$
Applying estimates (2.1), (2.8), and (2.9) we obtain
$$
  |f_{zz}\circ\vpf_w^{-1}(z)|\le \frac{C}{1-|\vpf_w^{-1}(z)|}\le C\frac{|1-\bar wz|}{(1-|w|)(1-|z|)} ,
$$
and, hence,
$$
  |(f_{2w})_{zz}(z)|\le \frac{C}{(1-|z|)^4} ,\, |(h_{2w})_{zz}(z)|\le \frac{C}{(1-|z|)^4}    \eqno (2.56)
$$
Obviously we have analogous estimates for other derivatives of second order.      

Return now to the map $f_{\mu_{F_w}}$.
As in the proof of Proposition 10 b) we can represent it in the chart $z$ on $D$ as 
$\hat g_w =\hat h_w \circ h_w$, where $\hat h_w$ is a holomorphic univalent uniformly bounded
function with derivative at zero uniformly bounded from below and from above. For such functions we have 
estimates (2.34) and ([Pom])
$$
  |\hat h''_w (z)|\le \frac{C}{(1-|z|)^2} .   \eqno (2.57)
$$ 
Applying this estimate and (2.31), (2.35), (2.56), we obtain at $\xi\le 0$
$$
  |(f_{\mu_{F_w}})_{\zeta\zeta}|\le Ce^{\xi}[(|(\hat g_w )_z\circ\exp\zeta |^2 +|(\hat g_w )_{zz}\circ\exp\zeta |]\le
$$
$$
  \le Ce^{\xi}[(|\hat h'_w \circ h_w \circ\exp\zeta |^2 |(h_w )_z \circ\exp\zeta |^2 +|\hat h''_w \circ h_w \circ\exp\zeta ||(h_w )_z\circ\exp\zeta |^2 +
$$
$$
 +|\hat h'_w \circ h_w \circ\exp\zeta ||(h_w )_{zz}\circ\exp\zeta |)]\le Ce^{\xi}(1+|\xi |^{-6})   
$$ 
at $|\xi |\le 1$ and analogous estimates for other second derivatives of $f_{\mu_{F_w}}$. Obviously, for derivatives of $\tilde f_{\mu_{F_w}}$
we have the same estimates.

On $\bbC\setminus D$ the functions $\hat g_w$ are univalent holomorphic with uniformly bounded derivatives at infinity
and uniformly bounded away from zero. Hence, there is the uniform estimate $|\hat g_w (z)|_{zz}\le C (|z|-1)^{-2}$. Thus
$|(f_{\mu_{F_w}})_{\zeta\zeta}(\zeta )|\le C |\zeta |^{-2}$ at $0<|\xi |\le 1$.  
We obtain the estimates for other derivatives of second order analogously. $\Box$

In conclusion of this section we obtain some estimates for derivatives of the principal solutions.
Remind the construction of the normal solutions (see, for example, [Ah]).

Let $f_{\mu}$ be the principal solutions with the Beltrami coefficient $\mu$ and 
put $f^0_{\mu}(z)=f_{\mu}(z)-f_{\mu}(0)$. Put 
$\tilde f (z)=\overline{f^0_{\mu}(\bar z^{-1})})^{-1}$ and  
$$
 \lambda =\left(\mu \frac{\tilde f_z}{\overline{\tilde f_z}}\right)\circ\tilde f^{-1}  ,
$$
Let $f_{\lambda}$ be the corresponding principal solution and $f^0_{\lambda}=f_{\lambda}-f_{\lambda}(0)$. 
Define 
$$
  f_c=f^0_{\lambda}\circ\tilde f  .    
$$   
That is
$$
  f_c (z) =f^0_{\lambda}\circ (\overline{f^0_{\mu}(\bar z^{-1})})^{-1}  .\eqno (2.58)
$$
Then
$$
  f=f_c /f_c (1)   \eqno (2.59)  ,
$$  
where $|f_c (1)|=1$.   

\begin{proposition} At assumptions of Lemma 1 we have the estimates
$$
  |f^0_{\mu}(z)|\ge c|z| ,   \eqno (2.60)
$$  
$$
  c|1-|z||\le |(f^0_{\mu})_z (z)|\le C|1-|z||^{-1}     \eqno (2.61)
$$
with $c,C$ depending only on $d$ and $b_1$.
$$
  |(f^0_{\mu})_{(k)}(z)|\le C|1-|z||^{-2},\,|(k)|=2  . \eqno (2.62) 
$$
with $C$ depending only on $d, b_1 ,b_2$. Analogous estimates hold for $\tilde f$.
$$
  c|1-|\tilde f^{-1}(z)||\le |(f^0_{\lambda})_z (z)|\le C|1-|\tilde f^{-1}(z)||^{-1}  \eqno (2.63)
$$  
also with $c,C$ depending only on $d$ and $b_1$.
\end{proposition}

{\bf Proof}. We have the representation on $D$: $f^0_{\mu}=H\circ f$, where $H$ is a holomorphic univalent function 
mapping zero to zero and $f$ is the normal mapping. Since $f^0_{\mu}$ is bounded on $D$ with a 
bound depending only on $d$ and we have estimate (2.1), it follows that the function $H$ is bounded and has the derivatives at
zero bounded from above by a constant depending only on $d$. Let show that its
$z$-derivative at zero is bounded also from below. It is enough to prove that $(f^0_{\mu})_z (0)$
is bounded from below, that is, to prove estimate (2.60).

Since $f_c$ differs from $f$ by the multiple equal to 1 by modulus, we see that
$|f_c (z)|\le C|z|$ with $C$ depending only on $d$ and $b_1$. Also,
$|f^0_{\lambda}(z) -z|\le C$ with $C$ also depending only on $d$ and $b_1$. It means 
that
$$
  |f_c (z)-\overline{f^0_{\mu}(\bar z^{-1})})^{-1}|\le C
$$
and, hence, 
$$
  |(f^0_{\mu}(\bar z^{-1}))^{-1}|\le C|z| .
$$  
It means that we obtain estimate (2.60). Thus $|H'(0)|$ is bounded also from below and we can
apply estimates (2.34) and (2.57) to $H$. Applying also (2.4) we obtain estimates (2.61), (2.62)
at $|z|\le 1$. 

On $\bbC\setminus D$ the function $f_c$ is univalent holomorphic with uniformly bounded derivatives at infinity.
Also, on this domain $f_c$ is uniformly bounded away from zero. We obtain estimates (2.61) and (2.62) as in Propositions
10 and 13. The estimates for $\tilde f$ follow by symmetry.

Solving the system
$$
  (f_c )_z =(f^0_{\lambda})_z \circ\tilde f \cdot\tilde f_z + (f^0_{\lambda})_{\bar z}\circ\tilde f \cdot\overline{\tilde f_{\bar z}} 
$$
$$
  (f_c )_{\bar z} =(f^0_{\lambda})_{\bar z}\circ\tilde f \cdot\overline{\tilde f_z} + (f^0_{\lambda})_z \circ\tilde f \cdot\tilde f_{\bar z} 
$$  
we obtain
$$
  (f^0_{\lambda})_z \circ\tilde f =\frac{(f_c )_z}{\tilde f_z}\cdot\frac{1-\bar\mu_{\tilde f}}{1-\mu_{\tilde f}\bar\mu_{\tilde f}}
$$
Estimates (2.63) follow from (2.61) and  from boundedness from below and from above of $(f_c )_z$  
$\Box$

\section{Estimates for differences of derivatives of the normal mappings with different complex dilatations.}

The proof of the lemma below is long and tedious but essentially simple. 
We adopt the notation $P(d)$ for the supremum of $p$ such that the series 
$1+\mS\mu +\mS\mu\mS\mu +...$ converge in $L^p$ if $\|\mu \|_C \le d$.
 
\begin{lemma} Let $\mu_1 ,\mu_2$ be functions on $D$ satisfying assumptions of Lemma 1
with the same $d,b_1 ,...,b_k$. Let $f_1 ,f_2$ be the corresponding normal mappings. Then
we have the estimates:

a) 
$$
  \|f_1-f_2 \|_C \le C\|\mu_1 -\mu_2 \|_p^{\alpha}  ,   \eqno (3.1)
$$
where $2<p<P(d)$ and $\alpha$ depends only on $d$ and $b_1$.

b) Let $(k)$ be a multi-index, $|(k)|\ge 1$. Fix some $0<R<1$. Then, for $z\in D$,
$$
  |(f_1 )_{(k)}(z)-(f_2 )_{(k)}(z)|\le \frac{C}{(1-|z|)^{|(k)|-1}}\min\left\{1,\left[\frac{\|\mu_1 -\mu_2 \|_p^{\alpha}}{1-|z|}+\right.\right.
$$
$$
 \left.\left.+\sup_{D_R}|\mu_1 -\mu_2 |^{\alpha}+[(1-R)/(1-|z|)]^{\alpha}+
 \sup_{0\le |(q)|\le |(k)|}\sup_{D_{R_z}}|(\mu_1)_{(q)}-(\mu_2)_{(q)}|^{\alpha}\right]\right\}  ,  \eqno (3.2)
$$ 
where $R_z <1$ is such that $1-R_z \ge a(1-|z|)$ for some $a$ depending only on $d$, $C$ is some 
uniform constant and $0<\alpha <1$ depends only on $d$, $2<p<P(d)$.
\end{lemma}

In the proof we shall use the terminology and the notations of the previous section.

{\bf Proof of estimate 3.1}. We put, as in Proposition 14,  
$f^0_{\mu_i}(z)=f_{\mu_i}(z)-f_{\mu_i}(0)$,  
$\tilde f_i (z)=\overline{f^0_{\mu_i}(\bar z^{-1})})^{-1}$,  
$$
 \lambda_i =\left(\mu_i \frac{\tilde f_{i,z}}{\overline{\tilde f_{i,z}}}\right)\circ\tilde f_i^{-1}  .
 \eqno (3.3)
$$
By definition, $f_{\lambda_i}$ is the corresponding principal solution, $f^0_{\lambda_i}=f_{\lambda_i}-f_{\lambda_i}(0)$, 
$f_{ci}=f^0_{\lambda_i}\circ\tilde f_i$, $f_i =f_{ci}/f_{ci}(1)$.   

By Proposition 14,
$$
  c|1-|z||\le |(f^0_{\mu_i})_z (z)|\le C|1-|z||^{-1}  ,   \eqno (3.4)
$$
$$
  c|1-|\tilde f_i^{-1}(z)||\le |f^0_{\lambda_i ,{(1)}}(z)|\le C|1-|\tilde f_i^{-1}(z)||^{-1}  ,    \eqno (3.5)
$$
$$
  |f^0_{\mu_i ,{(k)}}(z)|\le C|1-|z||^{-2},\,|(k)|=2  . \eqno (3.6) 
$$

We must estimate the value
$$
  |f_{c1}(z)-f_{c2}(z)||f_{c1}(1)|^{-1}+|f_{c2}(z)||f_{c1}(1)|^{-1}-f_{c2}(1)|^{-1}| .
$$
Since $|f_{c1}(1)|=1$, it is enough to estimate the difference  
$$
  |f_{c1}(z)-f_{c2}(z)|\le |f^0_{\lambda_1}\circ\tilde f_1 (z)-f^0_{\lambda_2}\circ\tilde f_1 (z)|+
  |f^0_{\lambda_2}\circ\tilde f_1 (z)-f^0_{\lambda_2}\circ\tilde f_2 (z)| .  \eqno (3.7)
$$

\begin{proposition}  Let $f_{\nu_i}$, $i=1,2$ be the principal solutions corresponding to 
compactly supported Beltrami coefficients $\nu_i$. Then
$$
  |f_{\nu_1}(z)-f_{\nu_2}(z)|\le \|(\nu_1 -\nu_2 )(|(f_{\nu_1})_z -1|+|(f_{\nu_2})_z -1|)\|_p   
$$  
$$
  \|f_{\nu_1 ,z}-f_{\nu_2 ,z}\|_p\le \|(\nu_1 -\nu_2 )(|(f_{\nu_1})_z -1|+|(f_{\nu_2})_z -1|)\|_p   
$$ 
$$
  |f^{-1}_{\nu_1}(\zeta )-f^{-1}_{\nu_2}(\zeta )|\le C\|(\nu_1 -\nu_2 )(|(f_{\nu_1})_{\zeta}-1|+|(f_{\nu_2})_{\zeta}-1|\|_p  
$$
for $2<p<P(d)$.
\end{proposition} 

{\bf Proof}. Let $h_i$ be the solution to the equation $h_i -\mS\nu_i h_i =1$. 
(Remind that here and below $\mS$ is the Beorling transform and $\mC$ is the Cauchy transform).
We have
$$
  |f_{\nu_1}(z)-f_{\nu_2}(z)|=|\mC\nu_1 h_1 (z)-\mC\nu_2 h_2 (z)| \le 
 C_p (\|(\nu_1 -\nu_2 )h_1 \|_p + \|\nu_2 \|_C \|h_1 -h_2 \|_p  .   \eqno (3.8)
$$

Further,
$$
  h_1 -h_2 =\mS[\nu_1 (h_1 -h_2)] +\mS[(\nu_1 -\nu_2 )h_2 ]  .
$$
That is,
$$
  h_1 -h_2 =(id - \mS\nu_1)^{-1}\mS[(\nu_1 -\nu_2 )h_2] ,
$$
and
$$  
  \|h_1 -h_2 \|_p \le C_p(\|(\nu_1 -\nu_2 )h_2 \|_p  .
$$   
with some constant $ C_p$ depending only on $p$. 
Since $h_i =(f_{\nu_i})_z -1$, we obtain the first two estimates of the proposition from (3.8). 

Prove the third estimate. 
Suppose $f_{\nu_1}(z_1 )=f_{\nu_2}(z_2 )$. Then
$$
  |z_1 -z_2 |=\mC\nu_1 h_1 (z_1 )-\mC \nu_2 h_2 (z_2 )|\le
C_p (\|(\nu_1 -\nu_2 )h_1 \|_p + \|\nu_2 \|_C \|h_1 -h_2 \|_p  .
$$
$\Box$

We can apply this proposition to $f_{\mu_i}$, $f_{\lambda_i}$ and, hence, to $f^0_{\mu_i}$, $f^0_{\lambda_i}$.
Also we have
$$
  |\tilde f_1 (z)-\tilde f_2 (z)|=|f^0_{\mu_1}(\bar z^{-1})f^0_{\mu_2}(\bar z^{-1})|^{-1}||f^0_{\mu_1}(\bar z^{-1})-f^0_{\mu_2}(\bar z^{-1})|\le
$$
$$
  \le C|z|^2 \|(\mu_1 -\mu_2 )(|(f_{\mu_1})_z -1|+|(f_{\mu_2})_z -1|)\|_p  . \eqno (3.9)
$$  
Here we applied estimate (2.60).

We must estimate the right part of inequality (3.7). We proceed in several steps.
In all inequalities below all constants such as $c$ or $\alpha$ depend only 
on $d$ and $b_1$.

1)
$$
  |f^0_{\mu_1}-f^0_{\mu_2}|_C \le C\|(\mu_1 -\mu_2 \|_p^{\alpha} ,
$$

{\bf Proof}. Fix some $p<p'<P(d)$ and some $r<1$. Applying Proposition 15 and estimate (3.4) we can
write
$$  
  \|(\mu_1 -\mu_2  )((f_{\mu_i})_z -1) \|_p \le C\left[r^{-1}\left(\int_{||z|-1|\ge r}|\mu_1 -\mu_2 |^p dS_z \right)^{1/p}+\right. 
$$
$$
  +\left.\left(\int_{||z|-1|\le r}|(f_{\mu_i})_z -1|^p dS_z \right)^{1/p}\right]
$$
The second integral we estimate by the Heolder inequality $\|fg\|_p \le \|f\|_{p'}\|g\|_{q'}$,  
$q'^{-1}+p'^{-1}=p^{-1}$. We obtain  
$$  
  \|(\mu_1 -\mu_2  )((f_{\mu_i})_z -1) \|_p \le  C[(1-r)^{-1}\|(\mu_1 -\mu_2 \|_p +
$$
$$  
  +\|(f_{\mu_i})_z -1\|_{p'} ({\rm mes}\{r\le |z|\le 1\})^{\frac{1}{p}-\frac{1}{p'}}]\le
  C[(1-r)^{-1}\|(\mu_1 -\mu_2 \|_p +r^{\frac{p'-p}{pp'}} .
$$     
If we define $r$ from the equation $(1-r)^{-1}\|(\mu_1 -\mu_2 \|_p =r^{\frac{p'-p}{pp'}}$, i.e., if
we put $1-r=\|(\lambda_1 -\lambda_2 \|_p^{\frac{pp'}{pp'+p'-p}}$, we get the estimate
$$
  |f_{\lambda_1}-f_{\lambda_2}|\le C\|(\lambda_1 -\lambda_2 \|_p^{\frac{p'-p}{pp'+p'-p}} .
$$
$\Box$

2) {\it At $|z|\le 1$}
$$
  |\tilde f_1 (z)-\tilde f_2 (z)|\le C\|\mu_1 -\mu_2 \|_p^{\alpha} ,
$$  

{\bf Proof}. It follows from step 1 and estimate (3.9). $\Box$

3)
$$
  \|(\tilde f_{1,z} -\tilde f_{2,z})\chi_D\|_p \le C\|\mu_1 -\mu_2 \|_p^{\alpha} ,
$$  
{\it Here $\chi_D$ is the characteristic function of $D$}.

{\bf Proof}. The proof in step 1 depends only on the right side of the first inequality 
of Proposition 15.
Since the second inequality has the same right part, we obtain our assertion from (3.9).$\Box$

4)
$$
  |f_{\lambda_1}-f_{\lambda_2}|_C \le C\|(\lambda_1 -\lambda_2 \|_p^{\alpha} ,
$$

{\bf Proof}. Again fix some $p<p'<P(d)$ and some $r<1$. Applying Proposition 15 and estimate (3.5) we can
write
$$  
  \|(\lambda_1 -\lambda_2  )((f_{\lambda_i})_z -1) \|_p \le C\left[r^{-1}\left(\int_{||\tilde f_i^{-1}(z)|-1|\ge r}|\lambda_1 -\lambda_2 |^p dS_z \right)^{1/p}+\right. 
$$
$$
  +\left.\left(\int_{||\tilde f_i^{-1}(z)|-1|\le r}|(f_{\lambda_i})_z -1|^p dS_z \right)^{1/p}\right]
$$
The second integral we again estimate by the Heolder inequality and obtain  
$$  
  \|(\lambda_1 -\lambda_2  )((f_{\lambda_i})_z -1) \|_p \le  C[(1-r)^{-1}\|(\lambda_1 -\lambda_2 \|_p +
$$
$$  
  +\|(f_{\lambda_i})_z -1\|_{p'} ({\rm mes}\{||\tilde f_i^{-1}(z)|-1|\le r\})^{\frac{1}{p}-\frac{1}{p'}}] .
$$
But 
$$
   {\rm mes}\{||\tilde f_i^{-1}(z)|-1|\le r\}=\int_{||z|-1|\le r} J(\tilde f_i )(z) dS_z  ,
$$
where $J(\tilde f_i )$ is the Jacobian of the transformation $z\mapsto\tilde f_i (z)$.   
Since $(\tilde f_i )_z -1$ belongs to $L_p$, we see that
$J(\tilde f_1 )$ restricted on the ring $\{||z|-1|\le r\}$ belongs to $L^{p/2}$. Hence,
by the Geolder inequality  
$$
   {\rm mes}\{||\tilde f_i^{-1}(z)|-1|\le r\}\le C r^{1-2/p} .
$$   
We obtain
$$  
  \|(\lambda_1 -\lambda_2  )((f_{\lambda_i})_z -1) \|_p \le
  C[(1-r)^{-1}\|(\lambda_1 -\lambda_2 \|_p +(1-r)^{\frac{(p'-p)(p-2)}{p^2 p'}} .
$$     
If we put $1-r=\|(\lambda_1 -\lambda_2 \|_p^{\frac{p^2 p'}{(p'-p)(p-2)+p^2 p'}}$, then we obtain the estimate
$$
  |f_{\lambda_1}-f_{\lambda_2}|\le C\|(\lambda_1 -\lambda_2 \|_p^{\frac{(p'-p)(p-2)}{(p'-p)(p-2)+p^2 p'}} .
$$
$\Box$

5)
$$
  \|\lambda_1-\lambda_2 \|_p \le C\|\mu_1 -\mu_2 \|_p^{\alpha} ,
$$  

{\bf Proof}. Since we have representation (3.3), we can see that
$$
\|\lambda_1-\lambda_2 \|_p \le C_p \left[ \|\mu_1 \circ \tilde f_1^{-1}-\mu_2 \circ \tilde f_1^{-1}\|_p +
 \|\mu_2 \circ\tilde f_1^{-1}-\mu_2 \circ\tilde f_2^{-1}\|_p +\right.
$$
$$
  \left.+\left\|\frac{\tilde f_{1,z}}{\overline{\tilde f_{1,z}}}\circ \tilde f_1^{-1}- 
 \frac{\tilde f_{2,z}}{\overline{\tilde f_{2,z}}}\circ \tilde f_1^{-1}\right\|_p+
 \left\|\frac{\tilde f_{2,z}}{\overline{\tilde f_{2,z}}}\circ \tilde f_1^{-1} -
 \frac{\tilde f_{2,z}}{\overline{\tilde f_{2,z}}}\circ \tilde f_2^{-1}\right\|_p \right] . \eqno (3.10)
$$
We obtain estimates for all terms in the right side by the same method as in the steps above. 

The first term we can write as
$$
\|\mu_1 \circ \tilde f_1^{-1}-\mu_2 \circ \tilde f_1^{-1}\|_p =\left(\int |\mu_1 -\mu_2 |^p J(\tilde f_1 )(z)dS_z \right)^{1/p}  ,
$$
Since $\mu_1 -\mu_2 \ne 0$ only on $D$, we can apply the estimate of
$J(\tilde f_1 )$ in $L^{p/2}(D)$. Also, from (3.4) follows the estimate $|J(\tilde f_1 )(\zeta )|\le C(1-|z|)^{-2}$
if $z\in D$. Fix some $r<1$. We obtain
$$
  \|\mu_1 \circ\tilde f_1^{-1}-\mu_2 \circ \tilde f_2^{-1}\|_p \le C (1-r)^{-2}\left(\int_{|z|\le r}|\mu_1 -\mu_2 |^p dS_z\right)^{1/p} +
C\left(\int_{1-r \le|z|\le 1}J(\tilde f_1 )(z)dS_z \right)^{1/p}\le 
$$  
$$
  \le C[(1-r)^{-2}\|\mu_1 -\mu_2 \|_p )+r^{1/p-2/p^2}]  .
$$  
We put $1-r=\|\mu_1 -\mu_2 \|_p^{\frac{p^2}{2p^2 +p-2}}$ and obtain 
$$
  \|\mu_1 \circ \tilde f_1^{-1}-\mu_2 \circ \tilde f_1^{-1}\|_p \le C\|\mu_1 -\mu_2 \|_p^{\frac{p-2}{2p^2 +p-2}} .
$$  
 
Consider the second term in (3.10). Remind that $|(\mu_2 )_{(k)}(z)|\le b_1 (1-|z|)^{-1}, |(k)|=1$.
Again fix some $r<1$. Applying step 2 and the third inequality of Proposition 15 we obtain
$$
  \|\mu_2 \circ \tilde f_1^{-1}-\mu_2 \circ \tilde f_2^{-1}\|_p \le  
$$
$$
  \le C\left((1-r)^{-p}\int_{r\le|\tilde f_1^{-1}(z)|\le 1}|\tilde f_1^{-1}(z)-\tilde f_2^{-1}(z)|^p dS_z \right)^{1/p} +
C[{\rm mes}\{|1-r\le |\tilde f_1^{-1}(z)|\le 1]\}]^{1/p} \le
$$
$$
  \le C[(1-r)^{-1}\|\mu_1-\mu_2 \|_p^{\alpha} +\left(\int_{1-r\le|z|\le 1\}}J(\tilde f_1 )(z)dS_z \right)^{1/p}\le .
$$
$$
  \le C[(1-r)^{-1}\|\mu_1 -\mu_2 \|_p^{\alpha}+(1-r)^{1/p-2/p^2}]
$$
Here $\alpha$ is the same as in step 2) and we again recall that $J(\tilde f_1 )$ belongs to $L^{p/2}$. 

We put $1-r=\|\mu_1 -\mu_2 \|_p^{\frac{\alpha p^2}{p^2 +p-2}}$ and obtain
$$
 \|\mu_2 \circ\tilde f_1^{-1}-\mu_2 \circ\tilde f_2^{-1}\|_p \le\|\mu_1 -\mu_2 \|_p^{\frac{\alpha (p-2)}{p^2 +p-2}} .
$$ 

Consider the third term in (3.10). We have
$$
  \frac{\tilde f_{1,z}}{\overline{\tilde f_{1,z}}}-\frac{\tilde f_{2,z}}{\overline{\tilde f_{2,z}}}=
\frac{\tilde f_{1,z}-\tilde f_{2,z}}{\overline{\tilde f_{1,z}}}-
\frac{\tilde f_{2,z}}{\overline{\tilde f_{2,z}}}\frac{\tilde f_{1,z}-\tilde f_{2,z}}{\overline{\tilde f_{1,z}}} .
$$
We see that we must estimate $|(\tilde f_{1,z}-\tilde f_{2,z})\circ\tilde f_1^{-1}||\tilde f_{1,z}\circ\tilde f_1^{-1}|^{-1}$. 

Applying left estimate (3.4) to $\tilde f_{1,z}$ and acting as above we see that
the third term in (3.10) is no greater than
$$
  C\left[(1-r)^{-1}\left(\int_{r\le|\tilde f_1^{-1}(z)]|\le 1}
|(\tilde f_{1,z}-\tilde f_{2,z})\circ\tilde f_1^{-1}(z)|^p dS_z \right)^{1/p} +
  (1-r)^{1/p-2/p^2}\right] . 
$$ 
for any $0<r<1$. Here we use the uniform boundedness of the third term in (3.12) and
again recall the estimate $|J(\tilde f_1 )(z)|\le C(1-|z|)^{-2}$.
We have
$$
   \int_{r\le|\tilde f_1^{-1}(z)]|\le 1}
|(\tilde f_{1,z}-\tilde f_{2,z})\circ\tilde f_1^{-1}(z)|^p dS_z=
$$
$$
  =\int_{r\le|z|\le 1}|\tilde f_{1,z}-\tilde f_{2,z}|^p |J(\tilde f_1 )dS_z
\le C(1-r)^{-2}\int_D |\tilde f_{1,z}-\tilde f_{2,z}|^p dS_{\zeta} .
$$
Applying step 3) we obtain for the third term the estimate
$$
  C(r^{-(1+2/p)}\|\mu_1 -\mu_2 \|_p^{\alpha} +r^{1/p-2/p^2}) , 
$$
where we set $\alpha$ as in step 3).
We put $r=\|\mu_1 -\mu_2 \|_p^{\frac{\alpha p^2}{p^2 +3p-2}}$ and obtain the estimate
$$
  \left\|\frac{\tilde f_{1,z}}{\overline{\tilde f_{1,z}}}\circ\tilde f_1^{-1}-\frac{\tilde f_{2,z}}{\overline{\tilde f_{2,z}}}\circ\tilde f_1^{-1}\right\|_p \le
C\|\mu_1 -\mu_2 \|_p ^{\frac{\alpha (p-2)}{p^2 +3p-2}}    .
$$ 

Consider the last term in (3.10). Denote $z_i =\tilde f_i^{-1}(z)$. We must estimate $p$-norm of the sum
$$
  \frac{\tilde f_{2,z}(z_1 )-\tilde f_{2,z}(z_2 )}{\overline{\tilde f_{2,z}}(z_1 )}+
\frac{\tilde f_{2,z}(z_2 )}{\overline{\tilde f_{2,z}}(z_2 )}\frac{\tilde f_{2,z}(z_2 )-\tilde f_{2,z}(z_1 )}{\overline{\tilde f_{2,z}}(z_1 )} .
$$
Analogously to the previous case we get for this $p$-norm the estimate
$$   
 C\left[r^{-1}\left(\int_{r\le|\tilde f_1^{-1}(z)]|\le 1}|\tilde f_{2,z}(z_1 )-\tilde f_{2,z}(z_2 )|^p dS_z\right)^{1/p} 
 +r^{1/p-2/p^2}\right]   
$$ 
for any $r>0$. 
Now at $|z_1 -z_2 |\le (1-r)/2$, applying estimate (3.6) and step 2), we have 
$$
  |\tilde f_{2,z}(z_1 )-\tilde f_{2,z}(z_2 )|\le
$$
$$
  \le\sup_{z\in [z_1 ,z_2 ],((k)|=1}||\tilde f_{2,z})_{(k)}(z)||z_1-z_2 |\le
 Cr^{-2}\|\mu_1 -\mu_2 \|_p^{\alpha}  .   \eqno (3.11)
$$   
Now we put  $1-r=\|\mu_1 -\mu_2\|_p^{\frac{\alpha p^2}{3p^2 +p-2}}$.   
In this case $\|\mu_1 -\mu_2\|_p^{\alpha}$ and,
hence, $|\tilde f_1^{-1}(z)-\tilde f_2^{-1}(z)|$ is small by comparison with $1-r$ and we,
indeed, can apply estimate (3.11). We obtain the estimate
$$
  \left\|\frac{\tilde f_{2,z}}{\overline{\tilde f_{2,z}}}(z_1 ) -
 \frac{\tilde f_{2,z}}{\overline{\tilde f_{2,z}}}(z_2 )\right\|_p \le
   C\|\mu_1 -\mu_2 \|_p^{\frac{\alpha (p-2)}{3p^2 +p-2}}  .  
$$
$\Box$

We obtained the estimate of the first difference in the right side of inequality (3.7).
We finish the proof of the proposition with estimation of the second difference.

6) {\it For $z\in D$}
$$
|f^0_{\lambda_2}\circ\tilde f_1 (z)-f^0_{\lambda_2}\circ\tilde f_2 (z)|\le \|\mu_1 -\mu_2 \|_p^{\alpha}
$$

{\bf Proof}. The estimate follows from step 2) and the well-known inequality
$|f_{\nu}(z_1 )-f_{\nu}(z_2 )|\le C|z_1 -z_2 |^{1-2/p}$, which holds for any compactly
supported $\nu$, $|\nu |\le d$ and $2< p<P(d)$ (see, for example [Al]).$\Box$
$\Box$

Now we begin the proof of estimate (3.2). Analogously to the notations $\mu_i ,f_i$
we adopt the notations $f_{wi},Z_{wi},h_{wi}$ and so on. At first we shall obtain estimates 
for the difference $\mu_{h_{w1}}-\mu_{h_{w2}}$ and its derivatives.

\begin{proposition} a) Suppose $|z|\le R<1$. There is the estimate
$$
  |(\mu_{h_{w1}}-\mu_{h_{w2}})_{(k)}(z)|\le 
$$
$$
  \le\frac{C}{(1-R)^{2|(k)|}}\left[\sup_{|(l)|\le k}\sup_{D_{R_w}}|(\mu_1 )_{(l)} -(\mu_2 )_{(l)}|
 +\frac{1}{(1-R)^2}|\mu_1 (w)-\mu_2 (w)|\right] ,  \eqno (3.12)
$$
where 
$$
  1-R_w^2 =b(1-|w|^2 )(1-R^2 )  \eqno (3.13)
$$ 
with some uniform $b$ independent of $R$. In fact, it depends only on the maximal dilatation of $f_i$ and we
can put $b=1/3$ if $\mu_1$, $\mu_2$ are small enough.

b) Let $p\ge 1$, $R<1$. Then
$$
  \|\mu_{h_{w1}}-\mu_{h_{w2}}\|_p \le C\left[\sup_{D_R}|\mu_1 (z)-\mu_2 (z)|+\left(\frac{1-R}{1-|w|}\right)^{1/p}+\right.
$$
$$  
  \left. +|\mu_1 (w)-\mu_2 (w)|^{1/(2p+1)}\right] .  \eqno (3.14)
$$
In all these inequalities $C$ is some uniform constant.
\end{proposition}

{\bf Proof}. a) Recalling (2.10) we see that any derivative $(\mu_{h_w})_{(k)}$
is a sum of items of the types 
$$
  P_l (\mu_{f_w} ,\bar\mu_0 ) \left(\frac{Z_{w,z}}{\overline{Z_{w,z}}}\right)_{(q)}(\mu_f)_{(l_1 )}^{k_1}...(\mu_f)_{(l_s )}^{k_s}\circ Z_w^{-1}
  (.)_{(i_1 )}^{j_1}...(.)_{(i_r )}^{j_r}   ,  \eqno (3.15)
$$
where $|(l_1 )| k_1 +...+|(l_s )| k_s \le |(k)|$, $|(q)|\le |(k)|$, $|(i_1 )|j_1+...+|(i_r |j_r \le |(k)|$,
$P_l$ is the derivative of order $l\le |(k)|$ of the fraction
$$
  \frac{\mu_{f_w}-\mu_0}{1-\bar\mu_0 \mu_{f_w}}
$$
considered as a function of $\mu_{f_w}$. That is $P_l$ is an uniformly bounded rational function of $\mu_f ,\bar\mu_0$,
and  $(.)$ can be either $Z_w^{-1}$ either $\overline{Z_w^{-1}}$. It follows that we can represent
any difference $(\mu_{h_{w1}})_{(k)}-(\mu_{h_{w2}})_{(k)}$ as a sum of terms of the types 
$$
  (P_l (\mu_{f_{w1}},\bar\mu_{1,0})-P_l (\mu_{f_{w2}},\bar\mu_{2,0})\circ Z_{w1}^{-1} [.],  \eqno (3.16)
$$ 
$$
  (P_l (\mu_{f_{w2}},\bar\mu_{2,0})\circ Z_{w1}^{-1}-P_l (\mu_{f_{w2}},\bar\mu_{2,0})\circ Z_{w2}^{-1}) [.]  , \eqno (3.17)
$$  
$$
  ((\mu_{f_{w1}})_{(l_i )}^{k_i}-(\mu_{f_{w2}})_{(l_i )}^{k_i})\circ Z_{w1}^{-1} [.],  \eqno (3.18)
$$
$$
  ((\mu_{f_{w2}})_{(l_i )}^{k_i}\circ Z_{w1}^{-1}-(\mu_{f_{w2}})_{(l_i )}^{k_i}\circ Z_{w2}^{-1} [.] ,  \eqno (3.19)
$$
$$
  \left(\left(\frac{Z_{w1,z}}{\overline{Z_{w1,z}}}\right)_{(l)}\circ Z_{w1}^{-1}-\left(\frac{Z_{w2,z}}{\overline{Z_{w2,z}}}\right)_{(l)}\circ Z_{w2}^{-1}\right) [.]
 \eqno (3.20)
$$     
$$
  ((Z_{w1}^{-1})_{(i_m )}^{j_m}- Z_{w2}^{-1})_{(i_m )}^{j_m}) [.]  ,  \eqno (3.21)
$$ 
where $[.]$ denotes each time a product of the multiples such as in (3.15) with omitted term corresponding to
the written difference. All these multiples are either derivatives of $Z_{wi}^{-1},i=1,2$, either derivatives of 
$\mu_{fi},i=1,2$ in the point $Z_{wi}^{-1}(z)$, either derivatives of $P_l$ with respect to $\mu_{f_i}$. The  
derivatives of $Z_{wi}^{-1}$ and of $P_l$ are uniformly bounded. 

$Z_{wi}(z)$ is a real analytic function of $\mu_{i,0}$
equal to $z$ identically at $\mu_{i,0}=0$ and, hence, we can estimate the differences in (3.20) and (3.21) as $C|\mu_{1,0}-\mu_{2,0}|$.  
From the other hand, the products in the square brackets in these cases have the estimates
$$
  C(1-|z|)^{-k_1 |(l_1 )|  +...+k_s |(l_s )|}(|1-wz|)^{-(k_1 |(l_1 )| +...+k_s |(l_s )|)}\le  C(1-|z|)^{-2|(k)|}  .
$$
It follows from estimate (2.5). We obtain for terms (3.20), (3.21) the estimate
$$
  \frac{C}{(1-R)^{2|(k)|}}|\mu_1 (w)-\mu_2 (w)|   \eqno (3.22) 
$$

The difference in (3.16) is no more than
$$
  C\sup_{D_R}|(\mu_{f_{w1}}-\mu_{f_{w2}})\circ Z_{w1}^{-1}|=C\sup_{D_R}|(\mu_1 -\mu_2 )\circ\vpf_w^{-1}\circ Z_{w1}^{-1}| .
$$ 

Now we have
$$
 1-|Z_{1,w}^{-1}(z)|^2 \ge C(1-|z|^2 )     \eqno (3.23)
$$
with $C$ depending only on $d$.
$$
  1-|\vpf_w^{-1}(z)|^2 =\frac{(1-|w|^2 )(1-|z|^2 )}{|1-\bar w z|^2}\ge \frac{1}{2}(1-|w|^2 )(1-|z|^2 ) \eqno (3.24)
$$
It implies
$$
  \sup_{D_R}|(\mu_{1,f_w}-\mu_{2,f_w})\circ Z_{1,w}^{-1}| \le C\sup_{D_{R_w}}|\mu_1 -\mu_2 |  ,
$$ 
where $R_w$ is such that
$$
  1-R_w^2 =b(1-|w|^2 )(1-R^2 )
$$
with some $b$ depending only $C$ in the right side of (3.23), i.e. , on maximal dilatations of $f_1$, $f_2$.
This constant $C$ can be made close to 1 if $\mu_1$, $\mu_2$ are small enough. In this case we can set
$b=1/3$. 

The product in the square brackets in (3.16) has the same estimate as in
(3.20), (3.21). We obtain for term (3.16) the estimate
$$
  \frac{C}{(1-R)^{2|(k)|}}\sup_{D_{R_w}}|\mu_1 -\mu_2 |  .  \eqno (3.25)   
$$

Consider now term (3.17). The function
$P_l$ has bounded derivatives as a function of $\mu_f$, and $Z_w^{-1}$ is an analytic function
of $\mu_0$ with a bounded derivative. Applying (2.5) and (3.23), we see that the difference 
in (3.17) is no more than
$$
  C\sup_{z\in D_R}\sup_{t\in [Z_{w1}^{-1}(z),Z_{w2}^{-1}(z)]}(|(\mu_{f_{w2}})_t (t)|+|(\mu_{f_{w2}})_{\bar t}(t)|)|\mu_{01}-\mu_{02}|\le
$$
$$
  \le \frac{C}{(1-R)^2}|\mu_1 (w)-\mu_2 (w)| .
$$
Since the product in the square brackets in (3.17) has the same estimate 
$C(1-|z|)^{-2k}$, we obtain for term (3.17) the estimate
$$
  \frac{C}{(1-R)^{2(k+1)}}|\mu_1 (w)-\mu_2 (w)|  .  \eqno (3.26)
$$

Considering term (3.18) we have
$$
|((\mu_{f_{w1}})_{(l)}^k-(\mu_{f_{w2}})_{(l)}^k )|\le C|p^{(k-1)}((\mu_{f_{w1}})_{(l)},(\mu_{f_{w2}})_{(l)})||\mu_{f_{w1}})_{(l)}-\mu_{f_{w2}})_{(l)}| ,
$$
where $p^{(k-1)}$ is the homogeneous polynomial in the variables $(\mu_{f_{w1}})_{(l)},(\mu_{f_{w2}})_{(l)}$ of degree $k-1$. 
Thus for the difference in (3.18) we have the estimate
$$
  \frac{C}{(1-R)^{2|(l_i )|(k_i -1)}}\sup_{D_R}|(\mu_{f_{w1}})_{(l_i )}-(\mu_{f_{w2}})_{(l_i )})\circ Z_{w1}^{-1}|  .
$$ 
Now remind  expression (2.7) for $\mu_{f_w}$. We have the obvious estimates
$$
  |(\vpf_w^{-1})_l (z)|\le C|1-z\bar w|^{-l}, |(r_w )_{(l)}(z)|\le C|1-z\bar w|^{-l} ,   \eqno (3.27)
$$
where we denote $r_w =\left(\frac{1-\bar z w}{1-z\bar w}\right)^2$.

We see that any derivative $(\mu_{f_w})_{(l)}$ is a sum of items of the types 
$$.
  \mu_{(m)}\circ\vpf_w^{-1}(.)_{(s_1 )}^{r_1}...(.)_{(s_q )}^{r_q}(r_w )_{(s)} ,
$$  
where $|(m)|\le |(l)|-|(s)|$, $s_1 r_1 +...+s_q r_q=|(l)|-|(s)|$, each $(.)$ can be either $\vpf_w^{-1}$ either
$\overline{\vpf_w^{-1}}$ and the corresponding derivative is either in $z$ either in $\bar z$. Applying (3.27), we see that we can represent
any difference $(\mu_{f_{w1}})_{(l)}(z)-(\mu_{f_{w2}})_{(l)}(z)$ as a sum of terms with the estimates
$$
  C|1-z\bar w|^{-|(l)|}\sup_{D_{R_w}}|(\mu_1 )_{(m)} -(\mu_2 )_{(m)}|,\,|(m)|\le |(l)|   .     
$$ 
Here we apply (3.24) and take into consideration that $s_1 r_1 +...+s_q r_q +|(s)|=|(l)|$. We obtain for the difference in (3.18) the estimate
$$
 C(1-R)^{-2|(l_i )|(k_i -1) -|(l_i )|}\sup_{m\le l_i}\sup_{D_{R_w}}|(\mu_1 )_{(m)} -(\mu_2 )_{(m)}| 
$$   

The product in the square brackets in (3.18) doesn't contain the term $(\mu_{f_{wj}})_{(l_i )}^{k_i}$ and, hence,
has the estimate: $C(1-R)^{-2(|(k)|+2k_i |(l_i |}$ . We obtain for term (3.18) the estimate
$$
  \frac{C}{(1-R)^{(2|(k)|-|(l_i )|}}\sup_{|(m)|\le |(k)|}\sup_{D_{R_w}}|(\mu_1 )_{(m)} -(\mu_2 )_{(m)}|   \eqno (3.28)
$$   

At last, the difference in (3.19) is no greater, than
$$
  C|p^{(k_i -1)}((\mu_{f_{w2}})_{(l_i )}\circ Z_{w1}^{-1},(\mu_{f_{w2}})_{(l_i )}\circ Z_{w2}^{-1})||(\mu_{f_{w2}})_{(l_i )}\circ Z_{w1}^{-1}-(\mu_{f_{w2}})_{(l_i )}\circ Z_{w2}^{-1}|,
$$
where $p^{(k_i -1)}$ is the homogeneous polynomial in the variables $(\mu_{f_{w1}})_{(l_i )}\circ Z_{w1}^{-1},(\mu_{f_{w2}})_{(l_i )}\circ Z_{w2}^{-1}$ of degree $k_i -1$.
Thus we can estimate this difference as
$$  
  \frac{C}{(1-R)^{2|(l_i )|(k_i -1)}}\sup_{D_R}|(\mu_{f_{w2}})_{(l_i )}\circ Z_{w1}^{-1}-(\mu_{f_{w2}})_{(l_i )}\circ Z_{w2}^{-1}|\le
$$
$$
  \le\frac{C}{(1-R)^{2|(l_i )|(k_i -1)}}\sup_{z\in D_R}\sup_{[Z_{w1}^{-1}(z),Z_{w2}^{-1}(z)]}(|((\mu_{f_{w2}})_{(l_i )})_t |+|((\mu_{f_{w2}})_{(l_i )})_{\bar t}||\mu_1 (w)-\mu_2 (w)|\le
$$
$$
  \le\frac{C}{(1-R)^{2|(l_i )|k_i -1)+2(|(l_i )|+1}}|\mu_1 (w)-\mu_2 (w)|=\frac{C}{(1-R)^{2(|(l_i )|k_i +1)}}|\mu_1 (w)-\mu_2 (w)| 
$$
Analogously to the previous case we obtain for term (3.19) the estimate
$$
  \frac{C}{(1-R)^{2(|(k)|+1)}}|\mu_1 (w)-\mu_2 (w)|   \eqno (3.29)
$$ 

From (3.22), (3.25), (3.26), (3.28) and (3.29) follows (3.12).

b) Applying (2.10), we have
$$
 \left(\int_D |\mu_{h_{w1}}(z)-\mu_{h_{w2}}(z)|^p dS_z \right)^{1/p} \le C\left[|\mu_{01}-\mu_{02}|+\right.
$$
$$
 +\left(\int_D |\mu_{f_{w1}}\circ Z_{w1}^{-1}-\mu_{f_{w2}}\circ Z_{w1}^{-1}|^p dS_z \right)^{1/p}+
\left(\int_D |\mu_{f_{w2}}\circ Z_{w1}^{-1}-\mu_{f_{w2}}\circ Z_{w2}^{-1}|^p dS_z \right)^{1/p}+
$$
$$
  +\left.\left(\int_D \left|\frac{Z_{w1,z}}{\overline{Z_{w1,z}}}\circ Z_{w1}^{-1}-\frac{Z_{w2,z}}{\overline{Z_{w2,z}}}\circ Z_{w2}^{-1}\right|^p dS_z \right)^{1/p}\right].  \eqno (3.30)
$$  
Consider the first integral. We have
$$
  J_1 =\int_D |\mu_{f_{w1}}\circ Z_{w1}^{-1}-\mu_{f_{w2}}\circ Z_{w1}^{-1}|^p dS_z \le
C\int_D (|\mu_{f_{w1}}(z)-\mu_{f_{w2}}(z)|)^p dS_z  ,   
$$
since the Jacobian of $Z_{w1}$ is uniformly bounded from above and from below. Also, applying (2.7), we have
$$
  \int_D |\mu_{f_{w1}}(z)-\mu_{f_{w2}}(z)|^p dS_z = \int_D |\mu_1 (z)-\mu_2 (z)|^p J_{\vpf_w}(z)dS_z  ,  \eqno (3.31)
$$
where $J_{\vpf_w}$ is the Jacobian of the map $\vpf_w$. This Jacobian is equal to $(1-|w|^2 )^2/|1-\bar w z|^4$ and there is the estimate
$$
  \int_{D\setminus D_R} |1-\bar w z|^{-4}dS_z \le 6\pi\frac{1-R^2}{(1-|w|^2 )^3}  .
$$
Indeed,
$$
  \int_{D\setminus D_R} |1-\bar w z|^{-4}dS_z =\int_R^1 \int_0^{2\pi}\frac{d\theta dr}{(1+|w|^2 r^2 -2|w|r \cos\theta )^2}=
$$
$$  
  =2\pi\int_R^1 \frac{r}{(1-|w|^2 r^2 )^2}\left(1+\frac{2}{1-|w|^2 r^2}\right)dr =
$$
$$
  =2\pi (1-R^2 )\left(\frac{1}{(1-|w|^2 )(1-R^2 |w|^2 )}+\frac{1}{(1-|w|^2 )(1-R^2 |w|^2 )^2}+\frac{1}{(1-|w|^2 )^2 (1-R^2 |w|^2 )}\right)\le
$$
$$
 \le 6\pi \frac{1-R^2}{(1-|w|^2 )^3} .
$$        
Further, the integral
$$
  \int_D J_{\vpf_w}(z)dS_z =(1-|w|^2 )^2 \int_D \frac{dS_z}{|1-\bar w z|^4}
$$
is uniformly bounded. Thus we obtain
$$
  (J_1 )^{1/p}\le C\sup_{D_R}|\mu_1 (z)-\mu_2 (z)|+
C\left[(1-|w|^2 )^2 \int_{D\setminus D_R}|1-\bar w z|^{-4} dS_z \right]^{1/p} \le 
$$
$$
  \le C\left[\sup_{D_R}|\mu_1 (z)-\mu_2 (z)|+\left(\frac{1-R}{1-|w|}\right)^{1/p}\right] \eqno (3.32) 
$$
for any $R<1$. From the other hand, if $|w|\le 1/2$, then from (3.31) follows the estimate
$$
  (J_1 )^{1/p}\le C\|\mu_1 -\mu_2 \|_p   .  \eqno (3.33)
$$ 

Consider the second integral in (3.30). Applying estimate (2.5), we have for any $r<1$
$$
  J_2 =\int_D |\mu_{f_{w2}}\circ Z_{w1}^{-1}-\mu_{f_{w1}}\circ Z_{w2}^{-1}|^p dS_z \le
$$
$$
  \le C\int_{D_r}\sup_{[Z_{w1}^{-1}(z),Z_{w2}^{-1}(z)]}(|(\mu_{f_{w2}})_{t}+(\mu_{f_{w2}})_{\bar t}||Z_{w1}^{-1}(z)- Z_{w2}^{-1}(z)|)^p dZ_z +C(1-r) \le
$$
$$
  \le C((1-r)^{-2p}|\mu_1 (w)-\mu_2 (w)|^p +(1-r) .
$$
That is,
$$
  (J_2 )^{1/p}\le C[(1-4)^{-2}|\mu_1 (w)-\mu_2 (w)|+(1-R_1 )^{1/p}] .  
$$   
If we put $(1-r)^{-2}|\mu_1 (w)-\mu_2 (w)|=(1-r)^{1/p}$, i.e., $1-r =|\mu_1 (w)-\mu_2 (w)|^{p/(2p+1)}$, then we obtain 
$$
  (J_2 )^{1/p}\le C|\mu_1 (w)-\mu_2 (w)|^{1/(2p+1)} .  \eqno (3.34)
$$

At last, we estimate the third integral in (3.30) analogously to the difference in (3.20): it is no greater than 
$$
   C|\mu_1 (w)-\mu_2 (w)|.                   \eqno (3.35)
$$
Collecting (3.30) and (3.33)-(3.35) we obtain (3.14).
$\Box$

We also need in an estimate for the norm of $\mu_{F_{w1}}-\mu_{F_{w2}}$ in $L_H^p$. Remind that
$H$ is the stripe $-\pi \le \vpf \le \pi$ and $L_H^p$ is the space of $\vpf$-periodic functions on $H$ with usual
$L^p$-norm. In what follows we don't write the index $H$ in our notations.

\begin{proposition} Let $p\ge 2$, $R<1$. There is the estimate
$$
  \|\mu_{F_{w1}}-\mu_{F_{w2}}\|_p \le C\left[\sup_{D_{r_w}}(|\mu_{1,z} -\mu_{2,z}|+|\mu_{1,\bar z} -\mu_{2,\bar z}|)+\right.
$$
$$
  \left.+\sup_{D_R}|\mu_1 (z)-\mu_2 (z)|+\left(\frac{1-R}{1-|w|}\right)^{1/p}+|\mu_1 (w)-\mu_2 (w)|^{1/(2p+1)}\right] , 
\eqno (3.36)
$$
where $1-r_w \ge b(1-|w|)$ with some $b<1$ depending only on $d$. If $\mu_1$, $\mu_2$
are small enough, then we can take $b=1/4$.  
\end{proposition}

{\bf Proof}. 
Remind equation (2.17). Since $\mu_{h_{wi}}(0)=0$, $i=1,2$, $p>2$ and the first derivatives of $\mu_{h_{wi}}$ are 
uniformly bounded on $D_{1/2}$ (see (2.12)), we have
$$
  \int_H |\mu_{F_{w1}} (\zeta ) -\mu_{F_{w2}} (\zeta ) |^p dS_{\zeta} \le\int_D\frac{|\mu_{1,h_w}(z)-\mu_{2,h_w}(z)|^p}{|z|^2}dS_z \le
$$
$$
\le \int_{|z|\le 1/2} \sup_{|t|\le |z|}(|\mu_{h_{w1},t}-\mu_{h_{w2},t}|+|\mu_{h_{w1},\bar t}-\mu_{h_{w2},\bar t}|)^p dS_z +
$$
$$
 +\frac{1}{4}\int_{|z|\ge 1/2}|\mu_{h_{w1}}(z)-\mu_{h_{w2}}(z)|^p dS_z .  \eqno (3.37)
$$
The first integral in the right side has the obvious estimate
$$
 C\sup_{|z|\le 1/2}(|\mu_{h_{w1},z}-\mu_{h_{w2},z}| +|\mu_{h_{w1},\bar z}-\mu_{h_{w2},\bar z}|  . 
$$
Applying estimate (3.12) we see that the integral is no greater than
$$
  C[\sup_{D_{r_w}}(|\mu_{1,z}-\mu_{2,z}|+|\mu_{1,\bar z}-\mu_{2,\bar z}|)+|\mu_1 (w)-\mu_2 (w)|]  .
$$
We apply estimates (3.14) to the second integral in (3.37). As a result we obtain estimate (3.36).
$\Box$
 
Our next step will be the proof of estimate (3.2) for the derivatives of first order. 
It is a consequence of estimate (3.36) and the next proposition.

\begin{proposition} Let $f_i$, $i=1,2$ be the normal solutions corresponding to the Beltrami coefficients $\mu_1$ and $\mu_2$
satisfying conditions of Theorem 1 with the same bounds $b,b_1 ,b_2$.
There is the estimate
$$
  |f_{1,z}(w)-f_{2,z}(w)|\le C\min\left[1,\frac{\|\mu_1 -\mu_2 \|_p^{\alpha}}{1-|w|}+ 
|\mu_1 (w)-\mu_2 (w)|^{\alpha}+\|\mu_{F_{1w}} -\mu_{F_{2w}}\|_p^{\alpha}\right]
$$
with some uniform $C$ and some $2<p<P(d)$, $\alpha <1$ depending only on $d,b_1, b_2$.
For the difference of $\bar z$-derivatives we have the analogous estimate.
\end{proposition}

{\bf Proof}.
By (2.15) and (2.28),
$$
  f_z (w)=\frac{1-|w'|^2}{1-|w|^2}(h_w )_z (0)(Z_w )_z (0)=\frac{\vpf_{w'}(1)}{\tilde h_w \circ Z_w^{-1}\circ\vpf_w (1)}\frac{1-|w'|^2}{1-|w|^2}(Z_w )_z (0)\tilde (h_w )_z (0) .
\eqno (3.38)
$$  
Remind that $w'=f(w)$ and $\tilde h_w =\exp\circ\tilde F_w \circ\log$, where 
$\tilde F_w =f_{\lambda_{F_w}}\circ\tilde f_{\mu_{F_w}}$ is defined in subsection 5)
of Section 2.

We can see that all multiples in (3.40) are uniformly bounded. For the multiple 
$\frac{\vpf_{w'}(1)}{\tilde h_w \circ Z_w^{-1}\circ\vpf_w (1)}$ it follows from the notion that $\vpf_w (1)=\frac{1-w}{1-\bar w}$
belongs to the unit circle and $|\tilde h_w \circ Z_w^{-1}|$ is uniformly bounded from below on the unit circle. 
According to (3.38) we can represent the difference $f_{1,z}(w)-f_{2,z}(w)$ as a sum of the terms
$$
  \frac{|w'_1 |^2-|w'_2 |^2}{1-|w|^2}O(1) ,  \eqno (3.39) 
$$
$$
  [(Z_{w1}^{-1})_z (0)-(Z_{w2})_z (0)]O(1) ,  \eqno (3.40)
$$
$$
  (\vpf_{w'_1}(1)-\vpf_{w'_2}(1))O(1) ,  \eqno (3.41)
$$
$$
  [\tilde h_{w1}\circ Z_{w1}^{-1}\circ\vpf_w (1)-\tilde h_{1w}\circ Z_{w2}^{-1}\circ\vpf_w (1)]O(1) ,  \eqno (3.42)
$$
$$
  [\tilde h_{1w}\circ Z_{w2}^{-1}\circ\vpf_w (1))^{-1}-\tilde h_{2w}\circ Z_{w2}^{-1}\circ\vpf_w (1)]O(1) ,  \eqno (3.43)
$$
$$
 [(\tilde h_{w1})_z (0)-(\tilde h_{w2})_z (0)]O(1) .  \eqno (3.44)
$$        
In fact, since the fractions $(1-|w'_i |)/(1-|w|)$ are uniformly bounded, it is enough, instead of (3.41), to estimate
$$
  C\min\left[1, \frac{|w'_1 -w'_2 |}{1-|w|^2}\right]  .  
$$
The difference $|w'_1 -w'_2 |=f_1 (w) -f_2 (w)$ we estimate by inequality (3.1). 

We can estimate term (3.40) as
$$
  C|\mu_1 (w)-\mu_2 (w)| .
$$
For term (3.41) we have the estimate
$$
  C|w'_1 -w'_2 |  
$$

To estimate term (3.42) we use the next proposition:

\begin{proposition} Let $F_{\nu}$ be the principal logarithmic solution corresponding to the
coefficient $\nu$, $|\nu |\le d<1$. Then
$$
  |F_{\nu}(\zeta_1 )-F_{\nu}(\zeta_2 )|\le C|\zeta_1 -\zeta_2 |^{1-2/p}
$$
for $p\le P(d)$.
\end{proposition}

The proof is analogous to the proof for the principal solution in the classical case.
$\Box$ 

Now, since $\tilde h_w =\exp\circ f_{\lambda_{F_w}}\circ\tilde f_{\mu_{F_w}}\circ\log$
and $|Z_{w1}^{-1}(z)-Z_{w2}^{-1}(z)|\le C|\mu_1 (w)-\mu_2 (w)|$, we obtain for term (3.42)
the estimate
$$ 
   C|\mu_1 (w)-\mu_2 (w)|^{(1-2/p)^2}  .
$$

For term (3.43) we have the estimate
$$
  C\max_{\xi =0}|\tilde F_{w1}(\zeta )-\tilde F_{w2}(\zeta )|
$$
and for term (3.44)
$$
  C\lim_{\xi \to -\infty}|\tilde F_{w1}(\zeta )-\tilde F_{w2}(\zeta  )| .
$$

We see that to prove Proposition 18 it remains to estimate the difference
$$
  \tilde F_{w1}-\tilde F_{w2} = (f_{\lambda_{F_{1w}}}\circ\tilde f_{\mu_{F_{1w}}}-f_{\lambda_{F_{2w}}}\circ\tilde f_{\mu_{F_{1w}}})+
 (f_{\lambda_{F_{2w}}}\circ\tilde f_{\mu_{F_{1w}}}-f_{\lambda_{F_{2w}}}\circ\tilde f_{\mu_{F_{2w}}}) .  \eqno (3.45)
$$
We proceed analogously to the proof of estimate (3.1) but we use now
the principal logarithmic solutions. Instead of Proposition 15 we have 

\begin{proposition}
Let $F_{\nu_1}$ and $F_{\nu_2}$ be the principal logarithmic solutions corresponding to the
coefficients $\nu_1$ and $\nu_2$, $|\nu_i |\le d<1$. Then, for $\zeta\in H$,
$$
  |F_{\nu_1}(\zeta )-F_{\nu_2}(\zeta )|\le C\|(\nu_1 -\nu_2 )(|(F_{\nu_1})_{\zeta}-1|+|(F_{\nu_2})_{\zeta}-1|\|_p  ,  
$$
$$
  \|(F_{\nu_1})_{\zeta}-(F_{\nu_2})_{\zeta}\|_p \le C\|(\nu_1 -\nu_2 )(|(F_{\nu_1})_{\zeta}-1|+|(F_{\nu_2})_{\zeta}-1|\|_p   , 
$$
$$
  |F^{-1}_{\nu_1}(\zeta )-F^{-1}_{\nu_2}(\zeta )|\le C\|(\nu_1 -\nu_2 )(|(F_{\nu_1})_{\zeta}-1|+|(F_{\nu_2})_{\zeta}-1|\|_p  
$$
for $2\le p\le P(d)$. 
\end{proposition}

{\bf Proof}. The proof of is completely analogous to the proof of
Proposition 15, we only change $\mC$ and $\mS$ to $P_h$ and $T_H$ correspondingly.
$\Box$

\begin{proposition} Suppose $|\mu_{F_{wi}}|\le d$, $|\lambda_{F_{wi}}|\le d$, $i=1,2$. Then
$$
  |f_{\lambda_{F_{w1}}}-f_{\lambda_{F_{w2}}}|\le C\|\lambda_{F_{w1}}-\lambda_{F_{2w}}\|_p^{\alpha}   ,
$$
$$
  \|(f_{\mu_{F_{w1}}})_{\zeta}-(f_{\mu_{F_{w2}}})_{\zeta}\|_p \le C\|\mu_{F_{w1}}-\mu_{F_{w2}}\|_p^{\alpha} ,
$$
$$
  |\tilde f_{\mu_{F_{w1}}}-\tilde f_{\mu_{F_{w2}}}|\le C\|\mu_{F_{w1}}-\mu_{F_{w2}}\|_p^{\alpha} ,
$$
$$
  |\tilde f^{-1}_{\mu_{F_{w1}}}-\tilde f^{-1}_{\mu_{F_{w2}}}|\le C\|\mu_{F_{w1}}-\mu_{F_{w2}}\|_p^{\alpha} ,
$$
where $2<p<P(d)$ and $\alpha$ depend only on $d$, and $C$ is uniform.
\end{proposition}

{\bf Proof}. Consider the first estimate. 
By Proposition 20 we must estimate $\|(\lambda_{F_{w1}}-\lambda_{F_{w2}})|(f_{\lambda_{F_{wi}}})_{\zeta}-1|\|_p$.
But $(f_{\lambda_{F_{wi}}})_{\zeta}-1$ belongs to $L^p$ with the norm bounded by a constant depending only
on the maximal dilatation of $\lambda_{F_{wi}}$ if $p<P(d)$. Also, we have estimate (2.33). Fix some $r<1$ and some $p'>p, p'<P(d)$. We have
$$
  \|(\lambda_{F_{w1}}-\lambda_{F_{w2}})|(f_{\lambda_{F_{wi}}})_{\zeta}-1|\|_p \le 
  C\left[r^{-3}\left(\int_{|{\rm Re}[\tilde f^{-1}_{\mu_{F_{wi}}}(\zeta )]|\ge r}|\lambda_{F_{w1}}-\lambda_{F_{w2}}|^p dS_{\zeta}\right)^{1/p}+\right.
$$
$$
  +\left.C(\|(f_{\lambda_{F_{wi}}})_{\zeta}-1\|_{p'}[{\rm mes}\{H\cap\{|{\rm Re}[\tilde f^{-1}_{\mu_{F_{w1}}}(\zeta )]|\le r\}\}]^{\frac{1}{p}(1-p/p')}\right]\le
$$
$$  
  \le C\left[ r^{-3}\|\lambda_{F_{w1}}-\lambda_{F_{w2}}\|_p +\left(\int_{H\cap\{|\xi |\le r\}}J(\tilde f_{\mu_{F_{1w}}},\zeta )dS_{\zeta}\right)^{\frac{p'-p}{pp'}}\right] .
$$
Since the restriction of the Jacobian $J(\tilde f_{\mu_{F_{wi}}})$ on the domain
$H\cap\{|\xi |\le 1\}$ belongs to $L^{p/2}$ with an uniform norm, we obtain for the last integral the estimate
$$
  \left(\int_H J^{p/2}\right)^{2/p}[{\rm mes}\{H\cap\{-r\le\xi\le 0\}\}]^{1-2/p}\le Cr^{1-2/p} . 
$$  
Now we put $r^{-3}\|\lambda_{F_{w1}}-\lambda_{F_{w2}}\|_p =r^\frac{(p-2)(p'-p)}{p'p^2}$, i.e,
$r=(\|\lambda_{F_{w1}}-\lambda_{F_{w2}}\|_p )^{\frac{p'p^2}{3p'p^2+(p-2)(p'-p)}}$. We obtain the estimate
$$ 
  \|(\lambda_{F_{w1}}-\lambda_{F_{w2}})|(f_{\lambda_{F_{w1}}})_{\zeta}|\|_p 
  \le C(\|\lambda_{F_{w1}}-\lambda_{F_{w2}}\|_p )^{\frac{(p-2)(p'-p)}{3p'p^2+(p-2)(p'-p)}} .
$$  

Analogously, to obtain the second estimate of the proposition we proceed applying (2.34) and the second inequality
of Proposition 20 
$$
  \|(f_{\mu_{F_{w1}}})_{\zeta}-(f_{\mu_{F_{w2}}})_{\zeta}\|_p 
  \le C\left[r^{-3}\left(\int_{|\xi |\le r}|\mu_{F_{w1}}-\mu_{F_{w2}}|^p dS_{\zeta}\right)^{1/p}+\right.
$$
$$
  \left.+\|(f_{\mu_{F_{w1}}})_{\zeta}-1\|_{p'}r^{\frac{1}{p}(1-p/p')}\right]\le C(r^{-3}\|\mu_{F_{w1}}-\mu_{F_{w2}}\|_p +r^{\frac{p'-p}{pp'}})
$$ 
for some $r<1$ and some $p'>p, p'<P(d)$. Putting $r=\|\mu_{F_{w1}}-\mu_{F_{w2}}\|_p^{\frac{pp'}{3pp'+p'-p}}$ we obtain
$$
  \|(f_{\mu_{F_{w1}}})_{\zeta}-(f_{\mu_{F_{w2}}})_{\zeta}\|_p \le C\|\mu_{F_{w1}}-\mu_{F_{w2}}\|_p^{\frac{p'-p}{3pp'+p'-p}} .
$$  
The prove of the last two estimates of the proposition is analogous. 
$\Box$

Return now to the proof of Proposition 18. To estimate the first difference in (3.45) we apply Propositions 20 and 21. 
We must estimate the norm $\|\lambda_{F_{w1}}-\lambda_{F_{w2}}\|_p$. We have
$$
  \|\lambda_{F_{w1}}-\lambda_{F_{w2}}\|_p \le C\left[\|\mu_{F_{w1}}\circ \tilde f^{-1}_{\mu_{F_{w1}}}-\mu_{F_{w2}}\circ \tilde f^{-1}_{\mu_{F_{w1}}}\|_p +
 \|\mu_{F_{w2}}\circ \tilde f^{-1}_{\mu_{F_{w1}}}-\mu_{F_{w2}}\circ \tilde f^{-1}_{\mu_{F_{w2}}}\|_p +\right.
$$
$$
  \left.+\left\|\frac{\overline{f_{\mu_{F_{w1}},\zeta}}}{f_{\mu_{F_{w1}},\zeta}}\circ (-\overline{\tilde f^{-1}_{\mu_{F_{w1}}}})- 
 \frac{\overline{f_{\mu_{F_{w2}},\zeta}}}{f_{\mu_{F_{w2}},\zeta}}\circ (-\overline{\tilde f^{-1}_{\mu_{F_{w1}}}})\right\|_p+
 \left\|\frac{\overline{f_{\mu_{F_{w2}},\zeta}}}{f_{\mu_{F_{w2}},\zeta}}\circ (-\overline{\tilde f^{-1}_{\mu_{F_{w1}}}})- 
 \frac{\overline{f_{\mu_{F_{w2}},\zeta}}}{f_{\mu_{F_{w2}},\zeta}}\circ (-\overline{\tilde f^{-1}_{\mu_{F_{w2}}}})\right\|_p\right] .
\eqno (3.46)
$$

Consider the first difference in the right side. Since $\exp\circ\tilde f_{\mu_{F_{1w}}}\circ\log$
is a homeomorphism of the plane, we can write
$$
  \|\mu_{F_{w1}}\circ \tilde f^{-1}_{\mu_{F_{w1}}}-\mu_{F_{w2}}\circ \tilde f^{-1}_{\mu_{F_{w1}}}\|_p^p =
 \int_{\tilde f^{-1}_{\mu_{F_{w1}}}(H)}|\mu_{F_{w1}}(\zeta )-\mu_{F_{w2}}(\zeta )|^p J(\tilde f_{\mu_{F_{w1}}},\zeta )dS_{\zeta}=
$$
$$
 =\int_{\bbC}|\mu_{F_{w1}}(\log z)-\mu_{F_{w2}}(\log z)|^p J(\tilde f_{\mu_{F_{w1}}},\log z)|z|^{-2}dS_z =
 \int_H |\mu_{F_{w1}}(\zeta )-\mu_{F_{w2}}(\zeta )|^p J(\tilde f_{\mu_{F_{w1}}},\zeta )dS_{\zeta} 
$$
From (2.32) follows the estimate 
$$
 |J(\tilde f_{\mu_{F_{w1}}},\zeta)|\le C(|\xi |^{-6}+1) .  \eqno (3.47)
$$ 
Also, restriction of $J(\tilde f_{\mu_{F_{w1}}})$ on the domain $H\cap\{|\xi |\le 1\}$ belongs to $L^{p/2}$ with an uniform estimate. 
We apply our usual method. We have for any $r<1$
$$
  \|\mu_{F_{w1}}\circ \tilde f^{-1}_{\mu_{F_{w1}}}-\mu_{F_{w2}}\circ \tilde f^{-1}_{\mu_{F_{w1}}}\|_p^p \le
C\left[r^{-6}\int_{H\cap\{|\xi |\ge r\}}|\mu_{F_{w1}}(\zeta )-\mu_{F_{w2}}(\zeta )|^p dS_{\zeta}+\right.
$$
$$
 \left.+\|\mu_{F_{w1}}-\mu_{F_{w2}}\|_C^p ({\rm mes}\{H\cap\{|\xi |\le r\}\})^{1-2/p}\right].
$$
We determine $r$ from the  equation $r^{-6}\|\mu_{F_{w1}}-\mu_{F_{w2}}\|_p^p =r^{1-2/p}$. I.e.,
we put $r=\|\mu_{F_{w1}}-\mu_{F_{w2}}\|_p^{\frac{p^2}{7p-2}}$ and obtain 
$$
  \|\mu_{F_{w1}}\circ \tilde f^{-1}_{\mu_{F_{w1}}}-\mu_{F_{w2}}\circ \tilde f^{-1}_{\mu_{F_{w2}}}\|_p^p \le
C\|\mu_{F_{w1}}-\mu_{F_{w2}}\|_p^{\frac{p(p-2)}{7p-2}}  .  \eqno (3.48)
$$

Consider the second term in the right side of (3.46). We can assume that $\|\mu_{F_{w1}}-\mu_{F_{w2}}\|_p$ is small
with some uniform estimate that we shall specify below.
Applying (2.11), (2.17) and the obvious estimate for $Z=Z_w^{-1}(z)$: $c(1-|z|)\le 1-|Z|\le C(1-|z|)$ with uniform $c,C$,  
we get the estimate $|(\mu_{F_{wi}})_{(l)}(\zeta )|\le Ce^{\xi}(1+|\xi |^{-2})$ if $|(l)|=1$, $i=1,2$. 
By the last inequality of Proposition 21,
$$
  |\tilde f^{-1}_{\mu_{F_{w1}}}(\zeta )-\tilde f^{-1}_{\mu_{F_{w2}}}(\zeta )|\le C\|\mu_{F_{w1}}-\mu_{F_{w2}}\|_p^{\alpha}  . 
$$
Suppose 
$$
  2C\|\mu_{F_{w1}}-\mu_{F_{w2}}\|_p^{\alpha}\le r<1 .   \eqno (3.49)
$$
Applying estimate (3.47), we can write 
$$
  \int_{{\rm Re}[\tilde f^{-1}_{\mu_{F_{wi}}}(\zeta )]|\le -r} |\mu_{F_{w2}}\circ \tilde f^{-1}_{\mu_{F_{w1}}}(\zeta )-\mu_{F_{w2}}\circ \tilde f^{-1}_{\mu_{F_{w2}}}(\zeta )|^p dS_{\zeta} \le
$$
$$
 \le C\int_{H\cap\{\xi\le -r\}}e^{p\xi}(1+|\xi |^{-2})^p (|\xi |^{-6}+1)|\tilde f^{-1}_{\mu_{F_{w1}}}(\zeta )-\tilde f^{-1}_{\mu_{F_{w2}}}(\zeta )|^p dS_{\zeta} \le
   Cr^{-2p-5}\|\mu_{F_{w1}}-\mu_{F_{w2}}\|_p^{p\alpha} , 
$$
where $\alpha$ is the same as in Proposition 21. If $r$ satisfies conditions (3.49),
we have 
$$
  \|\mu_{F_{w2}}\circ \tilde f^{-1}_{\mu_{F_{w1}}}-\mu_{F_{w2}}\circ \tilde f^{-1}_{\mu_{F_{w2}}}\|_p \le
$$
$$
  \le C[r^{-2-5/p}\|\mu_{F_{w1}}-\mu_{F_{w2}}\|_p^{\alpha}+C[{\rm mes}\{H\cap\{|{\rm Re}[\tilde f^{-1}_{\mu_{F_{w1}}}(\zeta )]|\le 3r/2 \}\}]^{1/p} .
$$
Since the restriction of the Jacobian $J(\tilde f_{\mu_{F_{1w}}})$ on the domain
$H\cap\{|\xi |\le 1\}$ belongs to $L^{p/2}$, we estimate the last integral as
$$
  \left(\int_H J^{p/2}\right)^{2/p}[{\rm mes}\{H\cap\{-r\le\xi\le 0\}\}]^{1-2/p}\le Cr^{1-2/p} . 
$$  
Now we put $r^{-2-5/p}(\|\mu_{F_{w1}}-\mu_{F_{w2}}\|_p^{\alpha} =r^{1/p-2/p^2}$, i.e., 
$r=\|\mu_{F_{w1}}-\mu_{F_{w2}}\|_p^{\alpha\frac{p^2}{2p^2 +6p-2}}$. We see that conditions
(3.49) are satisfied at $\|\mu_{F_{w1}}-\mu_{F_{w2}}\|_p$ small enough.
We obtain
$$
  \|\mu_{F_{w2}}\circ \tilde f^{-1}_{\mu_{F_{w1}}}-\mu_{F_{w2}}\circ \tilde f^{-1}_{\mu_{F_{w2}}}\|_p \le
C\|\mu_{F_{w1}}-\mu_{F_{2w}}\|_p^{\beta\frac{p-2}{2p^2 +6p-2}}  .  \eqno (3.50)
$$ 

Consider the third term in (3.46). We have
$$
  \frac{\overline{f_{\mu_{F_{w1}},\zeta}}}{f_{\mu_{F_{w1}},\zeta}}-\frac{\overline{f_{\mu_{F_{w2}},\zeta}}}{f_{\mu_{F_{w2}},\zeta}}=
\frac{\overline{f_{\mu_{F_{w1}},\zeta}}-\overline{f_{\mu_{F_{w2}},\zeta}}}{f_{\mu_{F_{w1}},\zeta}}-
\frac{\overline{f_{\mu_{F_{w2}},\zeta}}}{f_{\mu_{F_{w2}},\zeta}}\frac{f_{\mu_{F_{w1},\zeta}}-f_{\mu_{F_{w2}},\zeta}}{f_{\mu_{F_{w1}},\zeta}} .
$$
For $f_{\mu_{F_{1w}},\zeta}$ we have estimate (2.34). Acting as above, we can see that
the third term in (3.48) is no greater than
$$
  C\left[r^{-1}\left(\int_{H\cap\{|{\rm Re}[\tilde f^{-1}_{\mu_{F_{w1}}}(\zeta )]|\ge r\}}
|(f_{\mu_{F_{w1},\zeta}}-f_{\mu_{F_{w2},\zeta}})\circ (-\overline{\tilde f^{-1}_{\mu_{F_{w1}}}})|^p dS_{\zeta}\right)^{1/p} +
  r^{1/p-2/p^2}\right] . 
$$ 
for any $0<r<1$. 
Applying estimate (3.47), we obtain
$$
\int_{H\cap\{{\rm Re}[\tilde f^{-1}_{\mu_{F_{wi}}}(\zeta )]|\ge r\}}
|(f_{\mu_{F_{w1},\zeta}}-f_{\mu_{F_{w2},\zeta}})\circ (-\overline{\tilde f^{-1}_{\mu_{F_{w1}}}}|^p dS_{\zeta}=
$$
$$
  =\int_{H\cap\{|\xi |\ge r\}}|f_{\mu_{F_{w1},\zeta}}-f_{\mu_{F_{w2},\zeta}}|^p |J(\tilde f_{\mu_{F_{w1}}},\zeta )|dS_{\zeta}
\le Cr^{-6}\int_H |f_{\mu_{F_{1w},\zeta}}-f_{\mu_{F_{2w},\zeta}}|^p dS_{\zeta} .
$$
Applying the second inequality of Proposition 21, we see that the third term in (3.46) satisfies the estimate
$$
  C(r^{-(1+6/p)}\|\mu_{F_{w1}}-\mu_{F_{w2}}\|_p^{\alpha}+r^{1/p-2/p^2}) , 
$$
where $\alpha$ is the same as in Proposition 21.
If we put $r^{-(1+6/p)}\|\mu_{F_{w1}}-\mu_{F_{w2}}\|_p^{\alpha}=r^{1/p-2/p^2}$, i.e., 
$r=\|\mu_{F_{w1}}-\mu_{F_{w2}}\|_p ^{\alpha\frac{p^2}{p^2 +7p-2}}$, we obtain the estimate
$$
  \left\|\frac{\overline{f_{\mu_{F_{w1}},\zeta}}}{f_{\mu_{F_{w1}},\zeta}}\circ (-\overline{\tilde f^{-1}_{\mu_{F_{w1}}}})- 
 \frac{\overline{f_{\mu_{F_{w2}},\zeta}}}{f_{\mu_{F_{w2}},\zeta}}\circ (-\overline{\tilde f^{-1}_{\mu_{F_{w1}}}})\right\|_p \le
C\|\mu_{F_{w1}}-\mu_{F_{w2}}\|_p^{\beta\frac{p-2}{p^2 +7p-2}}   \eqno (3.51)  .
$$ 

Consider the last term in (3.46). Denote $\zeta_i =-\overline{\tilde f^{-1}_{\mu_{F_{wi}}}(\zeta )}$, $i=1,2$. We must estimate 
$p$-norm of the sum
$$
  \frac{\overline{f_{\mu_{F_{w2}},\zeta}}(\zeta_1 )-\overline{f_{\mu_{F_{w2}},\zeta}}\circ (\zeta_2 )}
{f_{\mu_{F_{w2}},\zeta}(\zeta_1 )} +
\frac{\overline{f_{\mu_{F_{w2}},\zeta}}\circ (\zeta_2 )}{f_{\mu_{F_{w2}},\zeta}(\zeta_2 )}
\frac{f_{\mu_{F_{w2}},\zeta}(\zeta_2 )-f_{\mu_{F_{w2}},\zeta}(\zeta_1 )}
{f_{\mu_{F_{w2}},\zeta}(\zeta_1 )} .   
$$
Let $r$ be such that
$$
  |\zeta_1 -\zeta_2 |\le r/2  .  \eqno (3.52)
$$  
Analogously to the previous case we obtain for this $p$-norm the estimate
$$   
 C\left[r^{-1}\left(\int_{|{\rm Re}[\tilde f^{-1}_{\mu_{F_{w1}}}(\zeta )]|\ge r}|f_{\mu_{F_{w2}},\zeta}(\zeta_1 )-f_{\mu_{F_{w2}},\zeta}(\zeta_2 )|^p dS_{\zeta}\right)^{1/p} 
 +r^{1/p-2/p^2}\right]  . \eqno (3.53)
$$ 
Now, by Proposition 21 and the estimate of Proposition 13,
$$
  |f_{\mu_{F_{w2}},\zeta}(\zeta_1 )-f_{\mu_{F_{w2}},\zeta}(\zeta_2 )|\le
$$
$$
  \le\sup_{\zeta\in [\zeta_1 ,\zeta_2 ],|(l)|=2}|(f_{\mu_{F_{w2}}}(\zeta ))_{(l)}||\zeta_1-\zeta_2 |\le
 Ce^{\xi}(1+r^{-6})\|\mu_{F_{w1}}-\mu_{F_{w2}}\|_p^{\alpha}   .
$$   
Applying (3.47), we obtain for the integral in (3.53) the estimate
$$
  Cr^{-7-6/p}\|\mu_{F_{w1}}-\mu_{F_{w2}}\|_p^{\alpha} .
$$  
Now we determine $r$ from the equation  $r^{-7-6/p}\|\mu_{F_{w1}}-\mu_{F_{w2}}\|_p^{\alpha} =r^{1/p-2/p^2}$, that is, 
$r=\|\mu_{F_{w1}}-\mu_{F_{w2}}\|_p^{\alpha\frac{p^2}{7p^2 +7p-2}}$.   
In this case $\|\mu_{F_{1w}}-\mu_{F_{w2}}\|_p$ is small by comparison with $r$ and
condition (3.52) is satisfied. We obtain the estimate
$$
  \left\|\frac{\overline{f_{\mu_{F_{w2}},\zeta}}}{f_{\mu_{F_{w2}},\zeta}}\circ (-\overline{\tilde f^{-1}_{\mu_{F_{w1}}}})- 
 \frac{\overline{f_{\mu_{F_{w2}},\zeta}}}{f_{\mu_{F_{w2}},\zeta}}\circ (-\overline{\tilde f^{-1}_{\mu_{F_{w2}}}})\right\|_p \le
 C\|\mu_{F_{w1}}-\mu_{F_{w2}}\|_p^{\alpha\frac{p-2}{7p^2 +7p-2}}  .  \eqno (3.54)
$$ 
Collecting (3.48), (3.50), (3.51), and (3.54), we obtain the estimate for the first difference in the
right side of (3.45).

It remains to estimate the second difference in (3.45). From Propositions 19 and 21 follows
$$
 |f_{\lambda_{F_{w2}}}\circ\tilde f_{\mu_{F_{w1}}}(\zeta )-f_{\lambda_{F_{w2}}}\circ\tilde f_{\mu_{F_{w2}}}(\zeta )|\le
$$
$$
  \le C|\tilde f_{\mu_{F_{w1}}}(\zeta )-\tilde f_{\mu_{F_{w2}}}(\zeta )|^{1-2/p}\le C(\|\mu_{F_{w1}}-\mu_{F_{w2}}\|_p^{\alpha (1-2/p)}  .
$$  
It finishes the proof of Proposition 18. $\Box$

Now we begin the prove of estimate (3.2) for derivatives of order higher than one.  
We use notations of step 7) of Section 2 with obvious modifications, in other words, 
$g_{wi}^{(k)}$, $i=1,2$ instead of $g_w^{(k)}$, $S_i^{(k)} =s_{0i}\circ s_{1i}\circ...\circ\tilde s_{ki}$ 
instead of $S$ and so on. Remind that $g_{wi}^{(k)}(z)=h_{wi}(z)$ at $|z|\ge 1/2$.

\begin{proposition} 
There are the estimates

a) If $|l|\le k+1$, $|z|\le 1/2$, then
$$
  |(\mu_{g_{w1}^{(k)}})_{(l)}(z)-(\mu_{g_{w2}^{(k)}})_{(l)}(z)|\le C(\sup_{D_{1/2}}\sup_{|(q)|\le |(l)|}|(\mu_{h_{w1}})_{(q)}-(\mu_{h_{w2}})_{(q)}|+  
$$
$$  
  +\max_{|(m)|\le k}|(h_{w1})_{(m)}(0)-(h_{w2})_{(m)}(0)|+\max_{|(m)|=k}|(\mu_{h_{w1}})_{(m)}(0)-(\mu_{h_{w2}})_{(m)}(0)|) .  
$$

b) If $|(l)|\le k+2$, then
$$
  |(S_1^{(k)})_{(l)}-(S_2^{(k)})_{l}|\le C(\max_{|(m)|\le k}|(h_{w1})_{(m)}(0)-(h_{w2})_{(m)}(0)|+\max_{|(m)|=k}|(\mu_{h_{w1}})_{(m)}(0)-(\mu_{h_{w2}})_{(m)}(0)|) . 
$$

c) For $p\ge 1$
$$
  \|\mu_{g_{w1}^{(k)}}-\mu_{g_{w2}^{(k)}}\|_p \le \| \mu_{h_{w1}}-\mu_{h_{w2}}\|_p +  
$$
$$
+C(\max_{|(l)|\le k}|(h_{w1})_{(l)}(0)-(h_{w2})_{(l)}(0)|+\max_{|(l)|=k}|(\mu_{h_{w1}})_{(l)}(0)-(\mu^{h_{w2}})_{(l)}(0)|) .   
$$
\end{proposition}

{\bf Proof}. a) and b) Remind that $g_w^{(k)}=h_w \circ S^{(k)}$. We have
$$
  \mu_{g_w^{(k)}}=\frac{\mu_{h_w}-\mu_{S^{(k)}}}{1-\bar\mu_{S^{(k)}} \mu_{h_w}}\cdot\frac{(S^{(k)})_z}{\overline{(S^{(k)})_z}}\circ (S^{(k)})^{-1} .
$$

Any derivative $(\mu_{g_w})_{(m)}$ is a sum of items of the types 
$$
  P_a (\mu_{h_w} ,\mu_{S^{(k)}} ,\bar\mu_{S^{(k)}} ) \left(\frac{(S^{(k)})_z}{\overline{(S^{(k)})_z}}\right)_{(l)}(\mu_{[.]})_{(l_1 )}^{k_1}...(\mu_{[.]})_{(l_s )}^{k_s}\circ (S^{(k)})^{-1}
  (.)_{(i_1 )}^{j_1}...(.)_{(i_r )}^{j_r}   ,  \eqno (3.55)
$$
where $|(l)|+|(l_1 )|k_1 +...+|(l_s )|k_s \le |(m)|$, $P_a$ is an uniformly bounded rational function of $\mu_{h_w}, \mu_S ,\bar\mu_S$,
$[.]$ can be any function from the tuple $h_w, S^{(k)}, \bar S^{(k)}$, any $(.)$ can be either $(S^{(k)})^{-1}$ either $\overline{(S^{(k)})^{-1}}$, 
and $|(1_1 )|j_1 +...+|(i_r )|j_r \le |(m)|$. It follows that we can represent
any difference $(\mu_{g_{w1}^{(k)}})_{(m)}-(\mu_{g_{w2}^{(k)}})_{(m)}$ as a sum of terms of the types
$$
  (P_a (\mu_{h_{w1}},\mu_{S^{(k)}_1},\bar\mu_{S^{(k)}_1})-P_a (\mu_{h_{w2}},\mu_{S^{(k)}_2},\bar\mu_{S^{(k)}_2}))\circ (S^{(k)}_1 )^{-1} [.],  \eqno (3.56)
$$ 
$$
  (P_a (\mu_{h_{w2}},\mu_{S^{(k)}_2},\bar\mu_{S^{(k)}_2})\circ (S^{(k)}_1 )^{-1}-P_a (\mu_{h_{w2}},\mu_{S^{(k)}_2},\bar\mu_{S^{(k)}_2})\circ (S^{(k)}_2 )^{-1}) [.]  , \eqno (3.57)
$$  
terms
$$
  ((\mu_{h_{w1}})_{(l_i )}^{k_i}-(\mu_{h_{w2}})_{(l_i )}^{k_i})\circ (S^{(k)}_1 )^{-1} [.],  \eqno (3.58)
$$
$$
  ((\mu_{h_{w2}})_{(l_i )}^{k_i}\circ (S^{(k)}_1 )^{-1}-(\mu_{h_{w2}})_{(l_i )}^{k_i})\circ (S^{(k)}_2 )^{-1} [.] ,  \eqno (3.59)
$$
$$
((\mu_{S^{(k)}_1})_{(l_i )}^{k_i}-(\mu_{S^{(k)}_2})_{(l_i )}^{k_i})\circ (S^{(k)}_1 )^{-1} [.],  \eqno (3.60)
$$
$$
  ((\mu_{S^{(k)}_2})_{(l_i )}^{k_i}\circ (S^{(k)}_1 )^{-1}-(\mu_{S^{(k)}_2})_{(l_i )}^{k_i}\circ (S^{(k)}_2 )^{-1}) [.] ,  \eqno (3.61)
$$
the analogous terms with $\bar\mu_{S^{(k)}_1},\bar\mu_{S^{(k)}_2}$, and
$$
  \left[\left(\frac{(S^{(k)}_1 )_z}{\overline{(S^{(k)}_1 )_z}}\right)_{(l)}\circ (S^{(k)}_1 )^{-1}-\left(\frac{(S^{(k)}_2 )_z}{\overline{(S^{(k)}_2 )_z}}\right)_{(l)}\circ (S^{(k)}_2 )^{-1}\right] [.]
 \eqno (3.62)
$$     
$$
  (((S^{(k)}_1 )^{-1})_{(i_q )}^{j_q}- ((S^{(k)}_2 )^{-1})_{(i_q )}^{j_q}) [.]  ,  \eqno (3.63)
$$ 
where $[.]$ denotes each time the product of the multiples such as in (3.55), where we omit the term corresponding to
the written difference. All these multiples are either derivatives of $(S^{(k)}_j )^{-1},j=1,2$, either derivatives of 
$\mu_{h_{wi}},\mu_{S^{(k)}_i},\bar\mu_{S^{(k)}_i},i=1,2$ in the point $(S^{(k)}_j )^{-1}(z)$, either derivatives of $P_l$ with respect to 
$\mu_{h_{wj}}$, $\mu_{S^{(k)}_j}$ or $\bar\mu_{S^{(k)}_j}$. All these derivatives are uniformly bounded. 
For derivatives of $\mu_{h_{wj}}$ it holds because $z$ belongs to the disk $D_{1/2}$ and for 
derivatives of $\mu_{S^{(k)}_i},\bar\mu_{S^{(k)}_i}$ it follows from Proposition 11.

For terms (3.56) -(3.63) we have the estimates:

For term (3.56)
$$
  \le C\sup_{D_{1/2}}(|\mu_{h_{w1}}-\mu_{h_{w2}}|+|\mu_{S^{(k)}_1}-\mu_{S^{(k)}_2}|) .  \eqno (3.56')
$$
For term (3.57)
$$
  \le C\sup_{D_{1/2}}|(S^{(k)}_1 )^{-1}(z)-(S^{(k)}_2 )^{-1}(z)|  .  \eqno (3.57')
$$
For term (3.58)
$$
  \le C\sup_{D_{1/2}}(|(\mu_{h_{w1}})_{(l_i )}-(\mu_{h_{w2}})_{(l_i )}|) .  \eqno (3.58')
$$
For term (3.59)
$$
  \le C\sup_{D_{1/2}}|(S^{(k)}_1 )^{-1}(z)-(S^{(k)}_2 )^{-1}(z)|  .  \eqno (3.59')
$$  
For term (3.60):
$$
  \le C\sup_{D_{1/2}}|(\mu_{S^{(k)}_1})_{(l_i )}-(\mu_{S^{(k)}_2})_{(l_i )}| .   \eqno (3.60')
$$
For term (3.61)
$$
  \le C|(S^{(k)}_1 )^{-1}(z)-(S^{(k)}_2 )^{-1}(z)|  .  \eqno (1.61')
$$  
For term (3.62)
$$
 \le C\sup_{D_{1/2}}(\sup_{|(q)|=|(l)|+1}|(S^{(k)}_1 )_{(q)}-(S^{(k)}_2 )_{(q)}|+|(S^{(k)}_1 )^{-1}(z)-(S^{(k)}_2 )^{-1}(z)|)  .  \eqno (3.62')
$$  
For term (3.63)
$$
 \le C\sup_{D_{1/2}}|((S^{(k)}_1 )^{-1})_{(i_q )}-((S^{(k)}_2 )^{-1})_{(i_q )}|  .  \eqno (3.63')
$$
Here we take into consideration that all derivatives of $S^{(k)}_i$, $(S^{(k)}_i )^{-1}$ and $\mu_{S^{(k)}_i}$, $\mu_{h_{wi}}$
are uniformly bounded on $D_{1/2}$. Terms (3.56') and (3.58') yield the term
$\sup_{D_{1/2}}\sup_{|(q)|\le |(l)|}|(\mu_{h_{w1}})_{(q)}-(\mu_{h_{w2}})_{(q)}|$ in the 
estimate of $|(\mu_{g_{w1}^{(k)}})_{(l)}(z)-(\mu_{g_{w2}^{(k)}})_{(l)}(z)|$.  

We must estimate $|((S^{(k)}_1 )^{-1})_{(l)}-((S^{(k)}_2 )^{-1})_{(l)}|, |(l)|\le k+1$, $|(S^{(k)}_1 )_{(l)}-(S^{(k)}_2 )_{(l)}|,|(l)|\le k+2$,
and $|(\mu_{S^{(k)}_1})_{(l)}-(\mu_{S^{(k)}_2})_{(l)}|, |(l)|\le k+1$ on $D_{1/2}$. Consider at first
the difference $|(S^{(k)}_1 )_{(l)}-(S^{(k)}_2 )_{(l)}|$.

The map $S^{(k)}_j$ is the composition
$$
  S_j = s_{0j}\circ s_{1j}...\circ s_{(k-1)j}\circ \tilde s_{kj} . 
$$
We can write the difference $(S^{(k)}_1 )_{(l)}-(S^{(k)}_2 )_{(l)}$ as a sum of terms of the types
$$
  ((s_{i1})_{(l_i )}-(s_{i2})_{(l_i )})\circ s_{(i+1)1}\circ...\circ\tilde s_{k1} [.]    \eqno (3.64)
$$
and
$$
  ((s_{i2})_{(l_i )})\circ s_{(i+1)1}\circ...\circ\tilde s_{k1}-(s_{i2})_{(l_i )})\circ s_{(i+1)2}\circ...\circ\tilde s_{k2}))[.] ,  (3.65)
$$
where $[.]$ denotes products having uniform estimates in $D_{1/2}$. We don't write analogous terms
containing differences of derivatives of $\tilde s_{kj}$, $j=1,2$.

Consider at first difference (3.65). It has the estimate
$$
 C|s_{(i+1)1}\circ...\circ\tilde s_{k1}(z)-s_{(i+1)2}\circ...\circ\tilde s_{k2}(z)|   .
$$
We can represent the last difference as a sum of terms of the types
$$
  s_{(i+1)1}\circ...\circ\tilde s_{k1}(z)-s_{(i+1)2}\circ...\circ\tilde s_{k1}(z) ,  \eqno (3.66)
$$
$$
  s_{(i+1)2}\circ...s_{j1}\circ...\circ\tilde s_{k1}(z)- s_{(i+1)2}\circ...s_{j2}\circ...\circ\tilde s_{k1}(z)|  \eqno (3.66')
$$
for $j=i+1,...k$. 

The following considerations are analogous to the reasoning leading to relations (2.42), (2.43). In more details,
let $\{c_{q\bar j}\},q+j=l+1$ be the coefficients of the form $R_{(l-1)l}$. By (2.37) - (2.40), we can see that $s_0$ is a function
of $a_0$; $s_l$, $l<k$ are functions of $c_{q\bar j}$, and $\tilde s_k$ is a function of $c_{q\bar j}, q+j=k+1$, $j\ge 1$. 
All these functions have uniformly bounded derivatives of any order at $|z|\le 1/2$. The same holds for $s_l^{-1}$.
When we consider $h_{wi}$, $i=1,2$ we denote by $\{c^i_{q\bar j}\},q+j=l+1, i=1,2$ the coefficients of the
corresponding forms $R_{(l-1)l,i}$. These coefficients are polynomial in $a_0^{-1}$
and in derivatives in zero of function $h_w$ and, hence, they are uniformly bounded.
It follows that differences (3.66), (3.66') have estimates
$$
  \le C\max\{|a_0^1 -a_0^2 |,|c^1_{q\bar j}-c^2_{q\bar j}|\} 
$$
for $q+j =m, m\le k$ or $q+j =k+1, j\ge 1$. 

 Further, the coefficients $c_{q\bar j}, q+j=l, l\le k$ or $q+j =k+1, j\ge 1$
are functions of $(h_w )_{z^p \bar z^q}(0)$ with $p+q\le k$ and $(\mu_{h_w})_{(l)}(0),|(l)|=k$.
All these functions have bounded derivatives. We get for terms of type (3.65) the estimate
$$
  \le C(\max_{|(l)|\le k}|(h_{w1})_{(l)}(0)-h_{w2})_{(l)}(0)|+\max_{|(l)|=k}|(\mu_{h_{w1}})_{(l)}(0)-(\mu_{h_{w2}})_{(l)}(0)|) . 
  \eqno (3.67)
$$
From the other hand, $z$-derivatives of $s_l$ are also functions of the same variables $(h_w )_{z^p \bar z^q}(0)$,
$p+q\le k$ and $\mu^{(k-1)}_{w,z^q \bar z^{k-q}}(0)$. We obtain for the terms of type (3.64) the same estimate (3.67).
Thus we proved estimate b).
  
Obviously, we get analogous estimate for $|((S^{(k)}_1 )^{-1})_{(l)}-((S^{(k)}_2 )^{-1})_{(l)}|, |l|\le k+1$. Also,
$|(\mu_{S^{(k)}_1})_{(l)}-(\mu_{S^{(k)}_2})_{(l)}|, |l|\le k+1$ is no greater, than $\sup_{|l|\le k+2}|(S^{(k)}_1 )_{(l)}-(S^{(k)}_2 )_{(l)}|$
and, hence, we also get for such term estimate (3.67). Since we have estimates (3.56') - (3.63'), we obtain a).

c) Since $\mu_{g_{wi}^{(k)}}(z)=\mu_{h_{wi}}(z)$ at $|z|\ge 1$, the estimate follows from a).$\Box$

Our next step will be estimates for $|h_{w1}(z)-h_{w2}(z)|$ and $|(h_{w1})_z (z)-(h_{w2})_z (z)|$.
Also we obtain estimates for $|g_{w1,z}^{(k)}(z)-g_{w2,z}^{(k)}(z)|$.

\begin{proposition} We have for some $0<\alpha <1$ depending only on $d$ and some uniform $C$ 
a)
$$
  |h_{w1}(z)-h_{2w}(z)| \le C\min\left[ 1,\frac{\|\mu_1 -\mu_2 \|_p^{\alpha}}{1-|w|}+\frac{|\mu_1 (w)-\mu_2 (w)|}{(1-|z|)^2}\right] ,
$$
b) 
$$
  |(h_{w1})_z (z)-(h_{w2})_z (z)| \le \frac{C}{(\min_{i=1,2}|1-\bar w Z_{wi}^{-1}(z)|)^2}\min\left[ 1,\frac{\|\mu_1 -\mu_2 \|_p^{\alpha}}{1-|w|}+\right.  
$$
$$ 
  \left.+|\mu_1 (w)-\mu_2 (w)|^{\alpha}+\frac{|\mu_1 (w)-\mu_2 (w)|}{1-|z|}+\|\mu_{F_{w1}} -\mu_{F_{w2}}\|_p^{\alpha}\right],
$$
c) for $|z|\le 1/2$
$$
  |g_{w1,z}^{(k)}(z)-g_{w2,z}^{(k)}(z)| \le C\min\left[1,\frac{\|\mu_1 -\mu_2 \|_p^{\alpha}}{1-|w|}+  
|\mu_1 (w)-\mu_2 (w)|^{\alpha}+\|\mu_{F_{w1}} -\mu_{F_{w2}}\|_p^{\alpha}+\right.
$$
$$
 \left. +\max_{m\le k-1}|(g_{w1}^{(m)})_{z^{m+1}}(0)-(g_{w2}^{(m)})_{z^{m+1}}(0)|+\max_{|(m)|\le k}|(\mu_{h_{w1}})_{(m)}(0)-(\mu_{h_{w2}})_{(m)}(0)|\right]  . 
$$ 
\end{proposition}

{\bf Proof}. a) Using the representation $h_w =\vpf_{w'}\circ f\circ\vpf_w^{-1}\circ Z_w^{-1}$, we see that we must estimate the
sum
$$
  |(\vpf_{w'_1}-\vpf_{w'_2})\circ f_1 \circ\vpf_w^{-1}\circ Z_{w1}^{-1}|+|(\vpf_{w'_2}\circ f_1 -\vpf_{w'_2}\circ f_2 )\circ\vpf_w^{-1}\circ Z_{w1}^{-1}|+
$$
$$
 +|\vpf_{w'_2}\circ f_2 \circ\vpf_w^{-1}\circ Z_{w1}^{-1}-\vpf_{2w'}\circ f_2 \circ\vpf_w^{-1}\circ Z_{w2}^{-1}| .  \eqno (3.68)
$$  
Considering $\vpf_w'$ as a function of $w'$ and applying inequality (2.4), we get that we 
can estimate the first term as
$$
  \sup_{[w'_1 ,w'_2]}(|(\vpf_w')_{w'}|+|(\vpf_w')_{\overline{w'}}||w'_1 -w'_2 |\le C\frac{|w'_1 -w'_2|}{1-|w|} .
$$
Applying estimate (3.1), we obtain
$$
  |(\vpf_{w'_1}-\vpf_{w'_2})\circ f_1 \circ\vpf_w^{-1}\circ Z_{w1}^{-1}(z)|\le C\min\left[ 1, \frac{\|\mu_1 -\mu_2 \|_p^{\alpha}}{1-|w|}\right] .  
$$
The second term in (3.68) we estimate as 
$$
  \sup |(\vpf_{2w'})_z||f_1 -f_2 |_C \le C(1-|w|)^{-1}|f_1 -f_2 |_C  .
$$
Again we obtain the estimate
$$
 |(\vpf_{w'_2}\circ f_1 -\vpf_{w'_2}\circ f_2 )\circ\vpf_w^{-1}\circ Z_{w1}^{-1}(z)|\le C\min\left[ 1, \frac{\|\mu_1 -\mu_2 \|_p^{\alpha}}{1-|w|}\right] .
$$
Further, we have the estimate $c\le |(f_w )_z (z)| \le C|1-\bar w z|^{-2}$, it is a particular case of
estimate (2.31). 
We see that the third term in (3.68) we can estimate as
$$
  |f_{w2}\circ Z_{w1}^{-1}(z)-f_{w2}\circ Z_{w2}^{-1}(z)|\le C\min\left[ 1,\frac{| Z_{w1}^{-1}(z)-Z_{w2}^{-1}(z)|}{(1-|z|)^2}\right]\le
$$
$$
  \le C_1 \min\left[ 1,\frac{|\mu_1 (w)-\mu_2 (w)|}{(1-|z|)^2}\right]
$$
with some uniform $C_1$. 

b) We must estimate the terms
$$
  |((\vpf_{w'_1})_z -(\vpf_{w'_2})_z )\circ f_1 \circ\vpf_w^{-1}\circ Z_{w1}^{-1}(z)| |(f_1 )_z \circ\vpf_w^{-1}\circ Z_{w1}^{-1}(z)\cdot
 (\vpf_w^{-1})_z \circ Z_{w1}^{-1}(z)\cdot (Z_{w1}^{-1})_z (z)| ,  \eqno (3.69)
$$ 
$$
  |(\vpf_{w'_2})_z \circ f_1 \circ\vpf_w^{-1}\circ Z_{w1}^{-1}(z) -(\vpf_{w'_2})_z \circ f_2 \circ\vpf_w^{-1}\circ Z_{w2}^{-1}(z)||(f_1 )_z \circ\vpf_w^{-1}\circ Z_{w1}^{-1}(z)\cdot
 (\vpf_w^{-1})_z \circ Z_{w1}^{-1}(z)\cdot (Z_{w1}^{-1})_z (z)| ,  \eqno (3.70)
$$ 
$$
  |(\vpf_{w'_2})_z \circ f_2 \circ\vpf_w^{-1}\circ Z_{w2}^{-1}(z)||((f_1 )_z -(f_2 )_z )\circ\vpf_w^{-1}\circ Z_{w1}^{-1}(z)|
 |(\vpf_w^{-1})_z \circ Z_{w1}^{-1}(z)\cdot (Z_{w1}^{-1})_z (z)| ,  \eqno (3.71)
$$ 
$$
  |(\vpf_{w'_2})_z \circ f_2 \circ\vpf_w^{-1}\circ Z_{w2}^{-1}(z)||(f_2 )_z \circ\vpf_w^{-1}\circ Z_{w1}^{-1}(z) -(f_2 )_z \circ\vpf_w^{-1}\circ Z_{w2}^{-1}(z)|
 |(\vpf_w^{-1})_z \circ Z_{w1}^{-1}(z)\cdot (Z_{w1}^{-1})_z (z)| ,  \eqno (3.72)
$$   
$$
  |(\vpf_{w'_2})_z \circ f_2 \circ\vpf_w^{-1}\circ Z_{w2}^{-1}(z)\cdot (f_2 )_z \circ\vpf_w^{-1}\circ Z_{w2}^{-1}(z)| 
 |(\vpf_w^{-1})_z \circ Z_{w1}^{-1}(z)-(\vpf_w^{-1})_z \circ Z_{w2}^{-1}(z)|(Z_{w1}^{-1})_z (z)| ,  \eqno (3.73)
$$
$$
  |(\vpf_{w'_2})_z \circ f_2 \circ\vpf_w^{-1}\circ Z_{w2}^{-1}(z)\cdot (f_2 )_z \circ\vpf_w^{-1}\circ Z_{2w}^{-1}(z)| 
 (\vpf_w^{-1})_z \circ Z_{w2}^{-1}(z)||(Z_{w1}^{-1})_z (z)-(Z_{w2}^{-1})_z (z)| .  \eqno (3.74)
$$
We don't write the analogous terms containing $\bar z$-derivatives.

To estimate term (3.69) we proceed as in the case of the first difference in (3.68). The derivatives of $(\vpf_{w'})_z$
with respect to $w'$ have the estimate $C(1-|w|)^{-2}$. The derivative $(\vpf_w^{-1})_z \circ Z_{w1}^{-1}(z)$ has the 
estimate $C(1-|w|^2 )|1-\bar w Z_{w1}^{-1}(z) |^{-2}$. Other multiples are uniformly bounded. 
Applying estimate (3.1), we obtain for term (3.69) the estimate
$$
  \frac{C}{|1-\bar w Z_{w1}^{-1}(z) |^2}\min\left[ 1,\frac{\|\mu_1 -\mu_2 \|_p^{\alpha}}{1-|w|}\right] . \eqno (3.69')
$$ 
For term (3.70) we obtain analogous estimate because the second derivative of $\vpf_{2w'}$ has the estimate $C(1-|w|)^{-2}$. 

In all other terms there is the multiple $(\vpf_{w'_2})_z =O(1-|w|)^{-1}$.

Thus for term (3.71) we have the estimate
$$
  \frac{C}{1-|w|}\|(f_1 )_z -(f_2 )_z \|_C\frac{1-|w|^2}{|1-\bar w Z_{w1}^{-1}(z) |^2} \le
$$
$$
  \le \frac{C}{|1-\bar w Z_{w1}^{-1}(z) |^2}\min\left[ 1,\frac{\|\mu_1 -\mu_2 \|_p^{\alpha}}{1-|w|}+ 
|\mu_1 (w)-\mu_2 (w)|^{\alpha}+\|\mu_{F_{1w}} -\mu_{F_{2w}}\|_p^{\alpha}\right]  \eqno (3.71')
$$
by Proposition 18.

Consider term (3.72).
Denote $z_i =\vpf_w^{-1}\circ Z_{wi}^{-1}(z)$, $i=1,2$. We see that we have the estimate
$$
  \frac{C}{|1-\bar w  Z_{w1}^{-1}(z) |^2}\sup_{t\in [ Z_{w1}^{-1}(z) , Z_{w2}^{-1}(z) ],|(k)|=2}|(f_2)_{(k)}\circ\vpf_w^{-1}(z)\cdot (\vpf_w^{-1})_z (t)|| Z_{w1}^{-1}(z) - Z_{w2}^{-1}(z)|  .
$$
Now, if $|(k)|=2$, 
$$
 |(f_2)_{(k)} \circ\vpf_w^{-1}(z)| \le C(1-|\vpf_w^{-1}(z)|^2 )^{-1}\le 
 C\frac{|1-\bar w z|^2}{(1-|w|^2 )(1-|z|^2 )},
$$
$$ 
 \le |(\vpf_w^{-1})_z |=\frac{1-|w|^2}{|1-\bar w z|^2}
$$
Since $c(1-|z|)\le 1-|Z_{wi}^{-1}(z)|\le C(1-|z|)$ for some uniform $c$ and $C$, we obtain the estimate
$$
  \frac{C}{|1-\bar w Z_{w1}^{-1}(z)|^2}\min\left[ 1,\frac{|\mu_1 (w)-\mu_2 (w)|}{1-|z|}\right] .  \eqno (3.72')
$$
For term (3.73) we have the estimate
$$
  \frac{C}{1-|w|}\sup_{[Z_{w1}^{-1}(z) ,Z_{w2}^{-1}(z)]}|(\vpf_w^{-1})_{zz}||Z_{w1}^{-1}(z) -Z_{w2}^{-1}(z)|\le C\frac{|Z_{w1}^{-1}(z) -Z_{w2}^{-1}(z)|}{(\min_{i=1,2}|1-\bar w Z_{wi}^{-1}(z)|)^3}\le
$$
$$  
  \frac{C}{(\min_{i=1,2}|1-\bar w Z_{wi}^{-1}(z)|)^2}\frac{|\mu_1 (w)-\mu_2 (w)|}{\min_{i=1,2}|1-\bar w Z_{wi}^{-1}(z)|}
$$
at $|\mu_1 (w)-\mu_2 (w)|\le 1-|z|$. 
We obtain the estimate
$$
  \frac{C}{(\min_{i=1,2}|1-\bar w Z_{wi}^{-1}(z)|)^2}\min\left[ 1,\frac{|\mu_1 (w)-\mu_2 (w)|}{1-|z|}\right] .  \eqno (3.73')
$$
At last, for term (3.74) we have the estimate
$$
  \frac{C}{|1-\bar w Z_{w2}^{-1}(z)|^2}|\mu_1 (w)-\mu_2 (w)| .  \eqno (3.74')
$$
Collecting estimates (3.69') - (3.74') we obtain b).

c) From the representation: $g_{wi}^{(k)}=h_{wi}\circ s_{0i}\circ...\circ \tilde s_{ki} =h_{wi}\circ S^{(k)}_i$, $i=1,2$, 
follows that we must estimate the terms of the types
$$
  (h_{w1,z}-h_{w2,z})\circ S^{(k)}_1 (z)[.],
$$
$$
  (h_{w2,z}\circ S_1 (z) -h_{w2,z}\circ S^{(k)}_2 (z))[.]  ,
$$
$$
  ((s_{j1})_z-(s_{j2})_z)\circ s_{(j+1)1}\circ...\circ \tilde s_{k1}(z) [.] ,   
$$
$$
  ((s_{j2})_z \circ s_{(j+1)1}\circ...\circ \tilde s_{k1}(z)-(s_{j2})_z \circ s_{(j+1)2}\circ...\circ \tilde s_{k2}(z))[.] ,\,j\le k-1  /  
$$
We don't write the analogous term with the difference $\tilde s_{k1}(z)-\tilde s_{k2}(z)$.
Multiples in the square brackets are uniformly bounded because derivatives of $s_{ij}$ are uniformly bounded 
and derivatives of $h_{wi}$ are uniformly bounded on $D_{1/2}$.

For the first term we have estimate b)  
$$
  C\max\left[ 1,\frac{\|\mu_1 -\mu_2 \|_p^{\alpha}}{1-|w|}+  
|\mu_1 (w)-\mu_2 (w)|^{\alpha}+\|\mu_{F_{1w}} -\mu_{F_{2w}}\|_p^{\alpha}\right]   .
$$ 
The second term is no greater than
$$
  C|S_1 (z)-S_2 (z)|,
$$
and, by point b) of Proposition 22, we get the estimate
$$
  C(\max_{|(l)|\le k}|(h_{w1})_{(l)}(0)-h_{w2})_{(l)}(0)|+\max_{|(l)|=k}|(\mu_{h_{w1}})_{(l)}(0)-(\mu_{h_{w2}})_{(l)}(0)|) . 
$$
The same estimate holds for the terms of third and fourth types.
Indeed, these terms are the particular cases of terms (3.65) and (3.66). Applying (2.43) we obtain c).
$\Box$

{\bf Proof of estimate (3.2)}. At first we estimate
$|g_{1w,z^{k+1}}^{(k)}(0)-g_{2w,z^{k+1}}^{(k)}(0)|$. We have representations (2.53) and the 
decomposition analogous to (2.55) but cut on the term of order $z^{k+1}$. Thus,
$$
 |(g_{w1}^{(k)})_{z^{k+1}}(0)-(g_{w2}^{(k)})_{z^{k+1}}(0)|\le C\left(\max_{\pd D_{1/2}}|h_{w1}-h_{w2}|+
\int_{D_{1/2}}\frac{|(\gamma_1 (g_{w1}^{(k)})_z -\gamma_2 (g_{w2}^{(k)})_z )(z)|}{|z|}dS_z \right) .  \eqno (3.75)
$$
Here we take into consideration that $g_{wi}^{(k)}(z)=h_{wi}(z)$ for $|z|\ge 1/2$.

For the first difference in the right side we have the first estimate of Proposition 23
$$
  \max_{\pd D_{1/2}}|h_{w1}-h_{w2}|\le C\min\left[1, \frac{\|\mu_1 -\mu_2 \|_p^{\alpha}}{1-|w|}+|\mu_1 (w)-\mu_2 (w)|\right] . \eqno (3.76)
$$
We must estimate the integral in (3.75). We have
$$
  \int_{D_{1/2}} \frac{|(\gamma_1 (g_{w1}^{(k)})_z-\gamma_2 (g_{w1}^{(k)})_z )(z)|}{|z|}dS_z \le 
$$
$$  
 \le\int_{D_{1/2}}\frac{|\gamma_1 (z)||((g_{w1}^{(k)})_z-(g_{w2}^{(k)})_z )(z)|}{|z|}dS_z +\int_{D_{1/2}}\frac{|(\gamma_1 -\gamma_2 )(z)||(g_{w2}^{(k)})_z (z)|}{|z|}dS_z
\eqno (3.77)
$$  
Consider the first integral in the right side. Applying third estimate of Proposition 23,
we have
$$
  I_1 =\int_{D_{1/2}}\frac{|\gamma_1 (z)||((g_{w1}^{(k)})_z-(g_{w2}^{(k)})_z )(z)|}{|z|}dS_z \le\sup_{|z|\le 1/2}|(g_{w1}^{(k)})_z-(g_{w2}^{(k)})_z |\int_{D_{1/2}}\frac{|\gamma_1 (z)|}{|z|}dS_z \le 
$$
$$
  \le C\min\left[ 1,\frac{\|\mu_1 -\mu_2 \|_p^{\alpha}}{1-|w|}+ |\mu_1 (w)-\mu_2 (w)|^{\alpha}+\|\mu_{F_{1w}} -\mu_{F_{2w}}\|_p^{\alpha}+\right.
$$
$$
 \left. +\max_{|(m)|\le k}|(\mu_{h_{w1}})_{(m)}(0)-(\mu^{h_{w2}})_{(m)}(0)|+\max_{m\le k-1}|(g_{w1}^{(m)})_{z^m}(0)-(g_{w2}^{(m)})_{z^m}(0)|\right]  . \eqno (3.78)
$$ 

Consider the second integral in the right side of (3.77). Remind that $\gamma (z)=|z|^{-(k+1)}\mu_{g_w^{(k)}}(z)$. It means
that for any $z=re^{i\vpf}$ we can write
$$
  |\gamma_1 (z)-\gamma_2 (z)|\le\frac{1}{(k+1)!}\left|\frac{d^{k+1}}{dr^{k+1}}(\mu_{g_{w1}^{(k)}}(re^{i\vpf})-\mu_{g_{w2}^{(k)}}(re^{i\vpf}))\right|  .
$$
By Proposition 22, we obtain
$$
  I_2 =\int_{D_{1/2}}\frac{|\gamma_1 -\gamma_2 |(z) |(g_{w2}^{(k)})_z (z)|}{|z|}dS_z \le 
C \sup_{D_{1/2}}\sup_{|l|=k+1}|(\mu_{g_{w1}^{(k)}})_{(l)}-(\mu_{g_{w2}^{(k)}})_{(l)}|\int_{D_{1/2}}\frac{|g_{w2}^{(k)})_z (z)|}{|z|}dS_z \le
$$
$$
  \le C(\sup_{D_{1/2}}\sup_{|(q)|\le k+1}|(\mu_{h_{w1}})_{(q)}-(\mu_{h_{w2}})_{(q)}|  
  +\max_{|(m)|\le k}|(h_{w1})_{(m)}(0)-(h_{w2})_{(m)}(0)|) .  \eqno (3.79)
$$

Gathering estimates (3.76), (3.78) and (3.79), we obtain by induction
$$
 |(g_{w1}^{(k)})_{z^{k+1}}(0)-(g_{w2}^{(k)})_{z^{k+1}}(0)|\le C\min\left[ 1, \frac{\|\mu_1 -\mu_2 \|_p^{\alpha}}{1-|w|}+\right. 
$$
$$
  +|\mu_1 (w)-\mu_2 (w)|^{\alpha}+\|\mu_{F_{w1}} -\mu_{F_{w2}}\|_p^{\alpha}+\sup_{D_{1/2}}\sup_{|(q)|\le k+1}|(\mu_{h_{w1}})_{(q)}-(\mu_{h_{w2}})_{(q)}|+
$$
$$
  \left. +\max_{|(m)|\le k}|(h_{w1})_{(m)}(0)-(h_{w2})_{(m)}(0)|\right] .  
$$

Now, if $|(l)|=k+1$, we have once more by induction, applying relations (2.43),
$$
  |(h_{w1})_{(l)}(0)-(h_{w2})_{(l)}(0)|\le C\min\left[ 1,\frac{\|\mu_1 -\mu_2 \|_p^{\alpha}}{1-|w|}+\right. 
$$
$$
 \left.+|\mu_1 (w)-\mu_2 (w)|^{\alpha}+\|\mu_{F_{1w}} -\mu_{F_{2w}}\|_p^{\alpha}+\sup_{D_{1/2}}\sup_{|(q)|\le k+1}|(\mu_{h_{w1}})_{(q)}-(\mu_{h_{w2}})_{(q)}|\right]
\eqno (3.80) 
$$
Further, by (3.12) we obtain 
$$
  \sup_{D_{1/2}}\sup_{|q|\le k+1}|(\mu_{h_{w1}})_{(q)}-(\mu_{h_{w2}})_{(q)}|\le C\min [1,\sup_{|q|\le k+1}\sup_{D_{R_w}}|(\mu_1)_{(q)}-(\mu_2)_{(q)}|+
$$
$$
  +|\mu_1 (w)-\mu_2 (w)|] ,
$$
where $1-R_w^2 =b(1-|w|^2 )$ with some uniform $b$. Also, Proposition 17 yields 
$$
  \|\mu_{F_{w1}} -\mu_{F_{w2}}\|_p \le C\min [1,\sup_{|q|=1}\sup_{D_{R_w}}|(\mu_1)_{(q)}-(\mu_2)_{(q)}|+
$$
$$
  +\sup_{D_R}|\mu_1 (z)-\mu_2 (z)|+[(1-R)/(1-|w|)]^{1/p}+|\mu_1 (w)-\mu_2 (w)|^{1/(2p+1)}] , 
$$
where we can take any $R<1$. Thus we obtain
$$
 |(h_{1w})_{(k+1)}(0)-(h_{2w})_{(k+1)}(0)|\le C\min \left[1,\frac{\|\mu_1 -\mu_2 \|_p^{\alpha}}{1-|w|}+\right.
$$
$$
 +\sup_{|q|\le k+1}\sup_{D_{R_w}}|(\mu_1)_{(q)}-(\mu_2)_{(q)}|+|\mu_1 (w)-\mu_2 (w)|^{\alpha}+\sup_{D_R}|\mu_1 -\mu_2 |^{\alpha}+
$$
$$
 \left.+[(1-R)/(1-|w|)]^{\alpha}+\sup_{|q|=1}\sup_{D_{R_w}}|(\mu_1)_{(q)}-(\mu_2)_{(q)}|^{\alpha}\right]  . \eqno (3.81)
$$ 
with some new $0<\alpha <1$. 

Now suppose that estimate (3.2) holds for multi-indexes $(l)$ with $|(l)|=k$. Recalling representation (2.15 we see
that we can represent the derivative $f_{(l)}(w)$, $|(l)|=k+1$ as a sum
of products, in which multiples have the types $(\vpf_{w'}^{-1})_{l_1}(0) =c_{l_1}(1-|w'|^2 )\bar w'^{l_1 -1}$, 
$(h_w )_{(l_2 )}(0)$, $(Z_w )_{(l_3 )}(0)$, $(\vpf_w)_{(l_4 )}(w) =c_{l_4}\bar w^{|(l_4 )| -1}(1-|w|^2 )^{-|(l_4 )|}$, 
where $l_1 \le k+1,|(l_i )|\le k+1, i\ge 2$. It implies that we can represent the difference $(f_1 )_{(k+1)}(w)-(f_2 )_{(k+1)}(w)$ as a
sum of items of the types 
$$
  ||w'_1 |^2 -|w'_2 |^2 |O((1-|w|)^{-(k+1)}) ,   \eqno (3.82)
$$
$$
  |w'_1  -w'_2 | |O((1-|w|)^{-k})  ,   \eqno (3.83)
$$
$$
  |(Z_{w1})_{(l)}(0)-Z_{w2})_{(l)}(0)|O((1-|w|)^{-k}) ,\, |l|\le k+1  ,  \eqno (3.84)
$$
$$
 |(h_{w1})_{(l)}(0)-(h_{w2})_{(l)}(0)|O((1-|w|)^{-k}) ,\, |l|\le k+1  .  \eqno (3.85)
$$
Items of type (3.82) have the estimate
$$
  C\min\left[1,\frac{|f_1 (w)-f_2 (w)|}{1-|w|}(1-|w|)^{-k}\right] \le C\min\left[ 1,\frac{|\mu_1 -\mu_2 \|_p^{\alpha}}{1-|w|}\right](1-|w|)^{-k}) .
$$
Here we apply estimate (3.1) and recall that $f_{(l)}(w)=O((1-|w|)^{-k}$. Also, term (3.83) has the estimate
$$
  C|\mu_1 -\mu_2 \|_p^{\alpha}(1-|w|)^{-k}) .
$$
Term (3.84) we estimate as
$$
  C|\mu_1 (w) -\mu_2 (w)|(1-|w|)^{-k} .
$$   
At last, for term (3.85) we have estimate (3.81). We obtain
$$
  |(f_1 )_{(l)}(w)-(f_2 )_{(l)}(w)|\le \frac{C}{(1-|w|)^k}\inf\{1,\left[\chi\left(\frac{\|\mu_1 -\mu_2 \|_p^{\alpha}}{1-|w|}\right)+\right.
$$
$$
 +\sup_{|q|\le k+1}\sup_{D_{R_w}}|(\mu_1)_{(q)}-(\mu_2)_{(q)}|+|\mu_1 (w)-\mu_2 (w)|^{\alpha}+\sup_{D_R}|\mu_1 -\mu_2 |^{\alpha}+
$$
$$
 \left.+[(1-R)/(1-|w|)]^{\alpha}+\sup_{|q|=1}\sup_{D_{R_w}}|(\mu_1)_{(q)}-(\mu_2)_{(q)}|^{\alpha}\right]\}  
$$ 
if $|(l)|=k+1$. Estimate (3.2) follows immediately.
$\Box$

\section{Approximation of a family of normal mappings by a smooth family} 

In this section we suppose that $\mu$ depends on a vector parameter $t\in\bbC^n$ and satisfies
conditions of Theorem 2. In this section and in what follows the word "uniform" applied
to an estimate or to a constant means, also, that it doesn't depend on $t$. 
We denote by $z$ the chart on a fiber. We adopt the notations $\pd^{(k,l)}f$ or $f_{(k,l)}$ with 
double multi-indexes, as it was described in Section 1. Also, we denote by $f_{(k)}$ derivatives 
in $z$, $\bar z$ and by $f_{(0,l)}$ derivatives in $t$ with the multi-index $(l)$.
We use the notations $\mu_t (z)=\mu (z,t)$ or $f_t (z)=f(z,t)$ when we aren't in danger
to mix them with notations for derivatives.

The next proposition is a corollary of Lemma 2. It shows that for $z$-derivatives of the normal mappings 
we have a Holder continuity with respect to the parameters .

\begin{proposition} In conditions of Theorem 2 let $f_t$ be the $\mu_t$-quasyconformal
normal mapping. Then
$$
  |f_{(k)}(z,t)-f_{(k)}(z,t+\delta t)|\le \frac{C}{(1-|z|)^{|(k)|-1}}\min \{1,(\delta t)^{\beta}((1-|z|)^{-s}) . \eqno (4.1)
$$  
for uniform $C$ and some $0<\beta <1$ and $s>0$ depending only on $d$ and $|(k)|$.
\end{proposition}

{\bf Proof}. By estimate (3.2) of Lemma 2, we have
$$
  |f_{(k)}(z,t)-f_{(z)}(z,t+\delta t)| \le
  \frac{C}{(1-|z|)^{|(k)|-1}}\min\left\{ 1,\frac{\|\mu_t -\mu_{t+\delta t}\|_p^{\alpha}}{1-|z|}+\right.
$$
$$  
  \left.+\sup_{|(q)|\le |(k)|}\sup_{D_{R_z}}|(\mu_t)_{(q)}-(\mu_{t+\delta t})_{(q)}|^{\alpha}+
  \sup_{D_R}|\mu_t -\mu_{t+\delta t}|^{\alpha}+[(1-R)/(1-|z|)]^{\alpha /p}\right\}  , 
\eqno (4.2)
$$ 
where $R$ can be arbitrary radius less than 1 and $R_z$ is such that 
$$
  1-R_z \le b(1-|z|)   .  \eqno (4.3)
$$
for some uniform $b>0$.
From the other hand, by inequality (1.5), we have
$$
 |(\mu_t)_{(q)}(z)-(\mu_{t+\delta t})_{(q)}(z)|\le C\frac{\delta t}{(1-|z|)^N}, \eqno (4.4)
$$
where $N$ depends only on $|(k)|$ if $|(q)|\le |(k)|$ and $C$ is uniform. 

Consider the right side of inequality (4.2). Set some $r<1$. We have
$$
 \|\mu_t -\mu_{t+\delta t}\|_p \le C\left[\frac{\delta t}{(1-r)^N}+(1-r)^{1/p}\right] .
$$
Put $r$ such that $\delta t (1-r)^{-N}=(1-r)^{1/p}$, i.e., $1-r=(\delta t )^{p/(Np+1)}$. 
We obtain
$$
  \|\mu_t -\mu_{t+\delta t}\|_p \le C(\delta t )^{\frac{1}{Np+1}}  .
$$
Also,
$$
 \sup_{D_R}|\mu_t -\mu_{t+\delta t}|^{\alpha}\le C\frac{\delta t^{\alpha}}{(1-R)^{N\alpha}}  .
$$
Put $R$ such that $[\delta t /(1-R)^N ]^{\alpha}= [(1-R)/(1-|z|)]^{\alpha /p}$, i.e, 
$1-R=(\delta t)^{p/(Np+1)}(1-|z|)^{1/(Np+1)}$. We obtain
$$
  \sup_{D_R}|\mu_t -\mu_{t+\delta t}|^{\alpha}+[(1-R)/(1-|z|)]^{\alpha /p} \le 
C\delta t^{\frac{\alpha}{Np+1}}(1-|z|)^{-\frac{N\alpha}{Np+1}}  .
$$
At last, by (4.4) and (4.3),
$$
\sup_{|(q)|\le |(k)|}\sup_{D_{R_z}}|(\mu_t)_{(q)}-(\mu_{t+\delta t})_{(q)}|^{\alpha}\le
C\frac{\delta t ^{\alpha}}{(1-|z|)^{N\alpha}} 
$$
We obtain for the sum in the right part of inequality (4.2) the estimate
$(\delta t)^{\beta}((1-|z|)^{-s}$, where we can put $\beta =\alpha /(Np+1)$, $s=N\alpha$.
$\Box$

The main result of this section is the next lemma about approximations:

\begin{lemma} In the above assumptions for every $\varepsilon >0$ and natural $m$ there exists a family of mappings $f_{at}$ 
smoothly depending on $z$ and approximating the family $f_t$ up to derivatives of order $m$
$$
  |(f_{at})_{(l)}-(f_t)_{(l)}|\le \varepsilon ,\,|(l)|\le m     \eqno (4.5)
$$
on $D$. The maps $f_{at}$ are quasiconformal with complex dilatations $\mu_{at}$ and map $D$ homeomorphically
onto some domain $\Omega_t$. The mapping $B\to C^0 (D):t\mapsto f_{at}$ is continuous.

There exists the decomposition
$$
  f_{at} =h_{at}\circ f^{\mu_{at}} ,   \eqno (4.6)
$$
where $f^{\mu_{at}}$ is the normal map with with the complex dilatation $\mu_{at}$ and $h_{at}$ is a holomorphic univalent
function on $D$ satisfying the estimate
$$
  |h'_{at}-1|\le\varepsilon ,\, |h''_{at}|\le \varepsilon (1-|t|)^{-1}    \eqno (4.7)
$$
at $m\ge 2$. Derivatives of $f_{at}$ satisfy estimates analogous to (1.5)
$$
  |(f_{at})_{(k,l)}(z) \le C(1-|z|)^{N_{|(l)|+|(k)|}} ,   \eqno (5.8)  
$$
where $C$ is uniform and $N_{|(l)|+|(k)|}$ doesn't depend on $\varepsilon$ (though it can depend om $m$). 
\end{lemma}

{\bf Proof}. 
Let $h(z)$ be a "cap", ${\rm supp}h \subset\{|t|\le 1\}$, $\int h dV_t=1$. 
We shall consider the approximations
$$
  f_{(l),\delta}(z,t)=\delta^{-2n}\int f_{(l)}(z,\zeta )h\left(\frac{t-\zeta}{\delta}\right)dV_{\zeta}= \int f_{(l)}(z,t-\delta\zeta )h(\zeta )dV_{\zeta} .
 \eqno (4.9)
$$
Since $h(\zeta )=0$ at $|\zeta |\ge 1$, we obtain by Proposition 24
$$
  |f_{(l),\delta}(z)-f_{(l)}(z)|\le C\frac{\delta^{\beta_k}}{(1-|z|)^{s_k+|(l)|-1}} ,  \eqno (4.10)
$$
for some $\beta_k$ and $s_k$ if $|(l)|\le k$ and $\delta^{\beta_k}/(1-|z|)^{s_k} \le 1$.  

Now we make $\delta$ depending on $|z|$. Namely we pick some $0<b<1$ and define
$$
  \delta \le [bC^{-1}(1-|z|)^{s_{2m}+2m-1}]^{\frac{1}{\beta_{2m}}}  .
$$
Then estimate (4.10) yields
$$
  |f_{(l),\delta}(z)-f_{(l)}(z)|\le b      \eqno (4.11)
$$
if $|(l)|\le 2m$.

Introduce also approximations of the functions $f_{(l)}(0,t)$ for $|(l)|\le m$
$$
  G_{(l)_0 ,\delta}(t)= \int f_{(l)}(0,t-\delta\zeta )h(\zeta )dV_{\zeta} .
$$
As a particular case of (4.11) we have the estimates
$$ 
  |G_{(l)_0 ,\delta}(t)-f_{(l)}(0,t)|\le b .   \eqno (4.12)
$$ 
at $|(l)|\le m$.

Now we describe the construction of our approximation. In what follows $z=re^{-\theta}$. 
We adopt the notations $(l)=(j\bar k)$, $|(l)|=j+k$, $f_{j\bar k}=f_{z^j \bar z^k}$.

We define the functions $f_{(l),\delta}$ as the approximations of $f_{(l)}$ for $|(l)|=m$. Now,
if $j+k=m-1$, we put 
$$
  g_{(j\bar k),\delta}(z,t)=G_{(j\bar k )_0 ,\delta}(t)+\int_{[0,|z|]} f_{(j+1,\bar k ),\delta}dz +f_{(j\overline{k+1}),\delta}d\bar z =
$$
$$
   =G_{(j\bar k )_0 ,\delta}(t)+\int_0^r [e^{i\theta}f_{(j+1,\bar k ),\delta}(r,\theta ,t)+e^{-i\theta}f_{(j\overline{k+1}),\delta}(r,\theta ,t)]dt  .
$$
Analogously, define by induction at $g+k =q$, $0\le q\le m-2$, 
$$
  g_{(j\bar k),\delta}(z,t)=G_{(j\bar 0)_0 ,\delta}(z)+\int_{[0,|z|]} g_{(j+1,\bar k )\delta}dz +g_{(j,\overline{k+1}),\delta}d\bar z  .
$$
In particular,
$$
  f_{\delta}(z,t)=g_{(0\bar 0 ),\delta}(z,t)=\int_{[0,|z|]} g_{(1\bar 0 ),\delta}dz +g_{(0\bar 1 ),\delta}d\bar z .
$$
We prove that $f_{\delta}$ at small enough $d$ satisfies all conditions of the lemma. 

At first, applying (4.11) and (4.12), we see that, if at $j+k=q$ there holds the inequality $|g_{(j\bar k ),\delta}(z,t)-f_{(j\bar k )}(z,t)| \le m_q b$, then 
at $j+k=q-1$
$$
  |g_{(j\bar k ),\delta}(z,t)-f_{(j\bar k )}(z,t)|\le b+\int_{[0,|z|]} [|g_{(j+1,\bar k),\delta}-f_{j+1,\bar k}| +|g_{(j,\overline{k+1}),\delta}-f_{j,\overline{k+1}}|dr \le (2m_q +1)b. 
$$
It follows that $m_{q-1}\le 2m_q +1$, and by induction
$$
  |g_{(j\bar k ),\delta}(z,t)-f_{(j\bar k)}(z,t)|\le (3^{m-j-k}+1)b  .   \eqno (4.13)
$$ 
In particular,
$$
  |f_{\delta}-f|\le (3^m +1)b.  \eqno (4.14)    
$$

Show now that $z$-derivatives of $f_{\delta}$ approximate the corresponding $z$-derivatives of $f$ up to degree $m$.
In the calculations below it is
essential that $\delta$ depends only on $r$ and, hence, there don't appear "large" derivatives originating from
$\delta^{-2n}h\left(\frac{t-\zeta}{\delta}\right)$.
Consider at first $(f_{\delta})_z$. We have
$$
  (f_{\delta})_z (z,t)-f_z (z,t)=\frac{e^{-i\theta}}{2}\left(\frac{\pd}{\pd r}+\frac{1}{ir}\frac{\pd}{\pd\theta}\right)
\int_0^r [e^{i\theta}(g_{(1\bar 0),\delta}-f_z )+e^{-i\theta}(g_{(0\bar 1),\delta}-f_{\bar z})]dr =
$$
$$
  =\frac{1}{2}[g_{(1\bar 0),\delta}-f_z +e^{-2i\theta}(g_{(0\bar 1 ),\delta}-f_{\bar z})]+\frac{1}{2r}\int_0^r [g_{(1\bar 0),\delta}-f_z -e^{-2i\theta}(g_{(0\bar 1),\delta}-f_{\bar z})]dr +
$$
$$
 +\frac{1}{2ir}\int_0^r \left[\frac{\pd}{\pd\theta}(g_{(1\bar 0),\delta}-f_z )+e^{-2i\theta}\frac{\pd}{\pd\theta}(g_{(0\bar 1),\delta}-f_{\bar z})\right]dr .
\eqno (4.15)
$$
Thus, applying (4.13) at $j+k=1$, we obtain
$$
  |(f_{\delta})_z -f_z |\le 2(3^{m-1}+1)b+\sup[|(g_{(1\bar 0),\delta}-f_z )_{\theta}|+|g_{(0\bar 1),\delta}-f_{\bar z})_{\theta}|]  .  \eqno (4.16)
$$
Now,
$$
  (g_{(1\bar 0),\delta}-f_z )_{\theta}=i\int_0^r [e^{i\theta}(g_{(2\bar 0),\delta}-f_{2\bar 0})-e^{-i\theta}(g_{(1\bar 1),\delta}-f_{1\bar 1})]dr +
$$
$$
 +\int_0^r [e^{i\theta}(g_{(2\bar 0),\delta}-f_{2\bar 0})_{\theta}+e^{-i\theta}(g_{(1\bar 1),\delta}-f_{1\bar 1})-{\theta}]dr .  \eqno (4.17)
$$
Again applying (4.13), we see that
$$
  |(g_{(1\bar 0),\delta}-f_z )_{\theta}|\le 2(3^{m-2}+1)b+ \sup [|(g_{(2\bar 0),\delta}-f_{2\bar 0})_{\theta}|+|g_{(1\bar 1),\delta}-f_{1\bar 1})_{\theta}|]  . \eqno (4.18)
$$
Proceeding in the same way we obtain at last
$$
 |(g_{(m-1,\bar 0),\delta}-f_{m-1,\bar 0})_{\theta}|\le b+\sup [|(f_{(m\bar 0),\delta}-f_{m\bar 0})_{\theta}|+|f_{(m-1,\bar 1),\delta}-f_{m-1,\bar 1})_{\theta}|],  \eqno (4.19)
$$
and we have analogous estimates for other differences $(g_{(j\bar k),\delta}-f_{j\bar k})_{\theta}$, $j+k=m-1$. 
But, by definition of $f_{(m\bar 0),\delta}$,
$$
  [(f_{(m\bar 0),\delta}-f_{m\bar 0})_{\theta}](z,t)=\int [(f_{m\bar 0})_{\theta}(z,t-\delta\zeta )-(f_{m\bar 0})_{\theta}(z,t)]h(\zeta )dV_{\zeta}  .\,   \eqno (4.20)
$$
and, by (4.13) with $j+k=m+1$, 
$$
  |[(f_{(m\bar 0),\delta}-f_{m\bar 0})_{\theta}|\le b.  \eqno (4.21)
$$
Evidently, we have analogous estimates for other derivatives, appearing in the process.
Collecting (4.15) - (4.21), we obtain 
$$
     |(f_{\delta})_z -f_z |\le K_m b ,  
$$
where $K_m$ is some integer-valued function of $m$, which we don't specify here because
it is not essential for us. The analogous estimate we have for $|(f_{\delta})_{\bar z} -f_{\bar z}|$. 

Now consider the $z$-derivative of difference (4.15). We apply to the right part the operator $\frac{\pd}{\pd z}=\frac{e^{-i\theta}}{2}\left(\frac{\pd}{\pd r}+\frac{1}{ir}\frac{\pd}{\pd\theta}\right)$.
We have
$$
  g_{(1\bar 0),\delta}-f_z =\vpf (t)+\int_{[0,|z|]} (g_{(2\bar 0),\delta}-f_{2\bar 0})dw +(g_{(1\bar 1),\delta}-f_{1\bar 1})d\bar z =
$$
$$
  =\vpf (t)+\int_0^r [e^{i\theta}(g_{(2\bar 0),\delta}-f_{2\bar 0})+e^{-i\theta}(g_{(1\bar 1),\delta}-f_{1\bar 1})]dr ,
$$
where $\vpf (t)$ is the value at zero. When we apply the operator $\pd /\pd z$, we obtain the terms
$$
  e^{i\theta}(g_{(2\bar 0),\delta}-f_{2\bar 0})+e^{-i\theta}(g_{(1\bar 1),\delta}-f_{1\bar 1})] ,   \eqno (4.22)
$$
$$
  \frac{1}{2r}\int_0^r [g_{(2\bar 0),\delta}-f_{2\bar 0}-e^{-2i\theta}(g_{(1\bar 1),\delta}-f_{1\bar 1})]dr  \eqno (4.23)
$$
and
$$
  \frac{1}{2ir}\int_0^r [(g_{(2\bar 0),\delta}-f_{2\bar 0})_{\theta}+e^{-2i\theta}(g_{(1\bar 1),\delta}-f_{1\bar 1})_{\theta}]dr 
\eqno (4.24)
$$
We consider these terms exactly as we considered the right side of (4.15). We again reduce 
the problem to the estimate $|[(f_{(j\bar k,\delta}-f_{j\bar k})_{\theta}|\le Cb$ at $j+k=m+1$. 
In the same way we proceed with derivatives of the term $e^{-2i\theta}(g_{(0\bar 1),\delta}-f_{\bar z})$.
The derivative of the multiple $e^{-2i\theta}$ yields only the item $-2ie^{-2i\theta}(g_{(0\bar 1),\delta}-f_{\bar z})$.

At differentiation of the second term in (4.15)) (the first integral) we obtain the terms
$$
  \frac{e^{-i\theta}}{4r}[g_{(1\bar 0),\delta}-f_z -e^{-2i\theta}(g_{(0\bar 1),\delta}-f_{\bar z})] ,  \eqno (4.25)
$$
$$
\frac{e^{-i\theta}}{4r^2}\int_0^r [g_{(1\bar 0),\delta}-f_z \pm e^{-2i\theta}(g_{(0\bar 1),\delta}-f_{\bar z})]dr  \eqno (4.26)
$$
and
$$    
-\frac{e^{-i\theta}}{4ir^2}\int_0^r [(g_{(1\bar 0),\delta}-f_z )_{\theta} -e^{-2i\theta}(g_{(0\bar 1),\delta}-f_{\bar z})_{\theta}]dr 
\eqno (4.27)
$$
All our derivatives are regular at zero. It means that terms of the type $r^{-1}\vpf (t)$, originating from initial
values at zero, must annihilate. Subtracting these initial items we obtain for terms (4.25), (4.26) the estimates
$$
  \sup \left[\left|\frac{\pd}{\pd r}(g_{(1\bar 0),\delta}-f_z )\right|+\left|\frac{\pd}{\pd r}(g_{(0\bar 1),\delta}-f_z )\right|\right]\le
$$
$$
  \le C\sup (|g_{(2\bar 0),\delta}-f_{2\bar 0}| +|g_{(1\bar 1),\delta}-f_{1\bar 1}|+|g_{(0\bar 2),\delta}-f_{0\bar 2}|)\le C'b
$$
with some uniform $C'$ by definition of $g_{(j\bar k),\delta}$

Considering term (4.27) we must, analogously, estimate the mixed derivatives of $(g_{(1\bar 0),\delta}-f_z )_{r\theta}$ and
$(g_{(0,\bar 1),\delta}-f_{\bar z})_{r\theta}$. Estimate, for example, the first difference. By (4.17), we see that
we must estimate $g_{(2\bar 0),\delta}-f_{2\bar 0}$, $g_{(1\bar 1),\delta}-f_{1\bar 1}$,
$(g_{(2\bar 0),\delta}-f_{2\bar 0})_{\theta}$ and $(g_{(1\bar 1),\delta}-f_{1\bar 1})_{\theta}$. The first two differences
are of order $b$ by (4.13) and, proceeding as after (4.17), we can see that we can estimates these terms
through $|f_{(j,k),\delta}-f_{(j,k)}|$ with $j+k=m+1$, i.e., these terms also are of order $b$.  

When we apply operator $\pd /\pd w$ to the third term in (4.15) (the second integral) we obtain the terms
$$
  \frac{e^{-i\theta}}{4r}[(g_{(1\bar 0),\delta}-f_z )_{\theta}-e^{-2i\theta}(g_{(0\bar 1),\delta}-f_{\bar z})_{\theta}] ,  
$$
$$
\frac{e^{-i\theta}}{4r^2}\int_0^r [(g_{(1\bar 0),\delta}-f_z )_{\theta}\pm e^{-2i\theta}(g_{(0\bar 1),\delta}-f_{\bar z})_{\theta}]dr  ,
$$
and
$$    
-\frac{e^{-i\theta}}{4ir^2}\int_0^r [(g_{(1\bar 0),\delta}-f_z )_{\theta^2} -e^{-2i\theta}(g_{(0\bar 1),\delta}-f_{\bar z})_{\theta^2}]dr  .
$$
To estimate the first two terms we must estimate the $r$-derivatives of $(g_{(1\bar 0),\delta}-f_z )_{\theta}$ and
$(g_{(0\bar 1),\delta}-f_{\bar z})_{\theta}$. We already made it when we considered term (4.27). To estimate the last
term we must estimate the the $r$-derivatives of $(g_{(1\bar 0),\delta}-f_z )_{\theta^2}$ and $(g_{(0\bar 1),\delta}-f_{\bar z})_{\theta^2}$.
Differentiating (4.17) with respect to $r$ and $\theta$ we obtain the terms analogous to already considered and the term
$$
  e^{i\theta}((g_{(2\bar 0),\delta}-f_{2\bar 0})_{\theta^2}+e^{-i\theta}(g_{(1\bar 1),\delta}-f_{1\bar 1})_{\theta^2} 
$$
Again analogously to (4.19) - (4.21) we reduce estimation to the inequality 
$|[(f_{(j\bar k,\delta}-f_{j\bar k})_{\theta}|\le Cb$ at $j+k=m+2$.
 
The case of derivatives of higher degree is analogous.
Applying, for example, operator $(\pd /\pd z )^l$ to the right side of (4.15) we obtain the terms 
of the types
$$
  e^{is\theta}[g_{(j\bar k),\delta}-f_{(j\bar k )}],\,j+k=l+1,  \eqno (4.28)
$$  
$$
  e^{is\theta}r^{-q}[g_{(j\bar k),\delta}-f_{(j\bar k )}],\,q\le l,\,j+k=l+1-q,  \eqno (4.29)
$$
$$
  e^{is\theta}r^{-q}[g_{(j\bar k),\delta}-f_{(j\bar k )}]_{\theta^p},\,q\le l,\,p\le q,\,j+k=l+1-q,   \eqno (4.30)
$$
where $s$ is some integer depending on the term. Also, we obtain integrals of the types
$$
   e^{is\theta}r^{-q}\int_0^r [g_{(j\bar k),\delta}-f_{(j\bar k )}]dr,\,q\le l+1,\,j+k=l+2-q,  \eqno (4.31)
$$
and
$$
  e^{is\theta}r^{-q}\int_0^r [g_{(j\bar k),\delta}-f_{(j\bar k )}]_{\theta^p}dr,\,q\le l+1,\,p\le q,\,j+k=l+2-q  .  \eqno (4.32)
$$
Again the singular parts must annihilate and terms (4.28) -(4.30) have estimates
$$
  \sup\left|\frac{\pd^q}{\pd r^q}[g_{(j\bar k),\delta}-f_{(j\bar k )}]\right|\le\sum_{j+k=l+1}\sup|(g_{(j\bar k),\delta}-f_{j\bar k})|   \eqno (4.33)
$$
and 
$$
  \sum_{j+k=l+1}\sup|(g_{(j\bar k),\delta}-f_{j\bar k})_{\theta^p}|    \eqno (4.34)
$$
Integrals (4.31), (4.32) have the estimates analogous to (4.33), (4.34), only $q$ and $p$ could be equal to $l+1$.
Analogously to the considerations after (4.17) (see (4.18) - (4.21)) we reduce estimates for right sides of
(4.33), ( 4.34) to estimates (4.13) for $|f_{(j\bar k),\delta}-f_{j\bar k}|$ at $j+k\le m+q$, $q\le m$. The maximal degree $j+k=2m$
occurs for the terms
$$
  r^{-m}\int_0^r [g_{(j\bar k),\delta}-f_{(j\bar k )}]_{\theta^m}dr,\,j+k=m .
$$  
We proved estimate (4.5).
 
To show that the mapping $t\mapsto f_{at}$ is continuous it is enough to prove it for 
the mapping $t\mapsto f_{(l),\delta}(\cdot,t)$. By (4.9), this follows from continuity
of the mapping $t\mapsto f_{(l)}(\cdot,t)$. But the last follows from estimate (4.1) of
Proposition 24: for any $z\in D$ we have $f_{(l)}(z,t)\to f_{(l)}(z,t_0)$ as $t\to t_0$. 
 
Suppose now $m\ge 2$.
We have decomposition (4.6) with some holomorphic $h_{at}$. Consider the map 
$f_{at}\circ (f^{\mu_t})^{-1}=h_{at}\circ f^{\mu_{at}}\circ (f^{\mu_t})^{-1}$. Let $\mu_{ct}$ be the Beltramy coefficient 
of the map $f^{\mu_{at}}\circ (f^{\mu_t})^{-1}$. By uniqueness, the normal map $f^{\mu_{ct}}$ coincides with 
$f^{\mu_{at}}\circ (f^{\mu_t})^{-1}$. Denote $\tilde\mu_{at}=\mu_{at}\circ (f^{\mu_t})^{-1}\overline{((f^{\mu_t})^{-1})_z}/((f^{\mu_t})^{-1})_z$.
We have
$$
  \mu_{ct}=\frac{\tilde\mu_{at}+\mu_{(f^{\mu_t})^{-1}}}{1-\overline{\mu_{(f^{\mu_t})^{-1}}}\tilde\mu_{at}}=
  \frac{\mu_{at}\circ (f^{\mu_t})^{-1} -\mu_t \circ (f^{\mu_t})^{-1}}{1-\overline{\mu_{(f^{\mu_t})^{-1}}}\tilde\mu_{at}}\frac{\overline{((f^{\mu_t})^{-1})_z}}{((f^{\mu_t})^{-1})_z} 
$$
since
$$
 \mu_{(f^{\mu_t})^{-1}}=\mu\circ (f^{\mu_t})^{-1}\frac{f^{\mu_t}_z}{\overline{f^{\mu_t}_z}}\circ (f^{\mu_t})^{-1}=-\mu_t \circ (f^{\mu_t})^{-1}\frac{\overline{((f^{\mu_t})^{-1})_z}}{((f^{\mu_t})^{-1})_z}
$$ 
We see that, if we fix any $\varepsilon$, than at appropriate $b$ (i.e., at appropriate approximation of $f^{\mu_t}$)
$$
  |\mu_{ct}|\le\varepsilon ,\,|(\mu_{ct})_z |\le\frac{\varepsilon}{1-|z|},\,|(\mu_{ct})_{\bar z}|\le\frac{\varepsilon}{1-|z|},
$$
$$   
 |(\mu_{ct})_{(l)}|\le\frac{\varepsilon}{(1-|z|)^2},|(l)|=2   .
$$
By Lemma 1 for any $\varepsilon$ at appropriate $b$
$$
  |(f^{\mu_{ct}})_z -1|\le \varepsilon ,\,|(f^{\mu_{ct}})_ {(l)}|\le\frac{\varepsilon}{1-|z|},|(l)|=2  .
$$  
From the other hand, since $f_{at}$ approximate $f^{\mu_t}$ up to second derivatives, 
we obtain analogous estimates for $f_{at}\circ (f^{\mu_t})^{-1}=h_{at}\circ f^{\mu_{ct}}$. 
Thus we obtain estimates (4.7). 
Hence, $h_{at}$ is an univalent function ([Pom]) and $f_{at}$ is a homeomorphism.
     
At last consider estimate (4.8). It is enough to prove it for the derivatives of approximation (4.9). We have, for example,
$$
  (f_{(l),\delta})(z,t)_t=\frac{1}{\delta^{2n+1}}\int f_{(l)}(z,\zeta )h_t \left(\frac{t-\zeta}{\delta}\right)dV_{\zeta} .
$$
It follows that (see (4.10), (4.11))
$$
  |(f_{(l),\delta})(z,t)_t| \le\frac{C}{(1-|t|)^{\frac{(s_{2m}+2m-1)(2n+1)}{\beta_{2m}}}}
$$   
Also,
$$
   (f_{(l),\delta})(z,t)_z =-\frac{(2n+1)\delta_z}{\delta^{2n+2}}\int f_{(l)}(z,\zeta )h\left(\frac{t-\zeta}{\delta}\right)dV_{\zeta} +
$$
$$
  +\frac{1}{\delta^{2m+1}}\int \left[f_{(l)}(z,\zeta )_z h\left(\frac{t-\zeta}{\delta}\right)-h_z \left(\frac{z-\zeta}{\delta}\right)\frac{\delta_w}{\delta^2}\right]dV_{\zeta} .
$$
We again obtain an estimate of type (4.8). We obtain estimates for higher derivatives analogously.
$\Box$

\begin{corollary} The form $dz+\mu d\bar z$ after the change of the variable $z=f_{at}^{-1}(z_1 )$ and division by
$(f_{at}^{-1})_{z_1}$ transforms into the form $dz_1 +\tilde\mu d\bar z_1$ defined on $\cup\Omega_t$, where $\tilde\mu $ satisfies 
the estimates
$$
  |\tilde\mu |\le d,\,|(\tilde\mu )_{(k)}(z_1 )|\le\frac{b_{|(k)|}}{({\rm dist}(z_1 ,\pd\Omega_t))^{|(k)|}}, |(k)|\le m  ,
$$
where $d$ and $b_k$ can be made arbitrary small at appropriate approximation $f_{at}$ and we have the estimate
$$  
  |(\tilde\mu )_{(k,l)}(z_1)| \le \frac{C}{({\rm dist}(z_1 ,\pd\Omega_t))^{N_{|(k)|+|(l)|}}} 
$$
with some $N_{|(k)|+|(l)|}$, $|(k)|\le m$.
\end{corollary}

{\bf Proof}. It is obvious that there are the estimates
$$
  |\tilde\mu |\le d,\,|(\tilde\mu )_{(k)}|\le\frac{b_{|(k)|}}{(1-|f_{at}^{-1}(z_1 )|)^{|(k)|}}, |(k)|\le m  ,
$$
$$  
  |(\tilde\mu )_{(k,l)}(z_1)| \le \frac{C}{(1-|f_{at}^{-1}(z_1 )|)^{N_{|(k)|+|(l)|}}} .
$$
But
$$
  c_1 {\rm dist}(z_1 ,\pd\Omega_t)\le (1-|f_{at}^{-1}(z_1 )|\le c_2 {\rm dist}(z_1 ,\pd\Omega_t)
$$
for some uniform $c_1 ,c_2$. Indeed, there is decomposition $f_{at}=h_{at}\circ f_{\mu_{at}}$,
where for $f_{\mu_{at}}$ we have inequalities (2.4) and $h_{at}$ has the derivative close to 1.
$\Box$

As a result we obtained the important reduction. To prove Theorem 2 it is enough to prove
the next theorem:

{\bf Theorem 2'}. {\it Let $\Omega \subset B\times\bbC, B\subset\bbC^n$ be a domain fibered
by topological disks $\Omega_t, t\in B$ and suppose $\Omega_t =g_t (D)$, where $g_t$ is
a $\nu_t$-quasiconformal map with $|\nu_t|$ uniformly bounded away from 1, and the mapping
$t\mapsto g_t$ is continuous as a mapping from $B$ to $C^0 (D)$. 
We suppose that $g_t$ satisfies the estimates
$$
                    c\le |(g_t )_z (z)|\le C     \eqno (4.35)
$$
for some uniform $c,C$,
$$ 
   |(g_t )_{(k)}(z)|\le \frac{B}{(1-|z|)^{|(k)|-1}}     \eqno (4.36)
$$
at $|(k)|\le P$, $P\ge 4$,
$$
  |(g_t )_{(k,l)}(z) \le C(1-|z|)^{-N}    \eqno (4.37)  
$$
at $|(k)|\le P ,|(l)|\le L$ with some uniform $B$, $C$, and $N$. Also, we suppose
that there is a decomposition $g_t =h_t \circ f^{\nu_t}$, where $f^{\nu_t}$ is
the $\nu_t$-quasiconformal normal map and $h_t$ is a holomorphic univalent function
satisfying the estimates
$$
  |h'(z)-1|\le \varepsilon ,\,|h''(z)|\le \varepsilon (1-|z|)^{-1}  \eqno (4.38)
$$
with some uniform $\varepsilon$.

Let $\mu$ be a function on $\Omega$ satisfying the estimates  
$$
  |\mu |\le b <1 ,              \eqno (4.39)
$$
$$
   |\mu_{(k)}(z,t)| \le b({\rm dist}(z,\pd\Omega_t))^{-|(k)|} ,  \eqno (4.40)
$$   
at $|(k)|\le K$
$$
  |\pd^{(k,l)}\mu (z,t)|\le C({\rm dist}(z,\pd\Omega_t))^{-N}     ,    \eqno (4.41)
$$
at $|(k)|\le K , K\ge 4, |(l)|\le L$. The constant $C$ here and in (4.37) can depend on $B$
but the exponent $N$ doesn't depend.
 
Then, if $K\le P-L$ and the constant $b$ in (4.39), (4.40) is small enough, 
there exists a solution $f$ to the Beltramy equation
$$
  f_{\bar z}=\mu f_z  ,
$$  
which is continuous in $C_0 (\bbC )$ as a function of $t$, is finitely smooth with respect to 
all variable up some be-degree $\{p,q\}$, where $p$ and $q$ can be arbitrary large if $K$ and $L$ are large enough, 
at every $t$ maps $\Omega_t$ homeomorphically onto some bounded subdomain of $\bbC$,
and satisfies the estimates
$$
  c({\rm dist}(z,\pd\Omega_t))^{\alpha}\le|f_z (z,t)|\le C({\rm dist}(z,\pd\Omega_t))^{-\alpha},\,|f_{\bar z}(z)|\le C({\rm dist}(z,\pd\Omega_t))^{-\alpha},
$$  
$$  
  0\ge\alpha <1 , \eqno (4.42)
$$
$$
  |f_{(k)}(z)|/|f_z (z.t)|\le C({\rm dist}(z,\pd\Omega_t))^{1-|(k)|} ,  \eqno (4.43)
$$
at $|(k)|\le p$,
$$
   |\pd^{(k,l)}f (z,t)|\le C({\rm dist}(z\pd\Omega_t))^{-M}    \eqno (4.44)
$$
for some uniform $C$ and $M$ at $|(k)|\le p,|(l)|\le q$.} 

The conditions $P\ge 4$, $K\ge 4$ are of technical character, we shall use them in Section 7.

\section{Extension of quasiconformal mappings}

In this section we consider a family of quasiconformal mappings $g_t$ satisfying estimate (4.35) -(4.38) of
Theorem 2'. Mostly we have deal with an individual map $g$ and we shall omit dependence on $t$.
Thus $g$ is a $\nu$-quasiconformal map, mapping $D$ onto the domain $\Omega$ and   
$$
   |\nu (z)|\le b<1 ,   \eqno (5.1)
$$
$$ 
   |\nu_{(k)}(z)|\le \frac{B}{(1-|z|)^{|(k)|}}    \eqno (5.2)
$$
at $|(k)|\le P$,
$$
  |\nu_{(k,l)}(z,t)| \le C(1-|z|)^{-N}    \eqno (5.3) 
$$
at $|(k)|\le P ,|(l)|\le L$ with some uniform $B$ and $C$. The last two inequalities follow
from (4.35) -(4.37). There is the decomposition $g=h\circ f^{\nu}$ and, in addition to 
estimates (4.38), from lemma 1 and (4.36) follows
$$
  \left|\frac{d^k}{dz^k}h(z)\right| \le C(1-|z|)^{k-1}    \eqno (5.4)
$$
at $k\le P$ with some uniform $C$. Also, we have the obvious estimates 
$$  
  a(1-|z|)\le{\rm dist}(g(z), \pd\Omega )\le A(1-|z|)    \eqno (5.5)
$$
for some uniform $a,A$.

In what follows we shall need estimates for derivatives of the normal mappings in the 
particular case when all derivatives of $\mu$ are uniformly bounded.
]
\begin{proposition}
Suppose $\mu$ is smooth, has the support in $D$, and smoothly depends on a vector parameter $t\in B\subset\bbC^n$. 
Suppose the derivatives $\mu_{(k,l)}$ satisfy the estimate
$$
   |\mu_{(k,l)}|\le M  
$$
at $|(k)|\le K+L$, $|(l)|\le L$ with some constant $M$ uniform with respect to the parameters. 
Then the derivatives of the normal mapping $(f^{\mu})_{(k,l)}$ satisfy the estimate
$$
   |(f^{\mu})_{(k,l)}|\le CM^{10(|(k)|+|(l)|)^2} 
$$ 
at $|(k)|\le K, |(l)|\le L$ with some uniform $C$. 
\end{proposition}

{\bf Proof}. 
At first we estimate the derivatives of the principal solution $f_{\mu}$. Consider the 
derivative of first order in $t$.
Namely we must estimate $t$-derivatives of the function ${\mC} h_{\mu}$, where $h_{\mu}$ is the solution
to the equation
$$
  h-\mu{\mS}h=\mu  .    \eqno (5.6)
$$
Differentiating by $t$ we obtain
$$
  h_t -\mu{\mS}h_t =\mu_t {\mS}h+\mu_t .  
$$
The $p$-norm of the right side is no greater than $CM$ and for $h_t$ we obtain the estimate
$\|h_t \|_p \le C_p CM/(1-d)$. It follows that for $|f_{\mu}|$ we also have the estimate $CM$
with some uniform $C$.

At further differentiation we obtain the equation
$$
  h_{0,(l)} -\mu{\mS}h_{0,(l)} =F,  
$$
where for $F$ we obtain by induction the estimate $\|F\|_p \le CM^{|(l)|}$. Thus 
$|(f_{\mu})_{(0,l)}|\le CM^{|(l)|}$.

Now remind that for a smooth compactly supported $\mu$ we can represent the function
$(f_{\mu})_z$ as $e^h$, where $h$ satisfies the equation
$$
  h_{\bar z}=\mu h_z +\mu_z    \eqno (5.7)
$$  
and tends to zero at infinity, i.e., $h=\mC H$, where $H$ is the unique solution to the equation
$$
  H-\mu{\mS}H=\mu_z  .
$$
For $H_t$ we obtain the equation
$$
   H_t -\mu{\mS}H_t =\mu_t \mS H+\mu_{zt}  .     \eqno (5.8)
$$
Since $\|\mu_t \mS h\|_p \le CM^2$, we obtain the estimate $\|H_t \|_p \le CM^2$. 
For $H_{(0,l)}$ we obtain by induction the equation analogous to (5.8) with the right side $F$
such that $\|F\|_p \le CM^{|(l)+1}$ and, hence, the same estimate for 
$|((f_{\mu})_z )_{(0,l)}|=|(e^{\mC H})_{(0,l)}|$.

Now, $(f_{\mu})_{z^2}=e^h h_z =e^h g$, where $g$ satisfies the equation obtained by differentiation
of equation (5.7)
$$
  g_{\bar z}-\mu g_z =\mu_z g+\mu_{z^2} .   
$$  
But $f$ is holomorphic outside of $D$ and $f(z)-z$ tends to zero when $z\to\infty$.
Hence, all derivatives of order higher than one also tend to zero at infinity and the same
holds for the function $g$. It implies that $g=\mC G$, where $G$ is the unique solution to the equation
$$
  G-\mu\mS G=\mu_z g+\mu_{z^2}  .   
$$
For $G_t$ we obtain the equation
$$
  G_t-\mu\mS G_t=\mu_t \mS G +\mu_{zt}g+\mu_z g_t+\mu_{z^2 t}  \eqno (5.9)
$$
Since $\|\mu_z g\|_p =\|\mu_z h_z\|_p \le CM^2$, we see that $\|G\|_p \le CM^2$,
and the right side of equation (5.9) has the $L_p$-estimate $CM^3$. At further differentiation
by parameters we obtain for $G_{(0,l)}$ the equation with the right part estimated in $L^p$
as $CM^{|(l)|+2}$. 
 
When we consider derivatives with respect to $z$ of higher order we analogously can see
that we obtain equations with right sides having the estimate $CM^s$ with the exponent
rising by 1 at each differentiation. As a result, we obtain the estimate
$$
  |(f_{\mu})_{(k,l)}|\le CM^{|(k)|+|(l)|}  .  \eqno (5.10) 
$$    
 
Consider now the normal solution $f=f^{\mu}$. By (2.58), (2.59), we can represent $f_{(k,l)}$ as
a sum of terms of the types 
$$
  a(f_{\lambda})_{(0,p)}(1)-f_{lambda})_{(0,p)}(0))\tilde f_{(0,q)}(1) (f_{\lambda})_{(r,s)}\circ \tilde f
  \tilde f_{(k_1 ,l_1 )}^{r_1}...\tilde f_{(k_j, l_j )}^{r_j} ,   \eqno (5.11)
$$
where $r\le |(k)|+|(l)|$, $p+q+s\le |(l)|$, $r_1 (|(k_1 )|+|(l_1 )|)+...+r_j (|(k_j )|+|(l_j )|)\le |(k)|+|(l)|$. 
For the product $\tilde f_{(k_1 ,l_1 )}^{r_1}...\tilde f_{(k_j, l_j )}^{r_j}$
according (5.10) we have the estimate $CM^{|(k)|+|(l)|}$. From the other hand, $\lambda_{(k,l)}$ can be
represented as a sum of items of the types
$$
  \mu_{(p,q)}\circ\tilde f^{-1}(\tilde f_z )_{r,s}\circ\tilde f^{-1}(\tilde f^{-1}_{(k_1 ,l_1 )})^{r_1}...(\tilde f^{-1}_{(k_j, l_j )})^{r_j}[.] ,
$$  
where $[.]$ is a multiple bounded by a constant independent of $(k,l)$, 
$p+r\le |(k)|+|(l)|$, $q+s+r_1 (|(k_1 )|+|(l_1 )|)+...+r_j (|(k_j )|+|(l_j )|)\le |(k)|+|(l)|$.  
But $|(\tilde f^{-1}_{(k_1 ,l_1 )})^{r_1}...(\tilde f^{-1}_{(k_j, l_j )})^{r_j}|\le CM^{|(k)|+|(l)|}$,
$|\mu_{(p,q)}|\le CM$, $|(\tilde f_z )_{r,s}|\le CM^{r+s+1}$ according to (5.10). We obtain 
$|\lambda_{(k,l)}|\le CM^{3(|(k)|+|(l)|)+2}$. We apply estimate (5.10) to $f_{\lambda}$ and 
obtain the estimates $|(f_{\lambda})_{(r,s)}|\le CM^{2(|(k)|+|(l)|)[3(|(k)|+|(l)|+2]}\le CM^{7(|(k)|+|(l)|)^2}$.
$|f_{\lambda})_{(0,p)}|\le CM^{|(l)|[3(|(k)|+|(l)|)+2]}$ for the terms in (5.11). As a result, we obtain for product
(5.11) the estimate $CM^{10(|(k)|+|(l)|)^2}$. $\Box$ 
 
Now we shall motivate the following constructions of this section. When we defined the transforms
${\mC}_m$ and ${\mS}_m$ in Section 1, we introduced into the kernels the counter-items of the types
$$ 
   \frac{(\bar\zeta^{-1}-\zeta )^{l-1}}{(\bar\zeta^{-1}-z)^l}  . \eqno (5.12)
$$  
  
We want to solve the Beltrami equation on the domain $\Omega$, and we need to define analogous 
counter-items to neutralize the growth of derivatives near the boundary. Namely we must find
some map replacing the mapping $z\mapsto\bar z^{-1}$ when we deal with the domain $\Omega$
instead of the disk $D$. For that we define some extension of the map $g$ on the domain 
$\bbC\setminus\Omega$. It appears, we can find a sufficiently good extension if $\varepsilon$ 
in (4.38) is sufficiently small. 

\begin{lemma}. 
Let the family of maps $g_t$ satisfy conditions of Theorem 2'. If $\varepsilon$ is small
enough we can define for each $t$ an extension of $g_t$ to a quasyconformal homeomorphisms of the plane
$G_t$; we denote by $\hat g_t$ its restriction on $\bbC\setminus D$. For the map $\hat g_t$
we have the estimates
$$
  \hat c_1  \le |\hat g_z (z,t)|\le \hat c_2  ,  \eqno (5.13)
$$
$$
  |\hat g_{(k)}(z,t)|\le B(|z|-1)^{1-|(k)|}   \eqno (5.14)
$$
at $|(k)|\le P-L$,  
$$
  |\hat g_{(k),(l)}(z,t)|\le C(|z| -1)^{-M}    \eqno (5.15)
$$
with some $M$ depending only on $N$ in inequality (4.37) at $|(k)|\le P-L,|(l)|\le L$, 
$$
  \hat a(|z|^2 -1)\le{\rm dist}(\hat g_(z,t), \pd\Omega_t )\le  \hat A(|z|^2 -1)  , \eqno (5.16)
$$  
$$
  c(1-|z|^2 )\le |g(z,t)-\hat g (\bar z^{-1},t)|\le C(1-|z|^2 )  .  \eqno (5.17)
$$
All constants in these inequalities are uniform. 
\end{lemma}

{\bf Proof}.
In most part of the proof we fix some $t$ and omit $t$-dependence.

Let $r_n >1/2$, $n\ge 1$ be some sequence tending to 1 and consider the sequence $D_{r_n}$ of disks of radii $r_n$ 
centered at zero. Let $\eta$ be a smooth function on the real axis $\eta (x) =1 ,x\le 0$,
$\eta (x) =0 ,x\ge 1$. Let $\nu_n$, $h_n$  be the functions $\nu_n (z)=\nu (z)\eta ((|z|-r_n )/(1-r_n))$,  
$h_n (z)=h(r_n z)$. Let $f^{\nu_n}$ be the normal $\nu_n$-quasyconformal homeomorphism mapping 
$D$ onto itself and $g_n$ be the map $g_n =h_n \circ f^{\nu_n}$.
The homeomorphism $g_n$ maps $D$ onto some domain $\Omega _n$ with a smooth boundary
and the sequence $g_n$ converges to $g$ uniformly on compact subsets of $D$. 

On any disk $D_{r_n}$ the map $g|_{D_{r_n}}$ can be represented as the composition $g=f_n \circ f^{\nu_n}$,
where $f_n$ is some function holomorphic on $f^{\nu_n}(D_{r_n})$. 

From (5.1) - (5.4) and (4.38) immediately follow the estimates 
$$
  |\nu_n (z)|\le b<1 ,  \eqno (5.18)
$$
$$
   |(\nu_n )_{(k)}(z)|\le B(1-r_n |z|)^{-|(k)|}   \eqno (5.19)
$$
at $|(k)|\le m$,
$$
  |(\nu_n )_{(k,l)}(z,t)| \le C(1-|r_n |)^{-N}    \eqno (5.20) 
$$
at $|(k)|\le m ,|(l)|\le L$,
$$
  |h'_n (z) -1|\le 1-r_n +2\varepsilon ,\,|h''_n (z)| \le \varepsilon (1-r_n^2 |z|^2 )^{-1}, \eqno (5.21)
$$
$$
  \left|\frac{d^k}{dz^k}h_n (z)\right| \le C(1-r_n |z|)^{k-1}   \eqno (5.22)
$$
at $k\le m$. Also, we have the estimate analogous to (4.35), (4.36), and (5.5)
$$
  c_1 \le|(g_n)_z |\le c_2  ,   \eqno (5.23)
$$
$$
  | (g_n )_{(k)}(z)|\le  B(1-|z|)^{1-|(k)|}    \eqno (5.24)
$$
at $|(k)|\le m$,
$$
  a(1-|z|^2 )\le{\rm dist}(g_n (z), \pd\Omega_n )\le A(1-|z|^2 )  . \eqno (5.25)
$$  
All constants in these estimates are uniform and independent of $n$ and $r_n$. Below in this 
section "uniform" means, in particular, that estimate or constant is independent of $n$ and $r_n$.

From Proposition 25 we obtain 
$$
  |f^{\nu_n}_{(k,l)}| \le C(1-|r_n |)^{-M}  ,    \eqno (5.26) 
$$
$$
  |(f_n )_{(k,l)}\circ f^{\nu_n}| \le C(1-|r_n |)^{-M} \eqno (5.27)  
$$
with some $M$ independent of $n$ at $|(k)|\le m-L, |(l)|\le L$.

\begin{proposition} 
 There exists some uniform $B$ such that at $1-|z|\ge B(1-r_m)$ 

 a)
$$
 |f^{\nu_n}_z (z) -f^{\nu_m}_z (z)| \le C\frac{(1-r_m )^{\beta}}{1-|z|}  .   \eqno (5.28)
$$ 
$$
|f^{\nu}_z (z) -f^{\nu_m}_z (z)| \le C\frac{(1-r_m )^{\beta}}{1-|z|}  .   \eqno (5.29)
$$ 

b)
$$
  |h'_m \circ f^{\nu_m}(z)-f'_m \circ f^{\nu_m}(z)|\le C\frac{(1-r_m )^{\beta}}{1-|z|} ,  \eqno (5.30)
$$
$$
  |h''_m \circ f^{\nu_m}(z)-f''_m \circ f^{\nu_m}(z)|\le\frac{(1-r_m )^{\beta}}{(1-|z|)^2}  .  \eqno (5.31)
$$ 

c)
$$
  |f'_n \circ f^{\nu_n}(z)-f'_m \circ f^{\nu_m}(z)|\le C\frac{(1-r_m )^{\beta}}{1-|z|} ,  \eqno (5.32)
$$
$$
    |f''_n \circ f^{\nu_n}(z)-f''_m \circ f^{\nu_m}(z)|\le\frac{(1-r_m )^{\beta}}{(1-|z|)^2}   .  \eqno (5.33)
$$ 
Everywhere $\beta >0$ depends only on $b$ in (5.1).
\end{proposition}

{\bf Proof}.  a) The second inequality is a limit case of the first one. 
We apply inequality (3.2) of Lemma 2 to estimate $|f^{\nu_n}_z -f^{\nu_m}_z |$. 
We estimate different terms in the right-side. At first consider $\|\nu_n -\nu_m \|_p$.

The functions $\nu_n$ and $\nu_m$ differ only in the ring $R_{m,n}=\{r_m \le |z|\le r_n +(1-r_n )/2=(1+r_n )/2\}$. Hence,
$$
  \|\nu_n -\nu_m \|_p =\left(\int_{R_{m,n}} dS_z\right)^{1/p} \le C_p (1+r_n -2r_m )^{1/p} \le
$$
$$
  \le C_p |1-r_m |^{1/p} ,\,n>m  ,  \eqno (5.34) 
$$
with some $C_p$ depending on $p$. Thus,
$$
   \|\nu_n -\nu_m \|_p^{\alpha}\le C|1-r_m |^{\alpha/p}  
$$
Now, since $\nu_n -\nu_m =0$ on $D_{r_m}$, we obtain putting $R=r_m$ 
$$
  \sup_{D_R}|\nu_n -\nu_m |^{\alpha}+\left(\frac{1-R}{1-|z|}\right)^{\alpha /p}=
  \left(\frac{1-r_m}{1-|z|}\right)^{\alpha /p}   .
$$ 
At last, $\sup_{D_{R_z}}|(\nu_n )_{(l)}-(\nu_m )_{(l)}|^{\alpha}=0$ for $z$ such that $R_z \le r_m$.
Thus we obtain (5.28) with $\beta=\alpha /p$.    

b) Let prove inequality (5.30). 
$$
  |h'_m \circ f^{\nu_m}-f'_m \circ f^{\nu_m}|\le |h'_m \circ f^{\nu_m}(z)-h'\circ f^{\nu_m}(z)|+
  |h'\circ f^{\nu_m}-h'\circ f^{\nu}|+|h'\circ f^{\nu}-f'_m \circ f^{\nu_m}| .
$$  
But
$$
  |h'_m (z)-h'(z)|=|r_m h'(r_m z)-h'(z)|\le (1-r_m )|h'(r_m z)|+|h'(r_m z)-h'(z)|\le
$$
$$
  \le c(1-r_m )+\varepsilon\frac{|r_m z -z|}{1-|z|}\le C\frac{1-r_m}{1-|z|} .
$$
Here we applied estimates (4.38). Hence,
$$
  |h'_m \circ f^{\nu_m}(z)-h'\circ f^{\nu_m}(z)|\le C\frac{1-r_m}{1-|f^{\nu_m}(z)|}\le C\frac{1-r_m}{1-|z|}  \eqno (5.35)
$$  
according estimates (2.4).

Also, applying (4.38), (2.4), estimate (3.1) of Lemma 2, and (5.34),
we have
$$
  |h'\circ f^{\nu_m}(z)-h'\circ f^{\nu}(z)|\le \sup_{[f^{\nu_m}(z),f^{\nu}(z)]}|h''||f^{\nu_m}(z)-f^{\nu}(z)|\le
  \frac{c}{1-|z|}(1-r_m )^{\beta}   \eqno (5.36)
$$
for some $\beta$ depending only on $b$ in (5.1). 

At last, 
$$
  h\circ f^{\nu}=f_m \circ f^{\nu_m}   \eqno (5.37)
$$
on $D_{r_m}$ and, hence,
$h'\circ f^{\nu}\cdot (f^{\nu})_z =f'_m \circ f^{\nu_m}\cdot (f^{\nu_m})_z$. Thus,
$$
  |h'\circ f^{\nu}-f'_m \circ f^{\nu_m}|=|h'\circ f^{\nu}||1-(f^{\nu})_z /(f^{\nu_m})_z |\le 
 C(1-r_m )^{\beta}/(1-|z|) \eqno (5.38)
$$
by (5.29), if $1-|z|\ge B(1-r_m)$. Taking into consideration (5.35) and (5.36), we obtain (5.30).

The proof of inequality (5.31) is analogous. We have
$$
  |h''_m \circ f^{\nu_m}-f''_m \circ f^{\nu_m}|\le |h''_m \circ f^{\nu_m}(z)-h''\circ f^{\nu_m}(z)|+
  |h''\circ f^{\nu_m}-h''\circ f^{\nu}|+|h''\circ f^{\nu}-f''_m \circ f^{\nu_m}| .
$$  
Applying (4.38) and (5.4), we can see that
$$
  |h''_m (z)-h''(z)|=|r^2_m h''(r_m z)-h''(z)|\le (1-r^2_m )|h''(r_m z)|+|h''(r_m z)-h''(z)|\le
$$
$$
  \le c\varepsilon\frac{1-r_m}{1-r_m |z|}+C\frac{|r_m z -z|}{(1-|z|)^2}\le C\frac{1-r_m}{(1-|z|)^2} .
$$
Hence,
$$
  |h''_m \circ f^{\nu_m}(z)-h''\circ f^{\nu_m}(z)|\le C\frac{1-r_m}{(1-|f^{\nu_m}(z)|)^2}\le C\frac{1-r_m}{(1-|z|)^2}  \eqno (5.39) 
$$  

Also, applying (5.4), (2.4), estimate (3.1) of Lemma 2, and (5.34), we obtain
$$
  |h''\circ f^{\nu_m}(z)-h''\circ f^{\nu}(z)|\le \sup_{[f^{\nu_m}(z),f^{\nu}(z)]}|h'''||f^{\nu_m}(z)-f^{\nu}(z)|\le
  \frac{C}{(1-|z|)^2}(1-r_m )^{\beta}   \eqno (5.40)  
$$
for some $\beta$ depending only on $b$ in (5.1). 

As above, from (5.37) we have 
$h''\circ f^{\nu}\cdot ((f^{\nu})_z )^2 +h'\circ f^{\nu}\cdot (f^{\nu})_{z^2}=f''_m \circ f^{\nu_m}\cdot ((f^{\nu_m})_z )^2 +f'_m \circ f^{\nu_m}\cdot ((f^{\nu_m})_{z^2}$. Thus,
$$
  |h''\circ f^{\nu}-f''_m \circ f^{\nu_m}|\le |h''\circ f^{\nu}|\left|1-\frac{((f^{\nu})_z )^2}{((f^{\nu_m})_z )^2}\right| +
$$
$$  
  +|(f^{\nu_m})_z |^{-2}(|h'\circ f^{\nu}-f'_m \circ f^{\nu_m}||(f^{\nu})_{z^2}|+
  |f'_m \circ f^{\nu_m}||(f^{\nu})_{z^2}-(f^{\nu_m})_{z^2}|   .
$$
Consider the right side. 
$$
  |h''\circ f^{\nu}(z)|\left|1-\frac{((f^{\nu})_z )^2}{((f^{\nu_m})_z )^2}\right|\le C{\varepsilon}\frac{(1-r_m )^{\beta}}{(1-|z|)^2}
$$
by (4,38), (2.4) and (5.29).  
$$
  |h'\circ f^{\nu}(z)-f'_m \circ f^{\nu_m}(z)||(f^{\nu})_{z^2}|\le C(1-|z|)^{-2}(1-r_m )^{\beta}  
$$
by Lemma 1 and (5.38). At last, the difference $|(f^{\nu})_{z^2}-(f^{\nu_m})_{z^2}|$ we 
estimate by inequality (3.2) of Lemma 2. From, in fact, the same considerations as in the proof
of inequality (5.28), we can see that
$$  
 |(f^{\nu})_{z^2}-(f^{\nu_m})_{z^2}|\le C(1-r_m )^{\beta}(1-|z|)^{-2}
$$
Thus,
$$
  |h''\circ f^{\nu}-f''_m \circ f^{\nu_m}|\le C(1-r_m )^{\beta}(1-|z|)^{-2} .  \eqno (5.41)
$$
Collecting inequalities (5.39) - (5.41) we obtain (5.31). 

c) We have 
$$
  |f'_n \circ f^{\nu_n}-f'_m \circ f^{\nu_m}|\le |f'_n \circ f^{\nu_n}-h'\circ f^{\nu}|+
  |h'\circ f^{\nu}-f'_m \circ f^{\nu_m}|  ,
$$
$$
    |f''_n \circ f^{\nu_n}-f''_m \circ f^{\nu_m}|\le |f''_n \circ f^{\nu_n}-h''\circ f^{\nu}|+
  |h''\circ f^{\nu}-f''_m \circ f^{\nu_m}|   .    
$$
The required estimates follow from (5.38) and (5.41).
$\Box$.

Now we shall prove the lemma itself. We proceed in several steps. 

1) {\it First extension of the map $g_n$}. 

We use the modified construction of the Loewner chains from the theory of univalent functions.   
For a function $f(z,t),\,(z\in \bbC\setminus D, 0\le t<\infty )$ we write $\dot f =\pd f/\pd t$. 

{\it We say that the family of functions $f(z,t)=f_t (z)$ ($z\in \bbC\setminus D, 0\le t<\infty )$ is a Loewner chain if
there exists a function $p(z,t),\,(z\in \bbC\setminus D, 0\le t<\infty )$ such that ${\rm Re}p(z,t)>0$ and}
$$
  f_t (z,t)=[zf_z (z,t)-\bar z f_{\bar z}(z,t)]p(z,t)  .  \eqno (5.42)
$$
For any $z, |z|=r$ vector $zf_z (z,t)-\bar z f_{\bar z}(z,t)$ is the vector of the outer normal at a boundary point of the set 
$\{f_t (z):|z|\le r\}$. Condition  (5.42) means that 
$$
  |\arg f_t (z,t)-\arg [zf_z (z,t)-\bar z f_{\bar z}(z,t)]|=|\arg p(z,t)|<\pi /2  ,    \eqno (5.43)
$$ 
and this means that the velocity vector $f_t$ on the boundary points out of this set. 

Instead of conditions (5.42), (5.43) we shall use the equivalent condition
$$
  \left|\frac{p-1}{p+1}\right|=\left|\frac{f_t -zf_z +\bar z f_{\bar z}}{f_t +zf_z -\bar z f_{\bar z}}\right|<1  . \eqno (5.44)
$$
In what follows we denote by $D_{\pm}$ the operators
$$
  D_{\pm}=\frac{\pd}{\pd t}\pm \left(z\frac{\pd}{\pd z}-\bar z\frac{\pd}{\pd\bar z}\right)  .
$$  
For $z:|z|\ge 1$ define the function $\hat f^{\nu_n}(z)=1/\overline{f^{\nu_n}(1/\bar z )}$. It is the symmetrical extension      
of $f^{\nu_n}$. 

Now define for $z:|z|\ge 1$ 
$$
  \hat h_n (z,t) =h_n \circ f^{\nu_n}(e^{-t}/\bar z )+[\hat f^{\nu_n}(e^t /\bar z )-f^{\nu_n}(e^{-t}/\bar z )]h'_n \circ f^{\nu_n}(e^{-t}/\bar z )  . 
$$
Let check for this function condition (5.44) at $|e^t /\bar z|\le 2$. We denote by $w$ the chart in the image of maps $f^{\nu_n}$ and $\hat f^{\nu_n}$ and by $\omega$ the
chart in the preimage of this map. We have
$$
  D_- \hat h_n (z,t)=
$$
$$
=h_{n,w}\circ f^{\nu_n}\left(\frac{e^{-t}}{\bar z}\right)\left[-f^{\nu_n}_{\omega}\left(\frac{e^{-t}}{\bar z}\right)\frac{e^{-t}}{\bar z}-f^{\nu_n}_{\bar\omega}\left(\frac{e^{-t}}{\bar z}\right)\frac{e^{-t}}{z}+
f^{\nu_n}_{\bar\omega}\left(\frac{e^{-t}}{\bar z}\right)\frac{e^{-t}}{z}-f^{\nu_n}_{\omega}\left(\frac{e^{-t}}{\bar z}\right)\frac{e^{-t}}{\bar z}\right]+
$$
$$
+\left[\hat f^{\nu_n}_{\omega}\left(\frac{e^t}{\bar z}\right)\frac{e^t}{\bar z}+\hat f^{\nu_n}_{\bar\omega}\left(\frac{e^t}{\bar z}\right)\frac{e^t}{z}+
\hat f^{\nu_n}_{\bar\omega}\left(\frac{e^t}{\bar z}\right)\frac{e^t}{z}-\hat f^{\nu_n}_{\omega}\left(\frac{e^t}{\bar z}\right)\frac{e^t}{\bar z}\right. +
$$
$$
+\left.f^{\nu_n}_{\omega}\left(\frac{e^{-t}}{\bar z}\right)\frac{e^{-t}}{\bar z}+f^{\nu_n}_{\bar\omega}\left(\frac{e^{-t}}{\bar z}\right)\frac{e^{-t}}{z}
 -f^{\nu_n}_{\bar\omega}\left(\frac{e^{-t}}{\bar z}\right)\frac{e^{-t}}{z}+f^{\nu_n}_{\omega}\left(\frac{e^{-t}}{\bar z}\right)\frac{e^{-t}}{\bar z}\right]
h'_n \circ f^{\nu_n}\left(\frac{e^{-t}}{\bar z}\right)+ 
$$
$$
 +\left[\hat f^{\nu_n}\left(\frac{e^t}{\bar z}\right)-f^{\nu_n}\left(\frac{e^{-t}}{\bar z}\right)\right]
 h''_n \circ f^{\nu_n}\left(\frac{e^{-t}}{\bar z}\right)\times
$$
$$
\times\left[-f^{\nu_n}_{\omega}\left(\frac{e^{-t}}{\bar z}\right)\frac{e^{-t}}{\bar z}-f^{\nu_n}_{\bar\omega}\left(\frac{e^{-t}}{\bar z}\right)\frac{e^{-t}}{z}+
f^{\nu_n}_{\bar\omega}\left(\frac{e^{-t}}{\bar z}\right)\frac{e^{-t}}{z}-f^{\nu_n}_{\omega}\left(\frac{e^{-t}}{\bar z}\right)\frac{e^{-t}}{\bar z}\right]=
$$
$$
 =2h'_n \circ f^{\nu_n}\left(\frac{e^{-t}}{\bar z}\right)\hat f^{\nu_n}_{\bar\omega}\left(\frac{e^t}{\bar z}\right)\frac{e^t}{z}- .
$$
$$
 -2\left[\hat f^{\nu_n}\left(\frac{e^t}{\bar z}\right)-f^{\nu_n}\left(\frac{e^{-t}}{\bar z}\right)\right]h''_n \circ f^{\nu_n}\left(\frac{e^{-t}}{\bar z}\right)f^{\nu_n}_{\omega}\left(\frac{e^{-t}}{\bar z}\right)\frac{e^{-t}}{\bar z}
$$
$$
  D_+ \hat h_n (z,t)=-2h'_n \circ f^{\nu_n}\left(\frac{e^{-t}}{\bar z}\right)f^{\nu_n}_{\bar\omega}\left(\frac{e^{-t}}{\bar z}\right)\frac{e^{-t}}{z}+
$$
$$
+\left[2\hat f^{\nu_n}_{\omega}\left(\frac{e^t}{\bar z}\right)\frac{e^t}{\bar z}+2f^{\nu_n}_{\bar\omega}\left(\frac{e^{-t}}{\bar z}\right)\frac{e^{-t}}{z}\right]
h'_n \circ f^{\nu_n}\left(\frac{e^{-t}}{\bar z}\right)- 
$$
$$
-2\left[\hat f^{\nu_n}\left(\frac{e^t}{\bar z}\right)-f^{\nu_n}\left(\frac{e^{-t}}{\bar z}\right)\right]h''_n \circ f^{\nu_n}\left(\frac{e^{-t}}{\bar z}\right)f^{\nu_n}_{\bar\omega}\left(\frac{e^{-t}}{\bar z}\right)\frac{e^{-t}}{z}=
$$
$$
 =2\hat f^{\nu_n}_{\omega}\left(\frac{e^t}{\bar z}\right)\frac{e^t}{\bar z}h'_n \circ f^{\nu_n}\left(\frac{e^{-t}}{\bar z}\right)- 
$$
$$
 -2\left[\hat f^{\nu_n}\left(\frac{e^t}{\bar z}\right)-f^{\nu_n}\left(\frac{e^{-t}}{\bar z}\right)\right]h''_n \circ f^{\nu_n}\left(\frac{e^{-t}}{\bar z}\right)f^{\nu_n}_{\bar\omega}\left(\frac{e^{-t}}{\bar z}\right)\frac{e^{-t}}{z}
$$
Applying inequalities (2.4) and estimate (2.1) of lemma 1 to $f^{\nu_n}$, we obtain
$$
  |f^{\nu_n}(e^{-t}z)-f^{\nu_n}(e^{-t}/\bar z)|\le Ce^{-t}(|z|^2 -1)/|z|
$$
for some $C$ independent of $n$ at $|e^t /\bar z|\le 2$. Hence,
$$
  \left|\hat f^{\nu_n}\left(\frac{e^t}{\bar z}\right)-f^{\nu_n}\left(\frac{e^{-t}}{\bar z}\right)\right|\le
$$
$$
 \le\left|\frac{1}{\overline{f^{\nu_n}(e^{-t}z)}}-f^{\nu_n}(e^{-t}z)\right|+|f^{\nu_n}(e^{-t}z)-f^{\nu_n}(e^{-t}/\bar z)|\le
$$
$$
  \le C[(1-|f^{\nu_n}(e^{-t}z)|^2 )+e^{-t}(1-|z|^2 )/|z|]\le C'(1-|e^{-t}/\bar z|^2 )  \eqno (5.45)
$$
for some $C'$ independent of $n$. The last inequality holds since $1-|e^{-t}z|^2 \le 1-|e^{-t}/\bar z|^2$ at $|z|\ge 1$,
$|e^{-t}z|\le 1$ and $e^{-t}(|z|^2 -1)/|z|\le 1 -e^{-t}/|z|\le 1-|e^{-t}/\bar z|^2$.  

Also,
$$
  \hat f^{\nu_n}_z (z)=\frac{1}{\overline{f^{\nu_n}(1/\bar z)}^2}\overline{f^{\nu_n}_z (1/\bar z)}\frac{1}{z^2} ,  \eqno (5.46)
$$
$$
  \hat f^{\nu_n}_{\bar z}(z)=\frac{1}{\overline{f^{\nu_n}(1/\bar z)}^2}\overline{f^{\nu_n}_{\bar z}(1/\bar z)}\frac{1}{\bar z^2} .  \eqno (5.47)
$$
Hence,
$$
  \hat f^{\nu_n}_{\bar\omega}\left(\frac{e^t}{\bar z}\right)\frac{e^t}{z}=\overline{\nu_n}(e^{-t}z)\frac{z}{\bar z}\hat f^{\nu_n}_{\omega}\left(\frac{e^t}{\bar z}\right)\frac{e^t}{\bar z} .
$$
Also,
$$
  \left|\frac{f^{\nu_n}_{\omega}\left(\frac{e^{-t}}{\bar z}\right)\frac{e^{-t}}{\bar z}}{\hat f^{\nu_n}_{\omega}\left(\frac{e^t}{\bar z}\right)\frac{e^t}{\bar z}}\right|=
  \left|\bar z^{-2}(f^{\nu_n}(e^{-t}z))^2\frac{f^{\nu_n}_{\omega}\left(\frac{e^{-t}}{\bar z}\right)}{\overline{f^{\nu_n}_{\omega}(e^{-t}z)}}\right|\le M   
$$
at $|z|\ge 1 $ for some $M$ independent of $n$. Thus dividing $D_{\pm} \hat h_n (z,t)$ by 
$\hat f^{\nu_n}_{\omega}\left(\frac{e^t}{\bar z}\right)\frac{e^t}{\bar z}h_{n,w}\circ f^{\nu_n}\left(\frac{e^{-t}}{\bar z}\right)$
and applying estimates (5.18) and (5.21) we obtain
$$
  \left|\frac{D_- \hat h_n (z,t)}{D_+ \hat h_n (z,t)}\right|\le\frac{b+CM\varepsilon }{1-CMb\varepsilon}\le b+B\varepsilon 
$$
for some uniform $B$. We see that $\hat h_n$ is
a Loewner chain at small enough $\varepsilon$ if $|e^t /\bar z|\le 2$.

Now we define
$$
  \hat h_n (z)= h_n \circ f^{\nu_n}(1/\bar z )+[\hat f^{\nu_n}(z)-f^{\nu_n}(1/\bar z )] h'_n \circ f^{\nu_n}(1/\bar z )  .  \eqno (5.48)
$$ 
at $|z|\ge 1$. 

We prove that $\hat h_n$ is a quasyconformal homeomorphism of $\bbC$ extending $g_n$. 
At first we note that $\hat h_n (z)=\hat h_n (z',t)$ for $z$ represented in the
form $z=e^t z'$, $t\ge 1, |z'|=1$. Thus in some neighborhood
of the unite circle the function $\hat h_n$ maps the point $e^t z'$ into a point on the trajectory of the vector field 
$P(e^t z')=\hat h_{n,t}(z',t)$ starting at $z'$. As it follows from (5.18), (5.21), the map $g_n$ extends to a $C^1$-diffeomorphism
of the unit circle onto the boundary $\pd\Omega_n$. Our vector field is transversal to this boundary and we obtain
homeomorphism of some neighborhood of $\Omega_n$ extending $g$.

Now we prove that  $\hat h_n$ is a local homeomorphism at any point $|z|>1$ with the complex dilatation bounded by some constant
less than 1 depending only on $\nu$. We denote $\omega=1/\bar z$. We have
$$
  \hat h_{n,z}(z)=-h_{n,w}\circ f^{\nu_n}\left(\frac{1}{\bar z}\right)f^{\nu_n}_{\bar\omega}\left(\frac{1}{\bar z}\right)\frac{1}{z^2}+
$$
$$
  +\left[\hat f^{\nu_n}_z (z)+f^{\nu_n}_{\bar\omega}\left(\frac{1}{\bar z}\right)\frac{1}{z^2}\right]h'_n \circ f^{\nu_n}\left(\frac{1}{\bar z}\right)-
$$ 
$$
  -\left[\hat f^{\nu_n}(z)-f^{\nu_n}\left(\frac{1}{\bar z}\right)\right]h''_n \circ f^{\nu_n}\left(\frac{1}{\bar z}\right)f^{\nu_n}_{\bar\omega}\left(\frac{1}{\bar z}\right)\frac{1}{z^2}=
$$
$$
  =\hat f^{\nu_n}_z (z) h'_n \circ f^{\nu_n}\left(\frac{1}{\bar z}\right)-  
\left[\hat f^{\nu_n}(z)-f^{\nu_n}\left(\frac{1}{\bar z}\right)\right]h''_n \circ f^{\nu_n}\left(\frac{1}{\bar z}\right)f^{\nu_n}_{\bar\omega}\left(\frac{1}{\bar z}\right)\frac{1}{z^2} ,
\eqno (5.49)
$$
$$
  \hat h_{n,\bar z}(z)=\hat f^{\nu_n}_{\bar z}(z) h'_n \circ f^{\nu_n}\left(\frac{1}{\bar z}\right)-
$$
$$  
 -\left[\hat f^{\nu_n}(z)-f^{\nu_n}\left(\frac{1}{\bar z}\right)\right]h''_n \circ f^{\nu_n}\left(\frac{1}{\bar z}\right)f^{\nu_n}_{\omega}\left(\frac{1}{\bar z}\right)\frac{1}{\bar z^2} .
\eqno (5.50)
$$
We have analogously to (5.45)
$$
  \left|\hat f^{\nu_n}(z)-f^{\nu_n}\left(\frac{1}{\bar z}\right)\right|=
$$
$$
  =\left|\frac{1}{\overline{f^{\nu_n}(1/\bar z)}}-f^{\nu_n}(1/\bar z)\right|\le C[1-|f^{\nu_n}(1/\bar z)|^2 ] .
\eqno (5.51)
$$
Thus,
$$
 \left|\hat h_{n,z}(z)-\hat f^{\nu_n}_z (z) h'_n \circ f^{\nu_n}\left(\frac{1}{\bar z}\right)\right|\le C\varepsilon,
\eqno (5.52)
$$
$$
  \left|\hat h_{n,\bar z}(z)-\hat f^{\nu_n}_{\bar z}(z) h'_n \circ f^{\nu_n}\left(\frac{1}{\bar z}\right)\right|\le C\varepsilon,
\eqno (5.53)
$$
for some $C$ depending only on $\nu$. Taking into consideration (5.46), (5.47), we see that $|\hat h_{n,z}|$ is bounded away
from zero and $|\hat h_{n,\bar z}/\hat h_{n,z}|$ is bounded by some uniform constant $b' <1$.

We don't consider behavior of our mapping at infinity because on the next step we modify it
outside some $D_r ,r>1$.

2) {\it Modification of $\hat h_n$}.

At $|z|\ge 1/r_n$ we define
$$
  \hat f_n (z)
  = f_n \circ f^{\nu_n}(1/\bar z )+[\hat f^{\nu_n}(z)-f^{\nu_n}(1/\bar z )] f'_n \circ f^{\nu_n}(1/\bar z )  .  \eqno (5.54)
$$
Let $\rho (s), s\ge 0$ be some monotonic smooth function, $\rho (s)=1, s\le 1$, $\rho (s)=0, s\ge 2$, 
$|\rho '(s)|\le 2$. Set some $d_n \ge 1-r_n$ and define at $|z|\ge 1$
$$
  \tilde g_n (z)=\hat h_n (z)\rho ((d_n )^{-1}(|z|-1)) +\hat f_n(z)[1-\rho ((d_n )^{-1}(|z|-1))]  .  \eqno (5.55)
$$

{\begin{proposition} There exists some uniform $M$ such that at 
$$
   d_n \ge M(1-r_n )^{\beta}  ,     \eqno (5.56)
$$   
where $\beta$ is the exponent from Proposition 26, the map $\tilde G_n$ defined on $D$ as $g_n$ and on $\bbC\setminus D$ as $\tilde g_n$ 
is a quasiconformal homeomorphism of the plane with uniformly bounded dilatation and with 
derivative $(\tilde g_n )_z$ uniformly bounded from below and from above.
\end{proposition}

{\bf Proof}. At first we prove that we can find $M$ such that on the domain 
$B_M =\{|z|\ge 1+M(1-r_n )^{\beta}\}$ the map $\hat f_n$ is a local homeomorphism.

Analogously to (5.49), (5.50) and applying (5.51), we obtain 
$$
  \hat f_{n,z}=\hat f^{\nu_n}_z (z) f'_n \circ f^{\nu_n}\left(\frac{1}{\bar z}\right)- 
$$
$$  
  -\left[\hat f^{\nu_n}(z)-f^{\nu_n}\left(\frac{1}{\bar z}\right)\right]f''_n \circ f^{\nu_n}\left(\frac{1}{\bar z}\right)f^{\nu_n}_{\bar\omega}\left(\frac{1}{\bar z}\right)\frac{1}{z^2} ,
 \eqno (5.57) 
$$
$$
  \hat f_{n,\bar z}(z)=\hat f^{\nu_n}_{\bar z}(z) f'_n \circ f^{\nu_n}\left(\frac{1}{\bar z}\right)-
$$
$$  
 -\left[\hat f^{\nu_n}(z)-f^{\nu_n}\left(\frac{1}{\bar z}\right)\right]f''_n \circ f^{\nu_n}\left(\frac{1}{\bar z}\right)f^{\nu_n}_{\omega}\left(\frac{1}{\bar z}\right)\frac{1}{\bar z^2} .   
 \eqno (5.58) 
$$   
Applying estimates (5.30), (5.31) and (5.51), we see that
$$
  \left|\hat f_{n,z}(z)-\hat f^{\nu_n}_z (z) h'_n \circ f^{\nu_n}\left(\frac{1}{\bar z}\right)\right|\le
  c\frac{(1-r_n )^{\beta}}{1-|f^{\nu_n}(1/\bar z)|^2} ,  \eqno (5.59) 
$$  
$$
  \left|\hat f_{n,\bar z}(z)-\hat f^{\nu_n}_{\bar z}(z) h'_n \circ f^{\nu_n}\left(\frac{1}{\bar z}\right)\right|\le
  c\frac{(1-r_n )^{\beta}}{1-|f^{\nu_n}(1/\bar z)|^2}   \eqno (5.60)
$$ 
with some uniform $c$. Taking into consideration inequalities (2.4), we see that at appropriate uniform $M$ the map 
$\hat f_n$ is a local homeomorphism on the domain $B_M$. 

To consider behavior of our mapping at infinity we put $\zeta =1/z$ and consider the function
$\hat f^1_n (\zeta )=1/\hat f_n (1/\zeta )$. It is easy to see that $\hat f^1_n (\zeta )\to 0$
as $\zeta\to 0$ and $\hat f^1_{n,\zeta}\to \overline{f^{\nu_n}_z (0)}f_{n,z}(0)$,  
$\hat f^1_{n,\bar\zeta}\to \overline{f^{\nu_n}_{\bar z}(0)}f_{n,z}(0)$. It implies that 
$\hat f_n$ extends as a local homeomorphism on infinity. 

It remains to prove that $\tilde g_n$ is a local homeomorphism on the domain $\{d_n \le |z|-1 \le 2d_n\}$. 
Indeed, then the map $\tilde G_n$ extends to a a local homeomorphism of the sphere and, hence, is a homeomorphism 
by the monodromy theorem.  

On the domain $\{d_n \le |z|-1 \le 2d_n\}$
$$
  \tilde g_{n,z} =\hat h_{n,z}\rho_n +\hat f_{n,z}(1-\rho)+(\rho_n )_z (\hat h_n -\hat f_n )  ,   
$$ 
where $\rho_n (z)=\rho (d_n^{-1}(|z|-1)$. We have
$$
  |\hat h_{n,z}\rho_n +\hat f_{n,z}(1-\rho_n ) -\hat h_{n,z}|=(1-\rho_n )|\hat f_{n,z}-\hat h_{n,z}| .
$$
By (5.52) and (5.59),
$$
  |\hat h_{n,z}\rho_n +\hat f_{n,z}(1-\rho_n ) -\hat f^{\nu_n}_z (z) h'_n \circ f^{\nu_n}(1/\bar z )|\le c\frac{(1-r_n )^{\beta}}{1-|f^{\nu_n}(1/\bar z)|^2}  
$$
Also, by (5.53) and (5.60), 
$$
  |\hat h_{n, \bar z}\rho_n +\hat f_{n,\bar z}(1-\rho_n ) -\hat f^{\nu_n}_{\bar z}(z) h'_n \circ f^{\nu_n}(1/\bar z )|\le c\frac{(1-r_n )^{\beta}}{1-|f^{\nu_n}(1/\bar z)|^2}
$$
Thus in the domain $\{d_n \le |z|-1 \le 2d_n\}$
$$
 |\hat h_{n,z}\rho_n +\hat f_{n,z}(1-\rho_n ) -\hat f^{\nu_n}_z (z) h'_n \circ f^{\nu_n}(1/\bar z )|\le c/M .
 \eqno (5.61)
$$
$$
  |\hat h_{n, \bar z}\rho_n +\hat f_{n,\bar z}(1-\rho_n ) -\hat f^{\nu_n}_{\bar z}(z) h'_n \circ f^{\nu_n}(1/\bar z )|\le c/M
 \eqno (5.62)
$$
with some uniform $C$, if we set $d_n$ as in (5.56).

Consider now $|\hat h_n -\hat f_n |$. Since $f_n \circ f^{\nu_n}=h\circ f^{\nu}$, we have
$$
  |h_n \circ f^{\nu_n}(z) -f_n \circ f^{\nu_n}(z)|\le |h_n \circ f^{\nu_n}(z)-h\circ f^{\nu_n}(z)|+
  |h\circ f^{\nu_n}(z)-h\circ f^{\nu}(z)| .
$$
But 
$$
  |h_n \circ f^{\nu_n}(z)-h\circ f^{\nu_n}(z)|=|h[r_n f^{\nu_n}(z)]-h[f^{\nu_n}(z)]|\le c(1-r_n )
$$  
since $h'$ is uniformly bounded. Also, taking into consideration (5.34),
$$
  |h\circ f^{\nu_n}(z)-h\circ f^{\nu}(z)|\le \sup_{[f^{\nu_n}(z),f^{\nu}(z)]}|h'||f^{\nu_m}(z)-f^{\nu}(z)|\le
  c(1-r_m )^{\beta}  .
$$
Thus,
$$ 
  |h_n \circ f^{\nu_n}(z) -f_n \circ f^{\nu_n}(z)|\le c(1-r_m )^{\beta}  .
$$
According to (5.48) and (5.54) and applying also estimate (5.30), we obtain
$$
  |\hat h_n (z) -\hat f_n (z)|\le c(1-r_m )^{\beta} .
$$
Since $|\rho'(d_n^{-1}(|z|-1))|\le 2/d_n$, we see that if we set $d_n$ as in (5.56), then
$$
  |(\rho_n )_z (\hat h_n -\hat f_n )|\le c/M
$$
and analogous estimate holds for $|(\rho_n )_{\bar z}(\hat h_n -\hat f_n )|$.

Recalling also (5.61), (5.62) we see that in the domain $\{d_n \le |z|-1 \le 2d_n\}$ we have the estimates
$$
  |\tilde g_{n,z}-\hat f^{\nu_n}_z (z) h'_n \circ f^{\nu_n}(1/\bar z )|\le c/M  ,
$$
$$
  |\tilde g_{n,\bar z}-\hat f^{\nu_n}_{\bar z}(z) h'_n \circ f^{\nu_n}(1/\bar z )|\le c/M  .
$$
We see that at appropriate uniform $M$ $\tilde g_n$ is a local homeomorphism. The assertions
about the dilatation and the derivative $\tilde g_{n,z}$ follow from estimates (5.52), (5.53) for
$\hat h_n$, (5.59), (5.60) and the last estimates.
$\Box$ 

3) {\it The final extension of $g_n$ and the extension of $g$}.

Now we define new extensions $\hat g_n (z)$. Fix  $1/2< r_1 <1$ and put $\hat g_1 =\tilde g_1$.  
As on the previous step, $\rho (s), s\ge 0$ is some monotonic smooth function, $\rho (s)=1, s\le 1$, $\rho (s)=0, s\ge 2$, 
$|\rho '(s)|\le 2$. Let $a_n >0$, $a_n <a_{n-1}$ be some sequence of constants, which we
shall specify later. We define
$$
  \hat g_n (z)=\tilde g_n (z)\rho (a_n^{-1}(|z|-1)) +\tilde g_{n-1}(z)[1-\rho (a_n^{-1}(|z|-1))] 
$$
at $1\le |z| \le 1+a_{n-1}$,
$$
   \hat g_n (z)=\tilde g_{n-1}(z)\rho (a_{n-1}^{-1}(|z|-1)) +\tilde g_{n-2}(z)[1-\rho (a_{n-1}^{-1}(|z|-1))] 
$$ 
at $1+a_{n-1}\le |z| \le 1+a_{n-2}$,

.............

$$
   \hat g_n (z)=\tilde g_1 
$$
at $|z|\ge 1+a_1$. 

Now we specify the sequences $r_n$, $d_n$ and $a_n$. We put
$$
  1-r_n \le\frac{1}{5^{1/\beta}}(1-r_{n-1}),\, d_n =C(1-r_n )^{\beta},\, a_n =2d_{n-1} ,  
\eqno (5.63)
$$  
where $\beta$ and $C$ are the constants from Proposition 27.

\begin{proposition} If we define the sequences $r_n$, $d_n$ and $a_n$ according to (5.63)
with some appropriate uniform $C$, then the map $g_n$ extended on $\bbC\setminus D$
as $\hat g_n$ will be a quasiconformal homeomorphism of the plane, 
and the derivative $(\hat g_n )_z$ will be uniformly bounded from below and from above. 
On each domain $1+a_k\le |z| \le 1+a_{k-1}$, $k\le n$
$$
  \hat g_n (z)=\hat f_k (z)\rho (a_k^{-1}(|z|-1)) +\hat f_{k-1}(z)[1-\rho (a_k^{-1}(|z|-1))] . \eqno (5.64)
$$  
\end{proposition}

{\bf Proof}.  It is easy to see that, if we set $a_n$ according to (5.63), the map $\hat g_n$ will
be described by expression (5.64) on each ring $1+a_k\le |z| \le 1+a_{k-1}$, $k\le n$.  
We must check only that it is a local diffeomorphism on any such domain. By induction,
it is enough to check it for $k=n$. 

We have
$$
  \hat g_{n,z} =\hat f_{n,z}\rho_n +\hat f_{n-1,z}(1-\rho_n )+(\rho_n )_z (\hat f_n -\hat f_{n-1})  , 
$$ 
where $\rho_n (z)=\rho(a_n^{-1}(1-|z|)$. The analogous expression we have for $\hat g_{n,\bar z}$. 
We have
$$
  |\hat f_{n,z}\rho_n +\hat f_{n-1,z}(1-\rho_n ) -\hat f_{n,z}|=(1-\rho_n )|\hat f_{n,z}-\hat f_{n-1,z}| .
$$
By (5.57), (5.59), (5.31), and taking into consideration (2.4),
$$
  |\hat f_{n,z}\rho_n +\hat f_{n-1,z}(1-\rho_n ) -\hat f^{\nu_n}_z (z) h'_n \circ f^{\nu_n}(1/\bar z )|\le c\frac{(1-r_{n-1})^{\beta}}{|z|^2 -1}+
$$  
$$
  +|\hat f^{\nu_n}_z (z) f'_n \circ f^{\nu_n}(1/\bar z )-\hat f^{\nu_{n-1}}_z (z) f'_{n-1}\circ f^{\nu_{n-1}}(1/\bar z )|+
$$  
$$
  +\left|\left[\hat f^{\nu_n}(z)-f^{\nu_n}\left(\frac{1}{\bar z}\right)\right]f''_n \circ f^{\nu_n}\left(\frac{1}{\bar z}\right)f^{\nu_n}_{\bar\omega}\left(\frac{1}{\bar z}\right)\frac{1}{z^2}-\right.
$$
$$
  \left.-\left[\hat f^{\nu_{n-1}}(z)-f^{\nu_{n-1}}\left(\frac{1}{\bar z}\right)\right]f''_{n-1}\circ f^{\nu_{n-1}}\left(\frac{1}{\bar z}\right)f^{\nu_{n-1}}_{\bar\omega}\left(\frac{1}{\bar z}\right)\frac{1}{z^2}\right|  
\eqno (5.65)
$$  
Now, applying (5.28) and (5.32), we have 
$$
  |\hat f^{\nu_n}_z (z) f'_n \circ f^{\nu_n}(1/\bar z )-\hat f^{\nu_{n-1}}_z (z) f'_{n-1}\circ f^{\nu_{n-1}}(1/\bar z )|\le |\hat f^{\nu_n}_z (z)-\hat f^{\nu_{n-1}}_z (z)||f'_n \circ f^{\nu_n}(1/\bar z )|+
$$
$$
   +|\hat f^{\nu_{n-1}}_z (z)|(|f'_n \circ f^{\nu_n}(1/\bar z )- f'_{n-1}\circ f^{\nu_{n-1}}(1/\bar z )| \le c\frac{(1-r_{n-1})^{\beta}}{|z|-1} .
$$
with some uniform $c$. 

Now remind that $|f'_n (z)|$ is uniformly bounded from below and from above since the first derivatives of $f^{\nu_n}$ and 
$g|_{D_{r_n}}=f_n \circ f^{\nu_n}$ are uniformly bounded from below and from above, and $|f''_n (z)|\le c(1-|z|)^{-1}$ since
$|f''_n (z)/f'(z)|\le c(1-|z|)^{-1}$ by properties of univalent functions (see [Pom]). 

We see that the second difference in the right side of (5.65) is no greater
then the sum of the terms: first,
$$
 \left|\hat f^{\nu_n}(z)-f^{\nu_n}\left(\frac{1}{\bar z}\right)-\hat f^{\nu_{n-1}}(z)+f^{\nu_{n-1}}\left(\frac{1}{\bar z}\right)\right|\left|f''_{n-1}\circ f^{\nu_{n-1}}\left(\frac{1}{\bar z}\right)\right|O(1)\le
$$
$$
  \le c\frac{(1-r_{n-1})^{\beta}}{|z|-1}  
$$  
since $f^{\nu_n}(z)-f^{\nu_{n-1}}(z)|\le c(1-r_{n-1})^{\beta}$ (see (3.1) and (5.34)); second,
$$
  \left|\hat f^{\nu_{n-1}}(z)-f^{\nu_{n-1}}\left(\frac{1}{\bar z}\right)\right|\left|f''_{n-1}\circ f^{\nu_{n-1}}\left(\frac{1}{\bar z}\right)\right|\left|f^{\nu_{n-1}}_{\bar\omega}\left(\frac{1}{\bar z}\right)-f^{\nu_n}_{\bar\omega}\left(\frac{1}{\bar z}\right)\right|\le
$$
$$
  \le c\frac{(1-r_{n-1})^{\beta}}{|z|-1} 
$$
by (5.51) and (5.28); and third,
$$
  \left|\hat f^{\nu_{n-1}}(z)-f^{\nu_{n-1}}\left(\frac{1}{\bar z}\right)\right|\left|f''_{n-1}\circ f^{\nu_{n-1}}\left(\frac{1}{\bar z}\right)-f''_n \circ f^{\nu_n}\left(\frac{1}{\bar z}\right)\right|O(1)\le
$$
$$
  \le c\frac{(1-r_{n-1})^{\beta}}{|z|-1} 
$$
by (5.51) and (5.32).

Thus,
$$
  |\hat f_{n,z}\rho_n +\hat f_{n-1,z}(1-\rho_n ) -\hat f^{\nu_n}_z (z) h'_n \circ f^{\nu_n}(1/\bar z )|\le c\frac{(1-r_{n-1})^{\beta}}{|z|-1} . \eqno (5.66)
$$
Analogously we obtain
$$
|\hat f_{n,\bar z}\rho_n +\hat f_{n-1,\bar z}(1-\rho_n ) -\hat f^{\nu_n}_{\bar z} (z) h'_n \circ f^{\nu_n}(1/\bar z )|\le c\frac{(1-r_{n-1})^{\beta}}{|z|-1}  \eqno (5.67) 
$$
with some uniform $c$.
 
Now we shall estimate $|\hat f_n -\hat f_{n-1}|$. Since $f_n \circ f^{\nu_n}=g|_{D_{r_n}}$, we see that
$f_n \circ f^{\nu_n}(1/\bar z)-f_{n-1}\circ f^{\nu_{n-1}}(1/\bar z)=0$ if $|1/z|\le r_{n-1}$.  
We must only estimate the difference of the terms containing $f'_n$ and $f'_{n-1}$. We have
$$
  |[\hat f^{\nu_n}(z)-f^{\nu_n}(1/\bar z )-
\hat f^{\nu_{n-1}}(z)+f^{\nu_{n-1}}(1/\bar z )]f'_n \circ f^{\nu_n}(1/\bar z)|\le c(1-r_{n-1})^{\beta}
$$
by (3.1) and (5.34), and
$$  
  |[\hat f^{\nu_{n-1}}(z)-f^{\nu_{n-1}}(1/\bar z )][f'_n \circ f^{\nu_n}(1/\bar z )-
f'_{n-1}\circ f^{\nu_{n-1}}(1/\bar z )]|\le c(1-r_{n-1})^{\beta}.
$$
by (5.51) and (5.31).

Since $|\rho'_n |\le 2/a_n$, we see that
$$
  |(\rho_n )_z (\hat f_n -\hat f_{n-1})|\le c/M
$$
in the ring $1+a_n\le |z| \le 1+a_{n-1}$. For $|(\rho_n )_{\bar z}(\hat f_n -\hat f_{n-1})|$ we have
an analogous estimate. 

As a result, collecting (5.66), (5.67), and the last estimates, we obtain
$$
   |(\hat g_n )_z -\hat f^{\nu_n}_z (z) h'_n \circ f^{\nu_n}(1/\bar z )|\le c/M  ,
$$
$$   
   |(\hat g_n )_{\bar z}-\hat f^{\nu_n}_{\bar z}(z) h'_n \circ f^{\nu_n}(1/\bar z )|\le c/M , 
$$
where $M$ is the constant from inequality (5.56) and $c$ is some uniform constant independent of $M$.
At appropriate $C$ we obtain that $\hat g_n$ is a quasyconformal map with uniformly bounded dilatation
and the derivative $(\hat g_n )_z$ uniformly bounded from below and from above.
$\Box$

Now we define the extension $\hat g$ as the limit of $\hat g_n$. By definition, $\hat g(z)=\hat g_n (z)$
if $|z|\ge 1+a_n$. Thus $\hat g_n$ converge to a quasiconformal map on $\hat\bbC\setminus D$ with 
derivative uniformly bounded from below and from above. Also, the map defined as $g$ on $D$ and as
$\hat g$ on $\bbC\setminus D$ is a one-to-one mapping. Indeed, if $g(z_1 )=\hat g (z_2 )$,
then the domains $g_n (D)$ and $\hat g_n (\bbC\setminus D)$ must intersect at great enough $n$.
Thus we obtained a homeomorphism of the plane and proved estimate (5.13).

4) {\bf Proof of estimates (5.14) and (5.15)}. 

It is enough to prove that these estimates hold for $\hat g_n$ at $1+a_n \le |z|\le 1+a_{n-1}$.

The maps $g$ and $f^{\nu_n}$ satisfy the estimates $|g_{(k)}(z)|\le C(1-|z|)^{1-|(k)|}$,
$|f^{\nu_n}_{(k)}(z)|\le C(1-|z|)^{1-|(k)|}$. From the equation $g|_{D_{r_n}}=f_n \circ f^{\nu_n}$
by successive differentiation we obtain the estimate 
$$
  |(f_n )_{z^k}\circ f^{\nu_n}(z)|\le c(1-|f^{\nu_n}(z)|)^{1-k}\le C (1-|z|)^{1-k} .
$$  
Also, we have estimates (5.26), (5.27) for $|f^{\nu_n}_{(k,l)}|$ and $|(f_n )_{(k,l)}|$.

According to representations (5.64), (5.54) the derivative $(\hat g_n )_{(k)}$ is a sum of items
containing the multiples $z^{-j}$, $\bar z^{-j}$, which don't influence an order in $|z|-1$
and multiples of the types
$$ 
  (f_n ){z^j}\circ f^{\nu_n}(1/\bar z ),\,f^{\nu_n}_{(l)}(1/\bar z ),\,\hat f^{\nu_n}_{(m)}(z)  ,
\eqno (5.68)
$$
and analogous terms containing $f_{n-1}$, $f^{\nu_{n-1}}$. Also, there can be the multiples 
$$
  (\rho_n )_{(s)}(z)  .  \eqno (5.69) 
$$
At differentiation of each multiple of type (5.68) we increase in order in $(|z|-1)^{-1}$ by one.  
From the other hand, differentiation of multiple (5.69) results in the additional multiple
$a_n^{-1}$. On our ring it is a value of order $(|z|-1)^{-1}$. By induction, we obtain 
estimate (5.14).

Now derivatives with respect to the parameter of terms of types (5.67) all have estimates $C(1-r_n )^{-M}$ 
for some $M$ by (5.26) and (5.27). But, if we set $r_n$, $d_n$ and $a_n$ according to (5.63), we have
$$
  1-r_n =5^{1/\beta}(1-r_{n-1})=5^{1/\beta}\left(\frac{d_{n-1}}{C}\right)^{1/\beta}= 
 5^{1/\beta}\left(\frac{a_n}{2C}\right)^{1/\beta} . 
$$
Again differentiation of terms with $\rho_n$ results in multiples $a_n^{-j}$ with some $j$.
Thus on the ring $1+a_n \le |z|\le 1+a_{n-1}$ we obtain estimate (5.15) with some uniform $C$ and $M$.
$\Box$

5). {\bf Proof of estimates (5.16), (5.17)}.

We know that $\hat g_n$ is a quasiconformal homeomorphism mapping $\hat\bbC\setminus D$ onto $\hat\bbC\setminus\Omega$
with complex dilatation $\hat\nu$. We can represent $\hat g$ as a composition $\tilde h \circ f^{\hat\nu}$,
where $f^{\hat\nu_n}$ is a normal homeomorphism mapping $\hat\bbC\setminus D$ onto itself and $\tilde h$ is an univalent 
holomorphic function. From estimates ((5.13) and (5.14) at $|(k)|\le 2$ follows that $|\hat\nu_{(k)}(z)|\le c(|z|-1)^{-1}$
at $|(k)|=1$. By Lemma 1, $c_1 \le |f^{\hat\nu}_z| \le c_2$ for some constants $c_1$, $c_2$. 
Applying (5.13), we see that for $\tilde h' (z)$ we have analogous estimates.
But for any univalent function $\tilde h$ we have
$$
  a_1 |\tilde h_z (z)|(|z|^2 -1)\le{\rm dist}(\tilde h(z), \pd\Omega )\le a_2 |\tilde h_z (w)|(|z|^2 -1)
$$
with some uniform $a_1 ,a_2$ (see [Pom]). Since for $f^{\hat\nu}$ we have estimates (2.4), we obtain (5.16). 

Let prove (5.17). The left inequality follows from (5.16) because 
$|g(z)-\hat g(\bar z^{-1})|\ge {\rm dist}(g(z),\pd\Omega )+{\rm dist}(\hat g(\bar z^{-1},\pd\Omega )$.

Suppose now that $1+a_n \le |1/\bar z|\le 1+a_{n-1}$.  
Now, by (5.64), we must estimate $g(z)-\hat f_n (1/\bar z)$ and $\hat f_n (1/\bar z)-\hat f_{n-1}(1/\bar z)$.
Applying (5.54) and (5.51), we obtain
$$
  |\hat f_n (1/\bar z)-f_n \circ f^{\nu_n}(z)|\le c(1-|z|) .
$$
But $f_n \circ f^{\nu_n}(z)=g(z)$ at $|z|\le r_n$ and, hence, $|g(z)-\hat f_n (1/\bar z)|\le c(1-|z|)$.
Analogously, we can see that $g(z)-\hat f_{n-1}(1/\bar z)|\le c(1-|z|)$. Also,
$$
  \hat f_n (1/\bar z)-\hat f_{n-1}(1/\bar z)=[\hat f^{\nu_n}(z)-f^{\nu_n}(1/\bar z )]f'_n \circ f^{\nu_n}(1/\bar z)-
$$
$$  
-[\hat f^{\nu_{n-1}}(z)-f^{\nu_{n-1}}(1/\bar z )]f'_{n-1}\circ f^{\nu_{n-1}}(1/\bar z) .
$$
By (5.51), this difference has the estimate $c(1-|z|)$. As a result, we obtain right
inequality (5.17).
$\Box$   

\section{Integral operators. $L^p$-estimates.}

We adopt the notations of the previous section. For $z\in D$, $w=\hat g(z)\in \Omega$ we  define
$$
  \hat w =\hat g (1/\bar z) .\eqno (6.1)
$$ 

\begin{proposition}
There are the estimates
$$
  c(1-|z|^2 ) \le  |w -\hat w |\le C(1-|z|^2 ) , \eqno (6.2)
$$
$$
  \left|\frac{\omega -\hat\omega}{w-\hat\omega}\right|\le C   ,\eqno (6.3)
$$
$$
  |w -\omega |\ge c|z-\zeta |,    \eqno (6.4)
$$  
$$
   c|z-\bar\zeta^{-1}|\le |w -\hat\omega |\le C|z-\bar\zeta^{-1}|   \eqno (6.5) 
$$
with uniform  $c$, $C$. .
\end{proposition}

{\bf Proof}. The first estimate follows immediately from (5.17). The second one,
also, follows from (5.16) and (5.17), taking into consideration that
$$
  |w-\hat\omega |\ge {\rm dist}(\hat\omega,\pd\Omega ) . 
$$
Prove left estimate (6.4). It is enough to prove that the map $g^{-1}$ has uniformly bounded
derivatives. But $g^{-1}=(f^{\nu})^{-1}\circ h^{-1}$, where $h^{-1}$ has a bounded derivative 
and $(f^{\nu})^{-1}$ is the normal map
with the Beltramy coefficient $\nu_{-1}=-\nu\circ (f^{\nu})^{-1}$, and by Lemma 1 and estimate (2.4)
$|(\nu_{-1})_{(k)}(z)|\le c(1-|z|)^{-1}$ at $|(k)|=1$. Again by Lemma 1 we obtain that 
$(f^{\nu})^{-1}_z ,(f^{\nu})^{-1}_{\bar z}$ are uniformly bounded.
Analogously,
$$
  |\hat g(z_1 )-\hat g(z_2 )|\ge c|z_1 -z_2 |  .   \eqno (6.6)
$$ 

Now let $\omega'$ be a point on $\pd\Omega$ closest to $\hat\omega$ (there can be several such points but it isn't essential). 
Suppose at first that $|w-\omega'|\le |\omega'-\hat\omega |$. Then 
$$
 |w -\hat\omega |\ge |\omega'-\hat\omega |\ge c(|\bar\zeta |^{-1}-1|)\ge c'|z-\bar\zeta^{-1}| 
$$
for some uniform $c, c'$. From the other hand, if $|w-\omega'|\le |\omega'-\hat\omega |$, then,
applying (6.5) and (6.6), we obtain
$$
 |w -\hat\omega |\ge |\omega'-w|\ge c|z-g^{-1}(\omega' )|\ge c'|z-\bar\zeta^{-1}| 
$$
for some uniform $c, c'$. 
 
From the other hand, according to (4.35) and (5.13), the map $g$ extended on $\bbC\setminus D$ as $\hat g$
is a Lipshitz homeomorphism of the plane with an uniformly bounded Lipshitz constant. Thus we obtain
right estimate (6.4).
$\Box$ 

Now we define integral transforms, which allow as to find solutions to the Beltrami equation on $\Omega$
with required estimates on the boundary. We define
$$
  P_m f(w)=\frac{1}{\pi}\int_{\Omega}\frac{f(\omega )}{w-\omega}\left(\frac{\omega -\hat\omega}{w-\hat\omega}\right)^m dS_{\omega}=
$$
$$
  =\frac{1}{\pi}\int_{\Omega_n}f(\omega )\left[\frac{1}{w-\omega}-\frac{1}{w-\hat\omega}-...-\frac{(\omega -\hat\omega )^{m-1}}{(w-\hat\omega )^m}\right]dS_{\omega}  .
\eqno (6.7)
$$  
The last representation follows from the identity
$$
  \frac{1}{w-\omega}\left(\frac{\omega -\hat\omega}{w-\hat\omega}\right)^m=\left(\frac{1}{w-\omega}-\frac{1}{w-\hat\omega}\right)\left(\frac{\omega -\hat\omega}{w-\hat\omega}\right)^{m-1} .
$$
Differentiating $P_m f(w)$ in $w$ we obtain the transform
$$
  T_m f(w)=-\frac{1}{\pi}\int_{\Omega}\frac{f(\omega )}{w-\omega}\left(\frac{\omega -\hat\omega}{w-\hat\omega}\right)^m  .
\left(\frac{1}{w-\omega}+\frac{m}{w-\hat\omega}\right)dS_{\omega}=
$$
$$
  =-\frac{1}{\pi}\int_{\Omega}f(\omega )\left[\frac{1}{(w-\omega )^2}-\frac{1}{(w-\hat\omega )^2}-...-\frac{(m-1)(\omega -\hat\omega )^{m-1}}{(w-\hat\omega )^{m+1}}\right]dS_{\omega}  .
\eqno (6.8)
$$
Here we understand the integral in terms of its principal value. In the chart $z$ we have
$$
  P_m f(z)=\frac{1}{\pi}\int_D \frac{f(g_n (\zeta ))}{g_n (z)-g_n (\zeta )}\left(\frac{g_n (\zeta ) -\hat g_n (\zeta )}{g_n (z)-\hat g_n (\zeta )}\right)^m 
|(g_n )_{\zeta}(\zeta )|^2 (1- |\nu_n (\zeta )|^2 )dS_{\zeta}  ,  \eqno (6.9)
$$
$$
  T_m f(z)=-\frac{1}{\pi}\int_D \frac{f(g_n (\zeta )}{g_n (z)-g_n (\zeta )}\left(\frac{g_n (\zeta) -\hat g_n (\zeta )}{g_n (z)-\hat g_n (\zeta )}\right)^m 
\left(\frac{1}{g_n (z)-g_n (\zeta )}+\right.
$$
$$
   \left. +\frac{1}{g_n (z)-\hat g_n (\zeta )}\right)||(g_n )_{\zeta}(\zeta )|^2 (1- |\nu_n (\zeta )|^2 )dS_{\zeta}  =
$$
$$
=-\frac{1}{\pi}\int_D f(g_n (\zeta ))\left[\frac{1}{(g_n (z)-g_n (\zeta ))^2}-\frac{1}{(g_n (z)-\hat g_n (\zeta ))^2}-...\right.
$$
$$ 
  \left. -\frac{(m-1)(g_n (\zeta )-\hat g_n (\zeta  ))^{m-1}}{(g_n (z)-\hat g_n (\zeta ))^{m+1}}\right]|(g_n )_{\zeta}(\zeta )|^2 (1- |\nu_n (\zeta )|^2 )dS_{\zeta}   .  \eqno (6.10)
$$

\begin{definition} We say that a function $f$ on $\Omega$ belongs to $L^p_s (\Omega)$ if the function $f(w)({\rm dist}(w,\pd\Omega ))^s$
belongs to $L^p$. We denote by $\|f\|_{p,s}$ the $L^p$-norm of the function $f(w)({\rm dist}(w,\pd\Omega ))^s$. A function $f$ belongs to
$C^0_s$ if $f(w)({\rm dist}(w,\pd\Omega ))^s$ is uniformly bounded. We denote by $\|f\|_{0,s}$ the $C^0$-norm of the function
$f(w)({\rm dist}(w,\pd\Omega ))^s$.   
\end{definition}

The equivalent conditions are:  $f(g(z))(1-|z| ^2 )^s$ belongs to  $L^p (D)$ and  $f(g(z))(1-|z| ^2 )^s$ is uniformly bounded on $D$. Corresponding
norms are equivalent to the norms of the Definition. 

We will need in the following estimates:

\begin{proposition} Define the integrals
$$
  J(w)=\frac{1}{\pi}\int_{\Omega}\frac{dS_{\omega}}{|w-\omega | |w-\hat\omega |^k} ,
$$
$$
  \tilde J (\omega)=\frac{1}{\pi}\int_{\Omega}\frac{dS_w}{|w-\omega | |w-\hat\omega |^k} .
$$
Suppose $k\ge 2$. Then
$$   
 J(w)\le C(1-|z|^2 )^{k-1},\,\tilde J (\omega)\le C(1-|\zeta |^2 )^{k-1} ,  (6.11)
$$   
where $w=g(z)$, $\omega =g(\zeta )$.
\end{proposition}

{\bf Proof}. a) Applying estimates (6.4) and (6.5) we can see that we must show that the integral
$$
   I(z)=\int_D \frac{dS_{\zeta}}{|z-\zeta ||z-\bar\zeta^{-1}|^2}=\int_D \frac{|\bar\zeta |dS_{\zeta}}{|z-\zeta | |1-z\bar\zeta |^k} .
$$
has the estimate
$$
   I(z)\le C (1-|z|^2 )^{k-1}   б  \eqno (6.12)
$$
and that the analogous integral
$$
  \tilde I(\zeta)=\int_D \frac{dS_z}{|z-\zeta ||z-\bar\zeta^{-1}|^k}=\int_D \frac{|\bar\zeta |dS_z}{|z-\zeta | |1-z\bar\zeta |^2} 
$$
has the estimate
$$
   \tilde I(\zeta )\le C (1-|\zeta |^2 )^{k-1}   .  \eqno (6.13)
$$
  
Let $\vpf_z$ be the map
$$
  \vpf_z (\tau )=\frac{z-\tau}{1-\bar z\tau} .
$$
We have
$$ 
  |z-\zeta|=\frac{|\tau |(1-|z|^2 )}{|1-\bar z \tau |} ,
$$
$$
  |z-\bar\zeta^{-1}|=\frac{1-|z|^2}{|\bar z-\bar\tau |}  .
$$
The Jacobian of the change of the variable $\zeta =\vpf_z (\tau )$ is
$$
  \frac{(1-|z|^2 )^2}{|1-\bar z \tau |^4}  .
$$
After this change of the variable we can write integral (6.12) as
$$ 
   (1-|z|^2 )^{1-k}\int_D \frac{|z-\tau |^k}{|\tau ||1-\bar z \tau |^3}dS_{\tau}  . \eqno (6.14)
$$   
 
But $|z-\tau |\le c|1-\bar z \tau |$ for some uniform $c$, and we see that it is enough to show
that the integral
$$ 
  \int_D \frac{dS_{\tau}}{|\tau ||1-\bar z \tau |}   
$$
is uniformly bounded. But we can write this integral as a sum of the integral over
the disk $D_{1/2}$ and of the integral over the ring $D\setminus D_{1/2}$. The
first integral is no greater, than
$$
  C\int_D \frac{dS_{\tau}}{|\tau |}  
$$
And the second one is no greater, than
$$
  C\int_D \frac{dS_{\tau}}{|1-\bar z \tau |}   .
$$  
Both these integral are uniformly bounded and, hence, integral (6.14) has the estimate $C(1-|z|^2 )^{-1}$.
   
Now notice that $\tilde I$ also reduces to integral (6.14), and we obtain estimate (6.13) exactly
as (6.12) . 
$\Box$

The following estimate is a corollary of this proposition.

\begin{proposition} Let $f$ belongs to $C^0_2$. Then 
$$
 \|P_m f \|_{0,1}\le C\|f\|_{0,2} 
$$
with some uniform $C$.
\end{proposition}    

{\bf Proof}. By (6.3), we have
$$
   |P_m f (w)|\le \frac{1}{\pi}\int_{\Omega}\left|\frac{f(\omega )}{w-\omega}\left(\frac{\omega -\hat\omega}{w-\hat\omega}\right)^m \right|dS_{\omega}\le
  \frac{1}{\pi}\sup_{\omega\in\Omega}|f(\omega)(\omega -\hat\omega )^2 | \int_{\Omega}\frac{dS_{\omega}}{|w-\omega | |w-\hat\omega |^2}  . 
$$
We apply first estimate (6.11).
$\Box$

\begin{proposition} 
Define the operators:
$$
  L_{\beta}f(\omega )=\frac{1}{\pi}\int_{\Omega}f(\omega ) \frac{\beta (w,\omega )}{(w-\hat\omega )^2}dS_{\omega} ,  
$$
and
$$
  \tilde L_{\beta}f(\omega )=\frac{1}{\pi}\int_{\Omega}f(\omega )\frac{\beta (w,\omega )}{(\hat w-\omega )^2}dS_{\omega}  ,
$$
where $\beta$ is uniformly bounded function. 
These operators are bounded in $L^p (\Omega )$, $1< p<\infty$. As a consequence, 
the operators $T_m$ and $\tilde T_m$ are bounded in $L^p (\Omega )$, $1< p<\infty$.
\end{proposition}

{\bf Proof}. Consider the operator $L_{\beta}$.
From boundedness of the derivatives of $g$ and estimates (6.4), (6.5) follows that we can write
$$
  w-\omega =g(z)-g(\tau )=h_1 (z,\tau )(z-\tau ),\,w-\hat\omega =h_2 (z,\tau )(z-\hat\tau^{-1}) ,
$$
where $h_1$ and $h_2$ are uniformly bounded from below and from above.  

The operator $L_{\beta}$ obviously is bounded in any $L^p$ for functions supported in any domain  $|g^{-1}(\omega )||\le c<1$. 
Suppose that the support of $f$ is contained in the domain $1/2 \le |g^{-1}(\omega )|\le 1$.
Go to the chart $t$ on the disk $D$. Recalling that the Jacobian of the transformation  $\tau\mapsto\omega$ 
is uniformly bounded, we see that we must estimate the norm of the operator
$$
  L'_{\beta}f(z) =\frac{1}{\pi}\int_D f(\tau )\frac{\beta (z,\tau )}{(z-\bar\tau^{-1})^2}dS_{\tau}  .  \eqno (6.15)
$$
Introduce a new variable $t=\bar\tau^{-1}$. Then $\tau=\overline{t^{-1}}$. Integral (6.15) transforms to
$$
   \frac{1}{\pi}\int\tilde f(t)\frac{\tilde\beta (z,t)}{(z-t)^2}dS_t ,
$$
where $\tilde f(t)=f(\overline{t^{-1}})$, $\tilde\beta=\beta (z,\overline{t^{-1}})J(t)$ and $J$ is the 
Jacobian of the transformation  $t\mapsto\tau$. Here we take the integral over $\bbC$. Remind that
support $\tilde f$ is contained in the domain ${\bbC}\setminus D$, the Jacobian $J$ is uniformly bounded
on this domain, and $z\in D\setminus\pd D$. 
The below considerations follow the method described in [Ah, Ch. 5, D].

Let be $t-z =\rho e^{i\theta}$. 
Our integral is of the type
$$
  \int_0^{2\pi}\int_0^{\infty}\tilde f(z+\rho e^{i\theta})\frac{h(z,z+\rho e^{i\theta})}{\rho}e^{-2i\theta}d\rho d\theta ,
$$
where $|h(z,z+\rho e^{i\theta})|/\rho \le c$
We have
$$
  \|L'_{\beta}f(z)\|_p \le 2\pi\max_{\theta}\left\|\int_0^{\infty}\tilde f(z+\rho e^{i\theta})\frac{(h(z,z+\rho e^{i\theta})}{\rho}d\rho\right\|_p  
$$
 We can suppose that the functions $f$ and $h$ are real. 
We have
$$
  \left|\int_0^{\infty}\tilde f(z+\rho e^{i\theta})\frac{(h(z,z+\rho e^{i\theta})}{\rho}d\rho\right|\le
   c\int_0^{\infty}\frac{|\tilde f(z+\rho e^{i\theta})|}{\rho}d\rho 
$$
for some $c$ and, hence,
$$
  \|L'_{\beta}f(z)\|_p\le 2\pi c \max_{\theta}\left\|\int_0^{\infty}\frac{|\tilde f(z+\rho e^{i\theta})|}{\rho}d\rho \right\|_p
$$
Let this maximum corresponds to $\theta =\Theta$.
 The norm in the right side doesn't change if we replace $z$ by $ze^{i\Theta}$. If we denote $\tilde f_{\Theta}(t)=\tilde f(te^{i\Theta})$, 
the integral in the right side becomes
$$
  H\tilde f_{\Theta}(w)=\int_0^{\infty}\frac{|\tilde f_{\Theta}(z+\rho )|}{\rho}d\rho  .
$$
The point $z=x+iy$ doesn't belong to the support of $\tilde f_{\Theta}$ and the value of the integral doesn't change if we extend
the function of $\rho$: $\tilde f_{\Theta}(z+.)$ as zero on the domain $\rho <0$, extend the integral over this domain,
and take the principal value. Hence, 
$$
   \int |H\tilde f_{\Theta}(x+iy)|^p dx \le A_p \int |\tilde f_{\Theta}(x+iy)|^p dx
$$ 
for some $A_p$ by one-dimensional Calderon-Zigmund inequality. 
Now we get for the two-dimensional norm
$$
  \|H\tilde f_{\Theta}\|_p^p =\int\int|H\tilde f_{\Theta}(x+iy)|^p dxdy\le 
$$
$$
  \le A_p^p \int\int|\tilde f_{\Theta}(x+iy)|^p dxdy=A_p^p \|\tilde f_{\Theta}\|_p^p .
$$
Since the functions $\tilde f$ and $\tilde f_{\Theta}$ have equal $L^p$-norms and the $L^p$-norm of $\tilde f$ can be estimated through 
$L^p$-norm of $f$, we obtain the required estimate.

Consider now the operator $T_m$.
The first item under the integral in the right side of (6.8) is the Beurling transform bounded in $L^p$. 
The other items are the integrals
$$
  \frac{1}{\pi}\int_{\Omega}f(\omega )\frac{(\omega-\hat\omega )^k}{(w-\hat\omega )^{k+2}}  .
$$  
These integrals are of the type considered above for $\beta (w,\omega )=(\omega-\hat\omega )^k /(w-\hat\omega )^k$. 
 
The case of the operators $\tilde L_{\beta}$ and $\tilde T_m$ is analogous. 
It is even simpler because we don't need now in the change of the variables $\tau\mapsto t$. 
The proof with this exception repeats the proof for $L_{\beta}$.
$\Box$

\begin{proposition}
Suppose an operator $L$ satisfies the estimate
$$
  Lf (w)\le C|w-\hat w|\int_{\Omega} |f(\omega )|\frac{|\omega -\hat\omega |}{|w-\omega ||w-\hat\omega |^3}dS_{\omega} . 
$$
Then this operator is bounded in $L^p (\Omega )$ for $1\le p\le\infty$. 
\end{proposition}

{\bf Proof}
We shall prove the boundedness in $L^1$ and $L^{\infty}$. Then the general case will follow from the Riesz-Thorin
interpolation theorem (see, for example, [RS]). We have
$$
  \|Lf\|_1 \le C\int_{\Omega}|w-\hat w |\int_{\Omega}|f(\omega )|\frac{|\omega -\hat\omega |}{|w-\omega ||w-\hat\omega |^3}dS_{\omega}dS_w \le
$$
$$
  \le C\int_{\Omega}|f(\omega )||\omega -\hat\omega |\int_{\Omega}\frac{dS_w}{|w-\omega ||w-\hat\omega |^2}dS_{\omega}\le C\|f\|_1
$$
by the second estimate of Proposition 30.  

From the other hand,
$$
  |Lf (w)|\le C|w-\hat w |\|f\|_{\infty}\int_{\Omega}\frac{|\omega -\hat\omega |}{|w-\omega ||w-\hat\omega |^3}dS_{\omega}\le
$$
$$
  \le C\|f\|_{\infty}|w-\hat w |\int_{\Omega}\frac{dS_{\omega}}{|w-\omega ||w-\hat\omega |^2} \le C\|f\|_{\infty}
$$
by the first estimate of Proposition 30.  
$\Box$

\begin{proposition}
Let $I_{m,s}$, $\tilde I_{m,s}$, $J_{m,s}$ be the operators 
$$
  I_{m,s}f(w)=(w-\hat w )^s \frac{1}{\pi}\int_{\Omega}\frac{f(\omega )}{w-\omega}\frac{(\omega -\hat\omega )^{m-s}}{(w-\hat\omega )^{m+1} }dS_{\omega} ,
$$ 
$$
  \tilde I_{m,s}f(w)=(w-\hat w )^s \frac{1}{\pi}\int_{\Omega}f(\omega )\frac{\hat w -\hat\omega}{(w-\omega )^2}\frac{(\omega -\hat\omega )^{m-s}}{(\hat w-\omega )^{m+1} }dS_{\omega} ,
$$
$$
  J_{m,s}f(w)=(w-\hat w )^s \frac{1}{\pi}\int_{\Omega}\frac{f(\omega )}{(w-\omega )^2}\frac{(\omega -\hat\omega )^{m-s}}{(w-\hat\omega )^m}dS_{\omega} ,
$$  .
These operators are bounded in $L^p (\Omega )$ for $1<p<\infty$ at $m\ge 0,s\le m$.
\end{proposition}

{\bf Proof}. Suppose at first that $m\ge 2$, $1\le s\le m-1$. Than the operators satisfy the estimate
of the previous proposition.

Suppose that $m\ge 1$, $s=m$. We have
$$
  I_{m,m}f(w)=(w-\hat w )^m \frac{1}{\pi}\int_{\Omega}\frac{f(\omega )}{w-\omega}\frac{dS_{\omega}}{(w-\hat\omega )^{m+1}}=
$$
$$
  =(w-\hat w )^m \frac{1}{\pi}\int_{\Omega}f(\omega )\left(\frac{\omega -\hat\omega }{(w-\omega )(w-\hat\omega )^{m+2}}+
  \frac{1}{(w-\hat\omega )(w-\hat\omega )^{m+2}}\right)dS_{\omega}=
$$
$$
  =I_{m+1,m}f(w)+L_m f(w) , 
$$  
where $L_m$ is an operator of the type considered in Proposition 22, and we obtain the estimate in $L^p$. 

Now suppose $m\ge 1$, $s=0$. We have
$$
  I_{m,0}f(w)=\frac{1}{\pi}\int_{\Omega}\frac{f(\omega )}{w-\omega}\frac{(\omega -\hat\omega )^m}{(w-\hat\omega )^{m+1}}dS_{\omega}=
$$
$$
  =\frac{1}{\pi}\int_{\Omega}f(\omega )\frac{(\omega -\hat\omega )^m}{(w-\hat\omega )^{m+1}(\hat w-\omega )}\left(1-\frac{w-\hat w}{w-\omega}\right)dS_{\omega}  . 
$$
In the right side we have the sum of two integrals, where the first one is of the type considered in Proposition 22
and the second one satisfies estimate of Proposition 23.

Suppose now $m=0$. We have
$$
  I_{0,0}f(w)=\frac{1}{\pi}\int_{\Omega}\frac{f(\omega )}{(w-\omega )(w-\hat\omega )}dS_{\omega}=
\frac{1}{\pi}\int_{\Omega}\frac{f(\omega )}{(w-\hat\omega )^2}\left(1+\frac{\omega-\hat\omega}{w-\omega}\right)dS_{\omega} .
$$
Again we obtain the sum of two integrals, where the first one is of the type of Proposition 22 and the second one
is $I_{1,0}$.    

The case of $\tilde I_{m,s}$ is analogous. We prove estimate only for $\tilde I_{m,0}$, which we shall
apply below.
$$
 \tilde I_{m,0}f(w)=\frac{1}{\pi}\int_{\Omega}f(\omega )\frac{\hat w -\hat\omega}{(w-\omega )^2}\frac{(\omega -\hat\omega )^m}{(\hat w-\omega )^{m+1}}dS_{\omega}=
$$
$$
  =\frac{1}{\pi}\int_{\Omega}f(\omega )\frac{\hat w -\hat\omega}{w-\omega}\frac{(\omega -\hat\omega )^m}{(w-\hat\omega )^{m+1}(\hat w-\omega )}\left(1-\frac{w-\hat w}{w-\omega}\right)dS_{\omega} .
$$  
Again we have the term satisfying conditions of Proposition 22 and the term satisfying the estimate of
Proposition 23.

Consider now the operator $J_{m,s}$.
We can write $J_{m,s}$ in the form
$$
  J_{m,s}f(w)=\frac{1}{\pi}\int_{\Omega}\frac{f(\omega )}{(w-\omega )^2} dS_{\omega}+
$$
$$
  +(w-\hat w )^s \frac{1}{\pi}\int_{\Omega}\frac{f(\omega )}{(w-\omega )^2}\left[\frac{(\omega -\hat\omega )^{m-s}}{(w-\hat\omega )^m}-\frac{1}{(w-\hat w)^s}\right]dS_{\omega}
$$
The first integral is bounded in $L^p$ by the Kalderon-Zigmund inequality. The second integral
after multiplying by $(w-\hat w )^s$ is   
$$
  \frac{1}{\pi}\int_{\Omega}\frac{f(\omega )}{(w-\omega )^2}\frac{(\omega -\hat\omega )^{m-s}[(w-\hat w)^s -(w-\hat\omega )^s ]+(w-\hat\omega )^s [(\omega -\hat\omega )^{m-s}-(w-\hat\omega )^{m-s}]}
{(w-\hat\omega )^m}dS_{\omega} =
$$   
$$
 =\sum_{k=0}^{s-1}\frac{1}{\pi}\int_{\Omega}f\omega )\frac{\hat w-\hat\omega}{(w-\omega )^2}(\omega -\hat\omega )^{m-s}\frac{(w-\hat w )^{s-1-k}}{(w-\hat\omega )^{m-k}}dS_{\omega}+
$$  
$$
  +\sum_{k=0}^{m-s-1}\frac{1}{\pi}\int_{\Omega}\frac{f(\omega )}{w-\omega}\frac{(\omega -\hat\omega )^{m-s-1-k}}{(w-\hat\omega )^{m-s-k}}dS_{\omega}  .
$$
All items of the both sums are of the types considered in Proposition 34. For example, the last item of the first sum is 
$\tilde I_{m-s+1,0}$ and the last item of the second sum is $I_{0,0}$.

For the completeness we consider also the case of non-integer $s$, though it isn't very essential. 
It is easy to see that boundedness in $L^p$ of the operator $J_{m,s}$ is equivalent to  boundedness in $L^p$ of the operator $J^+_{m,s}$,
where
$$
  J^+_{m,s}f(w) =|w-\hat w |^s \frac{1}{\pi}\int_{\Omega}\frac{f(\omega )}{(w-\omega )^2}\left(\frac{\omega -\hat\omega}{w-\hat\omega}\right)^m \frac{dS_{\omega}}{|\omega -\hat\omega |^s} .
$$  .
Indeed, $\|J^+_{m,s}f\|_p =\|J_{m,s}f^+ \|_p$, where $f^+ (w)=f(w)(w-\hat w)^s |w-\hat w|^{-s}$. 

Define the family of operators $J^+_{m,s,z}$, $z=t+iy$, $0\le t\le 1$
$$
  J^+_{m,s,z}f(w)=|w-\hat w |^{sz}\frac{1}{\pi}\int_{\Omega}\frac{f(\omega )}{(w-\omega )^2}\left(\frac{\omega -\hat\omega}{w-\hat\omega}\right)^m \frac{dS_{\omega}}{|\omega -\hat\omega |^{sz}} .
$$  
It is an analytic family of operators. On the left and right boundaries of the strip $0\le t\le 1$ the $L^p$-norms of 
$J^+_{m,s,z}$ are uniformly bounded. Indeed, these operators differ from the cases already considered only by the multiple
$|\omega -\hat\omega |^{-isy}$ under the integral, which we can include in the function $f$, and by the multiple
$|w-\hat w |^{isy}$, which doesn't change the $L^p$-norm. Thus the conditions of the Stein interpolation theorem 
(see, for example, [RS]) for this
family of operators are satisfied, and we obtain $L^p$-estimates for all $2\le s\le m-1$.
$\Box$

\begin{proposition} 
The operators $T_m$ are bounded in $L^p_s (\Omega )$ for $2\le p<\infty$, $m\ge 3$, $2\le s\le m-1$.
\end{proposition}

{\bf Proof}. Denote $f_s (w)=f(w)(w-\hat w)^s$. It is enough to prove the estimate
$$
  \frac{1}{\pi}\int_{\Omega}|w-\hat w |^{sp}|T_m f(w)|^p dS_w \le \|f_s \|_p^p .
$$
Recalling (6.8), we see that the estimates for $T_m f$ in $L^p_s$ follow from the estimates 
for the integrals $I_{m,s}f_s$ and $J_{m,s}f_s$ obtained in the previous proposition.
$\Box$

\section{The operator $\tilde T_m$ and uniform estimates.} 

Now we return to the program described in Section 1. 
Our purpose is to prove Theorem 2'. We can find a solution to the Beltrami equation
with $\mu$ satisfying conditions of Theorem 2' analogously to the classical method replacing 
the transforms $\mC$ and $\mS$ by the transforms $P_m$ and $T_m$. In such a way
we obtain the solution with estimates of its derivatives with respect to the parameters,
but this solution isn't necessary a homeomorphism mapping of $\Omega$ onto its image.

But we can obtain a $\mu$-quasiconformal homeomorphism if we find a solution to the Beltramy 
equation $F$, which is $\mu$-quasiholomorpic and satisfies the estimate 
$$
  |F_{ww}(w)/F_w (w)|\le \delta /{\rm dist}(w,\pd\Omega )   \eqno (7.1)
$$
with sufficiently small $\delta$. Indeed, $\Omega =h(D)$, where $h$ satisfies estimates (4.38),
and $F\circ h$ is quasiholomorphic on $D$ with the complex dilatation 
$\tilde\mu (z)=\mu\circ h(z) h_z (z)/\overline{h_z (z)}$. It isn't difficult to check
that we can write $F\circ h$ as $\tilde h\circ f^{\tilde\mu}$, where $f^{\tilde\mu}$ is the
$\tilde\mu$-normal map and $h$ is holomorphic, and $|\tilde h''(z)/\tilde h'(z)|\le 1/(1-|z|)$
if the constants $\epsilon$ in (4.38) and $\delta$ in (7.1) are small enough. 
Here we don't give details because we shall return to these matters in the next section.
It follows that $\tilde h$ is univalent (see [Pom]) and, hence, $F$ is a homeomorphism. 

Now, if $F$ is a $\mu$-quasiholomorphic map, then $f=(\log F_w )_w=F_{ww}/F_w$ satisfies the equation
$$
  f_{\bar w}=\mu f_w +\mu_w f+\mu_{ww} .   \eqno (7.2)  
$$ 
We can solve this equation by iteration method. On the first step we find a function $f_1$
satisfying the equation $(f_1 )_{\bar w}=\mu (f_1 )_w +\mu_{ww}$. If $f_1$ is such that
$|(f_1 )_w (w)|\le b/{\rm dist}(w,\pd\Omega )$ with sufficiently small $b$, then we shall solve 
the equation $(f_2 )_{\bar w}=\mu (f_2 )_w +\mu_w f_1 +\mu_{ww}$ and so on. On each step
we must solve the equation
$$
   f_{\bar w}-\mu f_w =G ,
$$
where $|G (w)|\le C/({\rm dist}(w,\pd\Omega ))^2$ for some $C$. We can hope that there
exist solutions represented as 
$$
  P_m ({\rm Id}-\mu T_m )^{-1}G=P_m G+ P_m \mu T_m G+P_m \mu T_m \mu T_m G...   \eqno (7.3)
$$
with appropriate $m$ and that these solutions have the estimate $|f(w)|\le c/{\rm dist}(w,\pd\Omega )$. 
There is a difficulty, the transform $T_m$ hasn't good uniform estimates. The key observation is
that we can change the order of integration in each term of series (7.3) and write these series
as 
$$
  \frac{1}{\pi}\int_{\Omega}{\cal P}_w (\omega )G(\omega )dS_{\omega}+\frac{1}{\pi}\int_{\Omega}\tilde T_m \mu{\cal P}_w (\omega )G(\omega )dS_{\omega}+
$$
$$
  +\frac{1}{\pi}\int_{\Omega}\tilde T_m \mu\tilde T_m \mu{\cal P}_w (\omega )G(\omega )dS_{\omega}+... , 
$$
where ${\cal P}_w (\omega )$ is the kernel of the operator $P_m$
$$
   {\cal P}_w (\omega )= \frac{1}{w-\omega}\left(\frac{\omega -\hat\omega}{w-\hat\omega}\right)^m  \eqno (7.4)
$$  
and $\tilde T_m$ is the transform  
$$
   \tilde T_m f(w)=-\frac{1}{\pi}\int_{\Omega}\frac{f(\omega )}{\omega -w}\left(\frac{w -\hat w}{\omega -\hat w}\right)^m 
\left(\frac{1}{\omega -w}+\frac{m}{\omega -\hat w}\right)dS_{\omega}   .  \eqno (7.5)
$$
This operator has the same kernel as $T_m$ but with transposed $w$ and $\omega$. 

The main reason, why it is useful to change the order of integration is that $\tilde T_m$
contains the multiple $(w -\hat w )^m$, and $w$ isn't a variable of integration. As a result,
$\tilde T_m$ has better uniform estimates than $T_m$. For example, the integral
$$
  |w -\hat w |\int_{\Omega}\frac{dS_{\omega}}{|\omega -\hat w |^2 |\omega -w|}
$$
is bounded and the integral
$$
  \int_{\Omega}\frac{|w -\hat w |}{|\omega -\hat w |^2 |\omega -w|}dS_w 
$$
isn't.
  
Thus we can write the sum of series (7.3) (defined at this moment only formally) as
$$
   \frac{1}{\pi}\int_{\Omega} G(\omega )f_w (\omega )dS_{\omega},   \eqno (7.6)
$$
where $f_w$ satisfy the equation
$$
   f_w -\tilde T_m \mu  f_w ={\cal P}_w  . \eqno (7.7)
$$    
Suppose this equation has a solution that we can represent as 
$$
   \frac{1}{w-\omega}\left(\frac{\omega -\hat\omega}{w-\hat\omega}\right)^2 g_w (\omega )  ,   \eqno (7.8)
$$
where $g_w$ is an uniformly bounded function.
Than, applying the estimates of Proposition 29 and the first estimate of Proposition 30, we
obtain for integral (7.6) the estimate $C/{\rm dist}(w,\pd\Omega )$. Thus we must prove
that for the solution to equation (7.7) there exists representation (7.8).    

In fact, we shall prove a more general assertion. 

\begin{lemma} Suppose $\Omega$ and $\mu$ satisfy the conditions of Theorem 2'.
Suppose a function $P_w$ can be written in the form
$$
  P_w (\omega )=\frac{1}{w-\omega}\left(\frac{\omega -\hat\omega}{w-\hat\omega}\right)^l L_w (\omega ),\,l\le m-6 , \eqno (7.9)
$$  
where $L_w (\omega )$ is uniformly bounded and 
$$
|L_w (\omega )-L_w (\omega_0 )|\le c\frac{|\omega -\omega_0 |}{|\omega -\hat\omega_0 |}  \eqno (7.10)
$$
with some uniform $c$. Then at $m\ge l+6$ the equation
$$
  f_w -\tilde T_m \mu  f_w =P_w    \eqno (7.11)
$$
has a unique solution representable in the form
$$
 f_w (\omega )=\frac{1}{w-\omega}\left(\frac{\omega -\hat\omega}{w-\hat\omega}\right)^l g_w (\omega ) , \eqno (7.12)
$$
where $g_w (\omega )$ is uniformly bounded. 
\end{lemma}

\begin{proposition} The function ${\cal P}_w$ satisfies conditions (7.9), (7.10).
\end{proposition}

{\bf Proof}. The only nontrivial part is to prove estimate (7.10) for 
$$
 L_w (\omega )=\left(\frac{\omega -\hat\omega}{w-\hat\omega}\right)^{m-l}  .
$$ 
Note at first that $|(L_w )_{\omega}(\omega)|\le c|w-\hat\omega |^{-1}\le c'|\omega-\hat\omega |^{-1}$
for some $c'$ and we have an analogous estimate for $|(L_w )_{\bar\omega}|$. From estimate
(6.2) follows that $|\omega-\hat\omega |\ge \delta{\rm dist}(\omega ,\pd\Omega )$ for some $\delta$.
Hence, if $|\omega -\omega_0 |\le \delta|\omega_0 -\hat\omega_0 |/2$, then for some uniform $C$
the derivatives $|(L_w )_{\omega}|$ and  
$|(L_w )_{\bar\omega}|$ are no greater than $C|\omega-\hat\omega_0 |^{-1}$ on the segment 
$[\omega ,\omega_0 ]$, and we obtain estimate (7.10). 

Suppose now that $|\omega -\omega_0 |>\delta|\omega_0 -\hat\omega_0 |/2$. Then 
$|\omega -\hat\omega_0 |\le |\omega_0 -\hat\omega_0 |+|\omega -\omega_0 |\le C|\omega -\omega_0 |$ 
for some $C$ and inequality (7.10) is trivial because  $L_w (\omega )$ is uniformly bounded.
$\Box$.

We shall prove Lemma 5 in the next section. Now we shall obtain some estimates necessary to
the proof.

Apply the transform $\tilde T_m$ to a function $f_w$ of type (7.12). We have
$$
   \tilde T_m f_w (\omega )=\frac{1}{\pi}\int \frac{g_w (t)}{(w-t)(t-\omega )}\left(\frac{t-\hat t}{w-\hat t}\right)^l \left(\frac{1}{t-\omega}+\frac{m}{t-\hat\omega}\right)\left(\frac{\omega -\hat\omega}{t-\hat\omega}\right)^{m-l}dS_t .
$$
Since 
$$
   \frac{1}{(w-t)(t-\omega )}=\frac{1}{w-\omega }\left(\frac{1}{w-t}+\frac{1}{t-\omega}\right)  ,
$$
we see that, if $f_w$ satisfies equation (7.11), then $g_w$ must satisfy the equation
$$
   g_w -{\cal T}_w \mu g_w =L_w ,    \eqno (7.13)
$$
where ${\cal T}_w$ is the transform
$$
   {\cal T}_w f (\omega )=\frac{1}{\pi}\int f(t)\left(\frac{\omega -\hat\omega }{t-\hat\omega}\right)^{m-l}\left(\frac{t-\hat t}{t-\hat\omega}\right)^l \left(\frac{w-\hat\omega}{w-\hat t}\right)^l   
$$
$$
\left(\frac{1}{(t-\omega )^2}-\frac{1}{(t-\omega )(t-w)}+\frac{m}{(t-\omega )(t-\hat\omega )}-\frac{m}{(t-w)(t-\hat\omega )}\right)dS_t  .
\eqno (7.14)
$$
We denote by ${\cal K}_w$ the kernel of this operator. In the rest of this section we study the operator ${\cal T}_w$.

In what follows we denote by $\chi_{\Omega}$ the characteristic function of the domain $\Omega$.
As we shall show below, even at the action of ${\cal K}_w$ on $\chi_{\Omega}$ there appears the term
$(\bar\omega -\bar w)/(\omega -w)$ with singular derivatives. Therefore, when we consider
the action of the transform ${\cal K}_w$, we must distinguish the "bad" part and study 
the action of ${\cal K}_w$ on this part. The next proposition describes the situation.

\begin{proposition} Suppose $\Omega$ and $\mu$ satisfy the conditions of Theorem 2'. Let $r_{wk}$, $k\ge 0$  be the functions
$$
  r_{wk}(\omega )=\left(\frac{\bar\omega -\bar w}{\omega -w}\right)^k \chi_{\Omega}(\omega )  .  \eqno (7.15)
$$  
Suppose $m\ge l+6$. Then
$$
  {\cal T}_w \mu r_{wk}(\omega )=\mu (w)r_{w,k+1}(\omega )+F_{wk}(\omega ),\,|F_{wk}(\omega )|\le ckb  , \eqno (7.16)
$$
with some uniform $c$ independent of $k$ ($b$ is the constant from estimate (4.40)).  

Fix a point $\omega_0 \in\Omega$. Then
$$
  F_{wk}(\omega )-F_{wk}(\omega_0 )=(\omega -\omega_0 )F_{w\omega_0 }(\omega ) ,
$$
where  $F_{w\omega_0 }$ satisfies the estimate
$$
  |F_{w\omega_0 }(\omega )|\le k^2 \frac{cb}{|\omega -\hat\omega_0 |} ,   \eqno (7.17)
$$
with some uniform $c$ independent of $k$.
\end{proposition}

{\bf Proof}. In what follows we denote
$$
  K_w (t,\omega )=\left(\frac{\omega -\hat\omega }{t-\hat\omega}\right)^{m-l}\left(\frac{t-\hat t}{t-\hat\omega}\right)^l \left(\frac{w-\hat\omega}{w-\hat t}\right)^l
$$
We have the representations
$$
 K_w (t,\omega )=1+(t-\omega )R_{\omega}(t,\omega ),\,K_w (t,\omega )=a_w (\omega)+(t-w)R_w (t,\omega ) ,
 \eqno (7.18)
$$
where 
$$
  a_w (\omega )=K_w (w,\omega )=\left(\frac{\omega -\hat\omega}{w-\hat\omega}\right)^{m-l},\,a_w (w)=1  .
$$  

There exists some uniform $0<d_1 \le d_2$ such that the disks $D_{\omega}$ and $D_w$ centered at $\omega$ and $w$
of radii $d_1 |\omega -\hat\omega |\le r_{\omega}\le d_2 |\omega -\hat\omega |$ and $d_1 \le |w -\hat w|r_w \le d_2 |w -\hat w|$ correspondingly 
are contained in $\Omega$. In what follows we shall specify $r_{\omega}$ and $r_w$. 

We need in some estimates, which we collect in the next proposition.   

\begin{proposition} For $t\in D_{\omega}$ there are the estimates
$$
  |R_{\omega}(t,\omega )|\le \frac{C}{|\omega -\hat\omega |},\,
  \left|\frac{\pd}{\pd\omega}R_{\omega}(t,\omega )\right|\le \frac{C}{|\omega -\hat\omega |^2}.  \eqno (7.19)
$$  
$$
  |(a_w )_{(j)}(\omega )|\le C|\omega -\hat\omega |^{-|(j)} ,\, |(j)|\le 3 .  \eqno (7.20)
$$  
For $t\in D_w ,\omega\in D_w$ we have the estimates
$$
  |R_w (t,\omega )|\le \frac{C}{|w-\hat w |},\,\left|\frac{\pd}{\pd\omega}R_w(t,\omega )\right|\le \frac{C}{|w -\hat w |^2},\,
  \left|\frac{\pd}{\pd t}R_w (t,\omega )\right|\le \frac{C}{|w -\hat w |^2},  \eqno (7.21)
$$  
$$
  \left|\frac{\pd^2}{\pd t\pd\omega}R_w (t,\omega )\right|\le \frac{C}{|w -\hat w |^3},\,
  \left|\frac{\pd^2}{\pd t^2}R_w (t,\omega )\right|\le \frac{C}{|w -\hat w |^3},  \eqno (7.22) 
$$   
$$
  \left|\frac{\pd^3}{\pd t^3}R_w (t,\omega )\right|\le \frac{C}{|w -\hat w |^4},  \eqno (7.23) 
$$   
and analogous estimates we have for multi-indexes containing $\bar\omega$ and $\bar t$-derivatives.
All constants are uniform.
\end{proposition} 

{\bf Proof.} 
Denote by $\hat h (\omega)$ the function $\omega \mapsto\hat\omega$. By (5.13), (5.14), we have the estimates
$$
  c_1 \le |\hat h_{(k)}(\omega )|\le c_2,|(k)|=1,\,|\hat h_{(k)}(\omega )|\le \frac{C}{|\omega -\hat\omega |^{|(k)|-1}},|(k)|\ge 2 . 
\eqno (7.24)
$$
We prove only estimates (7.19) - (7.23). The cases of $\bar\omega$ and $\bar t$-derivatives are analogous.

To prove inequalities (7.19) we must estimate $(K_w )_t$ and $(K_w )_{t\omega}$ at $t\in D_{\omega}$.
In what follows we denote by $O$ fractions of the type
$$
  \left(\frac{\omega -\hat\omega }{t-\hat\omega}\right)^{j_1}\left(\frac{t-\hat t}{t-\hat\omega}\right)^{j_2}\left(\frac{w-\hat\omega}{w-\hat t}\right)^{j_3},
$$
where $j_1 \le m-l$, $j_2 \le l$, $j_3 \le l$.

The derivative $(K_w )_t$ contains the terms of the types
$$
  ({t-\hat\omega})^{-1}O,\,\frac{1-\hat h_t (t)}{t-\hat\omega}O,\,\frac{\hat h_t (t)}{w-\hat t}O.
$$ 
The first two terms obviously have estimates $C|\omega -\hat\omega |^{-1}$ (we apply first estimate (7.24)). 
The last term has estimate $C|t-\hat t |^{-1}$. At $t\in D_{\omega}$ it is equivalent to the estimate  
$C|\omega -\hat\omega |^{-1}$ .

Consider now $(K_w )_{t\omega}$. There appear the terms
$$
  \frac{\hat h_{\omega}(\omega )}{(t-\hat\omega )^2}O,\,\frac{1-\hat h_{\omega}(\omega )}{(t-\hat\omega )^2}O,\,
  \frac{\hat h_{\omega}(\omega )}{(t-\hat\omega )((w-\hat t )}O,  \eqno (7.25)
$$
$$
  \frac{(1-\hat h_t (t))\hat h_{\omega}(\omega )}{(t-\hat\omega )^2}O,\,\frac{(1-\hat h_t (t))(1-\hat h_{\omega}(\omega ))}{(t-\hat\omega )^2}O,\,    
  \frac{(1-\hat h_t (t))\hat h_{\omega}(\omega )}{(t-\hat\omega )(w-\hat t )}O,  \eqno (7.26)
$$
$$
  \frac{1-\hat h_{\omega}(\omega )}{(t-\hat\omega )((w-\hat t )}O .  \eqno (7.27)
$$  
Again applying first estimate (7.24) and taking into consideration that $t\in D_{\omega}$, we obtain for all
these terms the estimate $C|\omega -\hat\omega |^{-2}$.

Estimate (7.20) easily follows from (7.24). 

First two estimates (7.21) follow from (7.19) if we put $\omega =w$. To prove the third estimate
we consider the derivative $(t_w )_{t^2}$ at $t\in D_{\omega}$. Analogously to (7.25) - (7.27),
there appear the terms
$$
  \frac{1}{(t-\hat\omega )^2}O,\,\frac{1-\hat h_t (t)}{(t-\hat\omega )^2}O,\,\frac{\hat h_t (t)}{(t-\hat\omega )(w-\hat t )}O ,
 \eqno (7.28)
$$
$$
  \frac{\hat h_{t^2}(t)}{t-\hat\omega}O,\,\frac{(1-\hat h_t (t))\hat h_t (t)}{(t-\hat\omega )(w-\hat t )}O,\,
  \frac{(\hat h_t (t))^2}{(w-\hat t )^2}O,\,\frac{\hat h_{t^2}(t)}{w-\hat t}O .  \eqno (7.29)
$$
As above, we see that all these terms have the estimate $C|\omega -\hat\omega |^{-2}$ 
(here we must apply also second estimate (7.24)).

To prove inequalities (7.22) we estimate $\omega$ and $t$-derivatives of terms (7.28), (7.29).
Again, when we differentiate $O$, there appears the multiples
$$
  \frac{1}{t-\hat\omega},\,\frac{\hat h_{\omega}(\omega )}{t-\hat\omega},\,\frac{1-\hat h_{\omega}(\omega )}{t-\hat\omega},\,
  \frac{\hat h_t (t)}{t-\hat\omega},\,\frac{\hat h_t (t)}{w-\hat t},\,\frac{\hat h_{\omega}(\omega )}{w-\hat t} .  
$$
All these multiples are of order $|\omega -\hat\omega |^{-1}$. Other multiples appearing at
differentiation are
$$
  \frac{1}{(t-\hat\omega )^3},\,\frac{\hat h_{\omega}(\omega )}{(t-\hat\omega )^3},\,\frac{(1-\hat h_t (t))\hat h_{\omega}(\omega )}{(t-\hat\omega )^3},
$$
$$  
  \frac{1-\hat h_t (t)}{(t-\hat\omega )^3},\,\frac{\hat h_{t^2}(t)}{(t-\hat\omega )^2},\,\frac{\hat h_t (t)}{(t-\hat\omega )^2 (w-\hat t )},\,
$$
$$
  \frac{(\hat h_t (t))^2}{(t-\hat\omega )(w-\hat t )^2},\,\frac{\hat h_{t^2}(t)}{(t-\hat\omega )(w-\hat t )},\,\frac{\hat h_{t^3}(t)}{t-\hat\omega},
$$
$$
  \frac{\hat h_{t^2}(t)\hat h_{\omega}(\omega )}{(t-\hat\omega )^2},\,\frac{(1-\hat h_t (t))\hat h_t (t)}{(t-\hat\omega )^2 (w-\hat t )},\,
\frac{(1-\hat h_t (t))\hat h_t (t)\hat h_{\omega}(\omega )}{(t-\hat\omega )^2 (w-\hat t )},
$$
$$  
  \frac{(1-\hat h_t (t))(\hat h_t (t))^2}{(t-\hat\omega )(w-\hat t )^2},\,\frac{(1-\hat h_t (t))\hat h_{t^2}(t)}{(t-\hat\omega )(w-\hat t )},\,  
\frac{h_t (t)\hat h_{t^2}(t)}{(t-\hat\omega )(w-\hat t )},\,  
$$
$$
  \frac{(\hat h_t (t))^3}{(w-\hat t )^3},\,\frac{2\hat h_t (t)\hat h_{t^2}(t)}{(w-\hat t )^2},\,\frac{\hat h_{t^3}(t)}{w-\hat t} .
$$  
As above, we obtain the estimate $C|\omega -\hat\omega |^{-3}$ at $t\in D_{\omega}$. 

To prove inequality (7.23) we must estimate $(K_w )_{t^4}$. As above, at each differentiation we 
either obtain the additional multiples $(t-\hat\omega )^{-1}O(1)$ or $(w-\hat t )^{-1}O(1)$
either replace a derivative of the function $\hat h$ by a derivative of order higher by 1. In either case we 
obtain an expression of order $|\omega -\hat\omega |^{-4}$.
$\Box$

Remark. The fact that we must differentiate up to the forth order and use estimate (7.23) explains 
the conditions $P\ge 4$, $L\ge 4$ of Theorem 2'.

Return now to Proposition 37.  
We shall consider the cases corresponding to various items in the kernel of integral (7.14).  

1) Consider at first the integral
$$
  J_{k1}(\omega )=\frac{1}{\pi}\int_{\Omega}\mu (t) r_{wk}(t)K_w (t,\omega )\frac{dS_t}{(t-\omega )^2}  .
$$  

{\bf Proof of representation (7.16) for $J_{k1}$.}
Suppose at first that $\omega\in D_w$ and consider the integral over the domain $D_w \cup D_{\omega}$. 
We have representation (7.18) and the representation
$$
  \mu (t)=\mu (w)+(t-w)\mu_w (t)  ,
$$  
where $\mu_w$ satisfies the estimates   
$$
|\mu_w (t)|\le cb|w -\hat w|^{-1},
|(\mu_w )_{(j)}(t)|\le cb |w -\hat w|^{-2},\,|(j)|=1  \eqno (7.30) 
$$
for some uniform $c$. 

Thus we consider the integral 
$$
   \frac{1}{\pi}\int_{D_{\omega}\cup D_w} (\mu (w)+(t-w)\mu_w (t))(a_w (\omega )+(t-w)R_w (t,\omega ))r_{wk}(t) \frac{dS_t}{(t-\omega )^2} . \eqno (7.31) 
$$  
At first we estimate the integral
$$
  I_{kw}(\omega)=\frac{1}{\pi}\mu (w)a_w (\omega )\int_{D_w}r_{wk}(t) \frac{dS_t}{(t-\omega )^2}+\frac{1}{\pi}\mu (w)a_w (\omega )\int_{D_{\omega}\setminus D_w}r_{wk}(t) \frac{dS_t}{(t-\omega )^2} .  \eqno (7.32)
$$
The first integral in the right side up to the multiple $\mu (w)a_w (\omega )$ is the value of the Beurling transform of the function $r_{wk}\chi_{D_w}$ in the point $\omega$. 
But it is easy to obtain this transform in the explicit form. Indeed, consider the function equals to
$$
  \frac{(\bar t -\bar w)^{k+1}}{(k+1)(t -w)^k} 
$$
in $D_w$ and to      
$$
  \frac{r_w^{2k+2}}{(k+1)(t-w)^{2k+1}}  
$$
in $\bbC\setminus D_w$. This function is continuous, tends to zero at infinity, and its $\bar t$-derivative
(in the sense of distributions) equals to $r_{wk}\chi_{D_w}$. Hence, $t$-derivative of this function is the
Beurling transform of $r_{wk}\chi_{D_w}$. Thus,
$$
  {\mS}r_{wk}\chi_{D_w}(\omega )=-\frac{k}{k+1}\left(\frac{\bar\omega -\bar w}{\omega -w}\right)^{k+1}=-\frac{k}{k+1}r_{w,k+1}(\omega )    
$$  
at $\omega\in D_w$. We see that we can write the first integral in the right side of (7.32) as
$$
  -\frac{k}{k+1}\mu (w)a_w (\omega )r_{w,k+1}(\omega )=-\frac{k}{k+1}\mu (w) r_{w,k+1}(\omega )+\mu (w)A_w (\omega) ,
$$
where $A_w$ is uniformly bounded. Below we shall prove that $A_w$ satisfies estimate (7.17).

Now we can estimate the second integral in the right side of (7.32) as
$$
 \left|\int_{D_{\omega}}(r_{wk}\chi_{D_{\omega}\setminus D_w}-r_{wk}(\omega ))\frac{dS_t}{(t-\omega )^2}\right|\le
  \max_{D_{\omega}}(|(r_{wk})_t | +|(r_{wk})_{\bar t}|)\int_{D_{\omega}}\frac{dS_t}{|t-\omega |}
  \le C_1 k
$$
for some uniform $C_1$. We obtain for $I_{kw}$ the representation
$$
  I_{kw}(\omega )=-\pi\frac{k}{k+1}r_{w,k+1}(\omega )+kf_w (\omega )
$$  
with uniformly bounded $f_w $.  
 
Now to finish with integral (7.31) we must estimate the integral
$$
  \int_{D_{\omega}\cup D_w}g_w (t,\omega )\frac{(\bar t -\bar w)^k}{(t-w)^{k-1}}\frac{dS_t}{(t-\omega )^2} , \eqno (7.33)
$$  
where $|g_w (t,\omega )|\le cb|w-\hat w|^{-1}$, $|(g_w )_t (t,\omega )|\le cb|w-\hat w|^{-2}$,
and we have analogous estimate for $|(g_w )_{\bar t}$ according to (7.30), (7.20). Consider at first the integral over $D_{\omega}$. Under the integral
we have the function 
$$
  \frac{\tilde g_w (t,\omega )}{(t-\omega )^2}  ,
$$  
where $t$-derivatives of the function $\tilde g_w$ we can estimate as
$$   
  cb|w-\hat w|^{-2}r_w +kcb|w-\hat w|^{-1}\le Ckb|w-\hat w|^{-1}\le Ckb|\omega -\hat\omega |^{-1} .   \eqno (7.34)
$$  
Thus we can estimate integral (7.33) over $D_{\omega}$ as
$$
  \left|\int_{D_{\omega}}\frac{\tilde g_w (t,\omega )}{(t-\omega )^2}dS_t\right|\le 
  \frac{Ckb}{|\omega -\hat\omega |}\int_{D_{\omega}}\frac{dS_t}{|t-\omega |}\le C'kb
$$
for some $C'$ since $D_{\omega}$ is a disk of radii $d|\omega -\hat\omega |$. We estimate
integral (7.33) over $D_w \setminus D_{\omega}$ analogously to the second integral in (7.32)
applying estimate (7.34).

Now, let $\omega\notin D_w$, and consider the integral over $D_{\omega}$. Instead of (7.31) we use the representation    
$$
  \int_{D_{\omega}} (\mu (\omega )+(t-\omega )\mu_{\omega}(t))(1+(t-\omega )R_{\omega}(t,\omega ))r_{wk}(t) \frac{dS_t}{(t-\omega )^2} , 
$$
 Here for $\mu_{\omega}(t)$ and $R_{\omega}(t,\omega )$ we have estimates (7.30) and (7.18).
We act as above but this case is more simple, we don't need in estimates (7.20). We estimate the integral
$$
  \int_{D_{\omega}}r_{wk}(t) \frac{dS_t}{(t-\omega )^2}
$$
as 
$$
  \max_{D_{\omega}}(|(r_{wk})_t | +|(r_{wk})_{\bar t}|)\int_{D_{\omega}}\frac{dS_t}{|t-\omega |}\le kc|\omega -w|^{-1}r_{\omega}\le kC|\omega -\hat\omega |^{-1}r_{\omega}
\le kC   
$$
for some uniform $C$. 

Instead of (7.33) we have the integral
$$ 
 \int_{D_{\omega}}g'_w (t,\omega )r_{wk}(t)\frac{dS_t}{(t-\omega )^2} ,  
$$
where $|g'_w (t,\omega )|\le cb|\omega-\hat\omega |^{-1}$. Acting as above we obtain the estimate
$Ckb$.
 
Now we estimate the integral over $\Omega\setminus (D_w \cup D_{\omega}$.

We set $\omega =g(\zeta )$, $w=g(z)$.  Let $\vpf_{\zeta}$ be the function 
$$
  \vpf_{\zeta}(\tau )=\frac{\zeta -\tau}{1-\bar\zeta \tau} ,
$$
and define $g_{\zeta}=g\circ\vpf_{\zeta}$. Note that $\zeta =\vpf_{\zeta}(0)$. 

Since $|r_{wk}|=1$ and $|\mu | \le\delta|$, it is enough to estimate the integral 
$$
  {\cal I}_{k1}(\zeta )=(1-|\zeta |^2 )^2 \int_{D\setminus g_{\zeta}^{-1}D_{\omega}}  \left|\frac{g_{\zeta}(0)-\widehat{g_{\zeta}(0)}}{g_{\zeta}(\tau )-\widehat{g_{\zeta}(0)}}\right|^{m-l} 
\left|\frac{g_{\zeta}(\tau )-\widehat{g_{\zeta}(\tau )}}{g_{\zeta}(\tau )-\widehat{g_{\zeta}(0)}}\right|^l 
\left|\frac{g(z)-\widehat{g_{\zeta}(0)}}{g(z)-\widehat{g_{\zeta}(\tau )}}\right|^l 
$$
$$
\frac{J_g (\vpf_{\zeta}(\tau ))}{|g_{\zeta}(\tau )-g_{\zeta}(0)|^2}\frac{dS_{\tau}}{|1-\bar\zeta \tau|^4} ,
$$  
where $J_g (t)$ is the Jacobian of the change of the variable $\omega \mapsto t=g^{-1}(\omega )$. 

Now, by inequality (6.4), we have
$$
  |g_{\zeta}(\tau )-g_{\zeta}(0)|\ge c|\vpf_{\zeta}(\tau )-\vpf_{\zeta}(0)| . 
$$
for some uniform $c$. But
$$
  \vpf_{\zeta}(\tau )-\vpf_{\zeta}(0)=\frac{|\tau |(1-|\zeta |^2 )}{|1-\bar\zeta \tau|} .  \eqno (7.35)
$$  
From the other hand, the domain $g_{\zeta}^{-1}D_{\omega}$ contains some domain 
$$
   \{|\vpf_{\zeta}(\tau )-\zeta |\le d'(1-|\zeta |^2 )\}=\left\{ \frac{|\tau |}{|1-\bar\zeta \tau|}\le d'\right\}  ,
$$
where $d'$ is some uniform constant.
Thus $|\tau |\ge d'|1-\bar\zeta \tau|$ on $D\setminus g_{\zeta}^{-1}D_{\omega}$. In particular,
the domain $g_{\zeta}^{-1}D_{\omega}$ contains some disk of radius uniformly bounded from below. 
Applying (7.35), we see that on $D\setminus g_{\zeta}^{-1}D_{\omega}$
$$
  |g_{\zeta}(\tau )-g_{\zeta}(0)|\ge c(1-|\zeta |^2 )   \eqno (7.36)
$$  
for some uniform $c$. We see that it is enough to estimate
the integral
$$
   |g(\zeta )(0)-\widehat{g(\zeta )(0)}|^{m-l}|g(z)-\widehat{g_{\zeta}(0)}|^l \int_D \frac{|g_{\zeta}(\tau )-\widehat{g_{\zeta}(\tau )}|^l dS_{\tau}}
{|g_{\zeta}(\tau )-\widehat{g_{\zeta}(0)}|^m |g(z)-\widehat{g_{\zeta}(\tau )}|^l |1-\bar\zeta \tau|^4}.
$$
Applying (6.2) and (6.5), we see that we must estimate the integral
$$
  (1-|\zeta |^2 )^{m-l}|1-z\bar\zeta |^l |\zeta |^{-m}\int_D \frac{|\vpf_{\zeta}(\tau )-\overline{\vpf_{\zeta}(\tau )}^{-1}|^l dS_{\tau}}{|\vpf_{\zeta}(\tau )-\bar\zeta^{-1}|^m 
|z-\overline{\vpf_{\zeta}(\tau )}^{-1}|^l |1-\bar\zeta \tau|^4} .  \eqno (7.37)
$$
We have
$$
\vpf_{\zeta}(\tau )-\overline{\vpf_{\zeta}(\tau )}^{-1}=\frac{(1-|\zeta |^2 )(1-|\tau |^2 )}{(\bar\zeta -\bar\tau )(1-\tau\bar\zeta )}  ,
\eqno (7.38)
$$ 
$$
  \vpf_{\zeta}(\tau )-\bar\zeta^{-1}=\frac{|\zeta |^2 -1}{\bar\zeta (1-\tau\bar\zeta )} ,  \eqno (7.39)
$$
$$
  (z-\overline{\vpf_{\zeta}(\tau )}^{-1})^{-1}=-\frac{\bar\zeta -\bar\tau}{1-z\bar\zeta -\bar\tau (\zeta -z)}=
\frac{\bar\zeta -\bar\tau}{1-z\bar\zeta}\frac{1}{1-\bar\tau z_{\zeta}},  \eqno (7.40) 
$$
where we denote
$$
  z_{\zeta}=(\vpf_{\zeta})^{-1}(z)=\frac{\zeta -z}{1-z\bar\zeta}  .  \eqno (7.41)
$$

We see that integral (7.35) equals to
$$
  (1-|\zeta |^2 )^m |1-z\bar\zeta |^2 |\zeta |^{-m}\int_D \frac{(1-|\tau |^2 )^l}{|\zeta -\tau |^l |1-\tau\bar\zeta |^l}\frac{|\zeta |^m |1-\tau\bar\zeta |^m}{(1-|\zeta |^2 )^m}
\frac{|\zeta -\tau |^l}{|1-z\bar\zeta |^l (|1-\bar\tau z_{\zeta}|^l}\frac{dS_{\tau}}{|1-\tau\bar\zeta |^4} =
$$
$$
  =\int_D \frac{(1-|\tau |^2 )^l |1-\tau\bar\zeta |^{m-l-4}}{|1-\bar\tau z_{\zeta}|^l}dS_{\tau}  .                            
$$
This integral is uniformly bounded at $m\ge l+4$.$\Box$

For convenience we write here the "general term" that appears in various integrals when we make the change of the variable:
$t=g_{\zeta}(\tau )$ and neglect uniformly bounded multiples. 
We denote
$$
   G_z (\tau ,\zeta) =(1-|\zeta |^2 )^2 \left[\frac{g_{\zeta}(0)-\widehat{g_{\zeta}(0)}}{g_{\zeta}(\tau )-\widehat{g_{\zeta}(0)}}\right]^{m-l} 
\left[\frac{g_{\zeta}(\tau )-\widehat{g_{\zeta}(\tau )}}{g_{\zeta}(\tau )-\widehat{g_{\zeta}(0)}}\right]^l
\left[\frac{g(z)-\widehat{g_{\zeta}(0)}}{g(z)-\widehat{g_{\zeta}(\tau )}}\right]^l \frac{J_g (\vpf_{\zeta}(\tau ))}{|1-\bar\zeta \tau|^4} .
$$
We write our correspondence in the form
$$
  G_z (\tau ,\zeta) \sim (1-|\zeta |^2 )^2 \frac{(1-|\tau |^2 )^l (|1-\tau\bar\zeta |^{m-l-4}}{|1-\bar\tau z_{\zeta}|^l} .  \eqno (7.42)  
$$

{\bf Proof of estimate (7.17) for $J_{k1}$}.

Fix some $\omega_0 \in\Omega$. Note that it is enough to prove estimate (7.17) when $\omega \in D_{\omega_0}$. Indeed,
$|\omega -\omega_0 |>c|\omega -\hat\omega_0 |$ with some uniform $c$ if $|\omega -\omega_0 |>d|\omega_0 -\hat\omega_0 |$ and we obviously have estimate
(7.17) since $|F_{wk}|\le ckb$. 

Let $D_{w1}$, $D_{w2}$, $D_{w3}$ be three disks centered at $w$ of radii $r_{1w}$,  
$r_{2w}=2r_{1w}$, $r_{3w}=3r_{1w}$. There can be
the case when at least one of the points $\omega$ or $\omega_0$ belongs to $D_{w1}$. 
We set $r_{\omega}=r_{\omega_0}=r_{1w}$ in this case.
The other point then belongs to $D_{w2}$ and both disks $D_{\omega}$ and $D_{\omega_0}$
are contained in $D_w$. We put $D_w =D_{w3}$ in this case and $D_w =D_{w1}$ if both 
$\omega$ or $\omega_0$ don't belong to $D_{w1}$. Thus we can consider only the cases: 

a) $D_w$ contains $D_{\omega}$ and $D_{\omega_0}$, $r_{\omega}=r_{\omega_0}=r_w /3$,

b) both $\omega$ and $\omega_0$ don't belong to $D_w$, $|w -\omega |>r_{\omega}=r_{\omega_0}<|w-\omega_0 |$. 

In the following proof we shall consider cases a) and b) separately. 

{\it Case a). The difference of the integrals over $D_w$}. 

We estimate the difference of the integrals
$$
   \frac{1}{\pi}\int_{D_w} \mu (t)K_w (t,\omega )r_{wk}(t)\frac{dS_t}{(t-\omega )^2} -
\frac{1}{\pi}\int_{D_w} \mu (t)K_w (t,\omega_0 )r_{wk}(t)\frac{dS_t}{(t-\omega_0 )^2}  .
$$
We represent this difference as
$$
  \frac{1}{\pi}\int_{D_w} (\mu (w)+(t-w)\mu_w (t))(a_w (\omega )+(t-w)R_w (t,\omega ))r_{wk}(t) \frac{dS_t}{(t-\omega )^2}-
$$
$$
 -\frac{1}{\pi}\int_{D_w} (\mu (w)+(t-w)\mu_w (t))(a_w (\omega_0 )+(t-w)R_w (t,\omega_0 ))r_{wk}(t) \frac{dS_t}{(t-\omega_0 )^2} .
\eqno (7.43)
$$ 

We shall estimate the difference of these integrals in several steps. The integral
$$
  \frac{1}{\pi}\int_{D_w}r_{wk}(t)\left(\frac{a_w (\omega )}{(t-\omega )^2}-\frac{a_w (\omega_0 )}{(t-\omega_0 )^2}\right)dS_t 
$$
equals to $\pi k/(k+1)(a_w (\omega )r_{w,k+1}(\omega )-a_w (\omega_0 )r_{w,k+1}(\omega_0 ))$, as we already
saw. Now $a_w (\omega )=1+(\omega -w)b_w (\omega )$, where $|b_w (\omega )|\le c|\omega -\hat\omega|^{-1}$
and $|(b_w )_{(j)}(\omega )|\le c|\omega -\hat\omega|^{-2}$ at $\omega\in D_w$ and $|(j)|=1$. 
It follows
$$
  |(a_w (\omega )r_{w,k+1}(\omega )-a_w (\omega_0 )r_{w,k+1}(\omega_0 ))-(r_{w,k+1}(\omega )-r_{w,k+1}(\omega_0 ))|\le 
$$
$$
 \le |b_w (\omega )-b_w (\omega_0 )|\frac{|\bar\omega -\bar w|^{k+1}}{|\omega -w|^k}+
 |b_w (\omega_0 )|\left|\frac{(\bar\omega -\bar w)^{k+1}}{(\omega -w)^k}-\frac{(\bar\omega_0 -\bar w)^k}{(\omega_0 -w)^{k-1}}\right|\le
$$
$$   
  \le C|\omega -\omega_0 |\left(\frac{|\bar\omega -\bar w|}{|\omega -\hat\omega|^2} +\frac{k}{|\omega -\hat\omega|}\right)\le\frac{Ck}{|\omega -\hat\omega|} .
$$  
Thus we obtained the representation
$$
  \mu (w)a_w (\omega )\int_{D_w}r_{wk}(t)\frac{dS_t}{(t-\omega )^2}=\pi\frac{k}{k+1}\mu (w)r_{w,k+1}(\omega )+A_w (\omega ) ,
$$
where 
$$
  |A_w (\omega )- A_w (\omega_0 )|\le ck\frac{|\omega -\omega_0 |}{|\omega -\hat\omega|} .
$$  
   
Now consider the integral
$$
 \frac{1}{\pi}\int_{D_w} \frac{(\bar t -\bar w)^k}{(t-w)^{k-1}}[(\mu_w (t)a_w (\omega )+\mu (w)R_w (t,\omega )+(t-w)\mu_w (t)a_w (\omega )R_w (t,\omega ))\frac{1}{(t-\omega )^2} -
$$
$$
  -(\mu_w (t)a_w (\omega_0 )+\mu (w)R_w (t,\omega_0 )+(t-w)\mu_w (t)a_w (\omega_0 )R_w (t,\omega_0 ))\frac{1}{(t-\omega_0 )^2}]dS_t  .  \eqno (7.44)
$$  
We represent the expression in the square brackets under the integral as
$$
  [r_1 (w,\omega )+(t-w)r'_1 (t,\omega )]\left[\frac{1}{(t-\omega )^2}-\frac{1}{(t-\omega_0 )^2}\right]+
$$
$$  
  (\omega -\omega_0 )r_2 (t,\omega ,\omega_0 )\frac{1}{(t-\omega_0 )^2} ,
$$
where, by (7.20), we have the estimates
$$
  |r_1 (w,\omega )|\le\frac{Cb}{|\omega -\hat\omega |},\,|r'_1 (t,\omega )|\le\frac{Cb}{|\omega -\hat\omega |^2},\  
   |r_2 (t,\omega ,\omega_0 )|\le \frac{Cb}{|\omega -\hat\omega |^2}  \eqno (7.45),
$$   
and, by (7.21), (7.22), 
$$ 
 |\pd_{(j),t} r'_1 (t,\omega )|\le\frac{Cb}{|\omega -\hat\omega |^3}, \,|\pd_{(j),t}r_2 (t,\omega ,\omega_0 )|\le\frac{Cb}{|\omega -\hat\omega |^3}, |(j)|=1,
 \eqno (7.46)
$$ 
$$
  |\pd_{(j),t} r'_1 (t,\omega )|\le\frac{Cb}{|\omega -\hat\omega |^4},\, |(j)|=1, \eqno (7.47)
$$ 
where $\pd_{(j),t}$ means $t$-derivative with multi-index $(j)$.

Return now to integral (7.44). At first we consider the integral
$$
  r_1 (w,\omega )\int_{D_w}\frac{1}{\pi}\frac{(\bar t -\bar w)^k}{(t-w)^{k-1}}\left[\frac{1}{(t-\omega )^2}-\frac{1}{(t-\omega_0 )^2}\right]dS_t  .   \eqno (7.48)
$$  
We must calculate the difference in the points $\omega$ and $\omega_0$ of the Beurling transform
of the function $(\bar t -\bar w)^k /(t-w)^{k-1} \chi_{D_w}(t)$. But this transform equals to
$$
  \frac{k-1}{k+1}\frac{(\bar t -\bar w)^{k+1}}{(t-w)^k}  .
$$  
if $t\in D_w$. Indeed, the Cauchy transform of the last function equals to
$$
  \frac{(\bar t -\bar w)^{k+1}}{(k+1)(t -w)^{k-1}} 
$$
in $D_w$ and
$$
  \frac{r_w^{2k+2}}{(k+1)(t-w)^{2k}}  
$$
in $\bbC\setminus D_w$. Indeed, this function is continuous, tends to zero at infinity, and its $\bar t$-derivative
(in the sense of distributions) equals to $(\bar t -\bar w)^k /(t-w)^{k-1} \chi_{D_w}(t)$. Hence, $t$-derivative of this function is the
Beurling transform of $(\bar t -\bar w)^k /(t-w)^{k-1} \chi_{D_w}(t)$. We see that integral (7.48) equals to
$$
   r_1 (w,\omega )\frac{k-1}{k+1}\left[\frac{(\bar\omega -\bar w)^{k+1}}{(\omega -w)^k}-\frac{(\bar\omega_0 -\bar w)^{k+1}}{(\omega_0 -w)^k}\right] .
$$
We obtain the estimate 
$$
  Ck|r_1 (w,\omega )||\omega-\omega_0 |\le\frac{Ckb}{|\omega -\hat\omega |}|\omega-\omega_0 |  ,  
$$
since the first order derivatives of the function $(\bar\omega -\bar w)^{k+1}/(\omega -w)^k$ have 
the uniform estimate $\le k+1$, and for $r_1$ we have estimate (7.45).

Now consider the difference 
$$
  \int_{D_w}r'_1 (t,\omega )\frac{(\bar t -\bar w)^k}{(t-w)^{k-2}}\frac{1}{(t-\omega )^2}dS_t-\int_{D_w}r'_1 (t,\omega )\frac{(\bar t -\bar w)^k}{(t-w)^{k-2}}\frac{1}{(t-\omega_0 )^2}dS_t  .  
$$  
We make the change of the variable $t\mapsto t+\omega -\omega_0$ in the first integral and denote by
$D_{w0}$ the domain $\omega -\omega_0 +D_w$. We get the sum
$$ 
  \int_{D_w \cap D_{w0}}\left[r'_1 (t+\omega -\omega_0 ,\omega )\frac{(\overline{t+\omega -\omega_0} -\bar w)^k}{(t+\omega -\omega_0 -w)^{k-2}}-
  r'_1 (t,\omega )\frac{(\bar t -\bar w)^k}{(t-w)^{k-2}}\right]\frac{dS_t}{(t-\omega_0 )^2} +
$$
$$
  +\int_{D_w \setminus D_{w0}}r'_1 (t,\omega )\frac{(\bar t -\bar w)^k}{(t-w)^{k-2}}\frac{1}{(t-\omega )^2}dS_t +
  \int_{D_{w0}\setminus D_w}r'_1 (t,\omega )\frac{(\bar t -\bar w)^k}{(t-w)^{k-2}}\frac{1}{(t-\omega_0 )^2}dS_t  . \eqno (7.49)
$$  
We represent the first integral as
$$
  \int_{D_w \cap D_{w0}}[r'_1 (t+\omega -\omega_0 ,\omega )-r'_1 (t,\omega )]\frac{(\overline{t+\omega -\omega_0} -\bar w)^k}{(t+\omega -\omega_0 -w)^{k-2}}\frac{dS_t}{(t-\omega_0 )^2}+
$$
$$  
  +\int_{D_w \cap D_{w0}}r'_1 (t,\omega )\left[\frac{(\overline{t+\omega -\omega_0} -\bar w)^k}{(t+\omega -\omega_0 -w)^{k-2}}-\frac{(\bar t -\bar w)^k}{(t-w)^{k-2}}\right]\frac{dS_t}{(t-\omega_0 )^2}
\eqno (7.50)
$$  
We consider the first integral in right side analogously to integral (7.33). We see that we must estimate
$$
  |\omega -\omega_0 |[\sup_{t\in D_w \cap D_{w0}, |(j)|=2}(|\pd_{(j),t}r'_1 (t,\omega )||t+\omega -\omega_0 -w|^2 )+
$$
$$  
  +\sup_{t\in D_w \cap D_{w0},|(j)|=1}(|\pd_{(1),t}r'_1 (t,\omega )|k|t+\omega -\omega_0 -w|)]
  \int_{D_w \cap D_{w0}}\frac{dS_t}{|t-\omega_0 |}
$$ 
But $|t+\omega -\omega_0 -w|\le c|\omega -\hat\omega |$ for some $c$. We apply inequalities (7.46) and (7.47)
and obtain the estimate $Ckb|\omega -\hat\omega |^{-1}|\omega -\omega_0 |$. From the other hand, we can represent
the difference
$$
  \frac{(\overline{t+\omega -\omega_0} -\bar w)^k}{(t+\omega -\omega_0 -w)^{k-2}}-\frac{(\bar t -\bar w)^k}{(t-w)^{k-2}}
$$
as $(\omega -\omega_0 )\tilde R_w (t,\omega ,\omega_0 )$, where for $\tilde R_w$ we have the estimates
$|\tilde R_w (t,\omega ,\omega_0 )|\le Ck|\omega -\hat\omega |$, $|\pd_{(j),t}\tilde R_w (t,\omega ,\omega_0 )|\le Ck^2$, $|(j)|=1$. 
We see that we can estimate the second integral in the right side of (7.49) as
$$
  |\omega -\omega_0 |C[\sup_{t\in D_w \cap D_{w0},|(j)|=1}(|\pd_{(j),t}r'_1 (t,\omega )|k|\omega -\hat\omega |+
   \sup_{t\in D_w \cap D_{w0}}|r'_1 (t,\omega )|k^2 ]
$$
$$   
   \int_{D_w \cap D_{w0}}\frac{dS_t}{|t-\omega_0 |}\le \frac{Ck^2 b}{|\omega -\hat\omega |}|\omega -\omega_0 | .
$$
Here we again apply estimates (7.45), (7.46).

The second and third integrals in (7.49) are over the lunules ${D_w \setminus D_{w0}}$ and ${D_{w0}\setminus D_w}$
of width $\sim|\omega -\omega_0 |$. We estimate, for example, the second integral as
$$
  [\sup_{t\in D_w \setminus D_{w0},|(j)|=1}(|\pd_{(j),t}r'_1 (t,\omega )|||t-w|^2 )+k\sup_{t\in D_w \setminus D_{w0}}(|r'_1 (t,\omega )|||t-w|)]
  \int_0^{\pi}\int_{r_{1\theta}}^{r_{2\theta}}d\rho d\theta \le 
$$
$$  
  \le\frac{Ckb}{|\omega -\hat\omega |}|\omega -\omega_0 |  
$$ 
since $|t-w|\le c|\omega -\hat\omega |$ for some $c$ and we have estimates (7.45), (7.46). 

Now to finish with integral (7.44) we must estimate the integral
$$
  \int_{D_w} r_2 (t,\omega ,\omega_0 )\frac{(\bar t -\bar w)^k}{(t-w)^{k-1}}\frac{dS_t}{(t-\omega_0 )^2} .
$$  
Applying (7.45), (7.46), we obtain the estimate
$$ 
  [\sup_{t\in D_w ,|(j)|=1}(|\pd_{(j),t}r_2 (t,\omega ,\omega_0 )||t-w|)+k\sup_{t\in D_w}(|r_2 (t,\omega ,\omega_0 )|] 
   \int_{D_w}\frac{dS_t}{|t-\omega_0 |}\le \frac{Ckb}{|\omega -\hat\omega |} . 
$$

{\it Case a). Difference of the integrals over $\Omega\setminus D_w$.} 

We put
$\tau =g_w^{-1}(t)=\vpf_{z}^{-1}\circ g^{-1}(t)$. We must estimate the difference of the integrals
$$
  \int_{D\setminus g_w^{-1}D_w} G_w (\tau ,\omega )\mu\circ g_w (\tau )r_{wk}\circ g_w (\tau )\frac{dS_{\tau}}{(g_w (\tau )-\omega )^2}- 
$$
$$
 -\int_{D\setminus g_w^{-1}D_w} G_w (\tau ,\omega_0 )\mu\circ g_w (\tau )r_{wk}\circ g_w (\tau )\frac{dS_{\tau}}{(g_w (\tau )-\omega_0 )^2} ,
$$ 
where
$$
  G_w (\tau ,\omega )=K_w (g_w(\tau ),\omega )J_{g_w}(\tau )= 
  \left(\frac{\omega -\hat\omega}{g_w (\tau )-\hat\omega}\right)^{m-l}\left(\frac{g_w (\tau )-\widehat{g_w (\tau )}}{g_w (\tau )-\hat\omega}\right)^l 
  \left(\frac{w-\hat\omega}{w-\widehat{g_w (\tau )}}\right)^l J_{g_w}(\tau ) .
$$
Here $J_{g_w}$ is the Jacobian of the transformation $t\mapsto\tau = \vpf_z^{-1}\circ g^{-1}(t)$.  

It is enough to estimate the $\omega$-derivative of the integral over $\Omega\setminus D_w$.
Now we can differentiate under the integral, i.e., we must estimate the expression
$$
  \int\frac{((G_w )_{\omega}\mu\circ g_w r_{wk}\circ g_w )(\tau ,\omega )}{(g_w (\tau )-\omega )^2}dS_{\tau}-  
  2\int\frac{(G_w \mu\circ g_w r_{wk}\circ g_w )(\tau ,\omega )}{(g_w (\tau )-\omega )^3}dS_{\tau}
  \eqno (7.51)
$$ 
We have 
$$
  (G_w )_{\omega}(\tau ,\omega )=
(m-l)\left(\frac{1 -\hat h_{\omega}(\omega)}{g_w (\tau )-\hat\omega}+\frac{\omega -\hat\omega}{(g_w (\tau )-\hat\omega )^2}\hat h_{\omega}(\omega )\right)
$$
$$ 
\left(\frac{\omega -\hat\omega}{g_w(\tau )-\hat\omega}\right)^{m-l-1}\left(\frac{g_w (\tau )-\widehat{g_w (\tau )}}{g_w (\tau )-\hat\omega}\right)^l 
  \left(\frac{w-\hat\omega}{w-\widehat{g_w (\tau )}}\right)^l J_{g_w}(\tau )+
$$
$$
  +l\frac{g_w (\tau )-\widehat{g_w (\tau )}}{(g_w (\tau )-\hat\omega )^2}\hat h_{\omega}(\omega )  
  \left(\frac{\omega -\hat\omega}{g_w (\tau )-\hat\omega}\right)^{m-l}\left(\frac{g_w (\tau )-\widehat{g_w (\tau )}}{g_w (\tau )-\hat\omega}\right)^{l-1} 
  \left(\frac{w-\hat\omega}{w-\widehat{g_w (\tau )}}\right)^l J_{g_w}(\tau )-
$$
$$  
  -l\frac{\hat h_{\omega}(\omega )}{w-\widehat{g_w (\tau )}}  
  \left(\frac{\omega -\hat\omega}{g_{\omega}(\tau )-\hat\omega}\right)^{m-l}\left(\frac{g_{\omega}(\tau )-\widehat{g_{\omega}(\tau )}}{g_{\omega}(\tau )-\hat\omega}\right)^l 
  \left(\frac{w-\hat\omega}{w-\widehat{g_{\omega}(\tau )}}\right)^{l-1} J_{g_{\omega}}(\tau )+
$$
$$  
  =G_w (\tau ,\omega )\left((m-l)\frac{1 -\hat h_{\omega}(\omega)}{\omega -\hat\omega}+m\frac{\hat h_{\omega}(\omega )}{g_w (\tau )-\hat\omega}
  +l\frac{\hat h_{\omega}(\omega )}{w-\hat\omega}\right) .  \eqno (7.52)
$$
Now it is easy to see that for $G_w$ we have the same estimate (7.42) as for $G_z$
$$
  G_w (\tau ,\omega) \sim (1-|\zeta |^2 )^2 \frac{(1-|\tau |^2 )^l (|1-\tau\bar\zeta |^{m-l-4}}{|1-\bar\tau z_{\zeta}|^l} .  \eqno (7.53)  
$$ 
Applying (7.52) and (7.36), we see that difference (7.51) has the estimate $C|1-|\zeta |^2 |^{-1}$ at $m\ge l+4$.

{\it Case b). The difference of the integrals over $D_{\omega}$ and $D_{\omega_0}$}. 

This case is more simple.
We put $\tilde K_w (t,\omega )=K_w (t+\omega ,\omega )$. In the integrals over $D_{\omega}$ and
$D_{\omega_0}$ we change the variable $t\mapsto t+\omega$ and $t\mapsto t+\omega_0$ 
correspondingly. Thus we must estimate the integral
$$
 \int_{|t|\le r_{\omega}}[\tilde K_w (t,\omega )r_{wk}(t+\omega )\mu (t+\omega ) -\tilde K_w (t,\omega_0 )r_{wk}(t+\omega_0 )\mu (t+\omega_0 )]\frac{dS_t}{t^2}  .    
$$ 
Applying estimate (7.19) and estimates for derivatives of $r_{wk}$ we represent the difference
in the square brackets as
$$
  (r_{wk}(\omega )\mu (\omega )-r_{wk}(\omega_0 )\mu (\omega_0 ))+tR_1 (t,\omega ,\omega_0 ) ,
$$
where $|R_1 |\le cbk^2 |\omega -\omega_0 ||\omega_0 -\hat\omega_0 |^{-2}$ with some uniform $c$.   
Further,
$$
  (r_{wk}(\omega )\mu (\omega )-r_{wk}(\omega_0 )\mu (\omega_0 ))\int_{|t|\le r_{\omega}}\frac{dS_t}{t^2}=0
$$  
and
$$
\left|\int_{|t|\le r_{\omega}}R_1 (t,\omega , )\frac{dS_t}{t}\right|
  \le\frac{ck^2 b|\omega -\omega_0 |}{|\omega_0 -\hat\omega_0 |^2}\int_{D_{\omega}}\frac{dS_t}{|t|}\le 
  \frac{Cbk^2 |\omega -\omega_0 |}{|\omega_0 -\hat\omega_0 |} .  
$$
It proves the estimate in our case.
 
{\it Case b). The difference of the integrals over $\Omega\setminus D_{\omega_0}$ and 
$\Omega\setminus D_{\omega}$.} At first we unify the domains of integration. 

\begin{proposition} Suppose, as above, $\omega =g(\zeta )$, $r_{\omega}=d|\omega -\hat\omega |^{-1}$. 
If $d$ is small enough with some uniform estimate there exists a homeomorphism $g_{\omega}$
mapping the ring $1/2 \le |\tau |\le 1$ onto the domain $D\setminus D_{\omega}$. We can represent
the map $g_{\omega}$ as a composition
$$
  g_{\omega}=g\circ\vpf_{\zeta}\circ\tilde g_{\omega} ,
$$
where $\tilde g_{\omega}$ is a homeomorphism mapping the ring $1/2 \le |\tau |\le 1$ onto the domain 
$D\setminus \vpf_{\zeta}^{-1}\circ g^{-1}(D_{\omega})$. Denote $\tau'_{\omega}=\tilde g_{\omega}(\tau )$. 
For the map $g_{\omega}$ we get the estimates
$$
  |(g_{\omega})_{\tau}(\tau )|\le C\frac{1-|\zeta |^2}{|1-\bar\zeta \tau'_{\omega}|^2},\,
    |(g_{\omega})_{\omega}(\tau )|\le \frac{C}{|1-\bar\zeta \tau'_{\omega}|^2} ,  \eqno (7.54)
$$
$$	
  |(g_{\omega})_{\tau\omega}(\tau )|\le	C\left(\frac{1}{|1-\bar\zeta \tau'_{\omega}|^4}+\frac{1}{|1-\bar\zeta \tau'_{\omega}|^2 (1-|\tau'_{\omega}|^2 )}\right),\,
 \eqno (7.55)
$$
\end{proposition}

{\bf Proof}. The domain $D_{\omega}$ with the boundary $\{\omega +r_{\omega}e^{i\vpf}\},0\le\vpf \le 2\pi$ 
transforms under the mapping $g^{-1}$ onto some domain with the boundary, described by the equation
$$
  |g(\zeta +\rho e^{i\vpf})-g(\zeta )|^2 =d^2 |\omega -\hat\omega |^2  .  \eqno (7.56)
$$
If $d$ is small enough with some uniform estimate, then this domain is star-like with respect to $\zeta$.
Indeed, we have the estimates 
$$
 |g_z |\le b,\,|g_{\bar z}/g_z |\le\delta <1 ,\, |g_{(j)}(z)|\le B|g(z)-\hat g(z)|^{-1}, |(j)|=2 . 
$$ 
Hence, we can write
$$
  g(\zeta +\rho e^{i\vpf})-g(\zeta )=(g_z (\zeta )e^{i\vpf}+g_{\bar z}(\zeta )e^{-i\vpf})\rho +\Delta_1  ,
$$
where $|\Delta_1 |\le 4Bd\rho$ if $\zeta +\rho e^{i\vpf}$ belongs to the curve described by equation
(7.56). Also,
$$
  \frac{\pd}{\pd\rho}g(\zeta +\rho e^{i\vpf})=(g_z (\zeta )+\Delta_2 )e^{i\vpf}+(g_{\bar z}(\zeta )+\Delta_3 )e^{-i\vpf} ,
$$ 
where $|\Delta_2 |\le 4Bd$, $|\Delta_3 |\le 4Bd$. 
Thus we can write $\rho$-derivative of the left side of equation (7.56) in the form
$$    
  2{\rm Re}[((\overline{g_z (\zeta )}e^{-i\vpf}+\overline{g_{\bar z} (\zeta )}e^{i\vpf})\rho +\overline{\Delta_1})
 ((g_z (\zeta )+\Delta_2 )e^{i\vpf}+(g_{\bar z}(\zeta )+\Delta_3 )e^{-i\vpf})]=
$$
$$
  =2\rho(|g_z (\zeta )|^2 +|g_{\bar z}(\zeta )|^2 )+4\rho {\rm Re}[g_z (\zeta )\overline{g_{\bar z} (\zeta )}e^{i\vpf}]+
g_z (\zeta )\Delta  ,  \eqno (7.57)
$$
where  $|\Delta |\le cBd\rho$ with some uniform $c$.   
The right side of (7.57) by modulus is no less, than $2\rho |g_z (\zeta )|^2 (1-\delta)^2 -b\Delta$ and non-equal to zero
for sufficiently small $d$. It follows that the domain with boundary (7.56) is star-like.

Applying the change of the variable $\tau =\vpf_{\zeta}^{-1}\circ g^{-1}$ we obtain
the domain with the boundary described by the equation: 
$$
  |g\circ\vpf_{g^{-1}(\omega )}(\rho e^{i\vpf}) -\omega |^2 -d^2 |\omega -\hat\omega |^2 =0  .  (7.58)
$$
Here we use the same notations in the chart $\tau$: $\tau =\rho e^{i\vpf}$.
This domain is also star-like with respect to zero and we see that its boundary can be described also by
the equation $\rho =\tilde h(\vpf ,\omega )$. By differentiation of equation (7.58) we get the estimates for the derivatives
$$
  |\tilde h_{\vpf}|\le C,\,|\tilde h_{\omega}|\le C|\omega -\hat\omega |^{-1},\,|\tilde h_{\vpf\omega}|\le C|\omega -\hat\omega |^{-1}, 
  |\tilde h_{\omega^2}|\le C|\omega -\hat\omega |^{-2} ,
$$
$$  
  |\tilde h_{\vpf\omega^2}|\le C|\omega -\hat\omega |^{-2}  .
$$
By homothety along radii we define a diffeomorphism $\tilde g_{\omega}$ mapping the ring $1/2 \le |\tau |\le 1$ onto the domain
$D\setminus \vpf_{\zeta}^{-1}\circ g^{-1}(D_{\omega})$. We use the same notation $\tau$ for the chart in the preimage
of $\tilde g_{\omega}$. The estimates for $\tilde h$ yield following estimates for $\tilde g_{\omega}$
$$
  |(\tilde g_{\omega})_{\tau}|\le C,\,|(\tilde g_{\omega})_{\omega}|\le C|\omega -\hat\omega |^{-1},\,|(\tilde g_{\omega})_{\tau\omega}|\le C|\omega -\hat\omega |^{-1}, 
  \eqno (7.59)
$$
and we have analogous estimates for derivatives containing $\pd\bar\tau$ or $\pd\bar\omega$.

The map $g_{\omega}=g\circ\vpf_{\zeta}\circ\tilde g_{\omega}$ maps the chart $\tau$ onto the original chart $t$.
We obtain estimates (7.53), (7.54) applying (7.59), estimates for derivatives of $\vpf_{\zeta}$, and
the estimate
$$
  |g_{(k)}\circ\vpf_{\zeta}(\tau )|\le C(1-|\vpf_{\zeta}(\tau )|^2 )^{1-|(k)|}=C\left(\frac{|1-\bar\zeta \tau |^2}{(1-|\zeta |^2 )(1-|\tau |^2 )}\right)^{|(k)|-1}
$$
at $|(k)|\ge 1$. $\Box$

Return now to the integral $J_{k1}$ over the domain $\Omega\setminus D_{\omega}$. We can write it as
$$
  \int_{1/2 \le |\tau |\le 1}k_w (\tau ,\omega )\mu\circ g_{\omega}(\tau )r_{wk}\circ g_{\omega}(\tau )\frac{dS_{\tau}}{(g_{\omega}(\tau )-\omega )^2} ,
$$  
where
$$
  k_w (\tau ,\omega )=K_w (g_{\omega}(\tau ),\omega )J_{g_{\omega}}(\tau )= 
  \left(\frac{\omega -\hat\omega}{g_{\omega}(\tau )-\hat\omega}\right)^{m-l}\left(\frac{g_{\omega}(\tau )-\widehat{g_{\omega}(\tau )}}{g_{\omega}(\tau )-\hat\omega}\right)^l 
  \left(\frac{w-\hat\omega}{w-\widehat{g_{\omega}(\tau )}}\right)^l J_{g_{\omega}}(\tau ) .
$$
Here $J_{g_{\omega}}$ is the Jacobian of the transformation $t\mapsto\tau =\tilde g_{\omega}^{-1}\circ\vpf_{\zeta}^{-1}\circ g^{-1}(t)$.

As above, to estimate the difference of the integrals over $\Omega\setminus D_{\omega_0}$ and 
$\Omega\setminus D_{\omega}$ it is enough to estimate $\omega$-derivative of this integral.

Again we can differentiate under the integral, i.e., we must estimate the integral
$$
  \int\frac{((k_w )_{\omega}\mu\circ g_{\omega}r_{wk}\circ g_{\omega})(\tau ,\omega )}{(g_{\omega}(\tau )-\omega )^2}dS_{\tau}-
$$
$$  
  -2\int\frac{(k_w \mu\circ g_{\omega}r_{wk}\circ g_{\omega})(\tau ,\omega )}{(g_{\omega}(\tau )-\omega )^3}((g_{\omega})_{\omega}(\tau )-1)dS_{\tau}+
$$ 
$$  
  +\int\frac{[k_w (\mu\circ g_{\omega})_{\omega}r_{wk}\circ g_{\omega}
  +k_w \mu\circ g_{\omega}(r_{wk}\circ g_{\omega})_{\omega}](\tau ,\omega )}{(g_{\omega}(\tau )-\omega )^2}dS_{\tau} .
\eqno (7.60)
$$

We have
$$
  (k_w )_{\omega}(t,\omega )=
(m-l)\left(\frac{1 -\hat h_{\omega}(\omega)}{g_{\omega}(\tau )-\hat\omega}-\frac{\omega -\hat\omega}{(g_{\omega}(\tau )-\hat\omega )^2}((g_{\omega})_{\omega}(\tau )-\hat h_{\omega}(\omega ))\right)
$$
$$ 
\left(\frac{\omega -\hat\omega}{g_{\omega}(\tau )-\hat\omega}\right)^{m-l-1}\left(\frac{g_{\omega}(\tau )-\widehat{g_{\omega}(\tau )}}{g_{\omega}(\tau )-\hat\omega}\right)^l 
  \left(\frac{w-\hat\omega}{w-\widehat{g_{\omega}(\tau )}}\right)^l J_{g_{\omega}}(\tau )+
$$
$$
  +l\left(\frac{(g_{\omega})_{\omega}(\tau )-(\hat h \circ g_{\omega})_{\omega}(\tau )}{g_{\omega}(\tau )-\hat\omega}-
  \frac{g_{\omega}(\tau )-\widehat{g_{\omega}(\tau )}}{(g_{\omega}(\tau )-\hat\omega )^2}((g_{\omega})_{\omega}(\tau )-\hat h_{\omega}(\omega ))\right)
$$  
$$  
  \left(\frac{\omega -\hat\omega}{g_{\omega}(\tau )-\hat\omega}\right)^{m-l}\left(\frac{g_{\omega}(\tau )-\widehat{g_{\omega}(\tau )}}{g_{\omega}(\tau )-\hat\omega}\right)^{l-1} 
  \left(\frac{w-\hat\omega}{w-\widehat{g_{\omega}(\tau )}}\right)^l J_{g_{\omega}}(\tau )-
$$
$$  
  -l\left(\frac{\hat h_{\omega}(\omega )}{w-\widehat{g_{\omega}(\tau )}}-\frac{w-\hat\omega}{(w-\widehat{g_{\omega}(\tau )})^2}(\hat h \circ g_{\omega})_{\omega}(\tau )\right)
$$
$$  
  \left(\frac{\omega -\hat\omega}{g_{\omega}(\tau )-\hat\omega}\right)^{m-l}\left(\frac{g_{\omega}(\tau )-\widehat{g_{\omega}(\tau )}}{g_{\omega}(\tau )-\hat\omega}\right)^l 
  \left(\frac{w-\hat\omega}{w-\widehat{g_{\omega}(\tau )}}\right)^{l-1} J_{g_{\omega}}(\tau )+
$$
$$  
  +\left(\frac{\omega -\hat\omega}{g_{\omega}(\tau )-\hat\omega}\right)^{m-l}\left(\frac{g_{\omega}(\tau )-\widehat{g_{\omega}(\tau )}}{g_{\omega}(\tau )-\hat\omega}\right)^l 
  \left(\frac{w-\hat\omega}{w-\widehat{g_{\omega}(\tau )}}\right)^l (J_{g_{\omega}})_{\omega}(\tau ) =
$$
$$  
  =k_w (t,\omega )\left((m-l)\frac{1 -\hat h_{\omega}(\omega)}{\omega -\hat\omega}-(m-2l)\frac{(g_{\omega})_{\omega}(\tau )-\hat h_{\omega}(\omega )}{g_{\omega}(\tau )-\hat\omega}
+l\frac{(g_{\omega})_{\omega}(\tau )-(\hat h \circ g_{\omega})_{\omega}(\tau )}{g_{\omega}(\tau )-\widehat{g_{\omega}(\tau )}}-\right.
$$
$$
  -\left.l\frac{\hat h_{\omega}(\omega )}{w-\hat\omega}+l\frac{(\hat h \circ g_{\omega})_{\omega}(\tau )}{w-\widehat{g_{\omega}(\tau )}}+
  \frac{(J_{g_{\omega}})_{\omega}(\tau )}{J_{g_{\omega}}(\tau )}\right) .  \eqno (7.61)
$$
  
Now we note that $(k_w )(\tau ,\omega )=G_z (\tau'_{\omega},\zeta )$ 
and apply estimate (7.42).
Also, we can apply formulas (7.35) and (7.38) - (7.40)
with $\tau =\tau'_{\omega}$. For derivatives of $g_{\omega}$ we have estimates (7.53), (7.54). We obtain   
$$
  (k_w )_{\omega}(t,\omega ) \sim (1-|\zeta |^2 )^2 \frac{|1-|\tau'_{\omega}|^2 )^l (|1-\tau'_{\omega}\bar\zeta |^{m-l-4}}{|1-\bar\tau'_{\omega}z_{\zeta}|^l}
$$
$$
  \left[\frac{1}{1-|\zeta |^2}+\frac{|\bar\zeta -\bar z||1-\tau'_{\omega}\bar\zeta|}{(1-|\zeta |^2 )|(1-|\tau'_{\omega}|^2 )}+
  \frac{|\bar\zeta -\bar z|}{|1-z\bar\zeta ||1-\bar\tau'_{\omega}z_{\zeta}|}+\right.
$$
$$
\left.+\frac{1}{(1-|\zeta |^2 )}\left(\frac{1}{|1-\tau'_{\omega}\bar\zeta|^2}+\frac{1}{1-|\tau'_{\omega}|^2}\right)\right]  
  \le
$$
$$
  \le C(1-|\zeta |^2 )|1-\tau'_{\omega}\bar\zeta |^{m-l-4}\left(\frac{1}{|1-\bar\tau'_{\omega}z_{\zeta}|}+\frac{1}{|1-\tau'_{\omega}\bar\zeta |^2}\right)  
\eqno (7.62)
$$
Again 
$$
  |g_{\omega}(\tau )-\omega |\ge c(1-|\zeta |^2 ) 
$$
in our domain analogously to estimate (7.36).
We see that the first integral in (7.60) has the estimate
$$
  \frac{C}{1-|\zeta |^2}\int_D \left(\frac{|1-\tau'_{\omega}\bar\zeta |^{m-l-4}}{|1-\bar\tau'_{\omega}z_{\zeta}|}+|1-\tau'_{\omega}\bar\zeta |^{m-l-6}\right)dS_{\tau}  .
$$  
We obtain estimate $O(1-|\zeta |^2 )^{-1}$ at $m\ge l+6$. Analogously, the second integral in
(7.60) has the estimate (here we apply estimate (7.53))
$$
  \frac{C}{1-|\zeta |^2}\int_D |1-\tau'_{\omega}\bar\zeta |^{m-l-6}dS_{\tau}  .
$$   
We also obtain estimate $O(1-|\zeta |^2 )^{-1}$ at $m\ge l+6$.    
  
Consider the last integral in (7.60). Applying (7.53), we have
$$
  |(\mu\circ g_{\omega})_{\omega}(\tau )|\le \frac{cb}{1-|\vpf_{\zeta}(\tau'_{\omega})|^2}\frac{1}{|1-\zeta\tau'_{\omega}|^2}\le
 \frac{Cb}{(1-|\zeta |^2 )(1-|\tau'_{\omega}|^2 )}  .
$$ 
Again applying estimate (7.42) to $k_w (\tau ,\omega )=G_z (\tau'_{\omega},\omega )$, we obtain
$$
  \left|\int\frac{(k_w (\mu\circ g_{\omega})_{\omega}r_{wk}\circ g_{\omega})(\tau ,\omega )}{(g_{\omega}(\tau )-\omega )^2}dS_{\tau}\right|\le
  \frac{Cb}{1-|\zeta |^2}\int_D \frac{|1-\tau'_{\omega}\bar\zeta |^{m-l-4}}{|1-\bar\tau'_{\omega}z_{\zeta}|}dS_{\tau}\le
  \frac{Cb}{1-|\zeta |^2} 
$$  
at $m\ge l+4$. Also,
$$
(r_{wk}\circ g_{\omega})_{\omega}(\tau ,\omega )=k\left[\frac{(\overline{g_{\omega}(\tau )}-\bar w )^{k-1}}{(g_{\omega}(\tau )-w)^k}(\overline{g_{\omega}(\tau )})_{\omega}-
\frac{(\overline{g_{\omega}(\tau )}-\bar w )^k}{(g_{\omega}(\tau )-w)^{k+1}}(g_{\omega}(\tau )_{\omega}\right] .
$$
Applying (7.53), we get for this function the estimate
$$
  \frac{Ck}{|1-\tau'_{\omega}\bar\zeta ||z-\vpf_{\zeta}(\tau'_{\omega})|}=\frac{Ck}{|1-\bar\zeta z|}\frac{1}{|\tau'_{\omega}-z_{\zeta}|} .
$$  
Again applying (7.42), we obtain
$$
  \left|\int\frac{(k_w \mu\circ g_{\omega}(r_{wk}\circ g_{\omega})_{\omega})(\tau ,\omega )}{(g_{\omega}(\tau )-\omega )^2}dS_{\tau}\right|\le
  \frac{Ckb}{|1-\bar\zeta z|}\int_D \frac{|1-\tau'_{\omega}\bar\zeta |^{m-l-4}}{|\tau'_{\omega}-z_{\zeta}|}dS_{\tau}\le\frac{Ckb}{1-|\zeta |^2} .
$$
at $m\ge l+4$. It finishes the proof of estimate (7.17) for $J_{kw}$.$\Box$

2) 
$$
  J_{k2}(\omega )=\frac{1}{\pi}\int_{\Omega}\mu (t)r_{wk}(t)K_w (t,\omega )\frac{dS_t}{(t-\omega )(t-w)}  .
$$ 
We follow the same steps as in case 1).

{\bf Proof of representation (7.16) for $J_{k2}$.}
At first we consider the case $\omega\in D_w$ and the integral over $D_w \cup D_{\omega}$, i.e., we consider 
the integral
$$
   \frac{1}{\pi}\int_{D_{\omega}\cup D_w} (\mu (w)+(t-w)\mu_w (t))(a_w (\omega )+(t-w)R_w (t,\omega ))r_{wk}(t) \frac{dS_t}{(t-\omega )(t-w)} .  
$$
We set $r_{\omega}=r_w /2$ and consider at first the integral over $D_w$. 

The integral 
$$
  \frac{1}{\pi}\mu (w)a_w (\omega )\int_{D_w}r_{wk}(t)\frac{dS_t}{(t-\omega )(t-w)}
$$
up to the multiple $\mu (w)a_w (\omega )$ is the Cauchy transform of the function $r_{wk}(t)(t-w)^{-1}\chi_{D_w}(t)$. This
Cauchy transform equals to
$$
  \frac{1}{k+1}\frac{(\bar\omega -\bar w)^{k+1}}{(\omega -w)^{k+1}}=\frac{1}{k+1}r_{w,k+1}(\omega )
$$
at $\omega\in D_w$ and
$$
  \frac{1}{k+1}\frac{r_w^{2k+2}}{(\omega -w)^{2k+2}} 
$$
at $\omega\notin D_w$. Indeed, this function is continuous, tends to zero at infinity, and its $\omega$-derivative
(in the sense of distributions) equals to $r_{wk}(\omega )(\omega -w)^{-1}\chi_{D_w}(\omega )$. 

Analogously to case 1) we obtain the representation
$$
  \mu (w)a_w (\omega )\frac{1}{\pi}\int_{D_w}r_{wk}(t)\frac{dS_t}{(t-\omega )(t-w)}=\frac{\pi}{k+1}\mu (w)r_{w,k+1}(\omega )+\mu (w)A'_w (\omega ) ,
$$
where $A'_w$ is uniformly bounded and, as in case 1), we shall show that it satisfies estimate (7.17).

The integral
$$
  \mu (w)a_w (\omega )\frac{1}{\pi}\int_{D_{\omega}\setminus D_w}r_{wk}(t)\frac{dS_t}{(t-\omega )(t-w)}  
$$ 
we estimate as
$$
  \frac{\mu (w)}{r_w}\frac{1}{\pi}\int_{D_{\omega}}\frac{dS_t}{|t-\omega |}= 2\mu(w )\frac{r_{\omega}}{r_w}=\pi\mu(w )
$$
since $r_w =2r_{\omega}$.

To estimate the part of the integral $J_{k2}(\omega )$ over the domain $D_w \cup D_{\omega}$ it is enough
now to estimate the integral
$$
  \int_{D_w \cup D_{\omega}}g_w (t,\omega )r_{wk}(t)\frac{dS_t}{t-\omega } ,
$$  
where $|g_w (t,\omega )|\le cb|\omega -\hat\omega |^{-1}$. This integral obviously has the estimate
$Cb$.
 
If $\omega\notin D_w$ we estimate the integral
$$
  \int_{D_{\omega}}\mu (t)r_{wk}(t)K_w (t,\omega )\frac{dS_t}{(t-\omega )(t-w)} 
$$
essentially as in case 1). We obtain the estimate
$$
   \frac{cb}{r_w}\int_{D_{\omega}}\frac{dS_t}{|t-\omega |}\le Cb  .
$$   
 
Consider now the integral over $\Omega\setminus D_{\omega}$. Analogously to case 1)
we make the change of the variable $t=g_{\zeta}(\tau )=g\circ\vpf_{\zeta}(\tau )$. We have
$$
 |g_{\zeta}(\tau )-g(z)| \ge \inf_{[z,\vpf_{\zeta}],|(j)|=1}|g_{(j)}||\vpf_{\zeta}(\tau )-z| \ge C\left|\frac{\zeta -z-\tau (1-z\bar\zeta )}{1-\bar\zeta\tau}\right|= 
$$
$$
  =C\left|\frac{(1-z\bar\zeta )(\tau -z_{\zeta})}{1-\bar\zeta\tau}\right|  , \eqno (7.63)
$$
where $z_{\zeta}$ is defined by (7.41). Applying (7.36), (7.42) and (7.63),
we see that it is enough to estimate the integral
$$
  (1-|\zeta |^2 )\int_D \frac{(1-|\tau |^2 )^l |1-\tau\bar\zeta |^{m-l-3}}{|1-\tau z_{\zeta}|^l |1-z\bar\zeta ||\tau -z_{\zeta}|}dS_{\tau}\le
 C\int_D \frac{|1-\tau\bar\zeta |^{m-l-3}}{|\tau-z_{\zeta}|}dS_{\tau}. 
$$
We obtain an uniform boundedness at $m\ge l+3$. 
 
 We obtained the representation
$$
  J_{k2}(\omega )=\frac{\mu (w)}{k+1}r_{w,k+1}(\omega )+\mu (w)A'_w (\omega )+F'_{w2}(\omega ),
$$
where $|A'_w (\omega )|\le C, |F'_{w2}(\omega )|\le Cb$. Note that the sum of the coefficients
$k/(k+1)$ and $1/(k+1)$ of cases 1) and 2) equals 1.

{\bf Proof of estimate (7.17) for $J_{k2}$}.
We can assume, as above, $\omega\in D_{\omega_0}$. 

Let $D_{w1}$, $D_{w2}$, $D_{w3}$, $D_{w4}$ be four disks centered at $w$ of radii $r_{1w}$, 
$r_{2w}=2r_{1w}$, $r_3{w}=3r_{1w}$, $r_{4w}=4r_{1w}$. There can be
the case when at least one of the points $\omega$ or $\omega_0$ belongs to $D_{w2}$. 
We put $r_{\omega}=r_{\omega_0}=r_{1w}$ in this case.
The other point then belongs to $D_{w3}$ and both disks $D_\omega$ and $D_{\omega_0}$ are contained
in $D_w$. We put $D_w =D_{w4}$ in this case. If both 
$\omega$ or $\omega_0$ don't belong to $D_{w2}$ we put $D_w =D_{w1}$. In the last case we can set $r_{\omega}=r_{\omega_0}$
such that $|\omega_0 -w|\ge r_{\omega_0}$ and $|\omega -w|\ge r_{\omega_0}$. Thus we can consider two different cases: .

a) $\omega$ and $\omega_0$ belong to $D_w$, $r_{\omega}=r_{\omega_0}=r_w /4$,

b) both $\omega$ and $\omega_0$ don't belong to $D_w$, $|w -\omega_0 |>r_{\omega_0}$, $|\omega -w|\ge r_{\omega_0}$. 

Again in the proof we consider cases a) and b) separately. 

{\it Case a). The difference of the integrals over $D_w$}. 

We consider the difference
$$
   \frac{1}{\pi}\int_{D_w} \mu (t)K_w (t,\omega )r_{wk}(t)\frac{dS_t}{(t-\omega )(t-w)} -
    \frac{1}{\pi}\int_{D_w} \mu (t)K_w (t,\omega_0 )r_{wk}(t)\frac{dS_t}{(t-\omega_0 )(t-w)}  .
$$
We represent this difference as
$$
  \frac{1}{\pi}\int_{D_w} (\mu (w)+(t-w)\mu_w (t))(a_w (\omega )+(t-w)R_w (t,\omega ))r_{wk}(t) \frac{dS_t}{(t-\omega )(t-w)}-
$$
$$
 -\frac{1}{\pi}\int_{D_w} (\mu (w)+(t-w)\mu_w (t))(a_w (\omega_0 )+(t-w)R_w (t,\omega_0 ))r_{wk}(t) \frac{dS_t}{(t-\omega_0 )(t-w)} .
\eqno (7.64)
$$ 

Consider at first the difference of the integrals over $D_w$. 
As in case 1), there is the term 
$$
  \mu (w)(k+1)^{-1}[a_w (\omega )r_{w,k+1}(\omega )-a_w (\omega_0 )r_{w,k+1}(\omega_0 )]=
$$
$$  
  =\mu (w)(k+1)^{-1}(r_{w,k+1}(\omega )-r_{w,k+1}(\omega_0 ))+\mu (w)(A'_w (\omega )-A'_w (\omega_0 )) ,
$$
and $|A'_w (\omega )-A'_w (\omega_0 )|\le C|\omega -\omega_0 |\omega -\hat\omega |^{-1}$. 

It remains to estimate the integral (see (7.44) 
$$
  \int_{D_w}r_{wk} (t)\left[(r_1 (w,\omega )+(t-w)r'_1 (t,\omega ))\left(\frac{1}{t-\omega}-\frac{1}{t-\omega_0}\right)+\right.
$$
$$  
  \left. +(\omega -\omega_0 )r_2 (t,\omega ,\omega_0 )\frac{1}{t-\omega_0}\right]dS_t ,  \eqno (7.65)
$$
where $r_1$, $r'_1$, $r_2$ satisfy estimates (7.45) - (7.47). 

The integral
$$
  r_1 (w,\omega )\int_{D_w}r_{wk}(t)\left(\frac{1}{t-\omega}-\frac{1}{t-\omega_0}\right)dS_t
$$
up to the constant $\pi r_1 (w,\omega )$ is the difference in the points $\omega$ and $\omega_0$ 
of the Cauchy transform of the function $r_wk \chi_{D_w}$. As above, this difference equals to
$$ 
  \frac{1}{k+1}\left[\frac{(\bar\omega -\bar w)^{k+1}}{(\omega -w)^k}-\frac{(\bar\omega_0 -\bar w)^{k+1}}{(\omega_0 -w)^k}\right] .
$$
and, as in case 1), we obtain the estimate $Cb|\omega -\hat\omega |^{-1}|\omega-\omega_0 |$.  
 
Now consider the difference 
$$
  \int_{D_w}r'_1 (t,\omega )\frac{(\bar t -\bar w)^k}{(t-w)^{k-1}}\frac{dS_t}{t-\omega}dS_t-\int_{D_w}r'_1 (t,\omega )\frac{(\bar t -\bar w)^k}{(t-w)^{k-1}}\frac{dS_t}{t-\omega_0}dS_t  .  
$$  
As in case 1), we make the change of the variable: $t\mapsto t+\omega -\omega_0$ in the first integral and denote by
$D_{w0}$ the domain $\omega -\omega_0 +D_w$. We obtain the sum
$$ 
  \int_{D_w \cap D_{w0}}\left[r'_1 (t+\omega -\omega_0 ,\omega )\frac{(\overline{t+\omega -\omega_0} -\bar w)^k}{(t+\omega -\omega_0 -w)^{k-1}}-
  r'_1 (t,\omega )\frac{(\bar t -\bar w)^k}{(t-w)^{k-1}}\right]\frac{dS_t}{t-\omega_0} +
$$
$$
  +\int_{D_w \setminus D_{w0}}r'_1 (t,\omega )\frac{(\bar t -\bar w)^k}{(t-w)^{k-1}}\frac{dS_t}{t-\omega}dS_t +
  \int_{D_{w0}\setminus D_w}r'_1 (t,\omega )\frac{(\bar t -\bar w)^k}{(t-w)^{k-1}}\frac{dS_t}{t-\omega_0}  . 
\eqno (7.66)
$$  
We again represent the first integral as
$$
  \int_{D_w \cap D_{w0}}[r'_1 (t+\omega -\omega_0 ,\omega )-r'_1 (t,\omega )]\frac{(\overline{t+\omega -\omega_0} -\bar w)^k}{(t+\omega -\omega_0 -w)^{k-1}}\frac{dS_t}{t-\omega_0}+
$$
$$  
  +\int_{D_w \cap D_{w0}}r'_1 (t,\omega )\left[\frac{(\overline{t+\omega -\omega_0} -\bar w)^k}{(t+\omega -\omega_0 -w)^{k-1}}-\frac{(\bar t -\bar w)^k}{(t-w)^{k-1}}\right]\frac{dS_t}{t-\omega_0}
  \eqno (7.67)
$$  
Applying the estimate
$$
  \left|\frac{(\overline{t+\omega -\omega_0} -\bar w)^k}{(t+\omega -\omega_0 -w)^{k-1}}-\frac{(\bar t -\bar w)^k}{(t-w)^{k-1}}\right|\le
  Ck |\omega -\omega_0 |
$$  
and acting as in case 1), we see that integral (7.67) is no greater, than   
$$
  C|\omega -\omega_0 |[\sup_{t\in D_w \cup D_{w0},|(j)|=1}(|\pd_{(j),t}r'_1 (t,\omega )||t+\omega -\omega_0 -w|)+ 
  \sup_{t\in D_w \cup D_{w0}}(|r'_1 (t,\omega )|k|]
$$
$$  
  \int_{D_w \cap D_{w0}}\frac{dS_t}{|t-\omega_0 |}
$$ 
Applying (7.45), (7.46), we obtain for integral (7.67) the estimate
$$
  \frac{Ckb}{|\omega -\hat\omega |}|\omega -\omega_0 | .
$$  

The second and third integrals in (7.66) we estimate analogously to the corresponding terms in (7.49).  
We estimate, for example, the second integral as
$$
  C\sup_{t\in D_w \setminus D_{w0}}(|r'_1 (t,\omega )|||t-w|)|\omega -\omega_0 |
    \le\frac{Ckb}{|\omega -\hat\omega |}|\omega -\omega_0 |  .
$$ 

Now to finish with integral (7.65) we must estimate the integral
$$
  \int_{D_w} r_2 (t,\omega ,\omega_0 )r_{wk}(t)\frac{dS_t}{t-\omega_0} .
$$  
Applying (7.45), we obtain the estimate
$$ 
  \sup_{t\in D_w}(|r_2 (t,\omega ,\omega_0 )|) 
   \int_{D_w}\frac{dS_t}{|t-\omega_0 |}\le \frac{Ckb}{|\omega -\hat\omega |} . 
$$ 

{\it Case a). The difference of the integrals over $\Omega\setminus D_w$}.  

We act as in case 1). After the change of the variable $t=g_w (\tau )$ we differentiate under the integral and
analogously to (7.51) we get the integral
$$
  \int\frac{((G_w )_{\omega}\mu\circ g_w r_{wk}\circ g_w )(\tau ,\omega )}{(g_w (\tau )-\omega )(g_w (\tau )-w)}dS_{\tau}-  
  2\int\frac{(G_w \mu\circ g_w r_{wk}\circ g_w )(\tau ,\omega )}{(g_w (\tau )-\omega )^2 (g_w (\tau )-w)}dS_{\tau}  \eqno (7.68)
$$  
We apply (7.52), (7.53), (7.36) and take into consideration that $|g_w (\tau )-w|\ge r_{\omega}\ge d(1-|\zeta|^2)$ 
for some uniform $d$ if $\omega\in\Omega\setminus D_w$. We obtain the estimate $Cb|\omega -\hat\omega |^{-1}$  
for integral (7.68) at $m\ge l+4$.

{\it Case b)}.

There is a difference by comparison with case 1). We have an additional pole in $w$ and,
after transition to the chart $\tau$ on $D$, we can't differentiate under the integral.

We use the decomposition
$$
   \frac{1}{(t-\omega )(t-w)}=\frac{1}{\omega -w}\left(\frac{1}{t-\omega}-\frac{1}{t-w}\right)  .
$$   
We define
$$
  {\cal J}^1_{kw}(\omega )=\int_{\Omega}K_w (t,\omega )\mu (t)r_{kw}(t)\frac{dS_t}{t-\omega}  
$$
and
$$
  {\cal J}^2_{kw}(\omega )=\int_{\Omega}K_w (t,\omega )\mu (t)r_{kw}(t)\frac{dS_t}{t-w}  .
$$
We write 
$$
  J_{k2}(\omega )-J_{k2}(\omega_0 )=(\omega -w)^{-1}({\cal J}^1_{kw}(\omega )-{\cal J}^1_{kw}(\omega_0 ))+
$$
$$
   +(\omega -w)^{-1}({\cal J}^2_{kw}(\omega )-{\cal J}^2_{kw}(\omega_0 ))+
$$
$$
  +[(\omega -w)^{-1}-(\omega_0 -w)^{-1}]({\cal J}^1_{kw}(\omega_0 )-{\cal J}^2_{kw}(\omega_0 )) .
  \eqno (7.69)
$$  
In what follows we consider several cases corresponding to different terms in the right side.

{\it The difference of the integrals ${\cal J}^1_{kw}(\omega )$ and ${\cal J}^1_{kw}(\omega_0 )$ over 
$D_{\omega}$ and $D_{\omega_0}$}.

We must estimate the difference
$$
  \int_{D_{\omega}}K_w (t,\omega )\mu (t)r_{kw}(t)\frac{dS_t}{t-\omega}-
\int_{D_{\omega_0}}K_w (t,\omega_0 )\mu (t)r_{kw}(t)\frac{dS_t}{t-\omega_0} .  
$$
As in case 1), we change the variable 
$t\mapsto t+\omega$ and $t\mapsto t+\omega_0$ in the integrals over $D_{\omega}$ and
$D_{\omega_0}$ correspondingly and denote $\tilde K_w (t,\omega )=K_w (t+\omega ,\omega )$. 
We get the integral
$$ 
 \int_{|t|\le r_{\omega}} [\tilde K_w (t,\omega )r_{wk}(t+\omega )\mu (t+\omega )- 
  \tilde K_w (t,\omega_0 )r_{wk}(t+\omega_0 )\mu (t+\omega_0 )]\frac{dS_t}{|t|}  .  \eqno (7.70)
$$
If $|t|\le r_{\omega}$ and $s\in [\omega ,\omega_0 ]$, then at $|(j)|=1$ we use the estimates 
$$
  |(\tilde K_w)_{(j),s}(t,s)|\le c|\omega -\hat\omega |^{-1},\,|(\mu_{(j)}(t+s)|\le cb|\omega -\hat\omega |^{-1},  \eqno (7.71)
$$
$$
  |(r_{wk})_{(j)}(t+s)|\le ck|t+s -w |^{-1}\le Ck|\omega -\hat\omega |^{-1}  \eqno (7.72)
$$
Thus for integral (7.70) we obtain the estimate 
$$
   Ckb\frac{|\omega -\omega_0 |}{|\omega -\hat\omega |}\int_{|t|\le r_{\omega}}\frac{dS_t}{|t|}
 \le Ckb .  \eqno (7.73)
$$ 

{\it The difference of the integrals ${\cal J}^1_{kw}(\omega )$ and ${\cal J}^1_{kw}(\omega_0 )$ 
over $\Omega\setminus D_{\omega}$ and $\Omega\setminus D_{\omega_0}$}.-

We proceed analogously to case 1) and apply Proposition 39. 
After change of the variable $t=g_{\omega}(\tau )$ and differentiation we obtain the expression
$$
  \int\frac{((k_w )_{\omega}\mu\circ g_{\omega}r_{wk}\circ g_{\omega})(\tau ,\omega )}{g_{\omega}(\tau )-\omega }dS_{\tau}-
$$
$$  
  -\int\frac{(k_w \mu\circ g_{\omega}r_{wk}\circ g_{\omega})(\tau ,\omega )}{(g_{\omega}(\tau )-\omega )^2}((g_{\omega})_{\omega}(\tau )-1)dS_{\tau}+
$$ 
$$  
  +\int\frac{[k_w (\mu\circ g_{\omega})_{\omega}r_{wk}\circ g_{\omega}
  +k_w \mu\circ g_{\omega}(r_{wk}\circ g_{\omega})_{\omega}](\tau ,\omega )}{g_{\omega}(\tau )-\omega }dS_{\tau} .
$$

The distinction from integrals (7.60) is that we replace the multiple $(g_{\omega}(\tau )-\omega )^{-2}$ with the multiple
$(g_{\omega}(\tau ) -w)^{-1}$ in the first and third integrals and replace the multiple $(g_{\omega}(\tau )-\omega )^{-3}$ 
with the multiple $(g_{\omega}(\tau )-\omega )^{-2}$ in the second integral. As a result, these
integrals have the uniform estimate $Ckb$ instead of $Ckb(1-|\zeta |^2 )^{-1}$. Together with 
estimate (7.73) it yields the estimate 
$$
  |{\cal J}^1_{kw}(\omega )-{\cal J}^1_{kw}(\omega_0 )| \le Ckb|\omega -\omega_0 | .
$$
Since $|\omega -w|\ge c|\omega -\hat\omega |$ for some $c$ in our case, we obtain the    
estimate $Ckb|\omega -\omega_0 ||\omega -\hat\omega |^{-1}$ for the first term in the
right side of identity (7.69).

{\it The difference ${\cal J}^2_{2w}(\omega )-{\cal J}^2_{2w}(\omega_0 )$}.

 The singular term $(t-w)^{-1}$ doesn't depend on $\omega$ and we can
differentiate under the integral. Namely we must show that the integral
$$
  \int_{\Omega} (K_w )_{\omega}(t,\omega )\mu (t)r_{kw}(t)\frac{dS_t}{t-w}
\eqno (7.74)
$$
is uniformly bounded. After obvious calculations we can write
$$
  (K_w )_{\omega}(t,\omega )=m\frac{\hat h_{\omega}(\omega )}{t-\hat\omega}K_w (t,\omega)+
$$
$$  
  +(m-l)\frac{1-\hat h_{\omega}(\omega )}{t-\hat\omega}K'_w (t,\omega )-
  l\frac{\hat h_{\omega}(\omega )}{w-\hat\omega}K_w (t,\omega),  \eqno (7.75)
$$
where 
$$
  K'_w (t,\omega )=\left(\frac{\omega -\hat\omega }{t-\hat\omega}\right)^{m-l-1}\left(\frac{t-\hat t}{t-\hat\omega}\right)^l \left(\frac{w-\hat\omega}{w-\hat t}\right)^l .    
$$
The last term in (7.75) by modulus is no greater than $C|t-\hat\omega |^{-1}|K'_w (t,\omega )|$.
After the change of the variable $t=g\circ\vpf_{\zeta}(\tau )$ applying (7.63), (7.42) and  
(7.39) we obtain the estimate for integral (7.74)
$$
  \sim\int_D \frac{|1-\tau\bar\zeta |^{m-l-2}}{|\tau -z_{\zeta}|}dS_{\tau}+
\int_D \frac{|1-\tau\bar\zeta |^{m-l-3}}{|\tau -z_{\zeta}|}dS_{\tau} .
$$
This integral is uniformly bounded at $m\ge l+3$.
   
{\it The last term in identity (7.69)}.

We have   
$$
  [(\omega -w)^{-1}-(\omega_0 -w)^{-1}]({\cal J}^1_{kw}(\omega_0 )-{\cal J}^2_{kw}(\omega_0 )) =
$$
$$
  =\frac{\omega -\omega_0}{((\omega -w)(\omega_0 -w)}\int_{\Omega}(K_w )_{\omega}(t,\omega_0 )\mu (t)r_{kw}(t)\frac{\omega_0 -w}{(t-\omega_0 )(t-w)}dS_t=
  \frac{\omega -\omega_0}{\omega -w}J_{k2}(\omega_0 ).  
$$
We already proved the uniform estimate for $J_{k2}$. Thus we estimate (7.17) for all terms in(7.69).
$\Box$

3)
$$
  J_{k3}(\omega )=\int_{\Omega}K_w (t,\omega )\mu (t)r_{wk}(t)\frac{dS_t}{(t-\omega )(t-\hat\omega )} . 
$$ 

After two "difficult" cases this case is "simple". We estimate the integral introducing our
usual chart $\tau =\vpf_{\zeta}^{-1}\circ g^{-1}(t)$. Applying (7.35), (7.39), and (7.42), we see
that we can estimate our integral as
$$
  Cb\int_D \frac{|1-\tau\bar\zeta |^{m-l-2}}{|\tau |}dS_{\tau} .
$$
We obtain the estimate $Cb$ at $m\ge l+2$.

{\bf Proof of estimate (7.17) for $J_{k3}(\omega )- J_{k3}(\omega_0 )$}.

If we set $r_w$ and $r_{\omega}=r_{\omega_0}$ as in case 1), we again can consider two cases:

a) $D_w$ contains $D_{\omega}$ and $D_{\omega_0}$, $r_{\omega}=r_{\omega_0}=r_w /3$,

b) both $\omega$ and $\omega_0$ don't belong to $D_w$, $|w -\omega |>r_{\omega}=r_{\omega_0}<|w-\omega_0 |$. 

{\it Case a) The difference of the integrals over $D_w$.} 

We represent the integral as
$$
  \mu (w)\left[\frac{1}{w-\hat\omega}\int_{D_w} \frac{r_{wk}(t)}{t-\omega}dS_t-
  \frac{1}{w-\hat\omega_0}\int_{D_w} \frac{r_{wk}(t)}{t-\omega_0}dS_t\right] +
$$
$$
  +\int_{D_w}r_{wk} (t)(t-w)\frac{r_1 (w,\omega )+(t-w)r'_1 (t,\omega )}{t-\hat\omega}\left(\frac{1}{t-\omega}-\frac{1}{t-\omega_0}\right)dS_t +
$$
$$  
  +(\omega -\omega_0 )\int_{D_w}r_{wk} (t)(t-w)\frac{r_2 (t,\omega ,\omega_0 )}{t-\hat\omega_0}\frac{dS_t }{t-\omega_0}dS_t ,  
  \eqno (7.76)
$$
where for $r$, $r'_1$ and $r_2$ we have estimates (7.45), (7.46). 
As above, the first difference equals to
$$
  \frac{\mu (w)}{k}\left[\frac{1}{w-\hat\omega}\frac{(\bar\omega -\bar w)^{k+1}}{(\omega -w)^k}-  
 \frac{1}{w-\hat\omega_0}\frac{(\bar\omega_0 -\bar w)^{k+1}}{(\omega_0 -w)^k}\right]
$$
and has the estimate $C|\omega -\hat\omega |^{-1}|\omega -\omega_0 |$. 
In the second integral the expression before the brackets has the estimate $cb|\omega -\hat\omega |^{-1}$ 
and its $t$-derivatives have the estimate $ckb|\omega -\hat\omega |^{-2}$. Following 
essentially the same considerations as for integral (7.66) we obtain the estimate $Ckb|\omega -\hat\omega |^{-1}|\omega -\omega_0 |$
The estimate for the last integral in (7.76) follows from the estimate for $|r_2 (t,\omega ,\omega_0 )|$.

{\it Case a). The difference of the integrals over $\Omega\setminus D_w$}.

Analogously to case 1) (integral (7.51)) we must estimate the expression
$$
  \int\frac{((G_w )_{\omega}\mu\circ g_w r_{wk}\circ g_w )(\tau ,\omega )}{(g_w (\tau )-\omega )(g_w (\tau )-\hat\omega )}dS_{\tau}-  
  \int\frac{(G_w \mu\circ g_w r_{wk}\circ g_w )(\tau ,\omega )}{(g_w (\tau )-\omega )^2 (g_w (\tau )-\hat\omega )}dS_{\tau}-
$$
$$
  -(\hat h )_{\omega}(\omega )\int\frac{(G_w \mu\circ g_w r_{wk}\circ g_w )(\tau ,\omega )}{(g_w (\tau )-\omega )(g_w (\tau )-\hat\omega )^2}dS_{\tau} ,
$$  
where the integrals are over $D\setminus g_w^{-1}(D_w )$. The difference with integral (7.51)
is that in the two first integrals we replace one multiple $g_w (\tau )-\omega$ in the denominator by 
the multiple $g_w (\tau )-\hat\omega$ and in the third integral the denominator is 
$(g_w (\tau )-\omega )(g_w (\tau )-\hat\omega )^2$ instead of $(g_w (\tau )-\omega )^3$.
Such replacing doesn't change the order in $1-|\zeta |^2$ and in $|1-\bar\zeta\tau|$ and we have
the same estimate as in case $C(1-|\zeta |^2 )^{-1}$ at $m\ge l+4$ as in case 1. 

{\it Case b)}.

 We consider the difference of the integrals over $D_{\omega}$ and $D_{\omega_0}$ 
analogously to cases 1) and 2). We must estimate the integral
$$
 \int_{|t|\le r_{\omega}} \left[\frac{\tilde K_w (t,\omega )r_{wk}(t+\omega )\mu (t+\omega )}{t+\omega -\hat\omega}- 
  \frac{\tilde K_w (t,\omega_0 )r_{wk}(t+\omega_0 )\mu (t+\omega_0 )}{t+\omega_0 -\hat\omega_0}\right]\frac{dS_t}{|t|}  .  
\eqno (7.77)
$$
Applying estimates (7.71), (7.72) and the estimate
$$
\left|\frac{1}{t+\omega -\hat\omega}-\frac{1}{t+\omega_0 -\hat\omega_0}\right|\le
 C\frac{|\omega -\omega_0 |}{|\omega -\hat\omega |^2}  ,
$$
we obtain for integral (7.77) the estimate $Ckb|\omega -\omega_0 ||\omega -\hat\omega |^{-1}$.

We consider, as in case 1), the difference of the integrals over $\Omega\setminus D_{\omega}$ and $\Omega\setminus D_{\omega_0}$. 
We apply Proposition 39 and again see that the
difference with case 1) is that we replace the multiples $g_{\omega}(\tau )-\omega $ or $(g_{\omega}(\tau )-\omega )^2$ 
by the multiple $g_{\omega}(\tau )-\hat\omega$ or $(g_{\omega}(\tau )-\hat\omega )^2$
in the denominators of the expressions under the integrals. As a result, we obtain the same estimate as
in case 1): $Ckb(1-|\zeta |^2 )^{-1}$ at $m\ge l+6$. $\Box$.

4) 
$$ 
  J_{k4}(\omega )=\int_{\Omega}\mu (t)r_{wk}(t)K_w (t,\omega )\frac{dS_t}{(t-w)(t-\hat\omega )} . 
$$ 
This case is most simple. We write the integral in the chart $\tau =\vpf_{\zeta}^{-1}\circ g^{-1}(t)$,
Applying (7.39), (7.42), and (7.63), we obtain the estimate
$$
  Cb\int_D \frac{|1-\tau\bar\zeta |^{m-l-2}}{|\tau -z_{\zeta}|}dS_{\tau} \le C'b
$$
at $m\ge l+2$.

Now in the our case the differentiation with respect to $\omega$ doesn't
lead to new singularities and for estimation of the difference $J_{k4}(\omega )-J_{k4}(\omega_0 )$
it is enough to estimate the integral from the derivative of the expression under the integral. 
It follows that we must estimate the integral
$$
  \int_{\Omega}\frac{(K_w )_{\omega}(t,\omega )\mu (t)r_{kw}(t)}{(t-\hat\omega )(t-w)}dS_t -
  \hat h_{\omega}(\omega) \int_{\Omega}\frac{(K_w )_{\omega}(t,\omega )\mu (t)r_{kw}(t)}{(t-\hat\omega )^2 (t-w)}dS_t .
$$  
Again introducing the chart $\tau$ and applying (7.39), (7.42), (7.63), and (7.75), we get the estimate
$$
  \sim b(1-|\zeta |^2 )^{-1}\int_D (|1-\tau\bar\zeta |^{m-l-1}+|1-\tau\bar\zeta |^{m-l-2})\frac{dS_{\tau}}{|\tau -z_{\zeta}|} .
$$ 
We obtain the estimate $Cb1-|\zeta |^2 )^{-1}$ at $m\ge 2$. $\Box$
 
\begin{proposition}
Let $f$ be a function defined on $\Omega$ and satisfying the estimates
$$
  |f(\omega )|\le a,\,|f(\omega )-f(\omega_0 )|\le \frac{a_1}{|\omega_0 -\hat\omega_0 |}|\omega -\omega_0 |.
$$
Then for every $0<\alpha <1$ $f$ satisfies the Holder estimate
$$
  |f(\omega )-f(\omega_0 )|\le\frac{a_1^{\alpha}(2a)^{1-\alpha}}{|\omega_0 -\hat\omega_0 |^{\alpha}}|\omega -\omega_0 |^{\alpha} .
$$  
\end{proposition}

{\bf Proof}. Define $b=2a(a_1 /2a)^{\alpha}=a_1^{\alpha}(2a)^{1-\alpha}\,r=(2a/a_1 )|\omega_0 -\hat\omega_0 |$. 
From the identity $a_1 |\omega_0 -\hat\omega_0 |^{-1}r =b|\omega_0 -\hat\omega_0 |^{-\alpha}r^{\alpha}=2a$ 
it follows that
$$
  |f(\omega )-f(\omega_0 )|\le \frac{a_1}{|\omega_0 -\hat\omega_0 |}|\omega -\omega_0 |\le 
  \frac{b}{|\omega_0 -\hat\omega_0 |^{\alpha}}|\omega -\omega_0 |^{\alpha} 
$$
if $|\omega -\omega_0 |\le r$. From the other hand, if $|\omega -\omega_0 |\ge r$, then 
$|f(\omega )-f(\omega_0 )|\le 2a \le b|\omega_0 -\hat\omega_0 |^{-\alpha}|\omega -\omega_0 |^{\alpha}$.
$\Box$

\begin{proposition} 
 Let $f$ be a function defined on $\Omega$ equal to zero in $w$ and satisfying the estimates
$$
  |f(\omega )|\le a, 
$$
$$  
  |f(\omega )-f(\omega_0 )|\le \frac{a_1}{|\omega_0 -\hat\omega_0 |^{\alpha}}|\omega -\omega_0 |^{\alpha}
  \eqno (7.78)
  $$
for some $\alpha >0$. Then
$$
  |{\cal T}_w f (\omega )|\le C(a+a_1 ) ,   
$$
$$
  |{\cal T}_w f (\omega )-{\cal T}_w f (\omega_0 )|\le \frac{C(a+a_1 )}{|\omega_0 -\hat\omega_0 |^{\alpha}}|\omega -\omega_0 |^{\alpha}
\eqno (7.79)
$$
for some uniform $C$.
\end{proposition} 
  
{\bf Proof}. We apply notations of Proposition 37, in particular we use the notations $D_w$, $D_{\omega}$,
$r_w$, $r_{\omega}$. Again we consider four cases.

1)
$$
 J_{1f}(\omega )=\int_{\Omega}f(t)K_w (t,\omega )\frac{dS_t}{(t-\omega )^2} .
$$
 
{\bf Proof of uniform boundedness}. 
Consider at first the integral over $D_{\omega}$. $K_w$ satisfies estimates of Proposition 40 with 
some uniform constants. We have
$$
  \left|\int_{D_\omega}f(t)K_w (t,\omega )\frac{dS_t}{(t-\omega )^2}\right|\le 
  \frac{C(a+a_1 )}{|\omega -\hat\omega |^{\alpha}}\int_{D_\omega}\frac{dS_t}{|t-\omega |^{2-\alpha}}=
$$
$$
  =\frac{2\pi C(a+a_1 )}{|\omega -\hat\omega |^{\alpha}}\int_0^{r_{\omega}}\frac{d\rho}{\rho^{1-\alpha}}\le C'(a+a_1 ) 
$$
with some uniform $C'$ since $r_{\omega}\le d|\omega -\hat\omega |$ for some uniform $d$.  
 
For the integral over $\Omega\setminus D_{\omega}$ we can apply the estimate of Proposition 37, 
case 1) since $f(t)$ is uniformly bounded. Thus we obtain $|J_{1f}(\omega )|\le C(a+a_1 )$.
$\Box$

{\bf Proof of estimate (7.79)}.
It is enough to check inequality (7.79) if $|\omega -\omega_0 |\le r_{\omega}/2$. Indeed, in the opposite case
$|{\cal T}_w f (\omega )-{\cal T}_w f (\omega_0 )|\le C(a+a_1 )=C(a+a_1 )|\omega_0 -\hat\omega_0 |^{-\alpha}|\omega_0 -\hat\omega_0 |^{\alpha}
\le C'|\omega_0 -\hat\omega_0 |^{-\alpha}|\omega -\omega_0 |^{\alpha}$ for some uniform $C'$. 

Let $D'_{\omega_0}$ be the disk of radius $2r_{\omega_0}$ centered at $\omega_0$.
Applying breaking of the identity we reduce the problem to estimations of integrals of the functions
$f_1$ and $f_2$ with the supports in $D'_{\omega_0}$ and $\Omega\setminus D'_{\omega_0}$ and satisfying
Holder estimates of type (7.79) with constants $c(a+a_1 )$.
 
We have the representation
$$
  \int_{D'\omega_0}f(t)K_w (t,\omega )\frac{dS_t}{(t-\omega )^2}=\int_{D'\omega_0}f(t)\frac{dS_t}{(t-\omega )^2}+
 \int_{D'\omega_0}f(t)R_w (t,\omega )\frac{dS_t}{t-\omega} ,  \eqno (7.80)
$$
where $R_w$ satisfies estimates (7.19), (7.20),
$$
  |R_w (t,\omega )|\le \frac{c}{|\omega_0 -\hat\omega_0 |},\,|(R_w )_{\omega}(t,\omega )|\le \frac{c}{|\omega_0 -\hat\omega_0 |^2}  .
$$  
The first integral in right side of (7.80) is the Beurling transform of the function $f_1$, which is bounded with the norm $c_{\alpha}$
in the Holder space $C^{\alpha}$ [As]. Consider the difference
$$
  \int_{D'_{\omega_0}}f(t)R_w (t,\omega )\frac{dS_t}{t-\omega}-\int_{D_{'\omega_0}}f(t)R_w (t,\omega_0 )\frac{dS_t}{t-\omega_0} .  
$$  
Following the same lines as in the proof of Proposition 38 we change the variable in the first integral
$t\mapsto t+\omega -\omega_0$. The disk $D'_{\omega_0}$ transforms into the disk $D'_{1\omega_0}$ and we get
the sum of the integrals
$$
  \int_{D'_{\omega_0}\cap D'_{1\omega_0}}[f(t+\omega -\omega_0)R_w (t+\omega -\omega_0,\omega )-f(t)R_w (t,\omega_0 )]\frac{dS_t}{t-\omega_0}+
$$
$$
  +\int_{D'_{1\omega_0}\setminus D'_{\omega_0}}f(t)R_w (t,\omega )\frac{dS_t}{t-\omega_0}
  +\int_{D'_{\omega_0}\setminus D'_{1\omega_0}}[f(t\omega -\omega_0)R_w (t+\omega -\omega_0,\omega )\frac{dS_t}{t-\omega_0} .  
  \eqno (7.81)
$$  
The first integral is no greater by modulus than
$$
  \int_{D'_{\omega_0}\cap D'_{1\omega_0}}|f(t+\omega -\omega_0)-f(t)|R_w (t+\omega -\omega_0,\omega )|\frac{dS_t}{|t-\omega_0 |}+
$$
$$
  +\int_{D'_{\omega_0}\cap D'_{1\omega_0}}|f(t)|||R_w (t+\omega -\omega_0,\omega )-R_w (t,\omega )|\frac{dS_t}{|t-\omega_0 |}\le
$$
$$
  \le \frac{C(a+a_1 )|\omega -\omega_0 |^{\alpha}}{|\omega_0 -\hat\omega_0 |^{1+\alpha}}\int_{D'_{\omega_0}\cap D'_{1\omega_0}}\frac{dS_t}{|t-\omega_0 |}+
  \frac{Ca|\omega -\omega_0 |}{|\omega_0 -\hat\omega_0 |^2}\int_{D'_{\omega_0}\cap D'_{1\omega_0}}\frac{dS_t}{|t-\omega_0 |}\le
$$
$$
  \le \frac{C(a+a_1 )|\omega -\omega_0 |^{\alpha}}{|\omega_0 -\hat\omega_0 |^{\alpha}}+\frac{Ca|\omega -\omega_0 |}{|\omega_0 -\hat\omega_0 |} .
$$  
But   
$$
   \frac{|\omega -\omega_0 |}{|\omega_0 -\hat\omega_0 |}\le\frac{|\omega -\omega_0 |^{\alpha}}{|\omega_0 -\hat\omega_0 |^{\alpha}}
$$
at $|\omega -\omega_0 |\le |\omega_0 -\hat\omega_0 |$ and we obtain the estimate for the first integral in (7.81)   
$C(a+a_1 )|\omega_0 -\hat\omega_0 |^{-\alpha}|\omega -\omega_0 |^{\alpha}$.  
  
We estimate the second and third integrals in (7.81) as in the proof of Proposition 37 taking into consideration
that the "width" of the lunules ${D'_{1\omega_0}\setminus D'_{\omega_0}}$ and ${D'_{\omega_0}\setminus D'_{1\omega_0}}$
is of order $|\omega -\omega_0 |$ and for $|R_w |$ we have estimate (7.19). We obtain the estimate 
$Ca|\omega_0 -\hat\omega_0 |^{-1}|\omega -\omega_0 |$ and, hence, as above, 
$Ca|\omega_0 -\hat\omega_0 |^{-\alpha}|\omega -\omega_0 |^{\alpha}$.  
  
Now consider the integral
$$
   \int_{\Omega\setminus D_{\omega_0}}f (t)\left[\frac{K_w (t,\omega )}{(t-\omega )^2}-\frac{K_w (t,\omega_0 )}{(t-\omega_0 )^2}\right]dS_t .
$$   
Applying Proposition 39 we see that we must estimate the difference
$$
  \int_{1/2 \le |\tau |\le 1}k_w (\tau ,\omega )f\circ g_{\omega}(\tau )\frac{dS_{\tau}}{(g_{\omega}(\tau )-\omega )^2}-
  \int_{1/2 \le |\tau |\le 1}k_w (\tau ,\omega_0 )f\circ g_{\omega_0}(\tau )\frac{dS_{\tau}}{(g_{\omega_0}(\tau )-\omega_0 )^2}=
$$   ,
$$ 
  =\int_{1/2 \le |\tau |\le 1}(k_w (\tau ,\omega )-k_w (\tau ,\omega_0 ))(f\circ g_{\omega}(\tau )\frac{dS_{\tau}}{(g_{\omega}(\tau )-\omega )^2} +
$$ 
$$
  +\int_{1/2 \le |\tau |\le 1}k_w (\tau ,\omega_0 )(f\circ g_{\omega}(\tau )-f\circ g_{\omega_0}(\tau ))\frac{dS_{\tau}}{(g_{\omega_0}(\tau )-\omega_0 )^2}  .
$$   ,
The first integral in the right side up to the bounded multiples is the integral considered
in case 1) of Proposition 37. We have for this term the estimate $Ca|\omega_0 -\hat\omega_0 |^{-1}|\omega -\omega_0 |\le
Ca|\omega_0 -\hat\omega_0 |^{-\alpha}|\omega -\omega_0 |^{\alpha}$. For the second integral we have the estimate
$$
  a_1 \int_{1/2 \le |\tau |\le 1}|k_w (\tau ,\omega_0 )|\frac{|g_{\omega}(\tau )- g_{\omega_0}(\tau )|^{\alpha}}{|g_{\omega_0}(\tau )-\widehat{g_{\omega_0}(\tau )}|^{\alpha}}\frac{dS_{\tau}}{|g_{\omega_0}(\tau )-\omega_0 |^2}\le
$$   
$$
   \le Ca_1\int_{1/2 \le |\tau |\le 1}\frac{|1-\tau'_{\omega_0}\bar\zeta |^{\alpha}}{(1-|\zeta |^2 )^{\alpha}(1-|\tau'_{\omega_0}|^2 )^{\alpha}}\frac{|\omega -\omega_0 |^{\alpha}}{|1-\tau'_{\omega_0}\bar\zeta |^{2\alpha}}
    \frac{(1-|\tau'_{\omega_0}|^2 )^l (|1-\tau'_{\omega_0}\bar\zeta |^{m-l-4}}{|1-\bar\tau'_{\omega_0} z_{\zeta}|^l}\le
$$
$$
   \le Ca_1 \frac{|\omega -\omega_0 |^{\alpha}}{|\omega_0 -\hat\omega_0 |^{\alpha}}\int_{1/2 \le |\tau |\le 1}\frac{|1-\tau'_{\omega_0}\bar\zeta |^{m-l-4-\alpha}}{|1-\bar\tau'_{\omega_0} z_{\zeta}|}	
  \le Ca_1 \frac{|\omega -\omega_0 |^{\alpha}}{|\omega_0 -\hat\omega_0 |^{\alpha}}
$$
at $m\ge l+5$. Here we applied estimates (7.54) and (7.39), (7.42) with $\tau =\tau'_{\omega_0}$.$\Box$   
   
2)
$$
 J_{2f}(\omega )=\int_{\Omega}f(t)K_w (t,\omega )\frac{dS_t}{(t-\omega )(t-w)} .
$$  

{\bf Proof of uniform boundedness.} As in the proof of Proposition 37 we can consider two cases:

a) $D_{\omega}\subset  D_w$, $r_{\omega}=r_w /2$,

b) $|\omega -w|\ge r_{\omega}$ .

Consider case a). We use the estimate
$$
  \int_{D_{\rho}}\frac{dS_t}{|t|^{1-\alpha}|t-\omega |}\le 10\pi\frac{\rho^{\alpha}}{\alpha}  \eqno (7.82)
$$
for some $C$. Here $D_\rho$ is a disk of radius $\rho$ centered at zero, $\omega\in D_\rho$. 
Indeed, if $t$ lies outside of the union of the disks $|t|\le|\omega |/2$ and $|t-\omega |\le |\omega |/2$, 
then $|t-\omega |\ge |t|/3$ and we see that integral (7.82) is no greater than
$$
  \frac{4\pi}{|\omega |}\int_0^{|\omega |/2}r^{\alpha}dr +\frac{4\pi}{|\omega |^{1-\alpha}}\int_0^{|\omega |/2}dr +  
  6\pi\int_{|\omega |/2}^{\rho}r^{\alpha -1}dr \le 
$$
$$
  \le 4\pi |\omega |^{\alpha} +6\pi{\rho}^{\alpha}/\alpha\le 10\pi{\rho}^{\alpha}/\alpha .
$$  

Now we consider integral $J_{2f}$ over $D_w$. Since $f(w)=0$, we get applying (7.82)
$$
  \left|\int_{D_w}f(t)K_w (t,\omega )\frac{dS_t}{(t-\omega )(t-w)}\right|\le 
  \frac{c(a+a_1 )}{|w-\hat w|^{\alpha}} \int_{D_w}\frac{dS_t}{|t-\omega ||t-w|^{1-\alpha}}\le
$$
$$  
   \le\frac{C(a +a_1 )r_w^{\alpha}}{|w-\hat w|^{\alpha}}\le C(a +a_1 ) .
$$  

Analogously, in case b) 
$$
  \left|\int_{D_{\omega}}f(t)K_w (t,\omega )\frac{dS_t}{(t-\omega )(t-w)}\right|\le
  \frac{ca}{|\omega -\hat\omega |}\int_{D_\omega}\frac{dS_t}{|t-\omega |} \le Ca .
$$  
 
From the other hand, the integral over $\Omega\setminus (D_{\omega}\cup D_w )$ has the same estimate as the integral
over $\Omega\setminus D_{\omega}$ in case 2 of Proposition 37, t.e., it is uniformly bounded.
$\Box$

{\bf Proof of estimate (7.79)}. As in the proof of Proposition 37 we consider two cases:

a) $D_w$ contains $D_{\omega}$ and $D_{\omega_0}$,

b) $|\omega -w|\ge r_{\omega}$ and $|\omega_w -w|\ge r_{\omega_0}=r_{\omega}$.

In case a) we must consider the difference of the integrals over $D_w$. Remind the representation
$$
  K_w (t,\omega )=a_w (\omega )+(t-w)R_w (t,\omega ) ,
$$
where $a_w (\omega )=((\omega -\hat\omega )/(w-\hat\omega ))^{m-2}$. We represent our difference as
$$
  (a_w (\omega )-a_w (\omega_0 ))\int_{D_w}f(t)\frac{dS_t}{(t-\omega )(t-w)}+
  a_w (\omega_0 )\int_{D_w}\frac{f(t)}{t-w}\left(\frac{1}{t-\omega}-\frac{1}{t-\omega_0}\right)dS_t +
$$
$$
  +\int_{D_w}f(t)\left(\frac{R_w (t,\omega )}{t-\omega}-\frac{R_w (t,\omega_0 )}{t-\omega_0}\right)dS_t .
\eqno (7.83)
$$  

Applying (7.82) we estimate the first integral as
$$ 
  Ca_1 \frac{|\omega -\omega_0 |}{|\omega -\hat\omega |^{1+\alpha})}\int_{D_w}\frac{dS_t}{|t-w|^{1-\alpha}|t-\omega |}\le 
  Ca_1 \frac{|\omega -\omega_0 |}{|\omega -\hat\omega |}\le Ca_1 \frac{|\omega -\omega_0 |^{\alpha}}{|\omega -\hat\omega |^{\alpha}} .
$$  
Passing to the second integral in (7.83) we suppose at first that $|\omega -\omega_0 |\ge\max\{|\omega -w|,|\omega_0 -w|\}/2$.
Suppose, for certainty, that this maximum equals to $|\omega -w|/2$. Applying (7.82) and taking into 
consideration that $|t-\omega |$ and $|t-\omega_0 |$ are no less than $|t-w|/2$ at
$|t-w|\ge 2|\omega -w|$, we get  
$$
  \left|\int_{D_w}\frac{f(t)}{t-w}\left(\frac{1}{t-\omega}-\frac{1}{t-\omega_0}\right)dS_t \right|\le  
   \frac{Ca_1}{|w-\hat w|^{\alpha}}\int_{|t-w|\le 2|\omega -w|}\left(\frac{1}{|t-\omega |}+\frac{1}{|t-\omega_0 |}\right)\frac{dS_t}{|t-w|^{1-\alpha}}+
$$
$$
  +\frac{Ca_1}{|w-\hat w|^{\alpha}}|\omega -\omega_0 |\int_{|t-w|\ge 2|\omega -w|}\frac{dS_t}{|t-w|^{1-\alpha}|t-\omega ||t-\omega_0 |} \le
$$
$$
  \le Ca_1 \frac{|\omega -w|^{\alpha}}{|w-\hat w|^{\alpha}}+Ca_1 \frac{|\omega -w|}{|w-\hat w|^{\alpha}}\int_{|\omega -w|}^{r_w}\rho^{\alpha -2}d\rho \le
  ca_1 \frac{|\omega -w|^{\alpha}}{|w-\hat w|^{\alpha}} \le
$$
$$
  \le C'a_1 \frac{|\omega -\omega_0 |^{\alpha}}{|\omega -\hat\omega |^{\alpha}} .
$$

Suppose now that $|\omega -\omega_0 |\le\max\{|\omega -w|,|\omega_0 -w|\}/2$. Suppose again that this maximum equals to $|\omega -w|$.
It also means that $|\omega -\omega_0 |\le\min\{|\omega -w|,|\omega_0 -w|\}=|\omega_0 -w|$. Denote by
$D_{0w}$ the disk $\{|t-w|\le |\omega_0 -w|/2\}$, by $D_{0\omega}$, and by $D_{0\omega_0}$ the disks 
$\{|t-\omega |\le |\omega_0 -w|/2\}$ and $\{|t-\omega_0 |\le |\omega_0 -w|/2\}$ correspondingly, We estimate 
separately the integrals over $D_{0w}$, $D_{0\omega}\cup D_{0\omega_0}$ and $D_w \setminus (D_{0w}\cup D_{0\omega}\cup D_{0\omega_0})$.
We have
$$  
  \left|\int_{D_{0w}}\frac{f(t)}{t-w}\left(\frac{1}{t-\omega}-\frac{1}{t-\omega_0}\right)dS_t \right|\le 
  \frac{Ca_1}{|w-\hat w|^{\alpha}}\int_{|t-w|\le |\omega_0 -w|/2}\left|\frac{1}{t-\omega}-\frac{1}{t-\omega_0}\right|\frac{dS_t}{|t-w|^{1-\alpha}} \le
$$  
$$
   \le Ca_1 \frac{|\omega -\omega_0 |}{|w-\hat w|^{\alpha}}\frac{1}{|\omega_0 -w|^2}\int_0^{|\omega_0 -w|/2}\rho^{\alpha}d\rho \le
   Ca_1 |\omega -\omega_0 |^{\alpha}\frac{|\omega -\omega_0 |^{1-\alpha}}{|\omega_0 -w|}\frac{|\omega_0 -w|^{\alpha}}{|w-\hat w|^{\alpha}}\le
$$
$$   
   \le Ca_1 |\omega -\omega_0 |^{\alpha}|\omega -\hat\omega |^{-\alpha}  .
$$   
Now consider the difference
$$
 \int_{D_{0\omega}\cup D_{0\omega_0}}\frac{f(t)}{t-w}\frac{dS_t}{t-\omega}-\int_{D_{0\omega}\cup D_{0\omega_0}}\frac{f(t)}{t-w}\frac{dS_t}{t-\omega_0} .
$$
We change the variable in the second integral $t\mapsto t+\omega_0 -\omega$ and represent this difference as the sum
of the integrals
$$
  \int_{D_{0\omega}}\left[\frac{f(t)}{t-w}-\frac{f(t+ \omega_0 -\omega )}{t+ \omega_0 -\omega -w}\right]\frac{dS_t}{t-\omega}+
$$
$$  
  +\int_{D_{0\omega_0}\setminus D_{0\omega}}\frac{f(t)}{(t-w)(t-\omega )}dS_t -  
  \int_{D_{0\omega}\setminus D_{0\omega_0}}\frac{f(t+ \omega_0 -\omega )}{(t+ \omega_0 -\omega -w)(t-\omega )}dS_t .
\eqno (7.84)
$$  
The first integral we estimate as
$$
  \int_{D_{0\omega}}|f(t)-f(t+ \omega_0 -\omega )|\frac{dS_t}{|t-w||t-\omega |}+
  \int_{D_{0\omega}}|f(t+ \omega_0 -\omega )|\left|\frac{1}{t-w}-\frac{1}{t+ \omega_0 -\omega -w}\right|\frac{dS_t}{|t-\omega |}\le
$$
$$  
  \le Ca_1 \frac{|\omega -\omega_0 |^{\alpha}}{|\omega -\hat\omega |^{\alpha}}\int_{D_{0\omega}}\frac{dS_t}{|t-w||t-\omega|}+
  C|\omega -\omega_0 |\int_{D_{0\omega}}\frac{|f(t+\omega_0 -\omega )|dS_t}{|t+ \omega_0 -\omega -w||t-w||t-\omega |}\le
$$
$$
  \le Ca_1 \frac{|\omega -\omega_0 |^{\alpha}}{|\omega -\hat\omega |^{\alpha}} \frac{1}{|\omega _0 -w|}\int_0^{|\omega _0 -w|}d\rho +
 Ca_1 \frac{|\omega -\omega_0 |}{|\omega_0 -w|}\frac{1}{|w-\hat w|^{\alpha}}\int_{D_{0\omega}}\frac{|t+ \omega_0 -\omega -w|^{\alpha} dS_t}{|t+ \omega_0 -\omega -w||t-\omega |}\le
$$
$$
  \le Ca_1 \frac{|\omega -\omega_0 |^{\alpha}}{|\omega -\hat\omega |^{\alpha}} +
  Ca_1 \frac{|\omega -\omega_0 |}{|\omega_0 -w|}\frac{1}{|w-\hat w|^{\alpha}}\int_{D_{0\omega}}\frac{dS_t}{|t+ \omega_0 -\omega -w|^{1-\alpha}|t-\omega |}\le
$$
$$
 \le Ca_1 \frac{|\omega -\omega_0 |^{\alpha}}{|\omega -\hat\omega |^{\alpha}}+
  Ca_1 \frac{|\omega -\omega_0 |}{|w-\hat w|^{\alpha}|\omega_0 -w|^{2-\alpha}}\int_0^{|\omega_0 -w|/2}d\rho \le 
$$ 
(since $|t+ \omega_0 -\omega -w|\ge |\omega_0 -w|/2$ if $t\in D_{0\omega}$)  
$$
  \le Ca_1 \frac{|\omega -\omega_0 |^{\alpha}}{|\omega -\hat\omega |^{\alpha}} +
  Ca_1 \frac{|\omega -\omega_0 |^{\alpha}}{|w-\hat w|^{\alpha}}\frac{|\omega -\omega_0 |^{1-\alpha}}{|\omega_0 -w|^{1-\alpha}}\le
  Ca_1 \frac{|\omega -\omega_0 |^{\alpha}}{|\omega -\hat\omega |^{\alpha}} .
$$  

The second and third items in (7.84) are 
integrals over the lunules $D_{0\omega_0}\setminus D_{0\omega}$ and $D_{0\omega}\setminus D_{0\omega_0}$
of width $\le c|\omega -\omega_0 |$. They have the estimates
$$
  Ca_1 \frac{|\omega -\omega_0 |}{|\omega_0 -w|^{1-\alpha}|w-\hat w|^{\alpha}}\le 
  Ca_1 \frac{|\omega -\omega_0 |^{\alpha}}{|w-\hat w|^{\alpha}}\frac{|\omega -\omega_0 |^{1-\alpha}}{|\omega_0 -w|^{1-\alpha}}\le
  Ca_1 \frac{|\omega -\omega_0 |^{\alpha}}{|\omega -\hat\omega |^{\alpha}} .
$$  

Now let estimate the integral over $D_w \setminus (D_{0w}\cup D_{0\omega}\cup D_{0\omega_0})$. Note
that if $t$ belongs to this domain, then $|t-\omega |\ge |t-w|/4$ and $|t-\omega_0 |\ge |t-w|/4$. We have
$$
  \left|\int_{D_w \setminus (D_{0w}\cup D_{0\omega}\cup D_{0\omega_0})}\frac{f(t)}{t-w}\left(\frac{1}{t-\omega}-\frac{1}{t-\omega_0}\right)dS_t \right|\le
$$
$$
  \le Ca_1 \frac{|\omega -\omega_0 |}{|w-\hat w|^{\alpha}}\int_{D_w \setminus (D_{0w}\cup D_{0\omega}\cup D_{0\omega_0})}\frac{dS_t}{|t-w|^{1-\alpha}|t-\omega ||t-\omega_0 |}\le
  Ca_1 \frac{|\omega -\omega_0 |}{|w-\hat w|^{\alpha}}\int_{|\omega_0 -w|/2}^{r_w}\rho^{\alpha -2}d\rho \le
$$
$$
  \le Ca_1 \frac{|\omega -\omega_0 |^{\alpha}}{|w-\hat w|^{\alpha}}\frac{|\omega -\omega_0 |^{1-\alpha}}{|\omega_0 -w|^{1-\alpha}}\le
   Ca_1 \frac{|\omega -\omega_0 |^{\alpha}}{|\omega-\hat\omega |^{\alpha}}  .
$$   

Consider now the third integral in (7.83). We represent it as the sum
$$
  \int_{D_w}(R_w (t,\omega )-R_w (t,\omega_0 ))f(t)\frac{dS_t}{t-\omega }+ 
  \int_{D_w}R_w (t,\omega_0 )f(t)\left(\frac{1}{t-\omega }-\frac{1}{t-\omega_0}\right)dS_t .  \eqno (7.85)
$$
The first integral has the estimate
$$
   Ca \frac{|\omega -\omega_0 |}{|\omega -\hat\omega |^2}\int_{D_w}\frac{dS_t}{|t-\omega |}\le
   Ca \frac{|\omega -\omega_0 |}{|\omega -\hat\omega |}\le Ca \frac{|\omega -\omega_0 |^{\alpha}}{|\omega -\hat\omega |^{\alpha}} .
$$   
The second integral in (7.85) we estimate by our usual method. After the change of the variable
$t\mapsto t+\omega_0 -\omega$ the domain $D_w$ transforms into the domain $D'_w$. We must estimate the sum    
$$
 \int_{D_w \cap D'_w}(R_w (t,\omega_0 )f(t)-R_w (t+\omega_0 -\omega ,\omega_0 )f(t+\omega_0 -\omega ))\frac{dS_t}{t-\omega}+
$$
$$
  +\int_{D_w \setminus D'_w}R_w (t,\omega_0 )f(t)\frac{dS_t}{t-\omega}-
  \int_{D'_w \setminus D_w}R_w (t+\omega_0 -\omega ,\omega_0 )f(t+\omega_0 -\omega )\frac{dS_t}{t-\omega } .
$$  
The first integral has the estimate
$$
  C\left(a\frac{|\omega -\omega_0 |}{|\omega -\hat\omega |^2}+a_1\frac{|\omega -\omega_0 |^{\alpha}}{|\omega -\hat\omega |^{1+\alpha}}\right)\int_{D_w \cap D'_w}\frac{dS_t}{|t-\omega |}\le
$$
$$
  \le C(a+a_1 )\frac{|\omega -\omega_0 |^{\alpha}}{|\omega -\hat\omega |^{\alpha}}  .
$$  
The second and third integrals are over the lunules of width of order $|\omega -\omega_0 |$ and have the estimates
$$
  Ca\frac{|\omega -\omega_0 |}{|\omega -\hat\omega |}\le Ca\frac{|\omega -\omega_0 |^{\alpha}}{|\omega -\hat\omega |^{\alpha}} .
$$  

Consider now case b) and estimate the difference of the integrals
$$  
  \int_{D_{\omega}}f(t)K_w (t,\omega )\frac{dS_t}{(t-\omega )(t-w)}-\int_{D_{\omega_0}}f(t)K_w (t,\omega_0 )\frac{dS_t}{(t-\omega_0 )(t-w)} .
$$  
We again apply the change of the variable $t\mapsto t+\omega_0 -\omega$ in the second integral. 
Denote $D'_{\omega_0}=D_{\omega_0}+\omega_0 -\omega$. Applying the estimates for derivatives of $K_w$ 
and recalling that now $|t-w|\ge d|\omega -\hat\omega |$ for some uniform $d$, we obtain
$$   
  \left|\int_{D_{\omega_0}\cap D'_{\omega_0}}\left(\frac{f(t)K_w (t,\omega )}{t-w}-\frac{f(t+\omega_0 -\omega )K_w (t+\omega_0 -\omega ,\omega_0 )}{t+\omega_0 -\omega -w}\right)\frac{dS_t}{t-\omega}\right|\le
$$
$$
  \le C\left(a\frac{|\omega -\omega_0 |}{|\omega -\hat\omega |^2}+a_1\frac{|\omega -\omega_0 |^{\alpha}}{|\omega -\hat\omega |^{1+\alpha}}\right)\int_{D_{\omega_0}\cap D'_{\omega_0}}\frac{dS_t}{|t-\omega |}\le
$$
$$
  \le C(a+a_1 )\frac{|\omega -\omega_0 |^{\alpha}}{|\omega -\hat\omega |^{\alpha}} ,
$$
$$
  \left|\int_{D_{\omega_0}\setminus D'_{\omega_0}}\frac{f(t)K_w (t,\omega )}{t-w}\frac{dS_t}{t-\omega}\right|\le  
  Ca\frac{|\omega -\omega_0 |}{|\omega -\hat\omega |}\le Ca\frac{|\omega -\omega_0 |^{\alpha}}{|\omega -\hat\omega |^{\alpha}} . 
$$
We have an analogous estimate for the integral over $D'_{\omega_0}\setminus D_{\omega_0}$.  
  
At last consider the integrals
$$
  \int_{\Omega\setminus D_w}\frac{f(t)}{t-w}\left[\frac{K_w (t,\omega )}{t-\omega }-\frac{K_w (t,\omega_0 )}{t-\omega_0 }\right]dS_t 
$$
in case a) and 
$$
 \int_{\Omega\setminus (D_{\omega}\cup D_{\omega_0})}\frac{f(t)}{t-w}\left[\frac{K_w (t,\omega )}{t-\omega }-\frac{K_w (t,\omega_0 )}{t-\omega_0 }\right]dS_t .
$$
in case b). As in case 1) we can write these integral as sums of the terms that up to the bounded multiples are the integrals considered in case 2) of Proposition 37
and the terms containing differences of values of the function $f$. We estimate these terms as in case 1) and obtain the required estimate. $\Box$

3)
$$
  J_{3f}(\omega )=\int_{\Omega}f(t)K_w (t,\omega )\frac{dS_t}{(t-\omega )(t-\hat\omega )}.
$$       
Let prove the uniform boundedness. The integral over $D_{\omega}$ has the estimate
$$
  \left|\int_{D_{\omega}}f(t)K_w (t,\omega )\frac{dS_t}{(t-\omega )(t-\hat\omega )}\right|\le 
  c(a+a_1 )\int_{D_{\omega}}\frac{dS_t}{|t-\omega |^{1-\alpha}|t-\hat\omega |}\le 
$$  
$$
  \le c(a+a_1 )\int_{D_{\omega}}\frac{dS_t}{|t-\omega |^{2-\alpha}}
$$
since $|t-\hat\omega |\ge c|t-\omega |$ for some uniform $c$. Hence, the integral is bounded. 
The integral over $\Omega\setminus D_{\omega}$ reduces to the integral of case 3) of Proposition 37. 
  
Let prove estimate (7.79) in our case. The difference
$$  
  \int_{D_{\omega_0}}f(t)K_w (t,\omega )\frac{dS_t}{(t-\omega )(t-\hat\omega )}-
  \int_{D_{\omega_0}}f(t)K_w (t,\omega_0 )\frac{dS_t}{(t-\omega_0 )(t-\hat\omega_0 )}
$$
we estimate by the usual method, applying the change of the variable in the second integral. Denote by
$D'_{\omega_0}$ the disk $D_{\omega_0}+\omega_0 -\omega$. In the usual way, applying estimates for
$|f(t)-f(t+\omega_0 -\omega )|$, $|K_w (t,\omega )-K_w (t+\omega_0 -\omega ,\omega_0 )|$ and   
$|(t-\hat\omega )^{-1}-(t+\omega_0 -\omega-\hat\omega )^{-1}|$ we obtain
$$
  \left|\int_{D_{\omega_0}\cap D'_{\omega_0}}\left[\frac{f(t)K_w (t,\omega )}{(t-\omega )(t-\hat\omega )}- 
  \frac{f(t+\omega_0 -\omega )K_w (t+\omega_0 -\omega ,\omega_0 )}{(t-\omega )(t+\omega_0 -\omega-\hat\omega )}\right]dS_t\right|\le
$$
$$
  \le C\left(a_1 \frac{|\omega -\omega_0 |^{\alpha}}{|\omega -\hat\omega |^{1+\alpha}}+a\frac{|\omega -\omega_0 |}{|\omega -\hat\omega |^2}\right]\int_{D_{\omega_0}\cap D'_{\omega_0}}\frac{dS_t}{|t-\omega |}\le
$$
$$
  \le C(a+a_1 )\frac{|\omega -\omega_0 |^{\alpha}}{|\omega -\hat\omega |^{\alpha}}  ,
$$  
$$
  \left|\int_{D_{\omega_0}\setminus D'_{\omega_0}}f(t)K_w (t,\omega )\frac{dS_t}{(t-\omega )(t-\hat\omega )}\right|\le
 Ca\frac{|\omega -\omega_0 |}{|\omega -\hat\omega |}\le Ca\frac{|\omega -\omega_0 |^{\alpha}}{|\omega -\hat\omega |^{\alpha}},
$$
and for the integral over $D'_{\omega_0}\setminus D_{\omega_0}$ we obtain an analogous estimate.

We obtain an estimate for the difference of the integrals over $\Omega\setminus D_{\omega_0}$ and  $\Omega\setminus D_{\omega_0}$ 
analogously to the previous cases.  
  
4)
$$
  J_{4f}(\omega )=\int_{\Omega}f(t)K_w (t,\omega )\frac{dS_t}{(t-w)(t-\hat\omega )}.
$$   
Again applying the inequality $|t-\hat\omega |\ge c|t-\omega |$ with some uniform $c$, we can see
that this integral is no greater by modulus than
$$
  C(a+a_1 )\int_{\Omega}\frac{dS_t}{|t-w|^{1-\alpha}|t-\hat\omega |} \le C(a+a_1 ) .
$$
  
Estimation of the difference $J_{4f}(\omega )-J_{4f}(\omega_0 )$ again presents no difficulties since
we can differentiate under the integral and we obtain the estimate analogously to case 4)
of Proposition 37.$\Box$.   
  
\section{Solutions to the Beltrami equation with estimates of derivatives}

{\bf Proof of Lemma 5}. Remind that to prove Lemma 5 means to show that equation (7.13) has an unique bounded solution. 
We represent the right side of this equation as
$$
  L_w (\omega )=A+R_w (\omega ) ,
$$
where $R_w (w) =0$ 

Remind the notation $r_{wk}$, $k\ge 1$ and put $r_{w0}=1$ identically. Introduce a linear space $X$ such that elements of $X$ are 
functions on $\Omega$ of the type
$$
  f(\omega )=\sum c_k r_{wk}(\omega ) +f_w (\omega ), 
$$
where the sum $\sum c_k r_{wk}(\omega )$ converges uniformly and the function $f_w (\omega )$ equals to zero at
$\omega =w$, has a finite $C^0$-norm, and satisfies the Holder condition 
$|f_w (\omega )-f_w (\omega_0 )|\le c|\omega-\omega_0 |^{\alpha}/|\omega_0 -\hat\omega_0 |^{\alpha}$ with some $\alpha <1$. 
The norm on this space is defined as 
$$
  \|\sum c_k (r_{wk}\|_C +\|f_w\|_{C^0}+\|f_w\|_{\alpha}  ,
$$
where $\|f_w\|_{\alpha}$ is defined as 
$$
\sup_{\omega\in\Omega ,\omega_0 \in\Omega}\frac{|f_w(\omega )-f_w(\omega_0 )|}{|\omega -\omega_0 |^{\alpha}}|\omega_0 -\hat\omega_0 |^{\alpha}  .
$$
By Proposition 40 the function $L_w =A+R_w$ belongs to $X$.

The operator ${\cal T}_w \mu$ isn't, in general, contracting on $X$ but, however, we can obtain
the solution to equation (7.13) by an iterations process if we put 
$$
  f_0 =L_w =A+R_w , 
$$
$$
   f_{k+1} =f_k +{\cal T}_w \mu f_k 
$$
at $k\ge 0$ and prove that $\|({\cal T}_w \mu )^k f_0\|_X \le a^k$ for some $a<1$.    

We shall prove by induction
$$
  ({\cal T}_w \mu )^k f_0 =\sum_{i=0}^k c_{ik} r_{wi}+f_{kw} ,
$$
where $f_{kw}(w)=0$ and
$$
  \sum_{i=0}^k|c_{ik}|\le (Cb)^k,\,\|f_{kw}\|_{C^0}\le (Cb)^k ,\,\|f_{kw}\|_{\alpha}\le (Cb)^k  \eqno (8.1)
$$  
for some uniform $C$.

Note at first that, by Proposition 37 at $i\ge 1$,
$$
  c_{i,k+1}=\mu (w)c_{i-1,k}    \eqno (8.2)
$$  
and 
$$
  c_{0,k+1}={\cal T}_w \mu f_{kw}(w)+\sum_{i=0}^k c_{ik}F_{wi}(w) ,\eqno (8.3)
$$
where, by Propositions 37 and 40,
$$
  \|F_{wi}\|_{C^0}\le c(i+1)b,\,\|F_{wi}\|_{\alpha}\le c(i+1)^2 b   \eqno (8.4)
$$
for some uniform $c$. Also, by Proposition 37,
$$
  f_{k+1,w}={\cal T}_w \mu f_{kw}-{\cal T}_w \mu f_{kw}(w)+\sum_{i=0}^k c_{ik}(F_{wi}-F_{wi}(w)) .
$$ 
Applying Propositions 37, 40, and 41, we get 
$$
  |c_{0,k+1}|\le cb(\|f_{kw}\|_{C^0}+\|f_{kw}\|_{\alpha})+cb\sum_{i=0}^k (i+1)|c_{ik}| ,  \eqno (8.5)
$$
$$ 
  \|f_{k+1,w}\|_{C^0}\le cb(\|f_{kw}\|_{C^0}+\|f_{kw}\|_{\alpha})+cb\sum_{i=0}^k (i+1)|c_{ik}| , \eqno (8.6)
$$
$$
  \|f_{k+1,w}\|_{\alpha}\le cb (\|f_{kw}\|_{C^0}+\|f_{kw}\|_{\alpha})+cb\sum_{i=0}^k (i+1)^2 |c_{ik}|  \eqno (8.7)
$$  
with some uniform $c$. In what follows we denote by $\beta$ the value $cb$.

From inductive relations (8.2) and (8.4) - (8.7) follows that the tuple $(\{|c_{ik}|\},\|f_{kw}\|_{C^0},\|f_{kw}\|_{\alpha})$ is majorated by the tuple  
$(\{d_{ik}\},A_k ,B_k$, for which we have the inductive relations
$$  
  d_{i,k+1}=\beta d_{i-1,k},\,i\ge 1 ,
$$
$$
  d_{0,k+1}=\beta\sum_{i=0}^k (i+1)d_{ik}+\beta (A_k +B_k ),
$$
$$
  A_{k+1}=\beta\sum_{i=0}^k (i+1)d_{ik}+\beta (A_k +B_k ),
$$
$$
  B_{k+1}=\beta\sum_{i=0}^k (i+1)^2 d_{ik}+\beta (A_k +B_k ),
$$
and we can put, replacing, if it is necessary, the constant $c$: $d_{00}=A_0 =B_0 =1$.  
By induction we easy get
$$
  d_{ik}=\beta^i d_{0,k-i}   \eqno (8.8)
$$
and, hence, we can rewrite our inductive relations as 
$$
  d_{0k}=\beta\sum_{i=0}^{k-1}(i+1)\beta^i d_{0,k-1-i}+\beta (A_{k-1} +B_{k-1}),
$$
$$
  A_k =\beta\sum_{i=0}^{k-1}(i+1)\beta^i d_{0,k-1-i}+\beta (A_{k-1}+B_{k-1}),
$$
$$
  B_k =\beta\sum_{i=0}^{k-1}(i+1)^2 \beta^i d_{0,k-1-i}+\beta (A_{k-1}+B_{k-1}). 
$$ 
Now, at $\beta$ small enough, $(i+1)\beta^i \le (2\beta )^i$ and $(i+1)^2 \beta^i \le (2\beta )^i$. 
We shall prove by induction that $d_{0k}\le (4\beta )^k$, $A_k \le (4\beta )^k$, $B_k \le (4\beta )^k$.
Indeed, we have
$$  
  d_{0k}\le\beta\sum_{i=0}^{k-1}2^i \beta^i 4^{k-1-i}\beta^{k-1-i}+2\beta 4^{k-1}\beta^{k-1}=
   4^{k-1}\beta^k (2+\sum_{i=0}^{k-1}1/2^i ) \le 4^k \beta^k .
$$
We obtain the estimates for $A_k$ and $B_k$ analogously.    

We get from (8.8) $d_{ik}\le 4^{k-i}\beta^k$ and, hence, 
$$
  \sum_{i=0}^k c_{ik}\le\sum_{i=0}^k d_{ik}\le \frac{4^{k+1}}{3}\beta^k .
$$  
We obtain estimates (8.1) with $C=5c$, where $c$ is the constant from (8.4) - (8.6). 
 
We obtained the bounded solution to equation (7.13) and, hence, a solution to equation (7.11) of type (7.12). 
$\Box$  

A corollary of Lemma 5 is 

\begin{proposition} Suppose a function $\psi$ defined on $\Omega$ satisfies the estimate
$$
  |\psi (w)|\le C|w-\hat w|^{-N}  .
$$
Then
$$
  |P_m ({\rm id}-\mu T_m )^{-1}\psi (w)|\le C_m C|w-\hat w|^{1-N}  \eqno (8.9)
$$  
at $m\ge N+6$, the constant $C_m$ depends on $m$.
\end{proposition}

{\bf Proof}. We follow the method explained before formulation of Lemma 5. 
We can write $P_m ({\rm id}-\mu T_m )^{-1}\psi$ as the sum
$$
  (P_m +P_m \mu T_m +P_m \mu T_m \mu T_m +...)\psi  .   \eqno (8.10) 
$$   
Recall that ${\cal P}_w (\omega )$ is the kernel of $P_m$
$$
  {\cal P}_w (\omega )= \frac{1}{w-\omega}\left(\frac{\omega -\hat\omega}{w-\hat\omega}\right)^m
$$ 
and denote by $K_m (w,\omega )$ the kernel of $T_m$
$$
  K_m (w,\omega )=\frac{1}{w-\omega}\left(\frac{\omega -\hat\omega}{w-\hat\omega}\right)^m \left(\frac{1}{w-\omega}+\frac{m}{w-\hat\omega}\right) .
$$
  In each term of sum (8.10) we can change the order of integration, for example,
$$   
  P_m \mu T_m \mu T_m \psi (w)=\int {\cal P}_w (\omega )\mu (\omega )\int K_m (\omega,\omega_1 )\mu (\omega_1 )\int K_m (\omega_1 ,t)\psi (t)dS_t dS_{\omega_1}dS_{\omega}=
$$   
$$
  =\int\psi (t)\int\mu (\omega_1 )K_m (\omega_1 ,t)\int{\cal P}_w (\omega )\mu (\omega )K_m (\omega,\omega_1 )dS_{\omega}dS_{\omega_1}dS_t =
$$
$$
  =\int\psi (t)\tilde T_m \mu \tilde T_m \mu {\cal P}_w (t)dS_t .  
$$   
We see that sum (8.10) equals to
$$
  \int\psi(t)({\rm Id} -\tilde T_m \mu )^{-1}{\cal P}_w (t)dS_t  .
$$   
By Lemma 5 at $m\ge N+6$, we can represent the expression $({\rm Id} -\tilde T_m \mu )^{-1}{\cal P}_w (t)$
as a function of type (7.12) for $l\ge N$, $m\ge l+6$. Therefore,
$$
  \frac{1}{\pi}\int_{\Omega}g_w (\omega )\psi (\omega )\left(\frac{\omega -\hat\omega}{w-\hat\omega}\right)^l \frac{dS_{\omega}}{w-\omega}  
$$
is another form of sum (8.10) and, by Proposition 30, this integral satisfies estimate (8.9). 
$\Box$ 

Now we have all instruments to prove Theorem 2'.

Let $F$ be a quasiconformal local homeomorphism with the complex dilatation $\mu$ defined on $\Omega_t$. 
Then $f_1 =(\log F_w )_z=F_{ww}/F_w$ satisfies the equation
$$
  (f_1 )_{\bar w}=\mu (f_1 )_w +\mu_w (f_1 )+\mu_{ww} .   \eqno (8.11)
$$   
From the other hand, if we have a solution to this equation we can find a  $\mu$-quasiholomorphic function $F$ 
such that $f_1 =F_{ww}/F_w$. Indeed, if we define $\tilde f_1 =\mu f_1+\mu_w$, then equation (8.11) implies  
$f_{\bar w}=\tilde f_w$ and we can define the function
$$
  g(w)=\int_0^w f_1 ds+\tilde f_1 d\bar s  .
$$
Here the integral doesn't depend on a way. Now we have
$$
   (e^g )_{\bar w}=\tilde f_1 e^g=(\mu f_1 +\mu_w )e^g= (\mu e^g)_w .
$$
Hence, we can define
$$
   F(w)=\int_0^w e^g ds+\mu e^g d\bar s .  \eqno (8.12)
$$ 
This map is $\mu$-quasiholomorphic and $F_{ww}/F_w =f_1$.

A generalization of equation (8.11) is the equation for $f_k =(f_1 )_{w^{k-1}}$:
$$
   (f_k )_{\bar w}=\mu (f_k )_w +P_k  ,\, P_k =\mu_w f_k +(P_{k-1})_w  .   \eqno (8.13)
$$
Here $P_j$, $1\le j\le k$ are defined by induction functions from $f_1,...,f_k $, where  
$$   
  f_j =(f_{j-1})_w ,\,j=2,...,k   .    \eqno (8.14)  
$$  

Suppose now that we have a solution to equation (8.13) with some functions $f_1,...,f_k $
not necessary satisfying relations (8.14). We define 
$$
  \tilde f_j =\mu f_j +P_{j-1},\,1\le j\le k   ,  \eqno (8.15)
$$  
From (8.13) follows $(\tilde f_k )_w =(f_k )_{\bar w}$ and the function 
$$
   \int_0^w f_k dt+\tilde f_k d\bar t 
$$   
is well-defined (i.e., it doesn't depend on a way of integration). Suppose now that we
have the relations
$$
   f_{j-1}=\int_0^w f_j dt+\tilde f_j d\bar t ,\,2\le j\le k .  \eqno (8.16)
$$
Than, by induction, these functions are well-defined and relations (8.14) are satisfied.   
Also, by induction, we obtain the relations
$$
  (f_j )_{\bar w}=\mu (f_j )_w +P_j , \,1\le j\le k  .  (8.17)
$$
Thus to obtain a solution to equation (8.13) with functions $f_j$ satisfying relations (8.14), 
(8.17) it is enough to satisfy equations (8.13), (8.15), and (8.16). We shall consider these equations 
as a system and we shall find a solution to this system satisfying estimates of Theorem 2'.

{\bf Remark}. It seems, the important case is $k=1$. The equations for higher derivatives
we consider mainly for completeness. 

\begin{proposition}. In conditions of Theorem 2' there exists a solution to system(8.13), (8.15), (8.16) satisfying estimates
$$
  |f_j (w)|\le Cb|w-\hat w|^{-j},\,1\le j\le k   .  \eqno (8.18)
$$
\end{proposition}

{\bf Proof}.
We can write $P_k$ as
$$
  P_k =\sum_{i=1}^k n_{ki}\mu_{w^i}f_{k+1-i}+\mu_{w^{k+1}}  ,    \eqno (8.19)
$$
where $n_{ik}$ are integer.

We solve system (8.13), (8.15), (8.16) by an iteration method. On the first step we solve the equation   
$$
  (f_{k1})_{\bar w}=\mu (f_{k1})_w +\mu_{w^{k+1}} .   \eqno (8.20)
$$
We have the solution to this equation
$$
  f_{k1}=P_m ({\rm id}-\mu T_m )^{-1}\mu_{w^{k+1}} .   \eqno (8.21)
$$  
Since we have the estimate $|\mu_{w^{k+1}}(w)|\le b|w-\hat w |^{-k-1}$, we obtain by Proposition 42
$$
  |f_{k1}(w)|\le Cb |w-\hat w |^{-k}  \eqno (8.22)
$$  
with some uniform $C$.   

Now we define the iterations
$$
  f_{k,i+1}=P_m ({\rm id}-\mu T_m )^{-1}P_k (f_{ki},...f_{1i}) ,   \eqno (8.23)
$$  
$$
  \tilde f_{j,i+1}=\mu f_{ji}+P_{j-1}(f_{ji},...f_{1i}),\,f_{j-1,i+1}(w)=\int_0^w f_{j,i+1}ds+\tilde f_{j,i+1}d\bar s,\,
 1\le j\le k .   \eqno (8.24)
$$ 

By induction, applying representation (8.19), we can see that
$$
  |f_{ji}(w)|\le Cb|w-\hat w|^{-j} .
$$
Here $C$, seems, can depend on $i$ but, in fact, iteration process (8.23), (8.24) converges
in $C^0_k \times C^0_{k-1}\times...\times C^0_1$ (see Definition 2). Indeed,
$$ 
  f_{k,i+1}-f_{ki} =P_m ({\rm id}-\mu T_m )^{-1}[P_k (f_{ki},...f_{1i}-P_k (f_{k,i-1},...f_{1,i-1})] .   
$$  
By (8.19),
$$
  \|P_k (f_{ki},...f_{1i}-P_k (f_{k,i-1},...f_{1,i-1})\|_{0,k+1}\le Cb(\|f_{ki}-f_{k,i-1}\|_{0,k}+...+
  \|f_{1i}-f_{1,i-1}\|_{0,1}) 
$$
and, hence, analogously to (8.22),
$$   
  \|f_{k,i+1}-f_{ki}\|_{0,k}\le Cb(\|f_{ki}-f_{k,i-1}\|_{0,k}+...+\|f_{1i}-f_{1,i-1}\|_{0,1}) .  \eqno (8.25)
$$
From (8.24) follows
$$
  \|f_{k,i+1}-f_{ki}\|_{0,j}\le Cb(\|f_{ki}-f_{k,i-1}\|_{0,k}+...+\|f_{1i}-f_{1,i-1}\|_{0,1}) .  \eqno (8.26)
$$
The constant $C$ in these inequalities doesn't depend on $i$. Thus at small enough $b$ the iterations converge.
$\Box$

In particular, $f_1$ is a solution to equation (8.11) satisfying the estimate $|f_1 (w)|\le C|w-\hat w|^{-1}$.
Formula (8.12) defines then a $\mu$-quasiholomorphic function $F$ such that 
$|F_{w^2}/F_w (w)|\le Cb/|w-\hat w |$ with some uniform $C$. It is easy to obtain estimates for other derivatives
of $F$. We have
$$
  \frac{F_{w\bar w}}{F_w}=\mu_w+\mu\frac{F_{w^2}}{F_w} 
$$
and we get the estimate $Cb/|w-\hat w |$. Also,   
$$
  \frac{F_{\bar w^2}}{F_w}=\frac{(\mu F_w )_{\bar w}}{F_w}=\mu_{\bar w}+\mu\frac{F_{w\bar w}}{F_w}
$$
and we again obtain the estimate $Cb/|w-\hat w |$. Also, applying estimates (8.18) we get
$$
  |(f_1 )_{(j)}(w)|\le Cb|w-\hat w |^{1-|(j)|},\,|(j)|\le k-1  .
$$
By integration we obtain the estimate
$$
   |\log F_w (w)|\le cb|\log (|w-\hat w |)|
$$
with some uniform $c$. I.e.,
$$
  |w-\hat w |^{cb}\le |F_w (w)|\le |w-\hat w |^{-cb}  .
$$
It is estimate (4.42).

\begin{proposition} $F$ is a homeomorphism.
\end{proposition}

{\bf Proof}. In conditions of Theorem 2' $\Omega =h (D)$, where $h$ is holomorphic and satisfies estimates (4.38) (we omit here the index $t$).
The map $F\circ h$ is quasiholomorphc on $D$ with the Beltramy
coefficient $\tilde\mu (z)=\mu\circ h(z) h_z (z)/\overline{h_z (z)}$. We have at $|(j)|=1$
$$
  |\pd_{(j)}\tilde\mu |\le |\pd_{(j)}\mu\circ h \cdot h_z \cdot h_z /\overline{h_z}+\mu\circ h (h_{z^2}/\overline{h_z}+  h_z \overline{h_{z^2}}/\overline{h_z}^2)|\le
$$
$$
  \le c(b+b\varepsilon )|z-\bar z^{-1}| .    
$$  
Analogously, at $|(j)|=2$
$$   
  |\pd_{(j)}\tilde\mu |\le c(b+b\varepsilon )|z-\bar z^{-1}|^{-2}   
$$
with some uniform $c$. We get for the normal map $f^{\tilde\mu}$ the estimates of Lemma 1 with small
coefficients. From the other hand,
$$
  F\circ h =\tilde h\circ f^{\tilde\mu}  
$$
with some holomorphic $\tilde h$. Denote $\tilde g=h\circ (f^{\tilde\mu})^{-1}$. From Lemma 1 follows
$$
  |\tilde g_z |\ge c,\,|(\tilde g_{\bar z}|\le Cb,\,\pd_{(j)}\tilde g )(z)|\le cb|z-\bar z^{-1}|^{-1},\,|(j)|=2 .
$$
for some uniform $c$, $C$. We see that $\tilde h =F\circ\tilde g$ satisfies the estimate
$$
  \left|\frac{\tilde h_{z^2}}{\tilde h_z}\right|(z)\le 10\sup_{|(j)|=1,|(s)|=2}\left|\frac{\pd_{(s)}F\circ\tilde g \cdot(\pd_{(j)}\tilde g )^2}{F_w \circ\tilde g\cdot\tilde g_z}(z)+
 \frac{\pd_{(j)}F\circ\tilde g \cdot\pd_{(s)}\tilde g}{F_w \circ\tilde g\cdot\tilde g_z}(z)\right|\le\frac{cb}{|z-\bar z^{-1}|} .
$$
By the criterium of univalence (for example [Pom]) $\tilde h$ is an univalent function if
$|(\tilde h_{z^2}/\tilde h_z )(z)|\le 1/(1-|z|^2 )$. This condition is satisfied and, hence, 
$F$ is a homeomorphism. $\Box$  

{\bf Proof of estimates for derivatives with respect to parameters}. 
In what follows some constants can depend on $b$ and we shall supply these constants by the subscript $b$.

We at first consider equation (8.20). We can write solution (8.21) as
$f_{k1} =P_m h$, where $h$ satisfies the equation
$$
  h=\mu T_m h +\mu_{w^k}  .
$$  
For a derivative with respect to a parameter $t$ we get the equation
$$
  h_t =\mu T_m h_t +\mu_t T_m h +\mu (T_m )_t \hat g_t h +\mu_{w^k,t} ,
$$
where  $\hat g (\omega)$ is the function $\omega\mapsto\hat\omega$ and $(T_m )_t$ is the operator with the kernel $(K_m )_t$ and $K_m$ is the kernel of $T_m$. 
$(K_m )_t$ is a sum of items of the types
$$
  \frac{1}{(w-\omega )^2}\frac{1}{w-\hat\omega}\left(\frac{\omega-\hat\omega}{w-\hat\omega}\right)^l ,\,
  \frac{1}{w-\omega}\frac{1}{(w-\hat\omega )^2}\left(\frac{\omega-\hat\omega}{w-\hat\omega}\right)^l ,
 \eqno (8.27)
$$  
where $l=m-1$ or $l=m$. 
According to (6.1) and (5.15), we have estimate $|\hat g_t (\omega )|\le M_b |\omega -\hat\omega |^{-n'}$ for some $c$ and $n'$. 
From the other hand, from the equation $h=({\rm Id} -\mu T_m )^{-1}\mu_{w^k}$ follows $\| h\|_{p,k}\le Cb$
for some uniform $C$. 

The operators with kernels (8.27) are of the types considered in Proposition 34. We obtain
$$
  \|(T_m )_t \hat g_t h \|_{p,n'+3} \le CM_b b
$$
for $m$ great enough, $p\ge 2$ and close enough to 2. From the other hand, $\mu_{wwt}$
and $\mu_t T_m h$ also belong to $L^p_{N}(\Omega )$ with some, maybe new, $N$. We get for $h_t$ the equation
$$
   h_t -\mu T_m h_t = H   , 
$$
where $H$ belongs to $L^p_N (\Omega )$ for some $N$ and has an uniform norm. Applying Proposition 35, we see that
the operator ${\rm Id} -\mu T_m$ is invertible in $L^p_N (\Omega )$ for $m$ great enough and $b=\|\mu \|_{C^0}$
small enough. Thus $h_t$ belongs to $L^p_N (\Omega )$. From the other hand,
$$    
  (f_{k1}(w))_t =P_m h_t +(P_m )_t \hat g_t h ,
$$
where $(P_m )_t$ is the operator with the kernel
$$
 \frac{m}{w-\omega}\left(\frac{\omega-\hat\omega}{w-\hat\omega}\right)^{m-1}\left(\frac{\omega-\hat\omega}{(w-\hat\omega )^2}-\frac{1}{w-\hat\omega}\right)  .
$$
Since $\hat g_t$ and $h_t$ belong to $L^p_N (\Omega )$, we obtain
$$  
  |f_{k1}(w)|\le C(\|h_t \|_{p,N}+\|h\hat g_t \|_{p,N})\|{\cal P}_{wN}\|_q  ,
$$
where $p^{-1}+q^{-1}=1$ and ${\cal P}_{wN}$ is the function
$$  
  {\cal P}_{wN}(\omega )= [({\cal P}_w )(\omega )+ (({\cal P}_w )_t (\omega )](\omega-\hat\omega )^{-N} .
$$  
But
$$
  \|{\cal P}_{wN}\|_q \le C|w-\hat w |^{-(N+1)}
$$
and we obtain  
$$    
  |(f_{k1}(w))_t|\le C(\|h_t \|_{p,N}+\|h\hat g_t \|_{p,N})|w-\hat w |^{-(N+1)}\le CM_b |w-\hat w |^{-(N+1)}  \eqno (8.28) 
$$   
with some uniform $C$.
   
We obtain estimates for higher derivatives of $h$ and $f_1$ by induction. At differentiation there appears
derivatives of the function $\hat g$ and operators with kernels of the types considered in Proposition 34.

Analogously to (8.28) we obtain the estimate
$$    
  |(f_{k1}(w))_{0,(l)}|\le CM_b |w-\hat w |^{-N},\,|(l)|\le L   
$$ 
with some new $N$ and $M_b$ if $m$ is great enough.

Consider now iterations (8.23), (8.24).

We see that $f_{k,i+1}-f_{ki} =P_m \Delta h_{ki}$, where $\Delta h_{ki}$
satisfies the equation
$$
  \Delta h_{ki} =\mu T_m \Delta h_{ki} +P_k (f_{ki},...f_{1i})-P_k (f_{k,i-1},...f_{1,i-1})   \eqno (8.29) 
$$
Also,
$$
\tilde f_{j,i+1}-\tilde f_{ji}=\mu (f_{ji}-f_{j,i+1})+P_{j-1}(f_{ji},...f_{1i})-P_{j-1}(f_{j,i-1},...f_{1,i-1}),  \eqno (8.30)
$$
$$
 f_{j-1,i+1}(w)-f_{j-1,i}(w)=\int_0^w (f_{j,i+1}-f_{ji})ds+(\tilde f_{j,i+1}-\tilde f_{ji})d\bar s,\, 
 1\le j\le k .   \eqno (8.31)  
$$
For the derivative $(\Delta h_{ki})_s$ we get the equation
$$
  (\Delta h_{ki})_t -\mu T_m (\Delta h_{kI})_t =\mu_t T_m \Delta h_{ki}+\mu (T_m )_t \Delta h_{ki}+
$$
$$  
  +[P_k (f_{ki},...f_{1i})]_t -[P_k (f_{k,i-1},...f_{1,i-1})]_t .
$$  
From (8.19) follows that we can write the difference $[P_k (f_{ki},...f_{1i})]_t -[P_k (f_{k,i-1},...f_{1,i-1})]_t$
as a sum of terms of the types
$$
  p_{kj}(f_{ji}-f_{j,i-1}) ,\,p'_{kj}(f_{ji}-f_{j,i-1})_s ,
$$
where $p_{kj}= n_{kj}(\mu_{w^{k+1-j}})_s$ and $p'_{kj}= n_{kj}\mu_{w^{k+1-j}}$  

We have the estimates
$$
 \|p_{kj}\|_{0,N}\le M_b   \eqno (8.32)
$$
for some $N$ and $M_b$ and
$$
  \|p'_{kj}\|_{0,j}\le C  \eqno (8.33)
$$
for some uniform $C$.  

By (8.30, (8.31), applying (8.25), (8.26), we get inductively
$$
  \|\Delta h_{ki}\|_{p,k+1} \le Cb^i .
$$
We also have from (8.25), (8.26)
$$ 
  \|f_{ki} -f_{k,i-1}\|_{0,k}\le C(ab)^i  
$$
for some uniform $C$ and $a$. 

Denote by $G_i$ the sum $\mu_t T_m \Delta h_{ki} +\mu (T_m )_t \Delta h_{ki} +\sum  p_{ki}(f_{ji}-f_{j,i-1})$. 
Then 
$$
  \|G_i \|_{p,N}\le M_b (ab)^i 
$$
for some uniform $N$, $M_b$ and $a$.  We get from (8.29)
$$ 
  (f_{k,i+1} -f_{ki})_t =(P_m )_t \Delta h_{ki}+P_m ({\rm id}-\mu T_m )^{-1}G_i +
$$
$$  
  +P_m ({\rm id}-\mu T_m )^{-1}\sum p'_{kj}(f_{ji}-f_{j,i-1})_t   \eqno (8.34)
$$

Analogously to (8.28) we obtain
$$
  |((P_m )_t \Delta h_{ki} +P_m ({\rm id}-\mu T_m )^{-1}G_i )(w)|\le M_b (ab)^i |w-\hat w |^{-N}  \eqno (8.35)
$$
for some uniform  $N$, $M_b$ and $a$. Estimate now the term
$$
 P_m ({\rm id}-\mu T_m )^{-1}\sum p'_{kj}(f_{ji}-f_{j,i-1})_t   \eqno (8.36)
$$
in (8.34). Suppose that 
$$
  \|(f_{ji}-f_{j,i-1})_t\|_{0,N-k+j}\le M_{bi}    \eqno (8.37)
$$ 
with the same $N$ as in (8.35). Then, applying (8.33), we get
$$ 
   \sum \|p'{kj}(f_{ji}-f_{j,i-1})_t \|_{0,N+1}\le cM_{bi}b  
$$
with the same $N$ and $M_{bi}$ and some $c$ independent of $i$. 
Applying Proposition 42. we obtain that term (8.36) has the estimate by modulus
$$
  cbM_{bi}|w-\hat w |^{-N}   \eqno (8.38)
$$  
with some uniform $c$ independent of $b$ and $i$.

By (8.34), (8.35), and (8.38), we get
$$
  |(f_{k,i+1} -f_{ki})_t |\le (M_b (ab)^i +cbM_{bi} )|w-\hat w|^{-N}  .  \eqno (8.39)
$$  

Now, by (8.24), applying estimates (8.32) and (8.33), it isn't difficult to obtain
$$
  |(f_{j,i+1} -f_{ji})_t |\le (M_b(ab)^i +cbM_{bi} )|w-\hat w|^{-N-j+k},\,1\le j\le k-1 .  \eqno (8.40)
$$
with some $c$ and $a$ independent of $i$. It is an estimate of type (8.37). Modifying, 
if necessary, $a$ and $c$ in (8.39) we can suppose that these constants are the same as in (8.40). 
We get the estimate
$$
  |(f_{j,i+1} -f_{ji})_t |\le (M_b(Ab)^i +AbM_{bi} )|w-\hat w|^{-N-j+k},\,1\le j\le k ,  \eqno (8.41)
$$
where $A=\max\{a,c\}$. 

Modifying $N$, if necessary, we can suppose that $f_{k1}\in C^0_N$ with the same $N$ as in (8.35).   
Then $f_{j1}\in C^0_{N-k+j}$. We put $f_{j,0}=0$, $1\le j\le k$ and $M_{b1}=\max \|f_{j1}\|_{0,N-k+j}$, $1\le j\le k$.
Using (8.41), we obtain inductively
$$
  |(f_{j,i+1} -f_{ji})_s |\le (Ab)^i (M_{b1}+iM_b )
$$
At small enough $b$ the iterations converge in $C^0_{N-k+j}$.
 
We can obtain estimates for higher derivatives with respect to the parameters by the same
method. Instead of (8.34) we obtain the equation for differences $(f_{k,i+1} -f_{ki})_{0,(l)}$ with 
multi-index $(l)$  
$$
  (f_{k,i+1} -f_{ki})_{0,(l)}=H_i +P_m ({\rm id}-\mu T_m )^{-1}\sum p'_{kj}(f_{ji}-f_{j,i-1})_{),(l)},
$$
where for $H_i$ we have the estimate $M_b (ab)^i$ with some $M_b$ and $a$. Other modifications are
obvious. $\Box$.
  
Thus we finished the proof of Theorem 2' and, hence, Theorem 2. 

Remark. It seems, we can prove the estimates for derivatives with respect to parameters
not applying the $L^p$-estimates of Section 6. For it we need in generalization of the estimates
of Section 7 on operators of the types $(P_m )_t$ and $(T_m  )_t$ and similar ones, 
containing $t$-derivatives. However the $L^p$-estimates
can be useful, and it is of some interest that we have estimates $({\rm dist}(w,\pd\Omega ))^{-N}$ for the growth 
 in $L^p$ also.  

\medskip
\centerline{REFERENCES}
\medskip

\begin{itemize}
\item{[Ah]} L. Ahlfors, {\em Lectures on Quasiconformal Mappings,} University Lecture Series 38. Providence,
(AMS), (2006). 
\item{[AhB]} L. Ahlfors, L. Bers,  {\em Riemann's Mapping Theorem for Variable Metrics,}
 Ann. Math. (2) 72 (1960), 385-404.
\item{[As]} K. Astala, T. Iwaniec, G. Martin, {\em Elliptic Partial Differrential equations and
Quasiconformal Mappings in the Plane} Prinston University Press, Princeton, New Jersey, (2009).
\item{[Br1]} M. Brunella, {\em Feuilletages Holomorphes sur les Surfaces Complex Compactes,} Ann. Scient. Ec. Norm. Sup.,
4e serie, 30 (1997), 569-594.
\item{[Br2]} M. Brunella, {\em Plurisubharmonic Variation of of the Leafwise Poincare Metric,} Internat. J. Math. 14 (2003),
139-151.
\item{[Br3]} M. Brunella, {\em Uniformisation of Foliations by Curves,} Holomorphic Dynamical Systems, Lecture Notes in Math., 
Springer-Verlag, (2012), 105-165.
\item{[CDFG]} S. Calsamiglia, B. Deroin, S. Frankel, A. Guillot, {\em Singular sets of holonomy maps for algebraic foliations,}
J. Eur. Math. Soc. 15 (2013), n. 3, 1067-1099.
\item{[DNS]} T.C. Dinh, V.A. Nguyen, N. Sibony, {\em Entropy for Hyperbolic Riemann Surface Laminations I,}
 Frontiers in complex dynamics. In celebration of John Milnor’s 80th birthday. NJ: Princeton University Press. 569-592 (2014).
\item{[Gl1]} A. A. Glutsyuk, {\em Hyperbolcty of the Leaves of a Generic One-dimensiopnal Holomrphic Foliations
on a Nonsingular Projective Algebraic Variety}, Proceedings of V.A.Steklov Iinst. of Math., V. 213, (1997), 90-111.
\item{[Gl2]} A. A. Glutsyuk, {\em On Simultaneous Uniformisation and Local Nonuniformisability,} C. R. Math. Acad. Sci. Paris,
334, No.6 (2002), 489-494. 
\item{[Il1]} Yu. S. Ilyashenko, {\em Foliations by Analytic Curves,} Mat. Sb, 88 (1972), 558-577.
\item{[Il2]} Yu. S. Ilyashenko, {\em Covering Manifolds for Analytic Families of Leaves of Foliations by Analytic Curves,}
Topological Meth. in Nonlinear Analysis, 11 (1998), 361-373.
\item{[Il3]} Yu. S. Ilyashenko, {\em Persistence Theorems and Simultaneous Uniformization,} Proceedings of V.A.Steklov Iinst. of Math.,
V.254, (2006), 184-200.
\item{[Le]} O. Lehto, K. Virtanen, {\em Quasiconformal Mapping in the Plane} Springer-Verlag,
Berlin-New York, (1971).
\item{[LN]} A. Lins-Neto, {\em Uniformization and the Poincare Metric on the leaves of a Foliation by Curves,}
Bol. Soc. Bras. Mat., Nova Ser. 31, No. 3 (2000), 351-366. 
\item{[Pom]} S. Pommerenke, {\em Univalent Functions} Vandenhoeck and Ruprecht, Gottingen, (1975).
\item{[RS]} M. Reed, B. Simon, {\em Methods of Modern Mathematical Physics. Volume 2} 
Academic Press, (1975).
\item{[Sh1]} A.A. Shcherbakov, {\em Metrics and Smooth Uniformisation of Leaves of Holomorphic Foliations,}
Moscow Math. J., 11, No.1 (2011), 157-178.
\item{[Sh1]} A.A. Shcherbakov, {\em Almost Complex Structure on Universal Coverings of Foliations}
Trans. of the Moscow Math. Soc., (2015), 137-179.
\item{[V]} A. Verjovsky, {\em A Uniformization Theorem for Holomorphic Foliations,} Contemporary Math., 58, part III (1987),
233-245.
\end{itemize}

\end{document}